\documentclass[10pt,twoside,openright]{report} 

\usepackage[phd]{edmaths}

\usepackage[ansinew]{inputenc} 
\usepackage[UKenglish]{babel}
\usepackage{textcomp}
\usepackage{mathrsfs}
\usepackage{amsfonts}
\usepackage{mathtools}
\usepackage{amssymb}
\usepackage{amsmath}
\usepackage{amsthm}
\usepackage{amscd}
\usepackage{tensor}
\usepackage{stmaryrd}
\usepackage{rotating}
\usepackage[text={6.8in,8.5in},centering]{geometry}
\usepackage{tikz}
\usetikzlibrary{arrows, calc, decorations.markings, positioning}
\usepackage{tikz-cd}
\usepackage{relsize} 
\usepackage[none]{hyphenat}
\usepackage{hyperref}

\usepackage{graphicx}
\graphicspath{ {./Images/} }

\usepackage{hyperref,xcolor}
\usepackage{makeidx}
\definecolor{deepred}{rgb}{0.5,0,0}
\definecolor{deepblue}{rgb}{0,0,0.5}
\definecolor{deepgreen}{rgb}{0,0.5,0}

\hypersetup{
    colorlinks=true,
    linkcolor=deepred,
    filecolor=deepred,      
    urlcolor=deepred,
    citecolor=deepblue,
}

\newenvironment{acknowledgements}{\chapter*{Acknowledgements}}{\addcontentsline{toc}{chapter}{Acknowledgements}}
\newenvironment{laysum}{\chapter*{Lay Summary}}{\addcontentsline{toc}{chapter}{Lay Summary}}
\newenvironment{preface}{\chapter*{Preface}}{\addcontentsline{toc}{chapter}{Preface}}

\newcommand{\Int}{\mathbb{Z}}
\newcommand{\Rat}{\mathbb{Q}}
\newcommand{\Real}{\mathbb{R}}
\newcommand{\dimm}{\text{dim}}
\newcommand{\Cin}[1]{\text{C}^\infty(#1)}
\newcommand{\Tan}{\text{T}}
\newcommand{\Cot }{\text{T}^*}
\newcommand{\Sec}[1]{\Gamma(#1)}
\newcommand{\Der}{\text{D}}
\newcommand{\Jet}{\text{J}}
\newcommand{\Proj}{\text{pr}}
\newcommand{\Id}{\text{id}}
\newcommand{\rkk}[1]{\text{rk}(#1)}
\newcommand{\Hom}[1]{\text{Hom}(#1)}
\newcommand{\Dim}[1]{\text{Dim}(#1)}
\newcommand{\Dr}[1]{\text{Der}(#1)}
\newcommand{\GL}[1]{\text{GL}(#1)}
\newcommand{\Ker}[1]{\text{ker}(#1)}
\newcommand{\Grph}[1]{\text{grph}(#1)}
\newcommand{\LGrph}[1]{\text{Lgrph}(#1)}

\newcommand{\Diff}{\text{Diff}}
\newcommand{\Obs}[1]{\text{Obs}(#1)}
\newcommand{\Dyn}[1]{\text{Dyn}(#1)}
\newcommand{\Lie}{\textsf{Lie}}
\newcommand{\Crnt}{\textsf{Crnt}}
\newcommand{\Dir}{\textsf{Dir}}
\newcommand{\LCrnt}{\textsf{LCrnt}}
\newcommand{\LDir}{\textsf{LDir}}
\newcommand{\Man}{\textsf{Man}}

\newcommand{\Ring}{\textsf{Ring}}
\newcommand{\Vect}{\textsf{Vect}}
\newcommand{\LVect}{\textsf{LVect}}
\newcommand{\Line}{\textsf{Line}}
\newcommand{\Symp}{\textsf{Symp}}
\newcommand{\Poiss}{\textsf{Poiss}}
\newcommand{\Cont}{\textsf{Cont}}
\newcommand{\Acts}{\mathbin{\rotatebox[origin=c]{-90}{$\circlearrowright$}}}
\newcommand{\dtimes}{\mathbin{\rotatebox[origin=c]{90}{$\ltimes$}}}
\newcommand{\utimes}{\mathbin{\rotatebox[origin=c]{-90}{$\ltimes$}}}
\newcommand{\ftimes}[2]{\tensor[_#1]{\times}{_#2}}
\newcommand{\LDer}{\mathcal{L}}

\newtheorem{prop}{Proposition}
\numberwithin{prop}{section}

\newtheorem{thm}{Theorem}
\numberwithin{thm}{section}

\newtheorem{conj}{Conjecture}
\numberwithin{conj}{section}

\numberwithin{equation}{section}

\makeatletter
\newenvironment{timeline}[6]{%

    \newcommand{\startyear}{#1}
    \newcommand{\tlendyear}{#2}

    \newcommand{\yearcolumnwidth}{#3}
    \newcommand{\rulecolumnwidth}{#4}
    \newcommand{\entrycolumnwidth}{#5}
    \newcommand{\timelineheight}{#6}

    \newcommand{\templength}{}

    \newcommand{\entrycounter}{0}

    \long\def\ifnodedefined##1##2##3{%
        \@ifundefined{pgf@sh@ns@##1}{##3}{##2}%
    }

    \newcommand{\ifnodeundefined}[2]{%
        \ifnodedefined{##1}{}{##2}
    }

    \newcommand{\drawtimeline}{%
        \draw[timelinerule] (\yearcolumnwidth+5pt, 0pt) -- (\yearcolumnwidth+5pt, -\timelineheight);
        \draw (\yearcolumnwidth+0pt, -10pt) -- (\yearcolumnwidth+10pt, -10pt);
        \draw (\yearcolumnwidth+0pt, -\timelineheight+15pt) -- (\yearcolumnwidth+10pt, -\timelineheight+15pt);

        \pgfmathsetlengthmacro{\templength}{neg(add(multiply(subtract(\startyear, \startyear), divide(subtract(\timelineheight, 25), subtract(\tlendyear, \startyear))), 10))}
        \node[year] (year-\startyear) at (\yearcolumnwidth, \templength) {\startyear};

        \pgfmathsetlengthmacro{\templength}{neg(add(multiply(subtract(\tlendyear, \startyear), divide(subtract(\timelineheight, 25), subtract(\tlendyear, \startyear))), 10))}
        \node[year] (year-\tlendyear) at (\yearcolumnwidth, \templength) {\tlendyear};
    }

    \newcommand{\entry}[2]{%

        \pgfmathtruncatemacro{\lastentrycount}{\entrycounter}
        \pgfmathtruncatemacro{\entrycounter}{\entrycounter + 1}

        \ifdim \lastentrycount pt > 0 pt%
            \node[entry] (entry-\entrycounter) [below of=entry-\lastentrycount] {##2};
        \else%
            \pgfmathsetlengthmacro{\templength}{neg(add(multiply(subtract(\startyear, \startyear), divide(subtract(\timelineheight, 25), subtract(\tlendyear, \startyear))), 10))}
            \node[entry] (entry-\entrycounter) at (\yearcolumnwidth+\rulecolumnwidth+10pt, \templength) {##2};
        \fi

        \ifnodeundefined{year-##1}{%
            \pgfmathsetlengthmacro{\templength}{neg(add(multiply(subtract(##1, \startyear), divide(subtract(\timelineheight, 25), subtract(\tlendyear, \startyear))), 10))}
            \draw (\yearcolumnwidth+2.5pt, \templength) -- (\yearcolumnwidth+7.5pt, \templength);
            \node[year] (year-##1) at (\yearcolumnwidth, \templength) {##1};
        }

        \draw ($(year-##1.east)+(2.5pt, 0pt)$) -- ($(year-##1.east)+(7.5pt, 0pt)$) -- ($(entry-\entrycounter.west)-(5pt,0)$) -- (entry-\entrycounter.west);
    }

        \begin{tikzpicture}
        \tikzstyle{entry} = [%
            align=left,%
            text width=\entrycolumnwidth,%
            node distance=10mm,%
            anchor=west]
        \tikzstyle{year} = [anchor=east]
        \tikzstyle{timelinerule} = [%
            draw,%
            decoration={markings, mark=at position 1 with {\arrow[scale=1.5]{latex'}}},%
            postaction={decorate},%
            shorten >=0.4pt]

        \drawtimeline
}
{
    \end{tikzpicture}
    \let\startyear\@undefined
    \let\tlendyear\@undefined
    \let\yearcolumnwidth\@undefined
    \let\rulecolumnwidth\@undefined
    \let\entrycolumnwidth\@undefined
    \let\timelineheight\@undefined
    \let\entrycounter\@undefined
    \let\ifnodedefined\@undefined
    \let\ifnodeundefined\@undefined
    \let\drawtimeline\@undefined
    \let\entry\@undefined
}
\makeatother

\usepackage[
backend=biber,
style=alphabetic,
sorting=nyt
]{biblatex}
\title{A Landscape of Hamiltonian Phase Spaces\\-\\On the foundations and generalizations of one of the most powerful ideas of modern science.}
\author{Carlos Zapata-Carratal{\'a}}
\date{2019}

\begin{document}

\parindent0pt

\pagenumbering{roman}

\maketitle

\declaration

\dedication{To my parents and teachers.\\ A mis padres y profesores.}

\dedication{In memory of Michael Atiyah.}

\begin{laysum}

Modern science and engineering uses mathematical quantities to represent parts of the world. These quantities need to be operated following some careful rules, known as dimensional analysis, if we want to obtain practical results from the theoretical knowledge we have about the world.\newline

One of the oldest and, to this day, most useful, bits of knowledge about the world is understanding how things move around and how some object's motion affects the motion of other objects. This is what we call mechanics.\newline

What we know today about how things move is based on a very long history of different ideas and theories that dates back $100$s of years. The invention of modern mathematics, at least modern from a historical perspective, has been motivated by the problems of mechanics. For example, the derivatives and integrals that we calculate today are a result of a group of stubborn physicists and mathematicians $300$s years ago trying to understand why the planets moved precisely in the way we see them move in the sky.\newline

The mathematics that we use to describe the motion of objects nowadays dates back roughly $200$ years. But the rules to operate with the quantities we measure in practice were only decided $100$ years ago. This difference means that today, when we want to, for example, send a rocket to Mars, we first measure the distance that the rocket needs to travel, let's say $50,000,000$ km, then the scientists and engineers come up with the right equations to calculate how long it will take to get there and then they introduce the digit $50,000,000$ without the units `km' for a computer or calculator to give back another digit, let's say $270$. If they used all the rules right, it is possible to give a unit to the new digit and conclude that the trip will last for $270$ days.\newline

This works in practical situations because scientists and engineers are always paying attention to what units their digits carry but when theoretical scientists and mathematicians study equations they don't worry about the units.\newline

In this thesis we try to understand two things that can be seen in the example of the rocket:
\begin{enumerate}
    \item What is the best mathematical way to represent the travel of a rocket to Mars? We need to know where the rocket is at all times but that's not enough, it is also important to know how quickly it is travelling. Imagine that in all our calculations to travel to Mars we only consider our position, then if we arrive on Mars, but we are traveling super fast, we would crash and all our efforts would have been for nothing!
    \item Is it possible that my mathematics take care of the units by themselves? I don't want to have to remove them manually and insert them back in every time I want to compute something. It is a hassle and there is more risk of making mistakes!
\end{enumerate}

\end{laysum}

\begin{abstract}

In this thesis we aim to revise the fundamental concept of phase space in modern physics and to devise a way to explicitly incorporate physical dimension into geometric mechanics. To this end we begin with a nearly self-contained presentation of local Lie algebras, Lie algebroids, Poisson and symplectic manifolds, line bundles and Jacobi and contact manifolds, that elaborates on the existing literature on the subject by introducing the notion of unit-free manifold.\newline

We give a historical account of the evolution of metrology and the notion of phase space in order to illustrate the disconnect between the theoretical models in use today and the formal treatment of physical dimension and units of measurement.\newline

A unit-free manifold is, mathematically, simply a generic line bundle over a smooth manifold but our interpretation, that will drive several conjectures and that will allow us to argue in favor of its use for the problem of implementing physical dimension in geometric mechanics, is that of a manifold whose ring of functions no longer has a preferred choice of a unit.\newline

We prove a breadth of results for unit-free manifolds in analogy with ordinary manifolds: existence of Cartesian products, derivations as tangent vectors, jets as cotangent vectors, submanifolds and quotients by group actions. This allows to reinterpret the notion of Jacobi manifold found in the literature as the unit-free analogue of Poisson manifolds. With this new language we provide some new proofs for Jacobi maps, coisotropic submanifolds, Jacobi products and Jacobi reduction.\newline

We give a precise categorical formulation of the loose term `canonical Hamiltonian mechanics' by defining the precise notions of theory of phase spaces and Hamiltonian functor, which in the case of conventional symplectic Hamiltonian mechanics corresponds to the cotangent functor. Conventional configuration spaces are then replaced by line bundles, what we call unit-free configuration spaces, and, after proving some specific results about the contact geometry of jet bundles, they are shown to fit into a theory of phase spaces with a Hamiltonian functor given by the jet functor. We thus find canonical contact Hamiltonian mechanics.\newline

Motivated by the algebraic structure of physical quantities in dimensional analysis, we redevelop the elementary notions of the theory of groups, rings, modules and algebras by implementing an addition operation that is only defined partially. In this formalism, the notions of dimensioned ring and dimensioned Poisson algebra appear naturally and we show that Jacobi manifolds provide a prime example. Furthermore, this correspondence allows to recover an explicit Leibniz rule for the Jacobi bracket, Jacobi maps, coisotropic submanifolds, and Jacobi reduction as their dimensioned Poisson analogues thus completing the picture that Jacobi manifolds are the direct unit-free analogue of ordinary Poisson manifolds.

\end{abstract}

\begin{acknowledgements}

The first person I would like to thank is Jose Figueroa-O'Farrill, who acted as my supervisor during all my years in Edinburgh and from whom I learnt a great deal. His help in resolving the many issues on which I got constantly stuck and his understanding during difficult times have de facto enabled this PhD thesis to be a reality. I thank Henrique Bursztyn for the opportunity he gave me to visit IMPA, for welcoming me in such a stimulating research environment and for the financial support that covered for my stay in Rio de Janeiro. I am indebted to Luca Vitagliano and Jonas Schnitzer, who very generously read my manuscripts and offered invaluable comments and remarks. I also thank Adolfo de Azcarraga for his continued encouragement and guidance and Jan Waterfield for helping me venture into the world of the baroque harpsichord, something that would provide much needed musical reliefs during my studies.\newline

A special mention goes to the late Michael Atiyah, with whom I had the privilege to engage in deep conversations, be present in a circle of academic all-stars, share Tarta de Santiago and weather media storms. Had it not been for the job he offered me as his assistant, the last few months of my PhD would have been financially troublesome. I will always dearly remember Michael's energy and passion about mathematics and science, even to his very final days, and how, in my position of his assistant, I was able to accompany him in the twilight of his life.\newline

I would like to thank the many members of the staff at the University of Edinburgh that have been immensely helpful during all this years, oftentimes going the extra mile just for my convenience. In particular Chris, Jill and Iain.\newline

I am very grateful to the people that accompanied me during my time in Edinburgh and that I now consider my friends. Xiling for welcoming me into my first office in JCMB for all the great times, from Scottish lakes to Mediterranean beaches. Tim for the very thought-provoking conversations about mathematics, physics, life and everything. Manuel for passing on the baton of PG rep and rightful title of Queen of the Postgraduates and for making our flat into our home. Matt for his precious time spent in my uninformed algebra questions and all the things we have come to enjoy together. Ruth for her musical collaboration. Scott for all the LAN party experiences and being my first gaming-fueled fried. Ross for the camaraderie of sharing supervisor and the many interesting discussions. Simon for creating a peaceful and quiet space for life in our flat. And very specially to Veronika, who helped me in complicated times and with whom I was lucky enough to share some of the most beautiful experiences I had in Edinburgh during my time there. Lastly, I cannot forget Eimantas, my gaming-friend-turned-flatmate during the last year of my PhD, for featuring proudly in the \emph{Chronicles of Oof} by my side, providing timely resupplies, helping me erect the olive monument assiduously joining the prayers to RNGesus, always keeping Ricardo in sight and gracing our flat with the distinguished presence of one true serendipitous genius Tommy Wiseau.\newline

Pau and I have grown academically together since our undergraduate in Valencia, us both being founding members of LEFT/SEMF and two of the ``Azcarraga's Three''. Good fortune wanted us in the same city for our respective PhDs, which has meant that, at a later time in life, we still look for enlightenment together and mature through early adult life. The passions that have driven me during the years of my doctoral studies (mathematics foundations, physics, music, strategy gaming, romantic paradigms...) have direct connection to concrete conversations with Pau and some form of mutual challenge or frustration for lack of understanding. My gratitude to him for all the very enlightning moments (although fun was not had) and for the long-lasting friendship that joins us.\newline

The first time I thought of the problems that ended up motivating the present thesis I was at the Deep Thinking room of Strohlehenalm, a mountain cabin situated in the snowy Austrian alps, where the Rethinking Foundations of Physics Workshop 2017 was taking place. I would like to thank the organizers of that particular event, and the members of the Basic Research Community for Physics as a whole, for encouraging me to think about these matters. I would like to mention Fede and Flavio, in particular, for the many conversations we have shared about science and life.\newline

Other people who I encountered along the way and with whom I shared small but beautiful portions of the last few years and I remember in gratitude: Didac, Silvia, Carlos, Tara, Jaime, Maria, David, Mariana, Abi...\newline

Carmen ha sido una companera de viaje imprescindible durante todos los a{\~n}os del doctordado. Ademas de las inumerables experiencias preciosas que hemos compartido, musicalmente, en otros hemisferios y con otras personas, su presencia en la distancia me ha permitido asipirar a mejores metas y a orientar mi vida de una forma mucho mas ordenada. Gran parte del esfuerzo que ha quedado plasmado en estas paginas es, sin duda, compartido con ella. Mi gratitud es infinita cuando pienso en el papel crucial que Carmen a jugado en mi maduraci{\'o}n y mi transici{\'o}n a la vida adulta, siendo una referencia de aceptaci{\'o}n y de cari{\~n}o incondicional.\newline

La SEMF ha sido un proyecto ilusionante que he llevado a cabo con la imprescindible ayuda de {\'A}lvaro, a quien menciono con especial afecto y gratitud por todo lo compartido en la distancia y por su dedicaci{\'o}n a la causa de la promoci{\'o}n del conocimiento y acervo intelectual de la sociedad en la que vivimos.\newline

Tambi{\'e}n querr{\'i}a recordar a mis amistades forjadas en la terreta d'Alacant, y por como nos hemos acompa{\~n}ado en una distacia que cuando se reduc{\'i}a, aun permit{\'i}a vivir cercanamente y celebrar nuestros proyectos conjuntamente: Camps(a), Badyn, Dieguette, Fit, Pascual, Rub{\'e}n, Borja, Aldo y Jaime.\newline

Con la perspectiva que me da mi historia acad{\'e}mica, querr{\'i}a agradecer a aquellos que supieron reconocer valor en mis aspiraciones y sue{\~n}os y que me apoyaron, de forma m{\'a}s o menos directa, en mi camino de crecimiento intelectual: Batiste Torregrosa, Francisco Carratala, Alberto Ferrer, Jose Reyes, Jose Gallego, David Grau, Fernando Parre{\~n}o, Fernando Moreno, Jaume Reus y Lydia Moll{\'a}.\newline

Y por {\'u}ltimo me gustaria recordar a mi familia. A mi Abuelita, que ya no vera estas l{\'i}neas, pero quien marc{\'o} fuertemente mis recovecos mas emocionales y a quien tendr{\'e} ligada la experiencia de la musica por el resto de mi vida. A mis t{\'i}os, t{\'i}as y abuelos por su apoyo y cari{\~n}o constante. A Luc{\'i}a por haber sido referencia sin saberlo. A Marcos y Alberto por dejarme ser con ellos un ni{\~n}o otra vez. A Ana por haber crecido juntos y por poder disfrutar ahora de una adultez temprana desde puntos alejados pero que a la vez se sienten muy cercanos. A Carmen (Potito) que pas{\'o} de ser el beb{\'e}, a mi compa{\~n}era de juegos y  a ser ahora una amiga {\'i}ntima con la que comparto mucho m{\'a}s que nuestra infancia. Y por {\'u}ltimo, a mis padres, sin quienes, evidentemente, no estar{\'i}a escribiendo estas l{\'i}neas. Pap{\'a} es el responsable {\'u}ltimo de esta tesis doctoral y su papel como gu{\'i}a, apoyo y est{\'i}mulo, infiri{\'o} en m{\'i} el momentum para continuar el resto de la trayectoria de mi vida con confianza e ilusion. Y Mam{\'a}, de quien nac{\'i} como persona ademas de ser quien me dio a luz. Sin el ba{\~n}o de amor que he sentido en casa todo este tiempo y sin los esfuerzos constantes por acercarme a mi centro emocional, nada de lo que he hecho en los {\'u}ltimos a{\~n}os habria sido posible.\newline

Y a Tibu, por haberme acompa{\~n}ado en cada hito de mi camino acad{\'e}mico.

\end{acknowledgements}

\tableofcontents

\pagenumbering{arabic}

\begin{preface}
\section*{\emph{or how a floating buoy took me from the Mediterranean Sea to the Firth of Forth, the Austrian Alps and Ipanema Beach}}

My father Paco, a hydraulics engineer (see \cite{zapata2003balsas}), introduced me to physics and mathematics as a hobby when I was 12 years old. During my high school years, I forged a close friendship with Alberto Balseyro, who would eventually follow a similar path to mine an become a researcher in mathematical physics (see \cite{alonso2018domain}). The initial years of our friendship mostly consisted in constantly challenging each other with physics problems of increasing difficulty. Since this was happening outside any formal education program, our focus was always the conceptual clarity of the solutions and the precise physical intuitions behind the mathematical tools necessary to reach them. We almost exclusively dealt with mechanical problems, but one day Alberto found a problem involving the buoyancy force of a submerged ball in some of the, to us, at the time, arcane, books that belonged to his uncle. Seeing that I was completely unfamiliar with hydrostatics, he was quick to challenge me with it. Around that time I had been particularly interested in multivariable calculus, as I was trying my best to master the rudiments of surface and volume integrals from the fantastic introductory textbook by Marsden and Tromba \cite{marsden2003vector}, and so I had a working understanding of pressure as a scalar field and the elements of infinitesimal force as a vector field. I felt I had all the ingredients to be able to directly compute the buoyancy force  as a surface integral on a sphere but, unfortunately, that relatively simple vector calculus exercise proved too difficult for me at the time and I had to use my last retort: asking my father.\newline

When I showed the problem to my father he immediately replied that he could hardly believe I was struggling with it and he quickly moved on to the blackboard that I had eagerly installed in my bedroom to show me the simple solution. He wrote in a single line that the buoyancy force was precisely the weight of the displaced volume of fluid. At that point I had never heard of Archimedes Principle and so I was completely baffled by the brevity of his answer. I then voiced my skepticism and demanded further justification. My father explained the intuitive arguments of the usual informal proof of Archimedes Principle but I was left deeply dissatisfied because I felt like I had all the tools to arrive to the same result from first principles solely by mathematical computations. This frustration gave me renewed energies to go back to the study of vector calculus and, after some time struggling with Jacobians and angular integrals, I eventually found the correct result. This is to this day, one of my happiest memories involving mathematics, since, to my teenage self, it seemed like a tremendous achievement. Years later I would learn about Stokes' theorem and how it can be used to prove Archimedes Principle for objects of arbitrary shape establishing the intuitions that had inspired me to postulate my explicit formula for the buoyancy force of a submerged sphere.\newline

This anecdote taught me a lesson that I have taken to heart ever since, and that is clearly reflected in my philosophy when approaching mathematical physics, in general, and the present thesis, in particular:
\begin{quote}
    \emph{If enough efforts are invested into capturing physical intuitions in a  mathematically precise way, even when this may seem superfluous and pedantic in a context of working models and effective theories, the mere structure of mathematics will yield physically insightful outcomes}.
\end{quote}

When I arrived in Edinburgh to start my PhD I didn't have any particular topic in mind, although I had been leaning towards differential geometric subjects for some time during the final years of my undergraduate and masters degree. Jose Figueroa-O'Farrill, my supervisor during my years in Edinburgh, pointed to some recent results on generalizations of symplectic reduction, a subject I had greatly enjoyed studying during my masters in Cambridge the year before, and so, motivated by some initial questions on how to extend notions that are well understood for Poisson manifolds such as reduction or a homological treatment of symmetries, we began an exploration of the field of Dirac geometry and Lie algebroids. The regular meetings with Jose, whose distinctive pedagogical style of presentation I greatly benefited from, helped me get deeper understanding of many topics on, and directly related to, Poisson geometry. After some time studying several generalizations of Poisson manifolds, I began to notice a pattern, what later became the Poissonization functors discussed in this thesis, and Jose encouraged me to investigate them further.\newline

I have had an interest in the foundations of classical and quantum mechanics since my early years of solving mechanics problems for fun. I was happy to discover that Poisson geometry and Lie algebroids offer the natural language in which many of the foundational problems involving phase spaces and quantization are formulated. In line with this interest in foundational questions, I applied to give a talk at the $2017$ Rethinking Foundations of Physics Workshop organized by the BRCP as a yearly event where a small selected group of young researchers meet in a mountain cabin near Dorfgastein, Austria, for an intense week of talks and focused discussion. The fruitful interaction with Federico Zalamea, one of the workshop attendees, led to the question of how to formally implement physical quantities in mathematical physics. After some serious thought during the workshop, this question was left unresolved and, when I was back to civilization, I asked colleagues and searched the literature for a clear answer to realize that almost no systematic research had been conducted in that direction. At the time, I just thought the question was not much more than a curiosity that I was personally itching to resolve, but that had no real physical or mathematical relevance, and so I moved on to continue with my PhD.\newline

During my doctoral visit to Henrique Bursztyn at IMPA, Rio de Janeiro, Brazil, during the northern spring of $2018$, I was lucky to meet Jonas Schnitzer, who was also visiting Henrique. One of the many generalizations that Jose and I had considered during our explorations were contact manifolds, that we regarded at the time as little more than the odd-dimensional analogues of symplectic manifolds. Jonas happened to be doing his PhD with Luca Vitagliano, who has been championing a line bundle approach to the more general Jacobi structures for some years and was responsible for the line bundle generalization of Dirac geometry, the ordinary version of which I had studied extensively during my years in Edinburgh. Jonas and I had many engaging interactions that eventually resulted in an increasing interest on the line bundle approach to Jacobi manifolds on my part. When I went back to Edinburgh I had come up with a more concrete research program, directly motivated by Jonas' explanations about the analogies between Poisson and Jacobi manifolds, that I presented to Jose, who agreed would be worth pursuing.\newline

A few months later, while I was working on this question, I was studying the natural structures that appeared on the jet bundle of the tensor product of two line bundles when I had a true moment of serendipity: the sections of the tensor powers of a line bundle seemed to behave algebraically precisely according to the axioms that we had decided characterized how physical quantities operate in conventional dimensional analysis back in the $2017$ Rethinking Workshop. Suddenly, what was, until then, a purely mathematical speculation gained a clear physical interpretation that demanded a much more precise answer. I was immediately thrown back to the question we had unsuccessfully tackled during the $2017$ Rethinking Workshop and, in light of some positive results I had already found on the interpretation of line bundles as generalized configuration spaces, the question about the role of line bundles in a canonical formalism of mechanics now became one that directly connected to my original question about the implementation of physical dimension into classical mechanics. My thesis was completed with the goal of finding a satisfactory answer to that question within the mathematical framework that I had been studying and developing during my doctorate.
\end{preface}

\chapter{Introduction} \label{Introduction}

\section{Of Phase Spaces and Units} \label{OfPhaseAndUnits}

This thesis is aimed at understanding the formal aspects of phase space kinematics and how one could implement the standard use of physical quantities carrying units within the mathematical framework of mechanics. In particular, we will be focusing on Hamiltonian mechanics, in which conservative physical systems are described by dynamical equations on Poisson or symplectic manifolds.\newline

The manipulation of physical quantities in conventional science and engineering is encapsulated in the rules and prescriptions of dimensional analysis. From a purely mathematical point of view we find two immediate questions:
\begin{quote}
    \emph{What is the natural mathematical framework for dimensional analysis and how can it be connected to a first-principles definition of physical quantities and measurement procedures?}
\end{quote}
\begin{quote}
    \emph{If dimensional analysis is given a precise algebraic formulation, what are its natural generalizations?}
\end{quote}
In relation to our interest in mechanics, we pose:
\begin{quote}
    \emph{Having identified the basic algebraic structure of physical quantities, how can we explicitly implement them in the mathematical formalisms of modern physics, and, in particular, in classical mechanics?}
\end{quote}

One of the central ideas argued in the body of the thesis is the identification of lines and line bundles as the appropriate objects to encapsulate the notion of freedom of choice of units, abstractly and in the context of mechanics respectively. In the mathematical preliminaries we shall present contact geometry under this light and thus the following question appears:

\begin{quote}
    \emph{If contact and Jacobi manifolds are the line bundle equivalents of Poisson and symplectic manifolds respectively, is it reasonable to expect that, if one replaces the configuration manifolds of classical mechanics by line bundles, canonical contact structures will appear playing a similar role to that of canonical symplectic forms on cotangent bundles in conventional classical mechanics?}
\end{quote}

Following from this, and connecting to the problem of incorporating physical quantities and units into mechanics, we wonder

\begin{quote}
    \emph{Can a formalism that takes line bundles, the unit-free analogues of conventional configuration spaces, into contact manifolds provide a theory of mechanics that is analogous to canonical Hamiltonian mechanics and that naturally incorporates the notion of physical dimension?}
\end{quote}

Throughout the thesis these questions shall be precisely formulated and, eventually, answered.

\section{The Categorical Imperative} \label{CategoricalImperative}

Throughout the course of my PhD, after numerous enlightning conversations with Pau Enrique (see \cite{moliner2017space}), I have come to know and love category theory. I was lucky to be able to attend some lectures of Tom Leinster's graduate course in Edinburgh, where I gained a firm understanding of the core concepts, and the online availability of the nLab collaborative project helped me form an overview of the subject with ease.\newline

I believe the language of category theory is the most effective way to connect and compare across different branches of mathematics so I have decided to use its notions and terminology throughout this thesis. However, my choice is not only one of convenience, but one of necessity, since the concepts of ``naturality'' or ``canonical construction'', which I think are essential in any serious foundational revision of physics, are best captured in the context of category theory.\newline

For this reason, many of the mathematical statements that will be used to make precise the questions introduced in Section \ref{OfPhaseAndUnits} above will be phrased in terms of functors and similarities between categories. The use of category theory has the added advantage that results proved within some set of assumptions are often easy to promote to precise conjectures in more general settings.

\section{\emph{Quo Vadimus}?} \label{QuoImus}

From the general motivating questions introduced in Section \ref{OfPhaseAndUnits} we extract here a more concrete list of goals to cover throughout the thesis:

\begin{enumerate}
    \item \textbf{\emph{Perspicuitati}}. Give a clear and concise account of the background material necessary to support the more conceptual discussion; aiming, whenever possible, for self-contained presentations that are accessible to graduates with a general knowledge of differential geometry. Succeeding in doing so would mean that this thesis could be taken as a short primer on the several subjects that are introduced as mathematical preliminaries: local Lie algebras, Lie algebroids, Poisson and symplectic manifolds, line bundles and Jacobi and contact manifolds.
    \item \textbf{\emph{Fundamenta}}. Conduct a detailed survey of the basic mathematical and physical ideas underlying the modern concept of phase space, as well as its historical origins, in particular in the context of classical Hamiltonian mechanics.
    \item \textbf{\emph{Exploratio}}. Having identified the defining mathematical features of phase space, proceed to systematically explore the generalizations that result from relaxing assumptions and discuss their physical significance.
    \item \textbf{\emph{Unitates}}. Revise the history and foundations of the notion of physical quantity and unit of measurement, as encountered in conventional dimensional analysis, and provide a precise mathematical characterization of their essential algebraic attributes.  
    \item \textbf{\emph{Mensura et Mechanica}}. Give a precise mathematical characterization to the term `canonical formalism' used loosely to refer to the theory of Hamiltonian mechanics on cotangent bundles of configuration spaces and argue whether replacing conventional configuration spaces by line bundles provides a similar framework.
    \item \textbf{\emph{Mensura et Mathematica}}. Explore the possibility of redefining the basic notions of algebra and geometry to accommodate for the algebraic structure of physical quantities found in \emph{Unitates}. Compare the notion of manifold and Poisson algebra that result from this exploration with the line bundle generalization of Poisson manifold found in \emph{Mensura et Mechanica}.
\end{enumerate}

\section{Structure of the Thesis, Notes on Style and Bibliographical Guide} \label{StructureOfThesis}

This thesis is comprised of four largely independent bodies of text:

\begin{enumerate}
    \item \textbf{Mathematical Preliminaries}. Chapter \ref{MathematicalPreliminaries} contains a mathematical exposition of a series of topics that has been made as self-contained as possible only assuming a firm background in graduate-level differential geometry and algebra. Due to our emphasis on providing a coherent narrative, we note that, contrary to what is customary in mathematical papers, our chapter \ref{MathematicalPreliminaries} on mathematical preliminaries contains a few original results. In particular, we provide a new, more focused, presentation of line bundles and Jacobi manifolds from the perspective of unit-free manifolds.
    \item \textbf{Foundations and Generalizations of Phase Spaces}. Sections \ref{BriefHistoryOfPhaseSpace} and \ref{FoundationsHamiltonianPhaseSpace} contain a historical account and a conceptual revision of the mathematical foundations of the modern notion of phase space. Chapters \ref{SymplecticPhaseSpaces} and \ref{ContactPhaseSpaces} present an exploration of several generalizations of phase spaces that are then compiled in chapter \ref{LandscapePhaseSpaces} where we present them in a pictorial representation that we call the Landscape of Phase Spaces.
    \item \textbf{Foundations and Axiomatization of Metrology}. Section \ref{BriefHistoryOfMetrology} contains a historical account of metrology and dimensional analysis. Section \ref{MeasurandFormalism} contains a proposal for an axiomatization of the usual structure of physical quantities in dimensioned analysis within the category of lines, first conceptually and physically motivated in that section, and then mathematically justified in Section \ref{MeasurandSpacesRevisited} by our original results on dimensioned rings and the potential functor presented in sections \ref{DimRingsModules}, \ref{DimAlgebras} and \ref{PotentialFunctor}.
    \item \textbf{Incorporating Physical Dimension into Geometry and Mechanics}. In sections \ref{UnitFreeRevisited}, \ref{JacobiManifoldsRevisited}, \ref{DimManifolds} and \ref{DimMechanics} results from all the previous sections are combined to argue that the notion of unit-free manifold, introduced in chapter \ref{MathematicalPreliminaries} via line bundles over smooth manifolds, proves effective in recovering a Hamiltonian theory of mechanics that is interpreted as to incorporate physical dimension naturally.
\end{enumerate}
Chapter \ref{SpeculationsDynIntQuant} contains some speculations and comments about lines of future research that open up from the results of this thesis. Chapter \ref{Conclusion} contains a few summarizing remarks.\newline

In the interest of keeping the discussion focused throughout the corpus of the thesis we have made the following stylistic choices:
\begin{itemize}
    \item[-] Most mathematical discussion will be conducted in a flowing, discursive manner, highlighting in \textbf{boldface typography} relevant newly defined terms and that will be mentioned multiple times throughout the thesis.
    \item[-] All the supporting references for each section are compiled in the bibliographical guide below and are largely omitted throughout the thesis.
    \item[-] Propositions, Theorems and Conjectures are occasionally highlighted to easily refer to them and to provide concise proofs.
    \item[-] If a proof is presented, it means that the claim is an original result of this thesis or that our original definitions and approach provide a new (often simpler) proof that further adds to our narrative or that illustrates a technical point.
\end{itemize}

What follows is a bibliographical guide to some of the supporting references for this thesis.

\subsection*{Chapter \ref{MathematicalPreliminaries}}

 General references: Category theory \cite{leinster2014basic,mac2013categories}, differential geometry \cite{lee2009manifolds} and commutative algebra \cite{kasch1982modules}.

\subsection*{Section \ref{Lines}}

J. Baez and J. Dolan introduced the notion of line objects and dimensional categories in \cite{dolan2009doctrines}. Our definition of the category of lines is a prime example of a dimensional category. They also noted that the structure physical quantities appears naturally in these categories although they did not study the possible implications this may have for the foundations of metrology, let alone any applications to mechanics.\newline

The work of J. Janyvska, M. Modgno and R. Vitolo \cite{janyvska2007semi,janyvska2010algebraic} has explored the the idea of finding a mathematical foundation for the theory of physical dimension in quite some detail. Our approach is compatible with theirs, although we work over vector spaces whereas they introduce the more refined notion of semi-vector spaces, and can be understood as a complementary development with a focus on classical mechanics and the formal generalizations of the algebraic properties of physical observables.

\subsection*{Section \ref{VectorBundles}}

General reference for vector bundles, modules of sections and differential operators \cite{nestruev2006smooth}.

\subsection*{Section \ref{LocalLie}}

Local Lie algebras were first studied by A. Kirillov in \cite{kirillov1976local} were a proof of the symbol part of proposition \ref{SymbolLieAlgebra} can be found. Propositions \ref{LocalLieRk1} and \ref{LocalLieRk2} are also proved indirectly in Kirillov's paper and they appear scattered within larger statements in the literature thereafter.\newline

General reference in Lie algebroids and Lie groupoids \cite{mackenzie2005general}. Integrability results of Lie algebroids \cite{crainic2003integrability}.

\subsection*{Section \ref{SymplecticGeometry}}

General reference in Poisson geometry \cite{fernandes2014lectures}. A first proof of proposition \ref{LieAlgebroidLinearPoisson} was given in \cite{vaintrob1997homological}. Presymplectic reduction \cite{echeverria1999reduction}. The Weinstein symplectic category is thoroughly discussed in \cite{weinstein2009symplectic}. Symplectic reduction and BRST cohomology \cite{figueroa1989brst}.

\subsection*{Section \ref{DiracGeometry}}

An effort has been made to present the linear theory of Dirac structures in detail, which appears only scattered in the literature. Introduction to Dirac manifolds \cite{bursztyn2013brief}, propositions \ref{IsotropicSplittingLinearCourant}, \ref{LinearDirac2Form} and \ref{LinearCourantComposition} can be found there. Proposition \ref{IsotropicSplittingLinearCourant} can be found in \cite{hu2009extended}. Proposition \ref{LinearDiracMorphism} can be found in \cite{alekseev2007pure}.\newline

Our presentation of Courant algebras is original but we were motivated by their introduction in \cite{bursztyn2007reduction}.\newline

T. Courant's thesis where Dirac manifolds were first defined \cite{courant1990dirac}. General reference on Courant algebroids and Lie bialgebroids: first two chapters of \cite{gualtieri2004generalized}. Courant morphisms \cite{bursztyn2008courant}. Proposition \ref{SeveraClass} was first found by Severa in communication with Weinstein, a first publication where this was used is \cite{vsevera2001poisson}, our proof focuses on the Leibniz cohomology perspective of the Severa class.\newline

Propositions \ref{AdjointMapCourant}, \ref{ExtendedActionsCourantAlgebroids}, \ref{EquivalenceExtendedActionsCourant} are all found in \cite{bursztyn2007reduction} but in considerable less detail and without fully exploiting the language of Leibniz algebras. This is the central reference on which we base our discussion of Courant and Dirac reduction.

\subsection*{Section \ref{LineBundles}}

Our approach to line bundles is strongly inspired by the work of L. Vitagliano and collaborators. Although our presentation on line bundles and related topics is largely original some important ideas, such as the product of line bundles of proposition \ref{LineProduct}, were obtained in private communications with J. Schnitzer and will appear soon in \cite{schnitzer2004thesis}.\newline

Proposition \ref{DerFunctor} can be found in \cite{tortorella2017deformations}, propositions \ref{DerBundleLineSubmanifold} and \ref{DerBundleGroupAction} are partially proved in that reference too.

\subsection*{Section \ref{ContactGeometry}}

General Jacobi geometry from a line bundle perspective in the first chapter of \cite{tortorella2017deformations}. Our main original addition to the exposition there is proposition \ref{ExtensionBySymbol}. Conditions 2, 3 and 4 in proposition \ref{CoisotropicSubmanifoldsJacobi} are proved to be equivalent in \cite{tortorella2017deformations} and our coisotropic reduction theorem for Jacobi manifolds is a reformulation of the similar result also in \cite{tortorella2017deformations}.\newline

The propositions mentioned already together with our results on the equivalence between coisotropic relations and Jacobi maps find their corresponding versions for trivial line bundles in \cite{ibanez1997coisotropic} and \cite{ibort1997reduction}.\newline

Jacobi geometry has seen a rise in popularity in recent years, both as a mathematical subject -  for the most part of the $1900$s, contact (Jacobi) geometry has been the neglected odd-dimensional sister of symplectic (Poisson) geometry - and for its applications to mathematical physics and other sciences. Just a few examples of recent developments are: the integrability of Jacobi manifolds to contact groupoids \cite{crainic2015jacobi}, the generalization of Dirac geometry for Jacobi manifolds \cite{vitagliano2015dirac,schnitzer2019normal}, a dissipative version of Liouville's theorem in contact manifolds \cite{bravetti2015liouville}, the identification of a contact structure in thermodynamics \cite{mrugala1991contact,grmela2014contact} or even applications to neuroscience \cite{petitot2017elements}. Most prominently for geometric mechanics, the work of M. de Le\'on \emph{et al.} \cite{deLeon2011methods,deLeon2017cosymplectic,valcazar2018contact} has shown how contact geometry is the natural framework for the dynamical formulation of mechanical systems subject to time-dependent forces and dissipative effects. Our formalism of canonical Hamiltonian contact mechanics presents yet another avenue for the application of Jacobi geometry to physics, this time tackling the issue of a systematic mathematical treatment of physical dimension and units in mechanics.

\subsection*{Section \ref{LDiracGeometry}}

We employ the terms L-Courant and L-Dirac to stay consistent with our naming conventions but our definitions are direct generalizations of the double der bundle of a line bundle and Dirac-Jacobi structures of Vitagliano and Wade \cite{vitagliano2015dirac}.

\subsection*{Sections \ref{BriefHistoryOfMetrology}, \ref{MeasurandFormalism}}

Human evolution \cite{ambrose2001paleolithic}. History of mathematics \cite{boyer2011history}. Units of measurement in ancient Egypt \cite{blanco2017atlas}. General history \cite{overy2010complete}. History of mechanics \cite{nolte2018galileo}. History of science and technology in the $1700$s \cite{wolf2019history}. Fourier's book containing the earliest version of dimensional analysis \cite{baron1822theorie}. A practitioner's account of the history of dimensional analysis \cite{macagno1971historico}. One of the papers that helped popularize the modern form of dimensional analysis \cite{buckingham1915principle}. A modern standard reference in dimensional analysis \cite{barenblatt1996scaling}. International Vocabulary of Metrology \cite{bimp2012metrology}.

\subsection*{Section \ref{BriefHistoryOfPhaseSpace}}

Overarching history of dynamics since Galileo \cite{nolte2018galileo}. The French enlightenment \cite{goodman1996republic}. A practitioner's detailed account of the papers by Lagrange and Poisson of $1809$ \cite{marle2009inception}. The images of figures \ref{duChatelet}, \ref{LagrangeSO(3)} and \ref{PoissonBracc} were obtained from \emph{Gallica}, the online archive of the Bibliotheque Nationale de France found at \url{https://gallica.bnf.fr/}. Hamilton's original paper from $1835$ \cite{hamilton1835vii}. Boltzmann's masterpiece \cite{boltzmann2012lectures}. An account of Lie's life and work \cite{fritzsche1999sophus}. Einstein's paper on gravitation \cite{einstein1915feldgleichungen}. An extensive overview of the development of geometric mechanics can be found in the historical notes of \cite{marsden2013introduction}. Paper where Lichnerowicz's isolation of the axioms of Poisson manifold are already apparent \cite{flato1976deformations}. Weinstein's account on the beginning of Lie groupoid theory can be found in \cite{weinstein1996groupoids}. Dirac's book containing his theory of constrained Hamiltonian systems \cite{dirac1966lectures}. One of Arnold's contributions to geometric mechanics \cite{arnold2013mathematical}. One of the centerpieces of the connection between string theory and non-commutative geometry \cite{seiberg1999string}. A practitioner's account of the recent history of Dirac geometry \cite{kosmann2013courant}. Generalizations of Lie algebras and in the number of arguments are extensively discussed here \cite{azcarraga2010nary}, an alternative route for the generalization of phase spaces that we do not explore in this thesis.

\subsection*{Section \ref{FoundationsHamiltonianPhaseSpace}}

Although no particular source was used in writing this section its contents have been strongly shaped by my personal communications with F. Zalamea whose thesis \cite{zalamea2016chasing} presents great contemporary account on the parallels and discrepancies between classical and quantum kinematics. In particular, part of the motivation to formulate the notion of theory of phase spaces as a commutative diagram of functors between categories of observables, phase spaces and dynamics came from conversations held during the $2017$ Rethinking Workshop.

\subsection*{Section \ref{CanonicalSymplectic}}

Canonical symplectic mechanics \cite{marsden1992lectures,abraham1978foundations}.

\subsection*{Section \ref{GeneralizationSymplectic} and \ref{PoissonizationFunctors}}

Propositions \ref{TwistedPreSymplecticDirac} and \ref{TwistedPoissonDirac} can be extracted from the proofs found in \cite{gualtieri2004generalized} but the explicit results are hard to pin-point in the literature. An overview of Poisson manifolds in relation to Lie algebroids and Lie groupoids, with a discussion on the Lie bracket on 1-forms is found in \cite{weinstein1998poisson}.

\subsection*{Section \ref{CanonicalContact}}

Proposition \ref{JetCanonicalContactStructure} is a standard result in the theory of jet bundles and a version of it for jets of higher orders can be found in \cite{saunders1989geometry}. Some of our theorems in this section can be found in \cite{ibort1997reduction} in the case of trivial line bundles. Our main contribution is to prove all of them in full generality, making extensive use of our formalism of unit-free manifolds.

\subsection*{Chapter \ref{DimAlgebraGeometry}}

Our notion of dimensioned ring was inspired by the typed families of fields of G.W. Hart introduced in \cite{hart2012multidimensional}. In fact, we recover the definition of typed family of fields precisely as a trivial dimensioned field and Hart's claim that ``all typed families of fields are isomorphic to the Cartesian product of a field and a group'' is a direct result of our proposition \ref{UnitsDimFields}.

\subsection*{Section \ref{CommentsObsSmoothDyn}}

A proof for the dimensionality of the space of $\text{C}^k$-differentiable derivations can be found in \cite{abraham2012manifolds}. The existence of compatible smooth atlases for $\text{C}^1$-differentiable manifolds is proved in \cite{munkres2016elementary}.

\subsection*{Section \ref{RecoveringHamiltonianDyn}}

A general Darboux theorem for precontact manifolds can be found in \cite{tortorella2017deformations}. Explicit coordinate expressions for the Hamiltonian dynamics on contact manifolds can be found in \cite{valcazar2018contact}.

\subsection*{Sections \ref{QuantumPhaseSpace}, \ref{QuantumUnitFree} and \ref{QuantizationIntegration}}

Von Newmann's original axiomatic presentation of quantum mechanics on Hilbert spaces \cite{von2018mathematical}. The convolution algebras of Lie groupoids as well as the C*-algebra approach to quantum mechanics can be found in \cite{landsman2012mathematical}, this an excellent compendium of the contemporary topics revolving classical and quantum phase spaces. The general theory of integration of Lie algebroids can be found in \cite{crainic2011lectures}. For the integration of Poisson manifolds to symplectic groupoids refer to \cite{crainic2004poiss}. For the integration of Dirac structures to presymplectic groupoids see \cite{bursztyn2004twdirac}. For the results on the integration of Jacobi manifolds by contact groupoids see \cite{crainic2015jacobi}. In \cite{landsman2006liephys} N.P. Landsman presents the connections Lie algebroids, Lie groupoids and non-commutative geometry and introduces some of the technical details of the potential quantization scheme enabled by the integration theory of Lie algebroids.

\chapter{Mathematical Preliminaries} \label{MathematicalPreliminaries}

Throughout this thesis, vector spaces and manifolds will be taken over the real numbers unless otherwise stated. We adopt the standard conventions of finite-dimensional smooth differential geometry; manifolds and maps between them will be assumed to be $\text{C}^\infty$-smooth unless otherwise stated. The category of real vector spaces will be denoted by $\Vect$ and the category of smooth manifolds will be denoted by $\Man$. We follow the standard conventions of commutative algebra for the occasional treatment of general rings, modules and algebras. We will also use the standard terminology of category theory.

\section{Lines} \label{Lines}

\subsection{The Category of Lines} \label{CategoryOfLines}

We identify \textbf{the category of lines}, $\Line$, as a subcategory of vector spaces $\Vect$. Objects are vector spaces over the field of real numbers $\mathbb{R}$ of dimension 1, a useful way to think of these in the context of the present work is as sets of numbers without the choice of a unit. An object $L\in\Line$ will be appropriately called a \textbf{line}. A morphism in this category $b\in\text{Hom}_{\Line}(L,L')$, usually simply denoted by $b:L\to L'$, is an invertible (equivalently non-zero) linear map. Composition in the category $\Line$ is simply the composition of maps. If we think of $L$ and $L'$ as numbers without a choice of a unit, a morphism $b$ between them can be thought of as a unit-free conversion factor, for this reason we will often refer to a morphism of lines as a \textbf{factor}. We consider the field of real numbers, trivially a line when regarded as real a vector space, as a singled out object in the category of lines $\Real\in \Line$.\newline 

Note that all the morphisms in this category are, by definition, isomorphisms, thus making $\Line$ into a groupoid; however, one should not think of all the objects in the category as being equivalent. As we shall see below, there are times when one finds factors between lines (invertible morphisms by definition) in a canonical way, that is, without making any further choices beyond the information that specifies the lines, these will be called \textbf{canonical factors}. In this manner, if given two lines $L,L'\in\Line$ there exists a canonical factor $b_{LL'}:L\to L'$ then we do regard $L$ and $L'$ as being equivalent. Our notation will reflect this fact by displaying a representative line in place of the equivalence class of lines for which there exist canonical factors mapping them to the representative and by using the equals sign ``='' between any of them. \newline

It is a simple linear algebra fact that any two lines $L,L'\in\Line$ satisfy
\begin{equation*}
    \dimm (L \oplus L')= \dimm L + \dimm L' = 2 > 1, \qquad \dimm L^* = \dimm L = 1, \qquad \dimm (L \otimes L') = 1.
\end{equation*}
Then, we note that the direct sum $\oplus$, which acts as a categorical product and coproduct in $\Vect$, is no longer defined in $\Line$, however, it is straightforward to check that $(\Line, \otimes, \mathbb{R})$ forms a symmetric monoidal category and that $*:\Line \to \Line$ is a duality contravariant autofunctor. The usual linear-algebraic definition of the tensor product via the universal property in $\Line$ implies, in particular, the facts listed below, which always hold in any symmetric monoidal category but that we mention here for convenience: 
\begin{itemize}
    \item For any two lines $L_1,L_2\in\Line$ there exists a canonical factor swapping their tensor product
    \begin{equation*}
        b_{12}:L_1\otimes L_2 \to L_2\otimes L_1.
    \end{equation*}
    \item The tensor unit $1_\otimes$ is precisely $\mathbb{R}$, this means that for any line $L\in\Line$ there is a canonical factor
    \begin{equation*}
        t_L:\Real\otimes L \to L.
    \end{equation*}
    \item For any three lines $L_1,L_2,L_3\in\Line$ there exists a canonical factor associating them
    \begin{equation*}
        a_{123}:L_1\otimes (L_2 \otimes L_3) \to (L_1\otimes L_2)\otimes L_3,
    \end{equation*}
    thus we can omit the brackets.
    \item More generally, given any collection of lines $L_1,\dots,L_k\in\Line$ and a permutation of $k$ elements $\sigma\in S_k$, there exists a canonical factor
    \begin{equation*}
        b_\sigma:L_1\otimes \cdots \otimes L_k \to L_{\sigma(1)}\otimes \cdots \otimes L_{\sigma(k)}.
    \end{equation*}
    \item Any two factors $b:L_1\to L_2$ and $b':L_1'\to L_2'$ have tensor product
    \begin{equation*}
        b\otimes b' :L_1\otimes L_1'\to L_2\otimes L_2'.
    \end{equation*}
    which is again a factor.
\end{itemize}
The autofunctor $*$ is in fact given by $\Real$-adjunction in the category of vector spaces $L^*:=\text{Hom}_{\Vect}(L,\Real)$. It is then easy to show that for any lines $L,L'\in\Line$ there are canonical factors $d_L:(L^*)^*\to L$ and $e_{LL'}:(L\otimes L)^*\to L^*\otimes L'^*$. Similarly, one can show that $L^*\otimes L'$ is canonically isomorphic to $\text{Hom}_{\Vect}(L,L')$ in the category of vector spaces and is indeed of dimension 1. In particular, $\text{Hom}_{\Vect}(L,L)\in\Line$ has a distinguished non-zero element, the identity id$_L$, thus we find that it is canonically isomorphic to $\Real$ as lines. Therefore, for any line $L\in\Line$ we find a canonical factor
\begin{equation*}
    p_L:L^*\otimes L\to \Real.
\end{equation*}
This last result, under the intuition of lines as numbers without a choice of unit, allows us to reinterpret the singled out line $\Real$ informally as the set of procedures common to all lines by which a number gives any other number in a linear way (preserving ratios). This interpretation somewhat justifies the following adjustment in terminology: we will refer to the tensor unit $\Real\in\Line$ as \textbf{the patron line}. In anticipation of our discussion about units of measurement and physical quantities of Section \ref{MeasurandFormalism}, the term \textbf{unit} is reserved for non-vanishing elements of a line $u\in L^\times$, where we have denoted $L^{\times}:= L\diagdown \{0\}$. \newline

We can see the role of the patron line $\Real$ more explicitly by considering the following map for any line $L\in\Line$:
\begin{align*}
\lambda: L\times L^{\times} & \to L^*\otimes L\\
(a,b) & \mapsto \lambda_{ab} \text{ such that } a=\lambda_{ab}(b).
\end{align*}
And, reciprocally, also define $\rho:L^{\times}\times L\to L^*\otimes L$. The following proposition gives mathematical foundation to the intuition of $L$ being a ``unit-free'' field of numbers and $\Real$ being the patron of ratios between ``unit-free'' numbers.
\begin{prop}[Ratio Maps]\label{RatMaps}
The maps $\lambda$ and $\rho$ are well-defined and linear in their vector arguments; furthermore, for any $a,b,c\in L^{\times}$, we have the following identities for the maps $l:=p_L\circ \lambda$ and $r:=p_L\circ \rho$
\begin{equation}\label{ratio1}
    l_{ab}\cdot l_{bc}\cdot l_{ca}=1 \qquad l_{ab}\cdot r_{ab}=1 \qquad r_{ab}\cdot r_{bc}\cdot r_{ca}=1.
\end{equation}
The maps $l$ and $r$ are called the \textbf{ratio maps} and the first and third equations will be called the \textbf{2-out-of-3 identity} for the maps $l$ and $r$ respectively.
\end{prop}
\begin{proof}
Since both $L$ and $\text{Hom}_\Vect(L,L)$ are 1-dimensional real vector spaces, it follows that there is a unique linear isomorphism mapping any two non-zero elements in $L$ and the zero element is only mapped by the zero linear map. This can be seen more explicitly with the use of the canonical factor $p_L:\Hom{L,L}\to \Real$, which allows us to define the map $l:=p_L\circ \lambda$ that when applied to two non-zero elements $a,b\in L$ gives the unique non-zero real number $l_{ab}$ acting as proportionality factor: $a=l_{ab}\cdot b$. To prove the 2-out-of-3 identity consider $a=\lambda_{ab}(b)$, $b=\lambda_{bc}(c)$ and $c=\lambda_{ca}(a)$. Combining the three equations we find $\lambda_{bc}\circ\lambda_{ca}(a)=\lambda_{ba}(a)$. Noting that $\lambda_{ab}=\lambda^{-1}_{ba}$ gives the desired result. Similarly for $\rho$ and the reciprocal identity follows by construction.
\end{proof}
Consider now two lines $L_1,L_2\in\Line$ with their corresponding maps $l^1$, $r^1$ and $l^2$, $r^2$. It follows from the definitions above that for any factor $B:L_1\to L_2$ we have $l^1_{ab}=l^2_{B(a)B(b)}$ and $r^1_{ab}=r^2_{B(a)B(b)}$. More generally, we can define functions on the space of factors between the lines $L_1$ and $L_2$ from pairs of line elements:
\begin{align}\label{ratio2}
    a_1\in L_1, b_2\in L_2^{\times} &\quad \mapsto \quad l^{12}_{a_1b_2}(B):=l^2_{B(a_1)b_2}\\
    b_1\in L_1^{\times}, a_2\in L_2 &\quad \mapsto \quad r^{12}_{b_1a_2}(B):=r^1_{b_1B^{-1}(a_2)}
\end{align}
for all $B:L_1 \to L_2$ a factor. It follows then by construction, that for any pair of non-zero line elements $b_1\in L_1^{\times}$ and $b_2\in L_2^{\times}$ the following identity holds
\begin{equation*}
    l^{12}_{b_1b_2}\cdot r^{12}_{b_1b_2} =1
\end{equation*}
as functions over the space of factors $\text{Hom}_\Line(L_1,L_2)$. Whenever there is no room for confusion, we will employ the following abuse of fraction notation for non-zero elements:
\begin{equation*}
    l_{ab}=\frac{a}{b}=\frac{1}{r_{ab}} \qquad \qquad l^{12}_{a_1b_2}=\frac{a_1}{b_2}=\frac{1}{r^{12}_{b_1a_2}}
\end{equation*}
which notationally reflect the identities (\ref{ratio1}).\newline

If we write the tensor product of $n$ copies of a line and $m$ copies of its dual as $\otimes^m L^* \otimes^n L$, composing the canonical factors $t_L$, $b_\sigma$ and $p_L$ iteratively we obtain the following canonical factor:
\begin{equation*}
    a_{nm}:\otimes^m L^*  \otimes^n L \to \begin{cases} 
      \otimes^{|n-m|}L & n>m \\
      \Real & n=m \\
      \otimes^{|n-m|}L^* & n<m 
   \end{cases}
\end{equation*}
This result then encourages the introduction of the following notation:
\begin{equation*}
    \begin{cases} 
      L^n:=\otimes^nL & n>0 \\
      L^n:= \Real & n=0 \\
      L^n:=\otimes^nL^* & n<0 
   \end{cases}
\end{equation*}
which is such that given two integers $n,m\in\Int$ and any line $L\in\Line$ the following equations hold
\begin{equation*}
    (L^n)^* = L^{-n} \qquad L^{n}\otimes L^{m} = L^{n+m}.
\end{equation*}
Thus we see how one single line and its dual $L,L^*\in\Line$ generate an abelian group with the tensor product as group multiplication, the patron $\Real\in\Line$ as group identity and the duality autofunctor as inversion. We call this group the \textbf{potential} of $L\in\Line$ and denote it by $L^\odot$. It trivially follows from the results above that $(\Real^\odot,\otimes) \cong 0$ and for $L\neq \Real$,  $(L^\odot,\otimes) \cong (\Int,+)$ as abelian groups.\newline

If we now consider a finite collection of lines $\{L_i\in\Line\}_{i=1}^k$ we define its \textbf{potential} by taking the union all the possible tensor powers:
\begin{equation*}
    (L_1\dots L_k)^\odot := \bigcup_{n_1,\dots,n_k\in\Int} L_1^{n_1}\otimes \cdots \otimes L_k^{n_k}.
\end{equation*}
An iterative use of the canonical factors $b_\sigma$, $e_{LL'}$ and $p_L$ allows us to write the following equations
\begin{align*}
    (L_1^{n_1}\otimes \cdots \otimes L_k^{n_k})^* &=L_1^{-n_1}\otimes \cdots \otimes L_k^{-n_k}\\
    (L_1^{n_1}\otimes \cdots \otimes L_k^{n_k}) \otimes (L_1^{m_1}\otimes \cdots \otimes L_k^{m_k})  &= L_1^{n_1+m_1}\otimes \cdots \otimes L_k^{n_k+m_k}
\end{align*}
and thus we see that the potential of a collection of lines forms an abelian group with the tensor product as multiplication, duality as inversion and the patron as group identity. The components of a potential $(L_1\dots L_k)^\otimes$ are called the \textbf{base lines} and we have
\begin{equation*}
    (L_1\dots L_k)^\odot\cong L_1^\otimes \oplus \dots \oplus L_k^\otimes \cong \Int^k
\end{equation*}
as abelian groups. Furthermore, the potential $(L_1\dots L_k)^\otimes$ has the structure of a $\Int$-module with the action given by:
\begin{equation*}
    m\cdot (L_1^{n_1}\otimes \cdots \otimes L_k^{n_k}):=\otimes^m(L_1^{n_1})\otimes \cdots \otimes \otimes^m(L_k^{n_k})=L_1^{m\cdot n_1}\otimes \cdots \otimes L_k^{m\cdot n_k}.
\end{equation*}

\subsection{The Category of L-Vector Spaces} \label{CategoryOfLVectorSpaces}

The identification of the category $\Line$ allows for a natural generalization of the category of real vector spaces where one replaces the field of scalars with a line, thought of informally as a ``unit-free'' field of scalars. The category of \textbf{line vector spaces} or \textbf{L-vector spaces} is defined formally as the product category
\begin{equation*}
    \LVect:= \Vect \times \Line.
\end{equation*}
Our notation for objects in this category will be $V^L:=(V,L)$ with $V\in\Vect$, $L\in\Line$, and similarly for morphisms $\psi^b:V^L\to W^{L'}$ with $\psi\in\text{Hom}_\Vect(V,W)$, $b\in\text{Hom}_\Line(L,L')$. Objects $V^L$ will be called L-vector spaces and morphisms $\psi^b$ will be called linear factors. The tensor structures of both $\Vect$ and $\Line$ induce a natural \textbf{tensor product of L-vector spaces} given by
\begin{align*}
    V^L\otimes W^{L'}:&=(V\otimes W)^{L\otimes L'}\\
    \psi^b\otimes \phi^{b'}:&=(\psi \otimes \phi)^{b\otimes b'}\\
    1_{\otimes\LVect}&=\Real^\Real,
\end{align*}
making $\textsf{LVect}$ into a symmetric monoidal category. We define a the \textbf{L-dual} of an L-vector space as
\begin{equation*}
    *(V^L):=(V^*\otimes L)^L.
\end{equation*}
We will often use the more compact notation $(V^*\otimes L)^L:=V^{*L}$. This construction is crucially idempotent
\begin{equation*}
    *(*(V^{L}))=(V^{**}\otimes L\otimes L^*)^L \cong (V^{**}\otimes \Real)^L \cong (V^{**})^L \cong V^L,
\end{equation*}
which justifies the chosen name of L-duality, in analogy with conventional duality of vector spaces. Given a linear factor $\psi^b:V^L\to W^{L'}$ we can define its L-dual as
\begin{align*}
\psi^{L'*L}: W^{*L'} & \to V^{*L}\\
\beta^{l'} & \mapsto \alpha^l
\end{align*}
with
\begin{equation*}
    \alpha = b^{-1}\circ \beta \circ \psi \qquad l=b^{-1}(l'),
\end{equation*}
where we regard L-dual vectors as linear maps $\alpha\in V^*\otimes L\cong \text{Hom}_\Vect(V,L)$. This construction clearly makes $*$ into a duality contravariant autofunctor. The \textbf{direct sum of L-vector spaces} is only defined for L-vector spaces sharing the same line:
\begin{equation*}
    V^L\oplus W^L := (V\oplus W)^L.
\end{equation*}

Subobjects in this category will be identified with linear factors whose vector component is injective, in particular, given any vector subspace $i:U\hookrightarrow V$ we have an inclusion linear factor given by $i^{\Id_L}:U^L\to V^L$. Quotients can then be taken in the obvious way $V^L/U:=(V/U)^L$. There is a natural notion of \textbf{L-annihilator} of a subspace
\begin{equation*}
    U^{0L}:=\{\alpha\in V^{*L} | \quad \alpha(u)=0\in L, \forall u\in U\}\cong U^0\otimes L.
\end{equation*}
It is then a direct result that there exists a canonical linear factor isomorphism
\begin{equation}\label{Lannquo}
    a^{\Id_L}:U^{*L}\to V^{*L}/U^{0L},
\end{equation}
giving the characterization of the L-dual of a subspace of a L-vector space, in direct analogy with the standard result for ordinary vector spaces.

\section{Vector Bundles} \label{VectorBundles}

\subsection{The Category of Vector Bundles} \label{CategoryVectorBundles}

Let $\epsilon: E\to M$ be a vector bundle over a smooth manifold $M$; the notations $E_M$, or simply $E$ at times, will be often used interchangeably. The typical fibre of a vector bundle $E_M$ is a vector space $V\in\Vect$ and its dimension is called the \textbf{rank} of the vector bundle, $\rkk{E_M}:=\dimm V=\dimm E_x$ for all $x\in M$. \textbf{Local trivializations} are maps of the form $T_U:E|_U\to U\times V$ for some open neighbourhood $U\subset M$ which overlap to give $\GL{V}$-valued transition functions. A local trivialization is equivalently given by a set of local sections $\{e_i:U\to E|_U, \quad \epsilon \circ e_i=\Id_U\}_{i=1}^{\rkk{E}}$ such that their values give a basis at every point $x\in U$. If such a basis exists globally we say that the vector bundle $E_M$ is \textbf{trivializable}, and the vector bundle structure is (non-canonically) isomorphic to that of the \textbf{trivial} vector bundle $\Proj_1:M\times V\to M$. A trivial vector bundle over the manifold $M$ with typical fibre $V$ will be usually denoted by $V_M$. A smooth map between two vector bundles $F: E_1 \to E_2$ is called a \textbf{vector bundle morphism} if there exists a smooth map between the bases $\varphi: M_1\to M_2$ such that the diagram:
\begin{equation*}
\begin{tikzcd}
E_1 \arrow[r, "F"] \arrow[d, "\epsilon_1"'] & E_2 \arrow[d, "\epsilon_2"] \\
M_1 \arrow[r, "\varphi"'] & M_2
\end{tikzcd}
\end{equation*}
commutes, which is referred to as $F$ covering $\varphi$, and $F$ restricts to a linear maps along the fibres. Equivalently, a vector bundle morphism is given by a smooth map between the bases $\varphi:M_1\to M_2$ together with an assignment of linear maps between the fibres:
\begin{equation*}
\{F_x: (E_1)_x \to (E_2)_{\varphi(x)} , x\in M_1\}.
\end{equation*}
Vector bundles over smooth manifolds with vector bundle morphisms form \textbf{the category of vector bundles} and we denote it by $\Vect_\Man$.\newline

If we fix a base manifold $M$, point-wise constructions allow to regard the subcategory of vector bundles over the same manifold $\Vect_M$ with morphisms covering the identity as an abelian symmetric monoidal category with duality. Under these conditions, we essentially recover the same categorical structure present in the category of real vector spaces $\Vect$. The \textbf{direct sum} in this category is given by the point-wise Whitney sum $A_M\oplus B_M$, the \textbf{tensor product} is similarly constructed with the fibre-wise tensor product $A_M\otimes B_M$ and hom sets, which are themselves naturally isomorphic to objects of the same category, are given by $\Hom{A_M,B_M}$, in particular the \textbf{duality functor} is $A^*_M:=\Hom{A_M,\Real_M}$. We will often abuse notation and use the same symbol for the bundle projection of any vector bundle canonically constructed by fibre-wise linear algebra from a given one $\epsilon: E\to M$ e.g. $\epsilon: E^*\to M$ or $\epsilon: E\otimes E \to M$.\newline

Let $\phi:N\to M$ be a smooth map and $\epsilon: E\to M$ a vector bundle, the \textbf{pull-back bundle} is defined as $\phi^*E_M:=N \ftimes{\phi}{\epsilon} E$ and will be sometimes also denoted by $E_N$. The category of vector bundles admits a categorical product: let two vector bundles $\epsilon_1: E_1\to M_1$ and $\epsilon_2: E_{M_2}\to M_2$ and let us denote the canonical projections of the Cartesian product of the bases:
\begin{equation*}
\begin{tikzcd}
M_1 & M_1\times M_2 \arrow[l,"\Proj_1"']\arrow[r,"\Proj_2"] & M_2.
\end{tikzcd}
\end{equation*}
then the \textbf{vector bundle product} is defined as
\begin{equation*}
    E_1\boxplus E_2 := \Proj_1^*E_1 \oplus \Proj_2^*E_2,
\end{equation*}
with the projection morphisms
\begin{equation*}
\begin{tikzcd}
E_1 & E_1 \boxplus E_2 \arrow[l,"P_1"']\arrow[r,"P_2"] & E_2
\end{tikzcd}
\end{equation*}
given fibre-wise by the projections of the vector direct sum in the obvious way. Although one is tempted to regard the above construction, together with the analogous one for tensor products of vector bundles with different bases, as the categorical operations of an abelian monoidal category with duality, several necessary conditions fail. For instance, the vector bundle product is not a coproduct, essentially since the Cartesian product of manifolds is not a coproduct generically, and the tensor product lacks a universal unit and even the partial units, trivial bundles with $\Real$ as typical fibre, are non-unique. Another crucial fact that will become relevant in sections below is the failure of the duality assignment $A_M \mapsto A^*_M$ to be a functor. As seen above, duality only becomes a functor in the restricted category $\Vect_M$, for a general vector bundle morphism $F: E_1 \to E_2$ covering a smooth map $\varphi: M_1\to M_2$ we have the following diagram
\begin{equation}\label{mappullback}
\begin{tikzcd}
E_1 \arrow[r, "F"] \arrow[d, "\epsilon_1"'] & \varphi^*E_2 \arrow[r, "\Id_{E_2}"] \arrow[d,"\Proj_2"] & E_2 \arrow[d, "\epsilon_2"] \\
M_1 \arrow[r, "\Id_{M_1}"'] & M_1 \arrow[r, "\varphi"'] & M_2
\end{tikzcd}
\end{equation}
and thus dual vector bundle maps can only be defined between pull-back bundles over the same base. These limitations will also become apparent, under a different guise, in our discussion about sections below.\newline

A \textbf{distribution} on a vector bundle $\epsilon: E\to M$ is a subset $D\subset E$ such that the subsets $D_x:= D \cap E_x$ are vector subspaces for all $x\in M$. If the rank of $D$ is independent of $x\in M$ we call $D$ a \textbf{regular distribution}. If for any $d_x\in D$ we can find a local section $s\in \Gamma(U)$, $x\in U$, such that $s(x)=d_x$ then we say that $D$ is a \textbf{smooth distribution}. Given a bundle map over the identity $F:A_M\to B_M$ the image $F(E)$ and the kernel $\ker{F}$ are smooth distributions that are not regular in general. A \textbf{subbundle} is a subset $D\subset E$ for which the restriction of the bundle projection makes $\epsilon|_D:D\to M$ into a vector bundle. These are the subobjects in the category of vector bundles. It is easy to see that $D$ being a subbundle is equivalent to $D$ being a smooth and regular distribution. The usual constructions of vector subspaces, such as spanning sums, intersections and quotients, apply to distributions equally by means of fibre-wise construction; notation analogous to the vector space case will be used for these.

\subsection{Modules of Sections} \label{Sections}

We define the space of \textbf{sections} of a vector bundle $\epsilon: E\to M$ as the set of maps that split the projection to the base:
\begin{equation*}
\Sec{E}:= \{s:M \to E |\quad  \epsilon \circ s = \Id_M\}.
\end{equation*}
The fibre-wise linear structure induces an addition operation and a multiplication by functions making $\Sec{E}$ into a $\Cin{M}$-module. Furthermore, using local trivializations we can see that modules of sections are locally free, generated by a number of basis elements equal to the rank of the vector bundle. It follows from our discussion about vector bundle morphisms in Section \ref{CategoryVectorBundles} that it is not possible to generically induce a well-defined map on sections from a given map between vector bundles. Thus we see that the section assignment $\Gamma:\Vect_\Man \to \textsf{RMod}$ fails to be a functor in general, in contrast with the similar contravariant functor for manifolds $\text{C}^\infty:\Man\to \Ring$. However, fixing a base manifold $M$ and considering the restricted category $\Vect_M$, the section assignment does become a functor:
\begin{equation*}
    \Gamma: \Vect_M\to \Cin{M}\textsf{Mod}.
\end{equation*}
This functor respects the abelian monoidal structure of $\Vect_M$ in the following sense: let two vector bundles $A,B\in\Vect_M$, using the fibre-wise linear structure one can prove the following list of canonical isomorphisms of $\Cin{M}$-modules:
\begin{align*}
    \Gamma(A\oplus B) &\cong \Sec{A} \oplus_{\Cin{M}} \Sec{B}\\
    \Sec{A^*} &\cong  \Sec{A}^*:=\text{Hom}_{\Cin{M}}(A,\Real_M)\\
    \Sec{A\otimes B} & \cong \Sec{A}\otimes_{\Cin{M}} \Sec{B}\\
    \Sec{\Hom{A,B}} &\cong \Sec{A^*\otimes B} \cong \text{Hom}_{\Vect_M}(A,B).
\end{align*}
Considering a distribution $D\subset A$, we denote by $\Sec{D}\subset \Sec{A}$ the set of all sections of $A$ taking values in the subspaces of the distribution at each fibre. When $D\subset A$ is a subbundle then $\Sec{D}$ is a $\Cin{M}$-submodule. It follows that spanning sums, intersections and quotients of subbundles induce the analogous constructions of $\Cin{M}$-modules under the section functor.\newline

The total space of a vector bundle $\epsilon:E\to M$ has two distinguished subsets of functions: the \textbf{fibre-wise constant functions} included by pull-back of the bundle projection $\epsilon^*:\Cin{M}\hookrightarrow \Cin{E}$ and the \textbf{fibre-wise linear functions} included by the fibre-wise action of sections of the dual vector bundle $l:\Sec{E^*}\hookrightarrow \Cin{E}$. The next proposition summarizes the nice properties of these special subsets of functions that makes them useful for our discussion about linear geometric structures in vector bundles in sections below.
\begin{prop}[Spanning Functions]\label{SpanningFunctions}
The fibre-wise constant and fibre-wise linear functions of a vector bundle $E_M$ form $\Real$-vector subspaces of $\normalfont\Cin{E}$ and the span
\begin{equation*}
    \normalfont \text{C}_s^\infty(E):=\epsilon^*\Cin{M}\oplus l_{\Sec{E^*}}
\end{equation*}
forms a $\normalfont \Cin{M}$-submodule of the ring of functions $\normalfont \Cin{E}$. Furthermore, the differentials of these functions span the cotangent bundle everywhere:
\begin{equation*}
    \normalfont d_e\text{C}_s^\infty(E)=\Cot_e E \quad \forall e\in E
\end{equation*}
and thus we call $\normalfont \text{C}_s^\infty(E)$ the \textbf{spanning functions} of the vector bundle $E$.
\end{prop}
\begin{proof}
First note that the inclusion maps $\epsilon^*$ and $l$ satisfy
\begin{equation*}
    \epsilon^*(f g)=\epsilon^*f\epsilon^*g \qquad \epsilon^*fl_\alpha=l_{f\cdot \alpha}
\end{equation*}
for all $f,g\in\Cin{M}$ and $\alpha\in\Sec{E^*}$, simply by construction using the fibre-wise linear structure of $E$. This clearly makes $\text{C}_s^\infty(E)$ into a $\Cin{M}$-module with the module multiplication being realized simply as point-wise multiplication on $E$. A quick coordinate computation using a local chart adapted to a trivialization shows that fibre-wise linear functions span cotangent spaces everywhere except for points on the zero section where they only span the vertical subspaces, however it is easily shown that the fibre-wise constants span the remaining horizontal subspaces at the zero section.
\end{proof}

Although, as we mentioned above, a general vector bundle morphism $F: E_1 \to E_2$ covering a smooth map $\varphi: M_1\to M_2$ does not induce a vector bundle morphism between the duals in a functorial way, several notions that capture duality can be defined on sections. Firstly consider diagram (\ref{mappullback}) and observe that we have the following induced maps for sections:
\begin{equation*}
\begin{tikzcd}
\Sec{E_1} \arrow[r, "F"] & \Sec{\varphi^*E_2} & \Sec{E_2} \arrow[l,"\varphi^*"'] \\
s_1 \arrow[r, mapsto] & F\circ s_1,s_2 \circ \varphi & s_2 \arrow[l, mapsto]
\end{tikzcd}
\end{equation*}
thus allowing to define a relation between sections on the intermediate object $\varphi^*E_2$. We say that sections $s_1$ and $s_2$ are $F$\textbf{-related} when their images under the induced maps above agree and we denote
\begin{equation*}
    s_1\sim_F s_2 \quad \Leftrightarrow \quad F  s_1 = \varphi^*s_2.
\end{equation*}
When $\varphi$ is a diffeomorphism this relation becomes a map
\begin{align*}
F_*: \Gamma(E_1) & \rightarrow \Gamma(E_2)\\
s& \mapsto F_*s := F\circ s \circ \varphi^{-1}
\end{align*}
that interacts with the module structures of sections in a natural way:
\begin{equation*}
F_*(s+r)=F_*s + F_*r \qquad F_*(f\cdot s)=(\varphi^{-1})^*f\cdot s
\end{equation*}
for all $s,r\in\Sec{E_1}$, $f\in\Cin{M_1}$. This map is called the \textbf{push-forward} of the vector bundle morphism $F$. When $F$ covers the identity map we will often use the same notation for the push-forward, omitting the $_*$ subscript. Another possibility for a notion of duality between sections is to consider the dual modules $\Sec{E}^*\cong \Sec{E^*}$ and define the \textbf{pull-back} of a general vector bundle morphism $F$ as
\begin{align*}
F^*: \Sec{E^*_2} & \rightarrow \Sec{E^*_1}\\
\eta & \mapsto F^*\eta.
\end{align*}
where
\begin{equation*}
(F^*\eta) |_x(e):= \eta|_{\varphi(x)}(F_xe)\quad \forall x\in M_1, e\in (E_1)_x.
\end{equation*}
This map extends to all the covariant tensors, i.e. sections of the tensor powers of the dual bundle $\otimes^k E_2^*$, with the natural properties
\begin{equation*}
F^*(\eta + \omega)=F^*\eta + F^*\omega \qquad F^*(\eta \otimes \omega)= F^*\eta \otimes F^*\omega \qquad F^*f=\varphi^*f
\end{equation*}
for all $\eta,\in\Sec{\otimes^k E_2^*}$, $\Sec{\omega\in\otimes^l E_2^*}$ and $f\in\Cin{M_2}\cong\Sec{\otimes^0 E_2^*}$. In turn, the notion of pull-back allows for an extension of the notion of $F$-relatedness and push-forward to contravariant tensors, i.e. sections of the tensor powers of the vector bundle $\otimes^k E_1$. The notion of pull-back of dual sections is the closest analogue to a contravariant duality functor in the category of vector bundles. Indeed, given two general vector bundle morphisms $F: E_1 \to E_2$ and $G:E_2\to E_3$ one readily checks
\begin{equation*}
    (G\circ F)^*=F^*\circ G^*.
\end{equation*}
In order to regard this construction as a contravariant functor one needs to introduce the \textbf{tensor functor} assigning vector bundles the tensor algebras of their dual bundles and sending vector bundle morphisms to their pull-backs:
\begin{align*}
\mathcal{T}: \Vect_\Man & \rightarrow \textsf{AssAlg}\\
E & \mapsto \bigoplus_{k=0}^\infty \otimes^k E^*\\
F & \mapsto F^*.
\end{align*}
Note that this is simply a specialization of the usual $\text{C}^\infty$ functor regarding the total spaces of vector bundles as manifolds and restricting to fibre-wise polynomial functions.

\subsection{Differential Operators} \label{DifferentialOperators}

Another approach to the question of maps between sections of vector bundles is to consider the module structures in their own right and somewhat forget that they come from an underlying geometric object. Given two vector bundles $A$ and $B$ over the same base $M$, we can identify a class of maps between the spaces of sections that only requires a compatibility with the module structures. We define the \textbf{differential operators of order at most} $k$, for $k$ some integer, as follows:
\begin{equation*}
    \Diff_k(A,B):=\{\Delta : \Sec{A} \to \Sec{B} |\quad c_{f_0}\circ c_{f_1}\circ \dots c_{f_k}\circ (\Delta)=0 \quad \forall f_0,f_1,\dots,f_k\in\Cin{M}\}
\end{equation*}
where we have simplified the notation for commutators of the module multiplications by using the maps $c_f:\text{Hom}_\Real(\Sec{A},\Sec{B})\to \text{End}_\Real(\Sec{A},\Sec{B})$ with explicit action on a $\Real$-linear map of sections $\Delta$ given by
\begin{equation*}
    c_f(\Delta):=[\Delta,f](s)=\Delta(f\cdot s )-f\cdot \Delta(s)
\end{equation*}
for $f\in\Cin{M}$, $s\in\Sec{A}$. Noting that the module product of the first term corresponds to the vector bundle $A$ whereas the one in the second term corresponds to $B$ we observe that both bundles sharing the same base manifold was crucial for the definition of differential operators. Also clear from the definition is that differential operators of order $0$ correspond to $\Cin{M}$-linear maps and thus are identified with vector bundle morphisms covering the identity:
\begin{equation*}
    \Diff_0(A,B)\cong \Hom{A,B}\cong \Sec{A^*\otimes B}.
\end{equation*}
We thus explicitly see that differential operators of order $0$ are the sections of some vector bundle. This is not an accident at order $0$ as spaces of general differential operators can be shown to be projective $\Cin{M}$-modules and thus, in virtue of the general theorem establishing the one-to-one correspondence between projective modules and vector bundles \cite[Theorem 11.32]{nestruev2006smooth}, they must correspond to the sections of some vector bundle. To see this explicitly, given a vector bundle $A$ we introduce its \textbf{$k$-jet bundle} $\Jet^k A$ as the vector bundle whose fibres are formed by the classes of local sections agreeing on derivatives of order up to $k$. The short exact sequence of vector bundles
\begin{equation*}
\begin{tikzcd}
0 \arrow[r] & \odot^k \Cot M\otimes A \arrow[r, "i"] & \Jet^kA \arrow[r, "\pi^k"]  & \Jet^{k-1}A \arrow[r]  & 0,
\end{tikzcd}
\end{equation*}
which follows structurally with $i$ being the inclusion of multidifferentials, is called the called the \textbf{Spencer sequence}. Since all the bundle morphisms cover the identity, we have a similar short exact sequence for the associated spaces of sections and there is a canonical right splitting given by the \textbf{$k$-jet prolongation} $j^k: \Sec{A}\to \Sec{\Jet^kA}$, which simply takes a section to the class of jets corresponding to its derivatives of order at most $k$. This map is sometimes called the \textbf{universal differential operator} since it realizes the isomorphism \cite[Proposition 11.51]{nestruev2006smooth} of $\Cin{M}$-modules:
\begin{equation*}
    \Diff_k(A,B)\cong \Sec{(\Jet^kA)^*\otimes B}
\end{equation*}
via the correspondence
\begin{equation*}
    \Hom{\Jet^k A,B}\ni\Phi \mapsto\Phi \circ j^k\in\Diff_k(A,B).
\end{equation*}

Note that the definition of differential operators above implies that we have $\Diff_{k-1}(A,B)\subset \Diff_k(A,B)$, furthermore, given any differential operator order $k$ we can associate a differential operator of order $0$ that captures its behaviour at the highest order. This is the \textbf{$k$-symbol map}:
\begin{align*}
\sigma: \Diff_k (A,B) & \to \Sec{\odot^k \Tan M\otimes A^*\otimes B}\\
\Delta & \mapsto \sigma_\Delta
\end{align*}
where
\begin{equation*}
    \sigma_\Delta(df_1,\dots,df_k):=c_{f_1}\circ \cdots \circ c_{f_k}(\Delta):A\to B
\end{equation*}
is a well-defined vector bundle morphism for all functions $f_1,\dots,f_k\in\Cin{M}$. The image of a given differential operator $\Delta$ under this map is usually called its $k$-symbol $\sigma_\Delta$. We can reformulate the definition of the symbol map alternatively by writing the following short exact sequence:
\begin{equation*}
\begin{tikzcd}
0 \arrow[r] & \Diff_{k-1}(A,B) \arrow[r ] & \Diff_{k}(A,B) \arrow[r, "\sigma"] & \Sec{\odot^k\Tan M\otimes A^*\otimes B}\arrow[r]  & 0.
\end{tikzcd}
\end{equation*}
 The symbol map can then be regarded as the natural quotient surjection induced by the inclusion $\Diff_{k-1}(A,B)\subset \Diff_k(A,B)$. It is clear by construction that this short exact sequence is recovered as the $\Sec{B}$-adjoint of the Spencer sequence above. Let us now consider a single vector bundle $\alpha: A\to M$, we will denote the space of differential operators from $A$ to itself as $\Diff_k(A):=\Diff_k(A,A)$ and the bundle of fibre-wise endomorphisms by End$(A)\cong A^*\otimes A$. Composition of differential operators gives $\bigcup_k^ \infty \Diff_k(A) \subset\text{Hom}_{\mathbb{R}}(\Sec{A})$ the structure of a filtered associative algebra since one can readily check:
\begin{equation*}
    \Diff_k(A) \circ \Diff_l(A) \subset \Diff_{k+l}(A).
\end{equation*}
Furthermore, as a direct consequence of our definition of general differential operators above, the commutator of differential operators satisfies
\begin{equation*}
    [\Diff_k(A), \Diff_l(A)] \subset \Diff_{k+l-1}(A).
\end{equation*}
Note that the symbol sequence at order 1 becomes
\begin{equation*}
\begin{tikzcd}
0 \arrow[r] & \Sec{\text{End}(A)} \arrow[r] & \Diff_{1}(A) \arrow[r, "\sigma"] & \Sec{\Tan M\otimes \text{End}(A)}\arrow[r]  & 0,
\end{tikzcd}
\end{equation*}
which can be obtained by adjuntion (and taking sections) of the Spencer sequence
\begin{equation*}
\begin{tikzcd}
0 \arrow[r] & \Cot M\otimes A \arrow[r, "i"] & \Jet^1A \arrow[r, "\pi^1"]  & A \arrow[r]  & 0.
\end{tikzcd}
\end{equation*}

A differential operator of degree at most 1, $\Delta\in\Diff_1(A,B)$, can be equivalently characterized by a pair of $\Real$-linear maps
\begin{align*}
\Delta: \Sec{A} & \to \Sec{B}\\
\delta: \Cot M & \to \Hom{A,B}
\end{align*}
such that $\delta$ is, furthermore, $\Cin{M}$-linear and they determine the interaction with the module structure via the Leibniz-like property
\begin{equation*}
    \Delta(f\cdot a) =f\cdot \Delta (a) + \delta (df)(a)
\end{equation*}
for all $a\in\Sec{A}$, $f\in\Cin{M}$. It is clear that this equation gives the defining condition for symbols of differential operators of degree at most 1; indeed, $\sigma_\Delta=\delta$ follows by construction. This is called the \textbf{Leibniz characterization} of the differtial operator $\Delta\in\Diff_1(A,B)$. Since throughout this thesis we will almost exclusively deal with differential operators of order at most 1, we will refer to them simply as \textbf{differential operators}. As shown above, any given differential operator factors through the jet bundle via the universal differential operator $j^1:\Sec{A}\to\Sec{\Jet^1 A}$ and the general characterization above gives, in this particular case, a redefinition of the inclusion map of the Spencer sequence:
\begin{equation*}
    j^1(f\cdot a)=f\cdot a + i(df\otimes a).
\end{equation*}

The endomorphisms of a vector bundle $\text{End}(A)$ contain a submodule of special vector bundle morphisms corresponding to $\Cin{M}$-multiples of the identity $\Id_A$, indeed the natural endomorphism induced by the $\Cin{M}$-module structure at the level of sections. This allows for the identification of a subset of differential operators whose symbol lies in this special submodule of endomorphisms. These are called the \textbf{derivations} of the vector bundle $A$:
\begin{equation*}
    \Dr{A}:=\sigma^{-1}(\Sec{\Tan M\otimes\Id_A})\subset \Diff_1(A).
\end{equation*}
The Leibniz characterization of differential operators allows us to explicitly define derivations as
\begin{equation*}
    \Dr{A}=\{D\in \Diff_1(A)|\quad D(f\cdot s)=f\cdot D(s) + X_D[f]\cdot s \quad \forall f\in \Cin{M},s\in\Sec{A}\},
\end{equation*}
where $X_D$ is the $\Cin{M}$-derivation, i.e. the vector field, corresponding to the symbol in the following sense:
\begin{equation*}
    D\in\Dr{A} \quad \Leftrightarrow \quad \sigma_D=X_D\otimes \Id_A.
\end{equation*}
It follows from this definition that $\Dr{A}\subset \Diff_1(A)$ is a $\Cin{M}$-submodule and thus we should expect derivations to be given by sections of some vector bundle. This is indeed the case, as we have the commutative diagram
\begin{equation*}
\begin{tikzcd}
\Der A \arrow[r, hookrightarrow] \arrow[d, "a"'] & (\Jet^1A)^*\otimes A \arrow[d, "\sigma"] \\
\Tan M \arrow[r, "-\otimes \Id_A"'] & \Tan M\otimes \text{End}(A)
\end{tikzcd}
\end{equation*}
that defines the map $a$ and the vector bundle $\Der A$, which we call the \textbf{der bundle} of $A$. The horizontal maps are injections and the vertical maps are surjections, all being vector bundle morphisms covering the identity. Note that the symbol map being $\Cin{M}$-linear in its arguments allowed us to regard $\sigma$ as a vector bundle morphism in the first place.  We can equivalently characterise the subbundle $\Der A\subset (\Jet^1A)^*\otimes A$ via the Spencer sequence by identifying elements $\delta \in(\Jet^1A)^*\otimes A$ that precompose with the injection to give a an element of the tangent bundle: $\delta\circ i =X_\delta\circ \text{id}_A$ for a (necessarily unique)  $X_\delta\in \Tan M$. With either characterization, it is then easy to show that $\Dr{A}=\Sec{\Der A}$. 

\subsection{Automorphisms and Group Actions} \label{VBundleAutomorphisms}

Invertible vector bundle morphisms of a vector bundle $A$ onto itself are called \textbf{automorphisms} and we denote them by $\text{Aut}(A)$. Clearly, the automorphisms of a vector bundle form a group with composition, this is indeed the isomorphism group of the vector bundle $A$ regarded as an object in the category $\Vect_\Man$. A vector bundle automorphism $F:A\to A$ necessarily covers a diffeomorphism of the base $\varphi:M\to M$. Similarly to how vector fields can be regarded as infinitesimal diffeomorphisms, 1-parameter families of automorphisms are identified with derivations of the vector bundle. We define \textbf{infinitesimal automorphisms} of $A$ as 1-parameter subgroups of automorphisms  $\{F^t:A\to A\}_{t\in I}$ covering  1-parameter subgroups of diffeomorphisms $\{\varphi^t:M\to M\}_{t\in I}$ with $F^0=\Id_E$, which, in turn, implies $\varphi^0=\Id_M$. Since automorphisms cover diffeomorphisms, push-forwards are always well-defined and we can give the following definition of the derivation induced by an infinitesimal automorphism:
\begin{equation*}
    D^F(s):=\tfrac{d}{dt}(F^t)_*(s)|_{t=0}
\end{equation*}
for all $s\in\Sec{A}$. This is indeed a derivation since we can directly compute
\begin{equation*}
\tfrac{d}{dt}(F^t)_*(f\cdot s)|_{t=0}=f\cdot \tfrac{d}{dt}(F^t)_*(s)|_{t=0} + \tfrac{d}{dt}(\varphi^t)^*(f)|_{t=0} \cdot s = f\cdot D^F(s) + X^\varphi[f]\cdot s
\end{equation*}
where $X^\varphi:=\tfrac{d}{dt}\varphi^t|_{t=0}$ is the vector field induced by the infinitesimal diffeomorphism, which clearly acts as the symbol of the derivation $D^F$.\newline

Let $G$ be a Lie group with lie algebra $\mathfrak{g}$ we say that $G$ acts on $A$ via a \textbf{vector bundle action}, $G\Acts A$, when there is a smooth map $F:G\times A\to A$ such that $F_g:A\to A$ is an automorphism for all $g\in G$ and the usual axioms of a group action are satisfied
\begin{equation*}
    F_g\circ F_h = F_{gh}\qquad F_{e}=\Id_A \qquad \forall g,h\in G.
\end{equation*}
It follows by construction that any such action induces a standard group action of $G$ on the base space $\phi:G\times M\to M$, which give the base smooth maps of the automorphisms and is thus called the base action. The infinitesimal counterpart of a vector bundle $G$-action is a morphism of Lie algebras
\begin{equation*}
    \Psi:\mathfrak{g}\to \Dr{A}.
\end{equation*}
Note that $G$ acts on both the domain and codomain of the infinitesimal action by the adjoint action and by push-forward respectively; the fact that $\Psi$ is defined as the infinitesimal counterpart of the action $F$ manifests as $G$-equivariance in the following sense
\begin{equation*}
    \Psi\circ \text{Ad}_g = (F_g)_*\circ \Psi \qquad \forall g\in G.
\end{equation*}
Denoting by $\psi:\mathfrak{g}\to \Sec{\Tan M}$ the infinitesimal counterpart of $\phi$, the infinitesimal vector bundle action satisfies 
\begin{equation*}
    \Psi(\xi)[f\cdot s]=f\cdot \Psi(\xi)[s]+\psi(\xi)[f]\cdot s \qquad\forall f\in\Cin{M},s\in\Sec{A},\xi\in\mathfrak{g}.
\end{equation*}
The \textbf{orbits} of a vector bundle $G$-action $F$ can be simply defined as the images of all group elements acting on a single fibre, and thus they are naturally regarded as the vector bundle restricted to the orbits of the base action $\phi$. In analogy with the case of smooth actions, we denote the set of orbits by $A/G$. Since a vector bundle action is comprised of fibre-to-fibre isomorphisms, any notion defined for usual group actions on smooth manifolds extends to a corresponding notion for vector bundle actions simply requiring the base action to satisfy the corresponding conditions. In particular, a \textbf{free and proper} vector bundle $G$-action gives a well-defined submersion vector bundle morphism
\begin{equation*}
\begin{tikzcd}
A \arrow[r, "Z"] \arrow[d] & A/G \arrow[d] \\
M \arrow[r,two heads, "\zeta"'] & M/G
\end{tikzcd}
\end{equation*}
where the vector bundle structure on the set of orbits $A/G$ is induced from the fact that the orbit space of the base action $M/G$ is a smooth manifold and from the fact that all pairs of fibres over the same orbit are mapped isomorphically by some group element. The linear structure of a fibre of $A/G$ is given by the linear structure of representatives on $A$ and it is well-defined because the $G$-action maps fibres to fibres as isomorphisms of vector spaces. Therefore we see that $A/G$ is a vector bundle over $M/G$ with the same typical fibre as $A$ over $M$. Informally, in the same sense that we say that $M/G$ is a manifold whose points are submanifolds that fit nicely together, we say that $A/G$ is a vector bundle whose total space points are vector subbundles restricted to submanifolds that fit nicely together.\newline

By construction, sections of the orbit bundle $\tilde{s}\in \Sec{A/G}$ are of the form $\tilde{s}(\varphi_G(x))=F_G(s(x))$ for some section $s\in \Sec{A}$. Conversely, for $s\in\Sec{A}$ to give a well-defined section $\tilde{s}\in\Sec{A/G}$ we require
\begin{equation*}
\begin{tikzcd}
F_G(s(\varphi_g(x))=F_G(s(x)) \quad \Leftrightarrow \quad F_g\circ s = s \circ \varphi_g
\end{tikzcd}
\end{equation*}
for all group elements $g\in G$. This then allows us to characterise the module of sections of the orbit bundle with the $G$\textbf{-invariant} sections under push-forward:
\begin{equation*}
\Sec{A/G}\cong \Sec{A}^G:=\{s\in\Sec{A}|\quad (F_g)_*s=s\quad  \forall g\in G\}.
\end{equation*}
Infinitesimally, following the definition of the derivation associated with a 1-parameter subgroup of automorphisms above, we identify the $\mathfrak{g}$\textbf{-invariant} sections as
\begin{equation*}
\Sec{A}^\mathfrak{g}:=\{s\in\Sec{A}|\quad \Psi(\xi)[s]=0 \quad \forall \xi\in \mathfrak{g}\}.
\end{equation*}
Then, if $G$ is a connected Lie group, we have:
\begin{equation*}
\Sec{A/G}\cong \Sec{A}^\mathfrak{g}.
\end{equation*}
We say that a subbundle $B\subset A$ is $G$\textbf{-invariant} if $F_g(B)=B$ for all $g\in G$, in this case the $G$-action $F$ on $A$ clearly restricts to a $G$-action on $B$. Taking the quotient vector bundle $A/B$ we can attempt a definition of a $G$-action on the quotient as follows:
\begin{align*}
\tilde{F} : G & \rightarrow \text{Aut}(A/B)\\
g & \mapsto \tilde{F}_g
\end{align*}
with
\begin{equation*}
F_g(a_x + B_x):= F_g(a_x)+ B_{\varphi_g(x)} \quad \forall x\in M, g\in G.
\end{equation*}
This is well-defined if and only if $B$ is a $G$-invariant subbundle of $A$. Note that if the space of sections of a subbundle $B\subset A$ is $G$-invariant then the subbundle itself will be clearly $G$-invariant; conversely, if the subbundle is $G$-invariant then the condition $F_g(b_x)\in B_{\varphi_g(x)}$ becomes $F_g(s(x))=r(\varphi_g(x))$ for any two $r,s\in \Sec{B}$ at the level of sections which is indeed equivalent to $G$-invariance of the submodule of sections of $B$. In summary:
\begin{equation*}
    F_G(B)\subset B \quad \Leftrightarrow \quad  F \text{ restricts to a $G$-action on } B.
\end{equation*}
We thus see that in a vector bundle $A$ equipped with a $G$-action and a $G$-invariant subbundle $D\subset A$ there are two different notions of quotient. It is a simple check to show that the two constructions commute in the sense that there is a natural isomorphism of vector bundles:
\begin{equation*}
\frac{A}{B}\big/G \cong\frac{A/G}{B/G}.
\end{equation*}
which, at the level of sections becomes
\begin{equation*}
\Gamma \left(\frac{A}{B}\big/G\right) \cong \frac{\Sec{A}^G}{\Sec{B}^G}.
\end{equation*}

\subsection{The Tangent Functor} \label{TangentFunctor}

To close this section we will give a brief account of some basic facts regarding tangent bundles of smooth manifolds, which are indeed the prime example of vector bundles one encounters in standard expositions of elementary differential geometry. The \textbf{tangent bundle} of a smooth manifold, $\Tan M$, can be defined as the manifold of tangent vectors, understood as equivalence classes of locally tangent smooth curves, at all points of $M$ or, equivalently, as the manifold of local derivations of smooth real-valued functions on $M$. These two equivalent characterizations are reflected in the isomorphism of $\Cin{M}$-modules
\begin{equation*}
    \Dr{\Cin{M}}\cong \Sec{\Tan M},
\end{equation*}
connecting an entirely algebraic object on the left, the derivations of the ring of functions, with a geometric one, the bundle of tangent vectors, on the right. The fact that smooth manifolds are locally modelled after vector spaces gives the fibre-wise linear structure of tangent bundles. For similar reasons, any map between manifolds $\varphi: M\to N$ induces a map of tangent vectors, this is indeed the \textbf{tangent map} or \textbf{differential} $\Tan \varphi:\Tan M\to \Tan N$, which is a vector bundle morphism covering $\varphi$. Furthermore, since we have the chain rule $\Tan (\psi\circ \varphi)=\Tan\psi \circ \Tan \varphi$ and $\Tan\Id_M=\Id_{\Tan M}$, the assignment of the tangent bundle to a given smooth manifold is functorial and is thus called the \textbf{tangent functor}
\begin{equation*}
    \Tan :\Man\to \Vect_\Man.
\end{equation*}
This functor interacts nicely with the natural constructions on smooth manifolds. In what follows we give a few details on some of the most relevant examples of such interaction.\newline

The tangent functor sends monomorphisms into monomorphisms. An monomorphism of smooth manifolds we take to be simply an embedding $i:S\hookrightarrow M$, indeed the datum of a submanifold. The fact that the tangent of the embedding $\Tan i:\Tan S \hookrightarrow \Tan M$ is an monomorphism follows by definition. We thus see that the tangent functor sends subobjects in $\Man$ to subobjects in $\Vect_\Man$. From an algebraic point of view, we can identify a submanifold with its vanishing ideal in the ring of functions $I_S\subset \Cin{M}$ so that there is a natural isomorphism of rings $\Cin{S}\cong \Cin{M}/I_S$. This, in turn, allows us to describe the vector fields on the submanifold $S$ solely from algebraic data on the vector fields on $M$ as follows: the vector fields which restrict tangentially to $S$ are defined algebraically by
\begin{equation*}
    \Gamma_S(\Tan M):=\{X\in\Sec{\Tan M}|\quad X[I_S]\subset I_S\},
\end{equation*}
and, in particular, those which restrict to the zero vector field are
\begin{equation*}
    \Gamma_{S_0}(\Tan M):=\{X\in\Sec{\Tan M}|\quad X[I_S]=0\},
\end{equation*}
then there is a natural isomorphism of $\Cin{S}$-modules
\begin{equation*}
    \Sec{\Tan S}\cong \Gamma_S(\Tan M)/\Gamma_{S_0}(\Tan M).
\end{equation*}
Conversely, we may consider subobjects in the vector bundle category lying within the image of the tangent functor. This is the categorical jargon for the familiar notion of $\textbf{tangent distributions}$. These give rise to the question of integrability: given a smooth distribution $D\subset \Tan M$ can we find submanifolds $\{N\subset M\}$ such that $D_x\subset \Tan _xN$ for all $x\in M$? If the answer is positive and the manifolds in the family $\{N\}$ partitioning $M$ are embedded, connected submanifolds we say that $D$ is an \textbf{integrable distribution}. The decomposition of $M$ into this kind of submanifolds is called a \textbf{foliation} on $M$ and each of the embedded submanifolds is called a \textbf{leaf}. If all the leaves $\{N\}$ are of the same dimension, we call it a \textbf{regular foliation}, otherwise we call it a \textbf{singular foliation}. The condition for the integrability of a smooth and regular tangent distribution $D\subset \Tan M$, that is, a subbundle, is given by \textbf{Frobenius' theorem}:
\begin{equation*}
D \text{ arises as the tangent bundles of a regular foliation } \quad \Leftrightarrow \quad [\Gamma(D),\Gamma(D)]\subset \Gamma(D).
\end{equation*}
Associated with a submanifold $i:S\hookrightarrow M$ there are other natural vector bundles in addition to its tangent bundle which are constructed via vector bundle pull-backs and fibre-wise quotients. The restriction of the ambient tangent bundle to the points of the submanifold is called the \textbf{restricted bundle} and defined formally as $\Tan|_S M:=i^*TM$. Since the tangent map of the embedding $i$ is injective, we can also define the \textbf{normal bundle} as the quotient $\Tan M/\Tan S:=\Tan M/(\Tan i(\Tan S))$.\newline

The tangent functor is compatible with the categorical product structures. It simply follows by construction of the fibre-wise product of vector bundles that we have the following natural isomorphism of vector bundles over the Cartesian product of manifolds:
\begin{equation*}
    \Tan(M_1\times M_2)\cong \Tan M_1 \boxplus \Tan M_2.
\end{equation*}
At the level of vector fields, this isomorphism implies that we can include vector fields of each factor into the Cartesian product
\begin{equation*}
\normalfont
\begin{tikzcd}
\Sec{\Tan M_1}\arrow[r,hook, "k_1"] & \Sec{\Tan (M_1\times M_2)} & \Sec{\Tan M_2}\arrow[l,hook',"k_2"'],
\end{tikzcd}
\end{equation*}
where the inclusions $k_i$, $i=1,2$, are defined by the following condition
\begin{equation*}
    k_iX_i[\Proj_j^* f_j] =  \begin{cases} 
      \Proj_i^*f_i & i=j \\
      0 & i\neq j
   \end{cases}
\end{equation*}
for all $f_i\in\Cin{M_i}$ and $X_i\in\Sec{\Tan M_i}$.\newline

The tangent functor sends epimorphisms into epimorphisms. Indeed, the tangent map of a surjective submersion is a surjective vector bundle morphism. More concretely, let us consider the case of a $G$-action on a smooth manifold $\phi:G\times M\to M$. We will often simply write $G\Acts M$ to symbolize a group action whenever the action map need not be made explicit. The tangent functor interacts nicely with group actions in several ways. Firstly, the tangent functor sends Lie groups into Lie groups: we can define the multiplication of elements of the tangent bundle of the Lie group $v_g\in \Tan _gG, v_h\in \Tan _hG$ in the following natural way:
\begin{equation*}
v_g\cdot v_h:=(\Tan _gR_h)v_g+(\Tan _hL_g)v_h \in \Tan _{gh}G
\end{equation*}
where $R,L$ denote the right and left actions of $G$ on itself. Given this multiplication rule, one can use the left-invariant canonical 1-form $\theta^L\in\Omega^1(G,\mathfrak{g})$ to make the usual diffeomorphism $\Tan G\cong G\times \mathfrak{g}$ into the morphism of Lie groups
\begin{equation*}
(\Tan G,\cdot) \cong G\times_{\text{Ad}}\mathfrak{g}
\end{equation*}
where $\mathfrak{g}$ is seen as an abelian group with the addition of vectors. The Lie algebra of this Lie group is then easily shown to be:
\begin{equation*}
\mathfrak{tg} \cong \mathfrak{g}\oplus_{\text{ad}}\mathfrak{g}
\end{equation*}
where the second component is taken as the abelian Lie algebra. Furthermore, given a smooth action on a manifold $\phi: G\times M \to M$ we can use the fact that $\Tan (G\times M)\cong \Tan G \times \Tan M$ to define an action as follows:
\begin{align*}
\Tan \phi: \Tan G\times \Tan M & \to \Tan M\\
v_g+w_x & \mapsto \Tan _g\phi(\cdot,x)v_g+\Tan _x\phi(g,\cdot)w_x.
\end{align*}
This is easily checked to be an action $\Tan G\Acts \Tan M$ with respect to the product of tangent vectors to $G$ defined above, we call this the \textbf{tangent lift} of the action $\phi$. Note that this is not a vector bundle action in the sense introduced in Section \ref{VBundleAutomorphisms} since the Lie group acting on the total space ($\Tan G$) is not the same as the group acting on the base space ($G$). When the $G$-action on $M$ is free and proper, so that the orbit space is a smooth manifold and there is a surjective submersion $q:M\to M/G$, the tangent lift of the action is also free and proper and there is a natural diffeomorphism:
\begin{equation*}
    \Tan(M/G)\cong \Tan M/\Tan G.
\end{equation*}
This quotient can be reformulated in terms of the \textbf{tangent action}, which is simply given by taking the differentials of the action diffeomorphisms $F_g=\Tan \phi_g$ and thus is clearly a vector bundle action. Note that a free and proper action defines a regular tangent distribution given simply by the image of the infinitesimal action map $\psi(\mathfrak{g})\subset \Tan M$. Recall that the infinitesimal action map is $G$-equivariant so it follows that for a free and proper action $\psi(\mathfrak{g})\subset \Tan M$ is a $G$-invariant subbundle. If $G$ is connected, then we find the isomorphism of vector bundles
\begin{equation*}
\Tan (M/G)\cong \frac{\Tan M}{\psi(\mathfrak{g)}}\big/ G.
\end{equation*}
Algebraically, the ring of functions of the orbit space is identified with the $G$-invariant functions in a natural way $\Cin{M/G}\cong \Cin{M}^G$. Then, by identifying the $G$-invariant vector fields as
\begin{equation*}
    \Gamma_G(\Tan M):=\{X\in\Sec{\Tan M}|\quad X[I_S]\subset I_S\}=\{X\in\Sec{\Tan M}|\quad X[\Cin{M}^G]\subset \Cin{M}^G\},
\end{equation*}
and, in particular, those which are tangent to the orbits
\begin{equation*}
    \Gamma_{G_0}(\Tan M):=\{X\in\Sec{\Tan M}|\quad X[\Cin{M}^G]=0\},
\end{equation*}
we can describe the vector fields of the orbit space via the natural isomorphism of $\Cin{M/G}$-modules
\begin{equation*}
    \Sec{\Tan M/G}\cong \Gamma_G(\Tan M)/\Gamma_{G_0}(\Tan M).
\end{equation*}

\section{Local Lie Algebras} \label{LocalLie}

\subsection{Locality and Derivative Algebras} \label{DerivativeAlgebras}

The presence of linear structures alongside the differentiable structure of a smooth manifold $M$, e.g. the vector spaces that locally model the space, the real numbers where smooth functions take values or the fibres of vector bundles over it, allows for a geometric characterization of the properties of algebraic structures associated with manifolds and vector bundles. One key notion is \textbf{locality}, which is understood broadly as the property of objects that take values over $M$ to sit nicely with the topological and differentiable structure, i.e. they only require the information of arbitrarily small neighbourhoods around a point to specify their values at that point. More concretely, locality appears as the compatibility of the notion of \textbf{support}, the subset where an object takes non-zero values, with the natural algebraic structures associated with a manifold. In particular, the ring multiplication of functions satisfies
\begin{equation*}
    \text{supp}(fg)\subset \text{supp}(f) \cap \text{supp}(g) \quad \forall f,g\in\Cin{M},
\end{equation*}
or, similarly, for the module multiplication of sections of any vector bundle $\alpha:A\to M$ we have
\begin{equation*}
    \text{supp}(f\cdot s)\subset \text{supp}(f) \cap \text{supp}(s) \quad \forall f\in\Cin{M}, s\in \Sec{A}.
\end{equation*}
As part of the general correspondence between projective modules and vector bundles, it can be shown that the natural notion of map between modules of sections compatible with the locality property of the module multiplication above is equivalent to that of a differential operator as defined in Section \ref{DifferentialOperators}. See \cite[Chapter 9]{nestruev2006smooth} for a detailed discussion on this topic. The discussion about the jet bundle as a means to recover differential operators as sections of some vector bundle found in that section precisely captures the idea that differential operators indeed only require local information (the partial derivatives of sections up to some finite order) to act on their arguments.\newline

A vector bundle $\alpha:A\to M$ is called a \textbf{local Lie algebra} when its module of sections is endowed with a ($\Real$-linear) Lie algebra structure $(\Sec{A},[,])$ such that
\begin{equation*}
    \text{supp}([a, b])\subset \text{supp}(a) \cap \text{supp}(b) \quad \forall  a,b\in \Sec{A}.
\end{equation*}
Note that this is the natural compatibility condition to impose on a Lie bracket that is introduced in the module of sections of a smooth vector bundle. As mentioned above, a map respecting the locality condition corresponds to a differential operator, thus we could define local Lie algebras as Lie brackets on sections $(\Sec{A},[,])$ that act as differential operators on each argument, in other words, a Lie algebra structure with respect to the $\Real$-vector space structure of $\Sec{A}$ such that its adjoint map is a differential operator of the form
\begin{equation*}
    \text{ad}_{[,]}: \Sec{A}\to \Diff_1(A).
\end{equation*}
The first source of natural examples of local Lie algebras are the bundles of differential operators of general vector bundles: recall from Section \ref{DifferentialOperators} that differential operators of a vector bundle $E$ correspond to the sections of the adjoint jet bundle
\begin{equation*}
    \Diff_1(E)\cong \Sec{(\Jet^1 E)^*\otimes E},
\end{equation*}
and the filtered associative algebra structure found in $(\Diff_1(E),\circ)$ allows us to define a commutator bracket satisfying
\begin{equation*}
    [\Diff_1(E),\Diff_1(E)]\subset \Diff_1(E),
\end{equation*}
thus giving a Lie bracket for sections of $A=(\Jet^1 E)^*\otimes E$. The fact that this bracket is local follows directly from the fact that differential operators are, themselves, local. We refrain from giving more examples at this point since any Lie algebra structure found throughout this thesis is, in fact, a local Lie algebra.\newline

In our discussion of Section \ref{DifferentialOperators} we identified the derivations of a vector bundle as a particular subset of differential operators that is closely tied to the geometry of the vector bundle via their correspondence with infinitesimal automorphisms. We are thus compelled to refine the general notion of local Lie algebra and define \textbf{derivative Lie algebra} as a vector bundle $\alpha:A\to M$ whose module of sections is endowed with a ($\Real$-linear) Lie algebra structure $(\Sec{A},[,])$ that acts as derivations on each of its arguments, that is, whose adjoint map is of the form
\begin{equation*}
    \text{ad}_{[,]}: \Sec{A}\to \Dr{A}.
\end{equation*}
Perhaps unsurprisingly, the first natural example of a derivative Lie algebra is the der bundle of a vector bundle with the commutator bracket restricted to the subspace of derivations, this case will be discussed in detail in Section \ref{LieAlgebroids} below. The given name to this subclass of local Lie algebras is motivated both by this example and the fact that the bracket of a derivative Lie algebra satisfies the \textbf{Leibniz identity} with respect to the $\Cin{M}$-module structure:
\begin{equation*}
    [a,f\cdot b]=f\cdot [a,b] + \lambda_a[f]\cdot b
\end{equation*}
for all $a,b\in\Sec{A}$, $f\in\Cin{M}$. Here $\lambda:\Sec{A}\to\Sec{\Tan M}$ denotes the asignment of the symbol of the derivation given by the second argument of the Lie bracket and we call it the \textbf{symbol} of the derivative Lie algebra $(\Sec{A},[,])$ for short. More precisely,
\begin{equation*}
    \lambda_a:=\sigma_{\text{ad}_{[,]}(a)}.
\end{equation*}
Note that we could similarly define $\lambda'$ as the symbol of the first argument of the bracket, but then it follows trivially that $\lambda'=-\lambda$.
\begin{prop}[Symbol and Squiggle of a Derivative Lie Algebra] \label{SymbolLieAlgebra}
Let $\alpha:A\to M$ be a vector bundle and $(\Sec{A},[,])$ a derivative Lie algebra structure, then the symbol map $\lambda:\Sec{A}\to \Sec{\Tan M}$ is a differential operator and a Lie algebra morphism. The symbol of $\lambda$ will be denoted by
\begin{equation*}
\normalfont
    \Lambda^\sharp \in\Sec{\Cot M\otimes A^* \otimes \Tan M}
\end{equation*}
and will be called the \textbf{squiggle} of the derivative Lie algebra $(\Sec{A},[,])$.
\end{prop}
\begin{proof}
That $\lambda$ is a differential operator follows simply from the fact that a local Lie bracket is a differential operator in each entry. To see this explicitly, consider the nested commutators $[[\lambda,f],g](a)$ for any pair of functions $f,g\in\Cin{M}$ and a section $a\in\Sec{A}$, then, by use of $\lambda'=-\lambda$, a direct computation shows that it vanishes when acting on an arbitrary function. The definition of the squiggle of the bracket is then simply the symbol of this map $\Lambda^\sharp=\sigma_\lambda$ when regarded as a differential operator $\lambda\in\Diff_1(A,\Tan M)$. In order to show that $\lambda$ is a Lie algebra morphism consider a bracket of the form $[[a,b],f\cdot c]$, then, by definition of the symbol of the bracket we can expand as
\begin{equation*}
    [[a,b],f\cdot c]=f\cdot [[a,b],c]+\lambda_{[a,b]}\cdot b.
\end{equation*}
On the other hand, by using the Jacobi identity of the Lie bracket first and then expanding by the definition of the symbol of the bracket we obtain, after cancellations:
\begin{equation*}
    [[a,b],f\cdot c]=[[a,f\cdot c],b]+[a,[b,f\cdot c]]=f\cdot ([[a,c],b]+[a,[b,c]])+\lambda_a[\lambda_b][f]\cdot c - \lambda_b[\lambda_a][f]\cdot c.
\end{equation*}
Since both expressions must agree for all $a,b,c\in\Sec{A}$ and $f\in\Cin{M}$, using the Jacobi identity of the Lie bracket once more, the desired result follows
\begin{equation*}
    \lambda_{[a,b]}=[\lambda_a,\lambda_b].
\end{equation*}
\end{proof}
It then follows that the bracket of a derivative Lie algebra $(\Sec{A},[,])$ satisfies the following general compatibility property with the $\Cin{M}$-module structure on sections:
\begin{equation*}
    [f\cdot a,g\cdot b]=fg\cdot[a,b]+f\lambda_a[g]\cdot b-g\lambda_b[f]\cdot a +\Lambda(df\otimes a,dg\otimes b)
\end{equation*}
for all $f,g\in\Cin{M}$, $a,b\in\Sec{A}$ and where we have used the fact that, by construction, we have
\begin{equation*}
    \Lambda^\sharp(df\otimes a)[g]\cdot b=-\Lambda^\sharp(dg\otimes b)[f]\cdot a
\end{equation*}
so we can regard the squiggle as a bilinear form $\Lambda\in\Sec{\wedge^2(\Tan M \otimes A^*)\otimes A}$. This identity will be called the \textbf{symbol-squiggle identity} of the local Lie algebra $(\Sec{A},[,])$.\newline

When the squiggle of a derivative Lie algebra $(\Sec{A},[,])$ vanishes, the symbol $\lambda:\Sec{A}\to \Sec{\Tan M}$ becomes a $\Cin{M}$-linear map and is thus induced from the push-forward of a vector bundle morphism covering the identity $\rho:A\to \Tan M$, $\lambda=\rho_*$. In this case, the symbol of the bracket is called an \textbf{anchor}.\newline

It turns out that local Lie algebra brackets are quite sensitive to the rank of the underlying vector bundle being 1 or greater than 1. This dichotomy, evidenced by the two propositions below, will motivate the definitions of the central objects of interest in this thesis: Jacobi manifolds and Lie algebroids.
\begin{prop}[Local Lie Algebras of Rank 1, \cite{kirillov1976local}] \label{LocalLieRk1}
A local Lie algebra structure $(\Sec{A},[,])$ with $\normalfont\rkk{A}=1$ is necessarily a derivative Lie algebra.
\end{prop}
\begin{proof}
This is a direct consequence of the fact that the bundle of fibre-wise endomorphism of a vector bundle $\alpha: A\to M$ with 1-dimensional fibres is trivial
\begin{equation*}
    \text{End}(A)=A^*\otimes A \cong \Real_M,
\end{equation*}
and thus all endomorphisms can be regarded as $\Cin{M}$-multiples of the identity map. It is clear then that the symbol of any differential operator $\Delta\in\Diff_1(A)$ will be uniquely determined by a derivation on the ring of functions, i.e. a vector field. The result then follows from the general fact about vector bundles of rank 1:
\begin{equation*}
    \rkk{A}=1 \quad \Rightarrow \quad \Diff_1(A)=\Dr{A}.
\end{equation*}
\end{proof}

\begin{prop}[Derivative Lie Algebras of Rank 2]\label{LocalLieRk2}
The symbol of any derivative Lie algebra $(\Sec{A},[,])$ with $\normalfont\rkk{A}\geq 2$ is necessarily an anchor.
\end{prop}
\begin{proof}
Restricting to a trivializing neighbourhood $U\subset M$ we can choose local sections $a,b\in \Gamma_U(A)$ that are $\Cin{U}$-linearly independent, this is guaranteed generically in sufficiently small open neighbourhoods since the rank of the vector bundle is 2 or greater: we are always able to choose two independent directions at any given fibre and then extend smoothly. Take two such sections $a,b\in \Gamma_U(A)$, two local functions $f,g\in \Cin{U}$ and consider the symbol-squiggle bracket $[f\cdot a, g\cdot b]$. Applying the Leibniz identity on each side in two different orders we get two expressions that must agree:
\begin{equation*}
    \lambda_{f\cdot a}[g]\cdot b - g\cdot\lambda_b[f]\cdot a + fg\cdot [a,b]=[f\cdot a, g\cdot b]=f\cdot \lambda_{a}[g]\cdot b - \lambda_{g\cdot b}[f]\cdot a + fg\cdot [a,b]
\end{equation*}
Since $a$ and $b$ are assumed to be $\Cin{U}$-independent, each factor accompanying them must vanish independently, thus giving
\begin{equation*}
    \lambda_{f\cdot a}=f\cdot\lambda_{a}
\end{equation*}
for all $a\in \Gamma_U(A)$, $f\in\Cin{U}$. This gives $\Cin{U}$-linearity of the symbol map restricted to the trivializing neighbourhood $\lambda|_U$, however this is clearly a trivialization-independent property as introducing other trivializations will give $\Cin{U}$-linear combinations of the local sections. Hence, $\lambda$ is globally $\Cin{M}$-linear, showing that it is an anchor, as desired.
\end{proof}

Defining a satisfactory general notion of \textbf{morphism of local Lie algebras} is no simple matter. This difficulty stems from the fact that local Lie algebras are algebraic structures on modules of sections of some vector bundle, where there are, of course, multiple ways in which one may require the Lie brackets and module maps to be compatible with the geometric structure of the underlying manifolds (see the discussion of Section \ref{Sections}). The general strategy to define an appropriate notion of morphism for specific classes of local Lie algebras will be to regard them in the context of a broader, well-defined category and demand that they form a subcategory. This will be exemplified in our definition of Lie algebroid morphism of Section \ref{LieAlgebroids} and Jacobi map of Section \ref{ContactGeometry}.

\subsection{Lie Algebroids} \label{LieAlgebroids}

A vector bundle $\alpha:A\to M$ with a derivative algebra structure $(\Sec{A},[,])$ is called a \textbf{Lie algebroid} if the squiggle of the bracket vanishes, i.e. the symbol of the bracket is $\Cin{M}$-linear and thus given by an anchor map $\rho:A\to \Tan M$. We will often denote the algebroid structure as the triple $(A,\rho,[,])$. It follows from proposition \ref{SymbolLieAlgebra} that the anchor map induces a morphism of Lie algebras $\rho_*:\Sec{A}\to \Sec{TM}$ which realizes the Leibniz identity of the local Lie bracket:
\begin{equation*}
    [a,f\cdot b]=f\cdot [a,b]+\rho_*(a)[f]\cdot b
\end{equation*}
for all $a,b\in\Sec{A}$, $f\in\Cin{M}$. Whenever there is no room for confusion we will drop the anchor map from our notation by writing $a[f]:=\rho_*(a)[f]$ to signify the action of sections of a Lie algebroid as derivations of functions on the base.\newline

Since the anchor is a vector bundle morphism $\rho: A\to M$, two distributions are naturally defined: the kernel $\Ker{\rho}\subset A$ and the image $\rho(A)\subset \Tan M$, both generically irregular. It follows from the Leibniz identity above that the vector spaces $\mathfrak{g}_x:=\Ker{\rho_x}\subset A_x$ inherit a Lie algebra structure from the Lie bracket on local sections extending elements of the fibre $A_x$. We call these $(\mathfrak{g}_x,[,]_x)$ the \textbf{isotropy Lie algebras} of the Lie algebroid $A$. The fact that the anchor is a Lie algebra morphism ensures that the image distribution $\rho(A)$ is involutive, then it will be generally integrable by a singular foliation on $M$. The image distribution is called the \textbf{characteristic distribution} of the Lie algebroid $A$.\newline

The bracket of a Lie algebroid $(A,\rho,[,])$ can be extended uniquely to a Gerstenhaber bracket on the graded algebra of multisections $(\Sec{\wedge^\bullet A},\wedge,\llbracket , \rrbracket)$ by setting
\begin{align*}
    \llbracket a,b \rrbracket &=[a,b]\\
    \llbracket a,f \rrbracket &=a[f]\\
    \llbracket f,g \rrbracket &= 0
\end{align*}
for all $a,b\in\Sec{\wedge^1A}=\Sec{A}$, $f,g\in\Sec{\wedge^0A}=\Cin{M}$ and defining the action on higher wedges of sections by graded derivations. This is called the \textbf{Gerstenhaber algebra} of the Lie algebroid $A$. Dually to this construction, we find the datum of a Lie algebroid structure on $A$ equivalently encapsulated in the graded algebra of dual forms $(\Sec{\wedge^\bullet A^*},\wedge)$ as a graded differential $d_A$ defined explicitly on a homogeneous element $\omega\in\Sec{\wedge^kA^*}$ as
\begin{equation*}
    d_A\omega(a_0,a_1,\dots,a_k):=\sum_{i< j} (-1)^{i+j-1}\omega([a_i,a_j],a_1,\dots,\hat{a_i},\dots,\hat{a_j},\dots,a_k) + \sum_{i=0}^k (-1)^{i-1} a_i[\omega(a_1,\dots,\hat{a_i},\dots,a_k)]
\end{equation*}
for all $a_0,a_1,\dots,a_k\in\Sec{A}$ and extended by linearity. We will call the differential graded algebra $\Omega^\bullet(A):=(\Sec{\wedge^\bullet A^*},\wedge ,d_A)$ the \textbf{exterior algebra} of the Lie algebroid $A$ and $d_A$ the \textbf{Lie algebroid differential}. We note that $\Omega^\bullet(A)$ is often also referred to as the \textbf{de Rham complex} of the Lie algebroid $A$ in the literature. The cohomology of this complex is called the \textbf{Lie algebroid cohomology (with trivial coefficients)} and will be denoted by $H^\bullet(A)$. It follows trivially from the definition of the exterior algebra above that, given a section of the Lie algebroid $a\in\Sec{A}$, there are natural notions of \textbf{interior product} $i_a$ and of \textbf{Lie derivative} $\LDer_a$ defined on homogeneous elements $\omega\in\Sec{\wedge^kA^*}$ by
\begin{align*}
    i_a\omega(a_0,a_1,\dots,a_{k-1}) &=\omega(a,,a_0,a_1,\dots,a_{k-1})\\
    \LDer_a\omega(a_0,a_1,\dots,a_k) &=\sum_{i=0}^k\omega(a_0,\dots,a_{i-1},[a,a_i],a_{i+1},\dots,a_k)-a[\omega(a_0,a_1,\dots,a_k)]
\end{align*}
and extended by linearity. A simple computation shows that $d_A,i_a,\LDer_a$ satisfy the usual identities of Cartan calculus
\begin{align*}
    \LDer_a &=i_a\circ d_A+d_A\circ i_a \\
    [\LDer_a,\LDer_b] &= \LDer_{[a,b]} \\
    [\LDer_a,i_b] &= i_{[a,b]} \\
    [i_a,i_b] &= 0,
\end{align*}
and since $i_a,\LDer_a,d$ are graded derivations of degree $-1,0+1$, respectively,  we will refer to these identities as the \textbf{Lie algebroid Cartan calculus (with trivial coefficients)} of $A$.\newline

The general problem of finding the right notion of morphism of local Lie algebras becomes more tractable in the case of Lie algebroids. Since the vector bundle of a Lie algebroid $(A,\rho,[,])$ is not subject to any general requirements, it will be desirable to find a notion of Lie algebroid morphism that involves, in particular, a general vector bundle morphism. Then, considering a vector bundle morphism between two Lie algebroids $F:A\to B$, we would like to demand that this map induces a Lie algebra morphism of the modules of sections. The issue is that, as was discussed in Section \ref{Sections}, vector bundle morphisms do not induce well-defined maps between the modules of sections in general. However, it was also shown that pull-backs of covariant tensors are always defined and thus we have a morphism of associative algebras $F^*:(\Sec{\wedge^\bullet B^*}, \wedge)\to (\Sec{\wedge^\bullet A^*},\wedge)$. Luckily, we can use the differential graded algebra characterization of a Lie algebroid on its exterior algebra introduced above to define a vector bundle morphism $F:A\to B$ covering a smooth map $\varphi:M\to N$ to be a \textbf{Lie algebroid morphism} when the pull-back is a morphism of differential graded algebras
\begin{equation*}
    F^*:\Sec{\wedge^\bullet B^*}, \wedge,d_B)\to (\Sec{\wedge^\bullet A^*},\wedge,d_A), \qquad d_A \circ F^* = F^*\circ d_B.
\end{equation*}
In particular, as a bundle map, a Lie algebroid morphism satisfies
\begin{enumerate}
    \item Compatibility with the anchor: $\rho_B\circ F=\Tan\varphi\circ\rho_A$.
    \item Compatibility with the Lie bracket: if $a\sim_F b$ and $a'\sim_F b'$ then $[a,a']_A\sim_F[b,b']_B$.
\end{enumerate}
Lie algebroids together with Lie algebroid morphisms then constitute the \textbf{category of Lie algebroids} $\Lie_\Man$, that can be naturally seen as as a subcategory of the category of vector bundles
\begin{equation*}
    \Lie_\Man\subset \Vect_\Man.
\end{equation*}

Lie algebroids should be regarded as generalizations of both tangent bundles and Lie algebras. Indeed, we clearly recover the usual notions of Lie algebra morphism, Lie algebra cohomology, de Rham differential, Schouten bracket, etc. when regarding Lie algebras as Lie algebroids with base manifold a point and tangent bundles as Lie algebroids with the identity anchor and the Lie bracket of vector fields. We note that, upon identifying the category of Lie algebroids, the tangent functor of Section \ref{TangentFunctor} becomes
\begin{equation*}
    \Tan: \Man\to \Lie_\Man.
\end{equation*}
Other, more interesting examples of Lie algebroids include regular involutive distributions on tangent bundles, with the anchor being simply the inclusion of the distribution into the tangent bundle, and infinitesimal actions of Lie algebras on manifolds, with the vector bundle given by the trivial Lie algebra bundle and the anchor by the action map. The datum of a Lie algebroid on a trivial vector bundle of rank 1 is equivalent to a smooth vector field on the base manifold: given a manifold $M$ and a vector field $X\in\Sec{\Tan M}$ there is a Lie algebroid structure $(A,\rho_X,[,]_X)$ given by
\begin{equation*}
    A=M\times \Real, \qquad \rho: f \mapsto f\cdot X, \qquad [f,g]_X=fX[g]-gX[f],
\end{equation*}
for all $f,g\in\Sec{A}\cong \Cin{M}$. Another important class of examples of Lie algebroids with surjective anchors are derivation bundles. Recall from Section \ref{DifferentialOperators} that given any vector bundle $\epsilon:E\to M$ its derivation bundle $\Der E$ is naturally equipped with an anchor map, the restricted symbol map, and its sections are equipped with the commutator bracket thus making it into a Lie algebroid. Another important example of Lie algebroid is the Atiyah algebroid $A_P$ of a principal $G$-bundle $\pi:P\to M$, as a vector bundle it appears in the Atiyah sequence constructed via the infinitesimal action and bundle projection
\begin{equation*}
\begin{tikzcd}
0 \arrow[r] & P\times_G \mathfrak{g} \arrow[r] & A_P \arrow[r]  & \Tan M \arrow[r]  & 0,
\end{tikzcd}
\end{equation*}
where $\mathfrak{g}$ is the Lie algebra of $G$, and the bracket is induced from the Lie bracket of left-invariant vector fields on $G$. Numerous other examples of Lie algebroids appear in the context of Poisson and Jacobi geometry, these will be discussed in detail in sections \ref{SymplecticGeometry}, \ref{DiracGeometry}, \ref{ContactGeometry} and \ref{LDiracGeometry} below.\newline

Having identified the derivation bundles of general vector bundles as Lie algebroids, we can recover the standard definition of a \textbf{connection} on a vector bundle $\epsilon: E \to M$ as a vector bundle morphism covering the identity $\nabla:\Tan M\to \Der E$ that is compatible with the anchors and whose curvature is the bilinear form measuring the degree to which $\nabla$ fails to be a Lie algebroid morphism. Given a Lie algebroid $\alpha:A\to M$ this can now be generalized in a direct manner by defining a $A$\textbf{-connection} as a vector bundle morphism $\nabla:A\to \Der E$ that is compatible with the anchors and whose \textbf{curvature} is defined as $R_\nabla(a,b):=[\nabla_a,\nabla_b]-\nabla_{[a,b]}$. When the curvature of $\nabla$ vanishes we say that the $A$-connection is \textbf{flat}, this is indeed equivalent to the map $\nabla:A\to \Der E$ being a Lie algebroid morphism. When $\nabla$ is a flat $A$-connection, the triple $(A,E,\nabla)$ is called a \textbf{representation} of the Lie algebroid $A$ on the vector bundle $E$.\newline

In the presence of a Lie algebroid representation $\nabla:A\to \Der E$ we can define a new complex of dual forms with values in the representation $E$:
\begin{equation*}
    \Omega_E^\bullet(A):=(\Sec{\wedge^\bullet A^*\otimes E},\wedge , d_A^E).
\end{equation*}
The wedge product is defined from the identity
\begin{equation*}
    (\omega \otimes e)\wedge(\eta\otimes e)=\omega\wedge\eta\otimes e
\end{equation*}
by $\Cin{M}$-linear extension. The differential is defined by setting $d_A^Ee:=\nabla_{-} e$ for $e\in\Sec{\wedge^0A^*\otimes E}\cong\Sec{E}$ and
\begin{equation*}
    d_A^E(\omega\wedge\xi):=d_A\omega\wedge\xi +(-)^k\omega\wedge d_A^E\xi
\end{equation*}
for $\omega\in\Omega^k(A)$ and $\xi\in\Omega^\bullet_E(A)$ and extending by $\Cin{M}$-linearity. This is called the \textbf{exterior algebra of forms with values in $E$} of the Lie algebroid $A$. The cohomology of this complex is denoted by $H^\bullet(A,E)$ and it is called the \textbf{Lie algebroid cohomology with coefficients in $E$}. Morally, in this construction $\Sec{E}$ plays the role the $\Cin{M}$ played in the definition of the exterior algebra of forms with trivial coefficients above, where the action of $\Sec{A}$ as derivations on functions has been replaced by the action of $\Sec{A}$ on $\Sec{E}$ via the representation $\nabla$. It is clear then that this complex carries a Cartan calculus, that we refer to as the \textbf{Cartan calculus of $A$ with values in }$E$, defined in an entirely analogous manner to the case of trivial coefficients by noting that the Lie derivative on sections of the vector bundle is given precisely by the representation map: $\LDer_ae=\nabla_a e$. A derivation bundle $\Der E$ carries a tautological Lie algebroid representation on itself $\Id_{\Der E}:\Der E\to \Der E$, the corresponding complex $\Omega_E^\bullet(DE)$ is called the \textbf{der complex} of $E$, and it is easily shown to be acyclic: $H^\bullet(DE,E)=0$. Dually to the der-complex we find the \textbf{multiderivations} $(\text{Der}^\bullet(E),\wedge,\llbracket , \rrbracket)$, defined inductively from $\text{Der}^0(E):=\Sec{E}$ and $\text{Der}^1(E):=\Dr{E}$. The multiderivations carry a Schouten-like bracket that is defined by setting:
\begin{align*}
    \llbracket D,D' \rrbracket &=[D,D']\\
    \llbracket D,a \rrbracket &=D(a)\\
    \llbracket a,b \rrbracket &= 0
\end{align*}
for $a,b\in\Sec{E}$, $D,D'\in\Der{E}$ and extending as graded derivations.\newline

\subsection{Lie Groupoids} \label{LieGroupoids}

In the context of category theory, the notion of symmetry within a class of objects is encapsulated by the concept of isomorphism. According to the standard definition, isomorphisms are morphisms between objects that can be reversed, i.e. they have inverses, and leave the object unchanged, i.e. they compose to the identity. We can identify this reversibility property of certain morphisms within a category as the essential feature of symmetry. We define a \textbf{groupoid} $G$ as a category whose morphisms are all isomorphisms. We regard groupoids as the natural categorical structures capturing the idea of symmetry, noting that this is a generalization of the common conception that groups, which are groupoids with a single object, are the natural mathematical representation of symmetry.\newline

In a similar vein to the formal introduction of Lie groups, we will be interested in groupoids that are also smooth manifolds. Before giving the precise definition of Lie groupoid, let us spell out the set-theoretic content of a general  groupoid (for the reminder of this section we assume any category to be small, so that the classes of objects and morphisms are sets):
\begin{itemize}
    \item The set of all morphisms between all pairs of objects will be called the \textbf{set of arrows} or, simply, the \textbf{groupoid}, $G$.
    \item The set of objects will be called the \textbf{base} of the groupoid, $M$.
    \item An arrow $g\in G$ between two objects $x,y\in M$ will be denoted indistinctly by $g:x\to y$, $g\in G(x,y)$ or $g_{xy}$. This gives two surjections from the groupoid onto the base $s,t:G\to M$ defined by $s(g)=x$ and $t(g)=y$, and aptly named \textbf{source} and \textbf{target} maps, respectively.
    \item As in any category, to each object $x\in M$ corresponds a unique identity arrow $1_x\in G(x,x)$, this assignment gives an injective map $u:M\to G$ that is called the \textbf{identity} map.
    \item From the general axiom of composition in a category we see that a groupoid carries a partial binary operation, since only arrows with matching sources and targets can be composed. This is the main feature that sets groupoids apart from groups. The \textbf{set of composable arrows} is defined as
    \begin{equation*}
        G  \ftimes{s}{t} G :=\{(g,h)\in G\times G|\quad t(g)=s(h)\},
    \end{equation*}
    then the categorical composition becomes a map $m:G _s\times_t G\to G$ called the \textbf{groupoid multiplication}.
    \item The structure described so far is common to all (small) categories, what characterizes a groupoid is that all the arrows have an inverse: for any $g\in G(x,y)$ there is a unique arrow $h\in G(y,x)$ such that $g\circ h = 1_x$ and $h\circ g = 1_y$, or, in terms of the structure maps introduced so far, $m(g,h)=u(x)$ and $m(h,g)=u(y)$. This gives an idempotent bijective map $i:G\to G$ sending an arrow to its inverse; we call $i$ the \textbf{inversion} map.
\end{itemize}
All of the maps above are collectively referred to as the \textbf{structure maps} of the groupoid $G$ with base $M$; they can be summarised in the following commutative diagram:
\begin{equation*}
\begin{tikzcd}
G  \arrow[loop left, "i"]\arrow[d,shift right, "s"'] \arrow[d,shift left, "t"] & G \ftimes{s}{t} G \arrow[l, "m"'] \\
M \arrow[u, "u"', bend right=50] & 
\end{tikzcd}
\end{equation*}

For any fixed object $x\in M$ we define the set of all the arrows departing from $x$ as the preimage by the source map $s^{-1}(x)\subset G$ and we call it the \textbf{$s$-fibre} over $x$; similarly, we define the subset of all the arrows arriving at $x$ by $t^{-1}(x)\subset G$ and call it the \textbf{$t$-fibre} over $x$. The inversion map clearly restricts to a bijection $i:s^{-1}(x)\to t^{-1}(x)$ by construction. For any given arrow $g\in G(x,y)$, composition with the partial groupoid multiplication allows us to define the \textbf{right action}
\begin{align*}
R_g: s^{-1}(y) & \to s^{-1}(x)\\
h & \mapsto m(h,g),
\end{align*}
and the \textbf{left action}
\begin{align*}
L_g: t^{-1}(x) & \to t^{-1}(y)\\
h & \mapsto m(g, h).
\end{align*}
Since all the arrows are invertible, the left and right action maps $L_g,R_g$ are bijections for all $g\in G$. The set of arrows sending $x$ to itself, $G_x:=G(x,x)=s^{-1}(x)\cap t^{-1}(x)$, is called the \textbf{isotropy} at $x$ and it becomes a group with multiplication given by the restriction of the partial groupoid multiplication, identity given by $1_x$ and inversion given by the restriction of the inversion of the groupoid. A groupoid defines an equivalence relation on the base as follows
\begin{equation*}
    x\sim_G y \quad \Leftrightarrow \quad \exists g\in G \quad | \quad g:x\to y.
\end{equation*}
An equivalence class of this relation is called an \textbf{orbit}. The orbit trough a point of the base $x\in M$ is denoted as $O_x\subset M$ and we readily check:
\begin{equation*}
    O_x=\{t(g)\in M,\quad g\in s^{-1}(x)\}.
\end{equation*}
The quotient by this equivalence relation gives the set of orbits denoted by $M/G$ and called the \textbf{orbit set} of the groupoid $G$.\newline

A groupoid $G$ whose underlying sets $(G,M,G\ftimes{s}{t}G)$ are smooth manifolds and whose structure maps  $(s,t,u,i,m)$ are smooth is called a \textbf{Lie groupoid}. Note that the requirement that $G\ftimes{s}{t}G$ is a smooth manifold is tantamount to demanding that the source and target maps $s,t$ are submersions. There are a few immediate consequences of this definition for the groupoid structure:
\begin{itemize}
    \item The source and target maps are fibrations, thus making the groupoid manifold into a doubly fibred bundle $s,t:G\to M$.
    \item The inversion $i$ and the left and right action maps $L_g,R_g$ become diffeomorphisms. 
    \item The isotropy groups are Lie groups sitting as embedded submanifolds in $G_x\subset G$.
    \item Orbits are immersed submanifolds $O_x\subset M$.
    \item The identity map is an embedding $u:M\to G$.
    \item The fibres restrict to orbits as principal $G_x$-bundles, with the bundle projections $t:s^{-1}(x)\to O_x$ or $s:t^{-1}(x)\to O_x$.
\end{itemize}

A morphism between two general groupoids $F:G\to H$ is simply a functor of the categorical structures. At the set-theoretic level, this is equivalent to a morphism of doubly fibred sets, i.e. a pair of maps $F:G\to H$ and $\varphi:M\to N$ such that the two following diagrams commute
\begin{equation*}
\begin{tikzcd}
G  \arrow[r, "F"]\arrow[d, "s_G"'] & H \arrow[d, "s_H"]\\
M \arrow[r,"\varphi"']  & N
\end{tikzcd}
\qquad
\begin{tikzcd}
G  \arrow[r, "F"] \arrow[d, "t_G"'] & H \arrow[d, "t_H"]\\
M \arrow[r,"\varphi"']  & N
\end{tikzcd}
\end{equation*}
that is compatible with the rest of structure maps in the following sense:
\begin{equation*}
    F(m_G(g,h))=m_H(F(g),F(h)), \qquad F\circ u_G = u_H\circ \varphi, \qquad F\circ i_G= i_H \circ F,
\end{equation*}
for all composable $g,h\in G$. A \textbf{morphism of Lie groupoids} is a functor between two Lie groupoids $F:G\to H$ such that the maps induced between the groupoid and base manifolds are smooth.\newline

The generalization of the notion of section in  multifibred manifolds is that of a map that splits the projections up to bijection. We then naturally define the \textbf{bisections} of a groupoid $s,t:G\to M$ as:
\begin{equation*}
    \Sec{G}:=\{b:M\to G|\quad \tau_b:=t\circ b:M\to M,\quad \sigma_b:=s\circ b:M\to M \text{ bijections }\}.
\end{equation*}
The identity map $u:M\to G$ gives a trivial example of bisection. We can exploit the groupoid structure on $G$ and define a point-wise operation of bisections $\cdot :\Sec{G}\times \Sec{G}\to \Sec{G}$ via
\begin{equation*}
    (a\cdot b)(x):=m(a(\sigma_a^{-1}(x)),b(\tau_b^{-1}(x))).
\end{equation*}
Note that the inversion induces a map of bisections $i:\Sec{G}\to \Sec{G}$ by setting $i(a)(x)=i(a(x))$. It is then straightforward to check that the bisections form a  group with the point-wise multiplication, with identity given by the identity map and with inverses obtained by the induced inversion map. We call $(\Sec{G},\cdot)$ the \textbf{bisection group} of the groupoid $G$. When $G$ is a Lie groupoid, bisections are required to be smooth and to split the source and target maps up to diffeomorphism. In this case, the bisection group $(\Sec{G},\cdot)$ becomes a generically infinite-dimensional manifold.\newline

Lie groupoids should be regarded as a unifying generalization of Lie groups, smooth manifolds and Lie group actions, as will become clear in the following examples. A Lie groupoid over a point manifold $p,p:G\to\{x\}$ is simply a Lie group and all the group structures present in a general groupoid collapse to give a single one $G_x\cong G\cong \Sec{G}$. A smooth manifold $M$ can be trivially made into the \textbf{pair groupoid} $\Proj_1,\Proj_2:M\times M\to M$ with partial multiplication given by $(x,y)\circ (y,z)=(x,z)$. More generally, the \textbf{fundamental groupoid} $\Pi_1(M)$ is defined as the manifold of smooth paths up to homotopy with concatenation of paths as composition, when $M$ is simply connected it is easy to see that $\Pi_1(M)\cong M\times M$. The isotropies of the fundamental group are the first homotopy groups with a base point and the bisection group is isomorphic to the diffeomorphism group of the base. Given an action of a Lie group on a manifold $\phi:G\times M\to M$ we can define the \textbf{action groupoid} $s,t:G\ltimes M\to M$ as the Cartesian product $G\ltimes M:=G\times M$ with source map $s(g,x)=x$, target map $t(g,x)=\phi(g,x)$ and multiplication $(h,y)\circ (g,x)=(hg,x)$, only defined when $y=\phi(g,x)$. It clearly follows from this definition that the usual notions of isotropy group and orbit for a group action coincide with the groupoid isotropy and orbits. The bisections of the action groupoid are naturally identified with the $G$-equivariant $G$-valued functions on $M$, $\Sec{G\ltimes M}\cong \Cin{M,G}^G$. Lastly, as an example of particular relevance in this thesis, given a vector bundle $\epsilon:E\to M$ we define the \textbf{general linear groupoid} as the manifold of all linear isomorphisms between fibres:
\begin{equation*}
    \GL{E}:=\{B_{xy}:E_x\to E_y \text{ isomorphism }, x,y\in M\}
\end{equation*}
The smooth structure of this groupoid is induced from the bundle structure of $E$ since a local trivialization $E|_U\cong U\times \Real^k$ gives a diffeomorphism $s^{-1}(U)\cong U \times \GL{k,\Real}\times U$. The structure maps $s,t:\GL{E}\to M$ are given in the obvious way and the partial groupoid multiplication is given by composition. The isotropy groups are simply the general linear groups of each fibre $\GL{E}_x\cong \GL{E_x}$. Since linear isomorphisms always exist between vector spaces of the same dimension, the general linear groupoid has its base a the single orbit, this is an example of what is generally called a \textbf{transitive groupoid}. It is clear by construction that there is a one-to-one correspondence between bisections of the general linear groupoid and the automorphisms of the vector bundle covering diffeomorphisms, it is then easy to see that the correspondence is, in fact, a group isomorphism
\begin{equation*}
    (\Sec{\GL{E}},\cdot)\cong (\text{Aut}{E},\circ).
\end{equation*}
This exemplifies a general feature of groupoids: they are finite-dimensional manifolds capturing the symmetry transformations of spaces which form infinite-dimensional groups with composition, particularly in the presence of a fibration.\newline

There is a natural generalization of the notion of action for groupoids: let $N$ be a manifold and $s,t:G\to M$ be a Lie groupoid, a (left) \textbf{groupoid action} of $G$ on $N$ is given by a pair of smooth maps $\phi:G\times N\to N$ and $p:N\to M$ such that
\begin{equation*}
    p(\phi(g,x))=t(g), \quad \phi(g,\phi(h,x))=\phi(g\circ h,x),\quad \phi(u(p(x)),x)=x,
\end{equation*}
for all composable $g,h\in G$ and $x\in N$. We will denote $G\Acts N$ and $N$ will be called a $G$-space. Note that for the case of $M$ being a single point the above definition recovers the usual notion of action of a Lie group $G$ on a manifold $N$. For any groupoid $s,t:G\to M$, an important family of a groupoid actions is when the $G$-space is a vector bundle over the same base $\epsilon: E\to M$ and $G$ acts via fibre-to-fibre linear isomorphisms. We call such an action $G\Acts E$ a \textbf{groupoid representation} and note that it is equivalently given by a groupoid morphism
\begin{equation*}
    \begin{tikzcd}[sep=small]
    G  \arrow[rr,"R"]\arrow[dr,shift right] \arrow[dr,shift left] & & \GL{E} \arrow[dl,shift right] \arrow[dl,shift left] \\
    & M & 
    \end{tikzcd}
\end{equation*}

In direct analogy with the invariant vector fields on a Lie group, we can identify the special subspace of vector fields on a Lie groupoid that are compatible with the double fibration and the partial multiplication structures. Let $s,t: G \to M$ be a Lie groupoid, then we see that the $s$-fibration defines a distribution spanned by the vector fields tangent to the $s$-fibres which we denote by $\Tan^s G:=\Ker{\Tan s}$. Composable elements of the groupoid act via diffeomorphisms on the $s$-fibres and thus also on tangent vector fields via the push-forward of the differential of the right action. It is then natural to define the \textbf{right-invariant vector fields} on the groupoid $G$ as
\begin{equation*}
    \Sec{\Tan G}^R:=\{X\in\Sec{\Tan^s G}|\quad X(hg)=\Tan_hR_g(X(h))\quad \forall(h,g)\in G\ftimes{s}{t} G\},
\end{equation*}
which are closed under the bracket of vector fields on $G$. The fact that elements of the groupoid act via push-forwards when restricted to the $s$-fibres implies that the values of a right-invariant vector field are fibre-wise determined from tangent vectors at the identity section $u(M)\subset G$. By defining $A_G:=\Tan^s G|_{u(M)}$ as a vector bundle over $M$, this observation shows that there is a isomorphism of $\Cin{M}$-modules
\begin{equation*}
    \Sec{G}^R\cong \Sec{A_G}
\end{equation*}
given by the correspondence $a\in\Sec{A_G}\mapsto a^R$ such that
\begin{equation*}
    a^R(g)=\Tan_{u(t(g))}R_g(a(t(g))).
\end{equation*}
The $\Cin{M}$-module structure is mapped via pull-back with the $s$ projection and thus, for $a,b\in\Sec{A_G}$ and $f\in\Cin{M}$, we can readily check
\begin{equation*}
    [a^R,s^*f\cdot b^R]=s^*f\cdot [a^R,b^R]+(\Tan t|_{\Tan^s G} a)[s^*f]\cdot b^R.
\end{equation*}
This shows that $A_G$ has the structure of a Lie algebroid $(A_G,\rho:=\Tan t|_{A_G},[,])$ with the anchor given by the tangent of the $t$-projection $\rho:=\Tan t|_{A_G}$ and Lie bracket induced from the bracket of right-invariant vector fields. We call $A_G$ the \textbf{Lie algebroid of the Lie groupoid $G$}. Had we chosen the distribution tangent to $t$-fibres and left-invariant vector fields, an entirely analogous construction will yield another Lie algebroid $A'_G$. These two constructions are, nonetheless, equivalent since the inversion diffeomorphism $i:G\to G$ induces an (anti-)isomorphism of Lie algebroids $\Tan i: A_G\to A_G'$ covering the identity on $M$. It follows structurally from this definition that the isotropy Lie algebras of $A_G$ are precisely the Lie algebras of the isotropy Lie groups $\mathfrak{g}_x\cong\Tan_eG_x$ and that the leaves of the singular foliation integrating the characteristic distribution of $A_G$ coincide with the Lie groupoid orbits $\Tan_xO_x=\rho((A_G)_x)$, for all $x\in M$.\newline

We recover all the basic examples of Lie algebroids from the examples of Lie groupoids given above. The Lie algebroid of a Lie group $p,p:G\to \{x\}$ regarded as Lie groupoid is simply the Lie algebra of the Lie group regarded as a vector bundle over a point $\mathfrak{g}\to \{x\}$. The Lie algebroid of the pair groupoid of a manifold $M$ is simply the tangent bundle $\Tan M$ with the identity anchor map. The Lie algebroid of an action groupoid $G\ltimes M$ is the action Lie algebroid $\mathfrak{g}\times M\to M$ with anchor map given by the infinitesimal action. Of particular relevance for this thesis is the case of the \textbf{Lie algebroid of the general linear groupoid} of a vector bundle $\GL{E}$, noting that elements of the groupoid naturally act on local sections of $E$ it is easy to show that there is an isomorphism of Lie algebroids
\begin{equation*}
    A_{\GL{E}}\cong \Der E.
\end{equation*}
What makes this isomorphism remarkable is that it allows us to regard the derivations of the vector bundle $\Dr{E}$, an infinite-dimensional space, as vector fields on a finite-dimensional manifold:
\begin{equation*}
    \Dr{E}\cong \Sec{A_{\GL{E}}}\cong\Sec{\Tan \GL{E}}^R\subset \Sec{\Tan \GL{E}}.
\end{equation*}
This illustrates yet another general point about groupoids: they are the spaces where the question of whether an infinitesimal set of symmetries, particularly in the presence of a fibration, has a global counterpart can be formulated in the conventional sense of integrability of vector fields. Directly connected to this is the integrability of Lie algebroids, which we discuss at the end of this section.\newline

Let $F:G\to H$ be a morphism of Lie groupoids covering the smooth map $\varphi:M\to N$, then it follows from the compatibility of $F$ and $\varphi$ with the structure maps of $G$ that $\Tan F|_{A_G}:A_G\to A_H$ is a well-defined map compatible with the Lie algebroid structure. We call this restriction of the tangent of a morphism of Lie groups $F:G\to H$ the \textbf{induced Lie algebroid morphism} and we denote it by $\Tan^sF:A_G\to A_H$. Recall that a representation of a groupoid $s,t:G\to M$ on a vector bundle $\epsilon:E\to M$ is characterized by a groupoid morphism $R:G\to \GL{E}$, then the induced Lie algebroid morphism $\Tan^s R:A_G\to \Der E$ is easily shown to give a representation of the Lie algebroid $A_G$ on the vector bundle $E$.\newline

Let $\alpha:A\to M$ be a Lie algebroid and $s,t:G\to M$ a Lie groupoid, we say that $G$ \textbf{integrates} $A$ if there is an isomorphism of Lie algebroids $A\cong A_G$. If given a Lie algebroid $A$ there exists at least one Lie groupoid that integrates it was say that $A$ is \textbf{integrable}. A Lie groupoid $s,t:G\to M$ is called $s$-simply connected if the $s$-fibres are all simply connected submanifolds. At this point the natural problem of \textbf{integrability} appears: given some data for Lie algebroids, can we find corresponding data on Lie groupoids integrating them? The known results about this problem are summarized below mirroring the classical theorems of Lie on the integrability of Lie algebras.
\begin{itemize}
    \item \textbf{Lie I.} \cite[Theorem 2.1]{crainic2003integrability} If a Lie algebroid $A$ is integrable then there exists a unique $s$-simply connected groupoid $G_A$ integrating it.
    \item \textbf{Lie II.} \cite[Proposition 2.2]{crainic2003integrability} Let $L:A\to B$ be a Lie algebroid morphism with $A$ and $B$ integrable Lie algebroids, then, if $H$ is a Lie groupoid integrating $B$, there exists a unique Lie groupoid morphism $F:G_A\to H$ such that $\Tan^sF=L$.
    \item \textbf{Lie III.} \cite[Theorem 4.1]{crainic2003integrability} In contrast with the case of Lie algebras, not all Lie algebroids are integrable. Integrability can be characterized by topological obstructions defined from the interaction of the foliation integrating the characteristic distribution of the Lie algebroid and the distribution of dimensions of the isotropy Lie algebras.
\end{itemize}

\section{Poisson, Presymplectic and Symplectic Geometry} \label{SymplecticGeometry}

A \textbf{Poisson algebra} is a triple $(A,\cdot,\{,\})$ with $A$ a $\mathbb{R}$-vector space, $(A,\cdot)$ a commutative (unital, associative) algebra and $(A,\{,\})$ a Lie algebra such that the following \textbf{Leibniz identity} holds
\begin{equation*}
    \{a,b\cdot c\}=\{a,b\}\cdot c + b\cdot \{a,c\}.
\end{equation*}
Equivalently, this condition can be rephrased as
\begin{equation*}
    \text{ad}_{\{,\}}:A\to \Dr{A,\cdot}
\end{equation*}
making the adjoint map of the Lie algebra a morphism of Lie algebras from $A$ to the derivations of both bilinear multiplications of the algebra, $\Dr{A,\{,\}}\cap \Dr{A,\cdot}$. Any element in this subspace of endomorphisms of $A$, i.e. both a derivation of the commutative product and the Lie bracket, will be called a \textbf{Poisson derivation}. In particular, the Poisson derivations given by the adjoint map of the Lie bracket are usually called \textbf{Hamiltonian derivations} and are denoted by $X_a:=\{a,-\}$. A linear map $\psi:A\to B$ is called a \textbf{morphism of Poisson algebras} if $\psi:(A,\cdot)\to (B,\cdot)$ is a commutative algebra morphism and $\psi:(A,\{,\})\to (B,\{,\})$ is a Lie algebra morphism. A subspace $I\subset A$ is called a \textbf{coistrope} if it is a multiplicative ideal and a Lie subalgebra.\newline

A \textbf{Poisson manifold} is a smooth manifold $P$ whose commutative algebra of smooth functions has the structure of a Poisson algebra $(\Cin{P},\cdot,\{,\})$, we also use the term \textbf{Poisson structure} to refer to this datum. A \textbf{Poisson map} is a smooth map between Poisson manifolds whose pull-back on functions gives a morphism of Poisson algebras. The Poisson derivations of $(\Cin{P},\cdot,\{,\})$ are called \textbf{Poisson vector fields} and correspond to infinitesimal Poisson diffeomorphisms. Hamiltonian derivations are called \textbf{Hamiltonian vector fields} in this context and we denote the adjoint map by $f\in\Cin{P}\mapsto X_f=\{f,-\}$. It follows directly by construction that the assignment of Hamiltonian vector fields is a $\mathbb{R}$-linear Lie algebra morphism,
\begin{equation*}
    X_{\{f,g\}}=[X_f,X_g].
\end{equation*}
The standard result that derivations on functions are isomorphic to vector fields allows for the more geometric definition of a Poisson manifold as a pair $(P,\pi)$ with $\pi\in\Sec{\wedge^2\Tan P}$ a bivector field with vanishing Schouten bracket $\llbracket \pi, \pi\rrbracket=0$. This condition is equivalent to the Jacobi identity of the natural bracket induced in on functions
\begin{equation*}
    \{f,g\}_\pi:=\pi(df,dg),
\end{equation*}
which satisfies the Leibniz rule automatically by construction. Such $\pi\in\Sec{\wedge^2\Tan P}$ is called a \textbf{Poisson bivector}. The bivector gives a musical map
\begin{equation*}
    \pi^\sharp:\Cot P\to \Tan P
\end{equation*}
which is a vector bundle morphism covering the identity. The Hamiltonian vector field of a function $f\in\Cin{P}$ can now be given explicitly as
\begin{equation*}
    X_f=\pi^\sharp(df).
\end{equation*}
The map $X_{-}=\pi^\sharp\circ d:\Cin{P}\to\Sec{\Tan P}$ will sometimes br referred to as the \textbf{Hamiltonian map} of the Poisson manifold $(P,\pi)$. The tangent distribution spanned by all the Hamiltonian vector fields at every point $\pi^\sharp(\Tan^* P)\subset \Tan P$ is called the \textbf{Hamiltonian distribution} of the Poisson manifold $(P,\pi)$. That the Hamiltonian distribution $\pi^\sharp(\Tan^* P)$ is involutive follows simply from the fact that the Hamiltonian map is a Lie algebra morphism.\newline

A family of Poisson structures of particular relevance in this thesis are those present on the total space of vector bundles. Let $\epsilon:E\to M$ be a vector bundle with the submodule of spanning functions $\text{C}_s^\infty(E)\cong \epsilon^*\Cin{M}\oplus l_{\Sec{E^*}}$, a Poisson bracket $(\Cin{E},\cdot, \{,\}$) is called a \textbf{spanning Poisson structure} on $E$ if
\begin{equation*}
    \{\text{C}_s^\infty(E),\text{C}_s^\infty(E)\}\subset \text{C}_s^\infty(E). 
\end{equation*}
More specifically, a \textbf{(fibre-wise) linear Poisson structure} on $E$ is one satisfying
\begin{align*}
    \{l_{\Sec{E^*}},l_{\Sec{E^*}}\} & \subset  l_{\Sec{E^*}}\\
    \{l_{\Sec{E^*}},\epsilon^*(\Cin{M})\} & \subset  \epsilon^*(\Cin{M})\\
    \{\epsilon^*(\Cin{M}),\epsilon^*(\Cin{M})\} & = 0.
\end{align*}
Note that the local nature of Poisson brackets, i.e. the fact that they are equivalently determined by a bivector field, implies that spanning Poisson structures are completely specified by their brackets on spanning functions, since from proposition \ref{SpanningFunctions} we see that spanning functions have differentials spanning cotangent spaces everywhere and thus fully specifying the bivector field.\newline

\begin{prop}[Lie Algebroids and Linear Poisson Structures] \label{LieAlgebroidLinearPoisson}
Let $(\alpha: A\to M,\rho,[,])$ be a Lie algebroid, then there exists a unique linear Poisson structure on the dual bundle $\normalfont (\Cin{A^*},\{,\})$ satisfying
\begin{align*}
    \{l_a,l_b\} & =  l_{[a,b]}\\
    \{l_a,\alpha^*f\} & =  \alpha^*\rho(a)[f].
\end{align*}
Furthermore, if $F:A_1\to A_2$ is a vector bundle morphism covering a diffeomorphism $\varphi:M_1\to M_2$, so that there is a well-defined pull-back bundle morphism $F^*:A_2^*\to A_1^*$, then $F$ is a Lie algebroid morphism iff $F^*$ is a linear Poisson map. This establishes a one-to-one correspondence between Lie algebroids and linear Poisson structures on vector bundles.
\end{prop}
\begin{proof}
Given a Lie algebroid $(\alpha:A\to M,\rho,[,])$ it is a routine computation to check that the following is a well-defined Poisson bracket on the dual vector bundle:
\begin{align*}
\{l_{a},l_{b}\} & = l_{[a,b]}\\
\{l_{a},\alpha^*f\} & = \alpha^*\rho(a)[f]\\
\{\alpha^*f,\alpha^*g\} & = 0,
\end{align*}
which is clearly linear. The converse construction, i.e. giving a Lie algebroid structure on $A$ when $A^*$ has linear Poisson structure, requires a bit more work. In order to construct the anchor map note that the Leibniz property of the linear Poisson bracket on $A^*$ implies
\begin{equation}
\alpha^*r(fa,g):=\{l_{f\cdot a},\alpha^*g\}=\{l_a,\alpha^*(fg)\}=\alpha^*f\alpha^*r(a,g).
\end{equation}
Then, since $\alpha^*$ is injective, we see that $r(a,f)\in\Cin{M}$ has to be of the form $\rho(a)[f]$ for some $\Cin{M}$-linear map $\rho$, thus giving the anchor. The fact that the natural inclusion map $l:\Sec{A}\to \Cin{A^*}$ is injective implies that the equation $l_{[a,b]}=\{l_a,l_b\}$ determines the bracket completely. The only thing left to verify is the Leibniz identity, which we check directly using the given definition of the anchor as a $\Cin{M}$-linear map:
\begin{equation*}
l_{[a,f\cdot b]} = \alpha^*f\{l_a,l_b\}+\{l_a,\alpha^*f\}l_b = l_{f\cdot[a,b]+\rho(a)[f]\cdot b}
\end{equation*}
where, again, we used the fact that $l$ is injective. The two constructions are trivially checked to be inverses of each other. The equivalence between Lie algebroid morphisms covering a diffeomorphism and linear Poisson maps easily follows from the observation that
\begin{equation*}
l_a\circ F^* = l_{F_*a}, \quad p^*f\circ F^* = p^*(f\circ \varphi^{-1})
\end{equation*}
for all $f\in\Cin(M)$ and $a\in\Sec{A}$, and the fact that the push-forward between Lie algebroid sections is a Lie algebra morphism.
\end{proof}

A submanifold $i:C\hookrightarrow P$ is called a \textbf{coisotropic submanifold} of the Poisson manifold $(P,\pi)$ if $\Tan_xC\subset (\Tan_xP,\pi_x)$ is a coisotropic subspace for all $x\in C$, that is
\begin{equation*}
    \pi_x^\sharp (T_xC^0)\subset T_xC.
\end{equation*}
Denoting by $I_C:=\Ker{i^*}\subset \Cin{P}$ the vanishing ideal of the submanifold $C$ we formulate the following result giving alternative equivalent definitions of coisotropic submanifolds.
\begin{prop}[Coisotropic Submanifolds of Poisson Manifolds, {\cite[Proposition 2.32]{fernandes2014lectures}}] \label{CoisotropicSubPoisson}
Let $i:C\hookrightarrow P$ be a closed submanifold of the Poisson manifold $(P,\pi)$, then the following are equivalent:
\begin{enumerate}
    \item $C$ is coisotropic,
    \item $I_C$ is a coisotrope in the Poisson algebra $(\Cin{P},\cdot,\{,\}_\pi)$,
    \item for all $c\in I_C$ the Hamiltonian vector field $X_c$ is tangent to $\normalfont C$: $X_c|_C\in \Sec{\Tan C}$.
\end{enumerate}
\end{prop}
Let two Poisson manifolds $(P_1,\pi_1)$ and $(P_2,\pi_2)$, we define the \textbf{product Poisson manifold} simply from the Cartesian product as $(P_1\times P_2,\pi_1+\pi_2)$, where the sum of bivectors is defined via the canonical isomorphism $\Tan(P_1\times P_2)\cong \Proj_1^*\Tan P_1\oplus \Proj_2^*\Tan P_2$. Equivalently, the product Poisson manifold structure on $P_1\times P_2$ is the Poisson bracket $\{,\}_{12}$ determined uniquely by the condition that the canonical projections
\begin{equation*}
\begin{tikzcd}
P_1 & P_1\times P_2 \arrow[l,"\Proj_1"']\arrow[r,"\Proj_2"] & P_2.
\end{tikzcd}
\end{equation*}
are Poisson maps, that is
\begin{equation*}
    \Proj_1^*\{f_1,g_1\}_1=\{\Proj_1^*f_1,\Proj_1^*g_1\}_{12} \qquad \Proj_1^*\{f_2,g_2\}_2=\{\Proj_2^*f_2,\Proj_2^*g_2\}_{12}
\end{equation*}
for all $f_1,g_1\in\Cin{P_1}$ and $f_2,g_2\in\Cin{P_2}$. For a given Poisson manifold $(P,\pi)$ let us denote the \textbf{opposite} Poisson manifold by $\overline{P}:=(P,-\pi)$. The proposition below gives a characterization of Poisson maps as \textbf{coisotropic relations}, defined as coisotropic submanifolds of the product Poisson manifold.
\begin{prop}[Product of Poisson Manifolds, {\cite[Proposition 2.35]{fernandes2014lectures}}] \label{ProductPoisson}
 Let two Poisson manifolds $(P_1,\pi_1)$ and $(P_2,\pi_2)$, and a smooth map $\varphi:P_1 \to P_2$, then $\varphi$ is a Poisson map iff its graph
\begin{equation*}
\normalfont
    \Grph{\varphi}\subset P_1\times \overline{P_2}
\end{equation*}
is a coisotropic submanifold in the product Poisson manifold.
\end{prop}

The general question of whether a submanifold of some manifold carrying a geometric structure inherits a similar structure, or to what degree it fails to do so, is generally known as the problem of \textbf{reduction}. In the context of Poisson geometry, a triple $((P,\pi),C,D)$, where $(P,\pi)$ is a Poisson manifold, $i:C\hookrightarrow M$ is an immersed submanifold and $D\subset \Tan P$ is a tangent distribution, will be called a \textbf{reduction triple} if they satisfy:
\begin{itemize}
\item $\Tan C\cap D$ is and integrable subbundle of $\Tan C$ that defines a foliation $\mathcal{F}$ on $C$,
\item the foliation $\mathcal{F}$ is regular, making the leaf space into a smooth manifold with submersion $q:C\to P':=C/\mathcal{F}$,
\item $D$ preserves the Poisson structure, that is, $\forall df,dg\in D^0$ $\Rightarrow$ $d\pi(df,dg)\in D^0$.
\end{itemize}
The triple is said to be \textbf{reducible} when the reduced manifold inherits a Poisson structure $(P',\pi')$ such that:
\begin{equation*}
q^*\{f,g\}'=i^*\{F,G\}
\end{equation*}
for all $f,g\in \Cin{P'}$ and $F,G\in \Cin{P}$ extensions of $q^*f,q^*g\in \Cin{C}$ satisfying $dF,dG\in D^0$.\newline

It follows from our discussion above that any coisotropic submanifold $C$ carries an involutive tangent distribution given by the images of the Hamiltonian vector fields of the (locally generating) elements of the vanishing ideal $I_C$. Let us denote this distribution by $X_{I_C}\subset \Tan C$ and note that $[X_{I_C},X_{I_C}]\subset X_{I_C}$ from the fact that $I_C$ is a coisotrope. This means that there will be a (singular) foliation in $C$, denoted by $\mathcal{X}_C$, integrating the tangent distribution of Hamiltonian vector fields of vanishing functions. The following result shows that, under some regularity assumptions, coisotropic manifolds provide reducible triples $((P,\pi),C,X_{I_C})$

\begin{prop}[Coisotropic Reduction of Poisson Manifolds, see \cite{marsden1986reduction}] \label{CoisotropicRedPoisson}
Let $(P,\pi)$ be a Poisson manifold and $i:C\hookrightarrow P$ a closed coisotropic submanifold. Assume that the tangent distribution $X_{I_C}$ integrates to a regular foliation with smooth leaf space $P:=C/\mathcal{X}_C$ in such a way that there is a surjective submersion $q$ fitting in the reduction diagram:
\begin{equation*}
    \begin{tikzcd}[sep=small]
        C \arrow[r, hookrightarrow, "i"] \arrow[d, twoheadrightarrow,"q"'] & (P,\pi) \\
        (P',\pi')
    \end{tikzcd}
    \end{equation*}
Then, the manifold $P'$ inherits a Poisson bracket on functions $(\normalfont \Cin{P'},\{,\}')$ uniquely determined by the condition
\begin{equation*}
\normalfont
    q^*\{f,g\}'=i^*\{F,G\} \quad \forall f,g\in\Cin{P'}
\end{equation*}
and for all $\normalfont F,G\in\Cin{P}$ leaf-wise constant extensions, i.e functions of the ambient Poisson manifold satisfying
\begin{equation*}
    q^*f=i^*F, \quad q^*g=i^*G.
\end{equation*}
\end{prop}
As a particular case of coisotropic reduction we find \textbf{Poisson submanifolds}, i.e. submanifolds $i:C\hookrightarrow P$ whose vanishing ideal $I_C$ is not only a Lie subalgebra but a Lie ideal, $\{I_C,\Cin{P}\}\subset I_C$. In this case the distribution of Hamiltonian vector fields clearly vanishes $X_{I_C}=0\subset \Tan C$ thus satisfying all the requirements of the proposition trivially and giving the reduced Poisson manifold simply as the submanifold itself carrying a Poisson structure $(\Cin{C},\{,\}')$. Another example of coistropic reduction is given by the presence of a \textbf{Hamiltonian group action}: a Lie group action on a Poisson manifold $G \Acts (P,\pi)$ via Poisson maps whose infinitesimal generators $\psi:\mathfrak{g}\to \Sec{\Tan P}$ given by Hamiltonian vector fields via a (co)moment map $\overline{\mu}:\mathfrak{g}\to \Cin{M}$ with the defining conditions:
\begin{equation*}
    \psi(\xi)=X_{\overline{\mu}(\xi)}, \qquad \{\overline{\mu}(\xi),\overline{\mu}(\zeta)\}_\pi=\overline{\mu}([\xi,\zeta]) \qquad \forall \xi,\zeta\in \mathfrak{g}.
\end{equation*}
The functional dual map $\mu:P\to \mathfrak{g}^*$ is called the \textbf{moment map}. When $0\in\mathfrak{g}^*$ is a regular value, $\mu^{-1}(0)\subset P$ is easily checked to be a coisotropic submanifold. If the Poisson action restricts to a free and proper action on $\mu^{-1}(0)$ then the orbit space $\mu^{-1}(0)/G$ is the coisotropic reduction of $(P,\pi)$. The procedure by which one obtains a new Poisson manifold $P'$ from a Hamiltonian group action $G\Acts P$ is called \textbf{Hamiltonian reduction}.\newline

A \textbf{presymplectic manifold} is a pair $(S,\omega)$ with $S$ a smooth manifold and $\omega\in\Sec{\wedge^2\Cot S}$ a closed 2-form, $d\omega=0$. We will also say that $S$ carries a \textbf{presymplectic structure} $\omega$. If the presymplectic form is exact, $\omega=d\theta$, we say that $(S,\theta)$ is an \textbf{exact presymplectic structure}. A smooth map $\varphi:S_1\to S_2$ between presymplectic manifolds is called a \textbf{presymplectic map} when $\varphi^*\omega_2=\omega_1$. Infinitesimal presymplectic diffeomorphisms then give the corresponding notion of \textbf{presymplectic vector field} as those $X\in\Sec{\Tan S}$ preserving the presymplectic form, $\LDer_X\omega=0$. The presymplectic form induces a musical map $\omega^\flat:\Tan S\to \Cot S$, the kernel of this vector bundle morphism is called the \textbf{characteristic distribution} of the presymplectic structure $\Ker{\omega^\flat}\subset \Tan S$. It follows from $d\omega=0$ that $\Ker{\omega^\flat}$ is an involutive tangent distribution. Let $f\in\Cin{S}$ be a function, we say that $X_f\in\Sec{\Tan S}$ is a \textbf{Hamiltonian vector field} for the function $f$ if $\omega^\flat(X_f)=df$. Hamiltonian vector fields need not exist for arbitrary functions and, when they do, they are not unique: adding any vector field tangent to the characteristic distribution will leave the Hamiltonian condition unchanged by definition. The subspace of all functions that admit Hamiltonian vector fields is called the space \textbf{admissible functions}:
\begin{equation*}
    \Cin{S}_\omega:=\{f\in\Cin{S}|\quad \exists \quad X_f\in\Sec{\Tan S}\}.
\end{equation*}
It follows by definition that $X_{fg}=f\cdot X_g+g\cdot X_f$, hence the admissible functions form, in fact, a subring $\Cin{S}_\omega\subset \Cin{S}$. Furthermore, we can define a bracket on admissible functions $f,g\in\Cin{S}_\omega$ by setting
\begin{equation*}
    \{f,g\}_\omega:=\omega(X_f,X_g)=X_f[g]=-X_g[f].
\end{equation*}
This bracket is clearly antisymmetric and satisfies the Leibniz identity by construction; a quick computation shows that the Jacobi identity of $\{,\}_\omega$ is equivalent to the presymplectic form being closed, $d\omega=0$. Note again that this bracket is well-defined because Hamiltonian vector fields are unique up to vector fields tangent to the kernel of the presymplectic form. We thus see that the admissible functions form a Poisson algebra $(\Cin{S}_\omega,\cdot, \{,\}_\omega)$.\newline

A submanifold of a presymplectic manifold $i:C \hookrightarrow S$ is called \textbf{isotropic} when the presymplectic form restricts to the zero form, i.e. when $i^*\omega=0$. An equivalent way to characterize isotropic manifolds is as the integral manifolds of the characteristic distribution.\newline

We will only consider the problem of reduction of presymplectic manifolds in the presence of symmetry. A Lie group action on a presymplectic manifold $\phi:G\times S\to S$ is called a \textbf{presymplectic action} if $\phi_g^*\omega=\omega$ for all $g\in G$. A \textbf{comoment map} for the presymplectic action is a linear map $\overline{\mu}:\mathfrak{g}\to \Cin{S}_\omega$ satisfying
\begin{equation*}
    \omega^\flat(\psi(\xi))=d\overline{\mu}(\xi), \qquad \{\overline{\mu}(\xi),\overline{\mu}(\zeta)\}_\omega=\overline{\mu}([\xi,\zeta]) \qquad \forall \xi,\zeta\in \mathfrak{g},
\end{equation*}
where $\psi:\mathfrak{g}\to \Sec{\Tan S}$ denotes the infinitesimal action. The functional dual of a comoment map $\mu:S\to \mathfrak{g}^*$ is easily checked to be $G$-equivariant with respect to the coadjoint action on $\mathfrak{g}^*$, $\mu$ is called a \textbf{moment map} for the presymplectic action.

\begin{prop}[Presymplectic Reduction, {\cite[Theorem 2]{echeverria1999reduction}}] \label{PresymplecticReduction}
Let $G\Acts (S,\omega)$ be a presymplectic action by a connected Lie group with moment map $\mu:S\to \mathfrak{g}^*$. Assume that $0\in\mathfrak{g}^*$ is a regular value so that $i:C:=\mu^{-1}(0)\hookrightarrow S$ is an embedded submanifold. Further assume that the group action restricts to a free and proper action so that to the orbit space $S':=C/G$ is a smooth manifold and there is a surjective submersion $q:C\to S'$. Then there is a presymplectic structure on the orbit space $(S',\omega')$, uniquely determined by the condition $q^*\omega'=i^*\omega$, with characteristic distribution given by $\normalfont \Ker{\omega'}\cong \Ker{\omega}/\psi(\mathfrak{g})|_C$.
\end{prop}

Note that both a Poisson manifold $(P,\pi)$ and a presymplectic manifold $(S,\omega)$ carry distributions that are involutive but generically irregular. Let us assume they are regular. In the case of a Poisson manifold, it is easy to check that an integral manifold of the Hamiltonian distribution $C$ becomes a Poisson submanifold with a bivector that is non-degenerate $(C,\pi|_C)$. A non-degenerate bivector has invertible musical map which defines a 2-form $\pi^{\sharp\sharp}$ which is, furthermore, closed from the fact that $\pi$ is a Poisson bivector. This makes $(C,\pi^{\sharp\sharp})$ into a non-degenerate presymplectic manifold. In the case of a presymplectic manifold, the foliation $\mathcal{K}$ integrating the regular characteristic distribution makes the leaf space into a smooth manifold $S':=S/\mathcal{K}$ with a surjective submersion $q:S\to S'$. Since the leaf space is obtained by removing the kernel of the presymplectic form, $\omega$ descends to a non-degenerate presymplectic form $(S',\omega')$. It is clear by construction that the admissible functions are then precisely the pull-back of general functions on the base, $\Cin{S}_\omega=q^*\Cin{S'}$. Since $\Ker{\omega'}=0$ by construction, all the functions on the leaf space $S'$ are admissible, and hence we see that $S'$ has a non-degenerate Poisson structure.\newline

These last remarks motivate the definition of a \textbf{symplectic manifold} as a non-degenerate Poisson manifold $(M,\pi)$ or, equivalently, a non-degenerate presymplectic manifold $(M,\omega)$. For historical reasons, we will adopt the usual definition of symplectic manifold as non-degenerate presymplectic manifold. All the notions introduced above for Poisson and presymplectic manifolds apply in a general sense to symplectic manifolds in particular. The distinctive feature of a symplectic manifold is that the musical map now becomes an isomorphism of vector bundles
\begin{equation*}
        \begin{tikzcd}
        \Tan M \arrow[r, "\omega^\flat",yshift=0.7ex] & \Cot M \arrow[l,"\pi^\sharp",yshift=-0.7ex]
        \end{tikzcd},
\end{equation*}
that relates the presymplectic 2-form and the Poisson bivector isomorphically via the equation
\begin{equation*}
    \omega(X,Y)=\pi(\omega^\flat (X),\omega^\flat(Y)) \qquad \forall X,Y\in\Sec{\Tan M}.
\end{equation*}
It follows from a simple point-wise linear algebra argument that non-degenerate 2-forms can only exist on even dimensional manifolds and that smooth maps preserving them, what we call \textbf{symplectic maps}, need to be local diffeomorphisms. In a symplectic manifold all functions are admissible and Hamiltonian vector fields are unique, furthermore
\begin{equation*}
\pi(df,dg)=\{f,g\}=\omega(X_f,X_g) \qquad \forall f,g\in\Cin{M}.
\end{equation*}
The characteristic distribution of a symplectic manifold vanishes everywhere $\Ker{\omega^\flat}=0$ whereas the Hamiltonian distribution spans the tangent bundle everywhere $\pi^\sharp(\Tan M)=\Tan M$. A submanifold in a symplectic manifold that is both coisotropic and isotropic it is called \textbf{Lagrangian}. Lagrangian submanifolds are the key entities that allow for the definition of a category of symplectic manifolds: \textbf{the Weinstein category of symplectic manifolds} $\Symp_\Man$ has objects smooth manifolds with a symplectic form $(M,\omega)$ and morphisms given by Lagrangian relations $R:(M_1,\omega_1)\dashrightarrow (M_2,\omega_2)$, i.e. Lagrangian submanifolds $R\subset M_1\times \overline{M_2}$. Note that symplectic maps are recovered as Lagrangian relations by taking their graphs, which are easily shown to be isotropic on $M_1\times \overline{M_2}$ and are clearly of maximal dimension. The term category is used here with some looseness since composition may not be strictly well-defined between all a priori composable Lagrangian relations. This is due to the composition of relations involving intersections of submanifolds, which can generally produce non-smooth spaces. Smoothness of the composition of Lagrangian relations will be discussed in the particular examples of subcategories of the Weinstein symplectic category that we will encounter in later sections.

\section{Dirac Geometry} \label{DiracGeometry}

\subsection{Linear Courant Spaces} \label{LinearCourant}

Let $V$ be a real vector space, a \textbf{Courant space} is a triple $(C,\langle, \rangle,\rho)$ with $C$ a real vector space, $\langle, \rangle$ a non-degenerate bilinear form and $\rho : C \to V$ a linear map that we will call the \textbf{anchor}. Given a subspace $N\subset C$ we say that $N$ is \textbf{isotropic} when $N\subset N^\perp$ and \textbf{coisotropic} when $N^\perp \subset N$. Note that non-degeneracy of the bilinear form implies $(N^\perp)^\perp=N$ for all subspaces $N\subset C$. When the anchor map is surjective and its kernel ker$(\rho)\subset C$ is maximally isotropic, the triple $(C,\langle, \rangle,\rho)$ is called an \textbf{exact Courant space}. This condition is equivalent to the anchor map fitting in the following short exact sequence:
\begin{equation*}
\begin{tikzcd}
0 \arrow[r] & V^* \arrow[r, "j"] & C \arrow[r, "\rho"] & V \arrow[r]  & 0
\end{tikzcd}
\end{equation*}
with $j= \sharp \circ \rho^*$ and where $\sharp: V^* \to V$, $\flat:V\to V^*$ is the musical isomorphism induced by the bilinear form. Note that, by construction, the injective map of the exact sequence is isotropic in the sense that $\langle j(V^*),j(V^*)\rangle=0$. The vector space $C$ is said to extend or to be a Courant extension of the vector space $V$. The bilinear form of an exact Courant space must have split signature, thus maximally isotropic subspaces, called \textbf{Lagrangian} subspaces, satisfy $L^\perp = L$. Two Lagrangian subspaces $L,L'\subset C$ are called \textbf{transversal} when $L\cap L'=\{0\}$ and $L + L'=C$, if this is the case we clearly have $C\cong L\oplus L'$ as vector spaces and the bilinear form in $C$ allows for the identifications $L^*\cong L'$ and $(L')^*\cong L$. We will be focusing our attention on exact Courant spaces for the remainder of this section. The first natural example of exact Courant space is the \textbf{standard Courant space} $\mathbb{V}:= V\oplus V^*$ with natural anchor $\rho_{\mathbb{V}}=$pr$_1$ and bilinear form given by:
\begin{equation*}
\langle v \oplus \alpha , w \oplus \beta \rangle_{\mathbb{V}}:= \tfrac{1}{2}(i_v\beta +i_w\alpha).
\end{equation*}
Given a Courant space $(C,\langle, \rangle_C,\rho)$ we can form the \textbf{opposite Courant space} as $(\overline{C},\langle, \rangle_{\overline{C}},\rho_{\overline{C}}):=(C,-\langle, \rangle_C,\rho_C)$. From any two given Courant spaces $(C,\langle, \rangle_C,\rho_C)$ and $(D,\langle, \rangle_D,\rho_D)$, we define their \textbf{direct sum} as $(C \oplus D,\langle, \rangle_{C \oplus D},\rho_D \oplus \rho_C)$ with
\begin{equation*}
\langle c_1 \oplus d_1 , c_2 \oplus d_2 \rangle_{C\oplus D}:= \langle c_1, c_2\rangle_C + \langle d_1, d_2\rangle_D.
\end{equation*}
Also note that given a vector space $C$ equipped with a non-degenerate bilinear form of split signature $\langle , \rangle$, a choice of a Lagrangian subspace $L\subset C$ makes the triple $(C,\langle , \rangle,q)$ into an exact Courant space where the anchor is given by the quotient map $q: C \to C/L$.\newline

An splitting of the exact sequence of an exact Courant space, $\nabla:V\to C$, is called an \textbf{isotropic splitting} if it embeds $V$ as an isotropic subspace $\nabla(V)\subset (C,\langle,\rangle)$.
\begin{prop}[Isotropic Splittings of Exact Courant Spaces] \label{IsotropicSplittingLinearCourant}
Let $\begin{tikzcd}[cramped, sep=small]0 \arrow[r] & V^* \arrow[r, "j"] & C \arrow[r, "\rho"] & V \arrow[r]  & 0 \end{tikzcd}$ be an exact Courant space, then isotropic splittings $\nabla:V\to C$ always exist; furthermore two isotropic splittings differ by antisymmetric bilinear forms $\wedge^2V^*$.
\end{prop}
\begin{proof}
Note that we can always find right splittings $\nabla : V\to C$ for the short exact sequence of vector spaces. Assume $\nabla_0,\nabla:V\to C$ are splittings (not necessarily isotropic), then it is clear that we can express their difference as maps:
\begin{equation*}
\nabla - \nabla_0 = \gamma
\end{equation*}
with $\gamma: V \to\Ker{\rho}$. From exactness we have $\Ker{\rho}=j(V^*)$ and so we can rewrite $\gamma=j\circ b^\flat$ for some linear map $b^\flat: V\to V^*$ or, equivalently, $b\in V^* \otimes V^*$. This means that, once a splitting $\nabla_0$ is fixed, any other choice of splitting is tantamount to choosing a bilinear form so that $\nabla = \nabla_0 + j\circ b^\flat$. Consider an arbitrary splitting $\nabla_0$ and define a second splitting $\nabla = \nabla_0 + j\circ b^\flat_{\nabla_0}$ to be given by the bilinear form $b_{\nabla_0}(v,w)=-\tfrac{1}{2}\langle \nabla_0(v),\nabla_0(w)\rangle$. A simple computation yields
\begin{equation*}
\langle\nabla(v),\nabla(w)\rangle = b^\flat_{\nabla_0}(\rho(j(v))(w) -\tfrac{1}{2}\langle \nabla_0(v),\nabla_0(w)\rangle-\tfrac{1}{2}\langle \nabla_0(w),\nabla_0(v)\rangle+\langle \nabla_0(v),\nabla_0(w)\rangle = 0
\end{equation*}
where the first term vanishes from exactness ($\rho\circ j = 0$), thus making $\nabla$ into an isotropic splitting. This shows that isotropic splittings always exist and that any two isotropic splittings will be related by an antisymmetric bilinear form $b\in \wedge^2 V^*$.
\end{proof}

A subspace of an exact Courant space $L\subset C$ is called a \textbf{Dirac space} if it is Lagrangian and fits into an exact sequence
\begin{equation*}
\begin{tikzcd}
0 \arrow[r] & \overline{W} \arrow[r, "j"] & L \arrow[r, "\rho"] & W \arrow[r]  & 0
\end{tikzcd}
\end{equation*}
for some $W\subset V$ and $\overline{W}\subset V^*$. Let $\nabla:V\to C$ be an isotropic splitting, the image of $L$ under the anchor $\rho(L)=W$ is called the \textbf{range} of the Dirac space and the subspace $\Ker{L}:=\rho(\nabla(W)\cap L)\subset W$ is called the \textbf{kernel} of the Dirac space. Note that the definition of $\Ker{L}$ does not depend on the choice of isotropic splitting as any other choice will give a subspace $\overline{\nabla}(W)\cap L$ that spans $\nabla(W) \cap L$ except, possibly, for components in the kernel of $\rho$.

\begin{prop}[The 2-form of a Dirac Space] \label{LinearDirac2Form}
Let $\begin{tikzcd}[cramped, sep=small]0 \arrow[r] & V^* \arrow[r, "j"] & C \arrow[r, "\rho"] & V \arrow[r]  & 0 \end{tikzcd}$ be an exact Courant space and $L\subset C$ a Dirac structure, then $L$ uniquely determines a 2-form $\omega_L$ on $\rho(L)\subset V$ and $\normalfont\Ker{L}=\Ker{\omega_L^\flat}$. Conversely, for each choice of isotropic splitting $\nabla:V\to C$, a subspace $W\subset V$ equipped with a 2-form $\omega\in\wedge^2 W^*$ uniquely defines a Dirac structure $L_W^\nabla\subset C$.
\end{prop}
\begin{proof}
Let $w_i\in \rho(L)=W$ for $i=1,2$ and take two arbitrary extensions in the Dirac space $e_i+j(\epsilon_i)\in L$ with $\epsilon_i\in \overline{W}$ and $e_i\in C$ such that $\rho(e_i)=w_i$. Then define the 2-form via the identity:
\begin{equation*}
\omega_L(w_1,w_2):=\epsilon_1(w_2)=-\epsilon_2(w_1).
\end{equation*}
Note that the second equality holds due to $L\subset C$ being Lagrangian. This is well-defined since any other choice of extensions will give the same values: suppose that $l=e+j(\epsilon)$ and $l'=e'+j(\epsilon')$ both extend $w\in W$ so $\rho(e)=w=\rho(e')$, this means that $e,e'$ can be seen as images of isotropic splittings that differ by an antisymmetric bilinear form. Then, choosing covector representatives given by $\epsilon'=\epsilon+\delta$ with $j(\delta)\in L$ as $\nabla(0)+j(\delta)\in L$, we clearly see that $\epsilon'_1(w_2)=-\epsilon'_2(w_1)=\epsilon_1(w_2)=-\epsilon_2(w_1)$. Choosing an arbitrary isotropic splitting, we clearly see that the kernel of $\omega_L$ or, equivalently, the $\omega_L$-orthogonal of $W$, will be given by the anchor image of elements that belong to $L^\perp$. Since a Dirac space is Lagrangian, $L^\perp = L$, according to our splitting-independent definition above those elements are precisely the kernel of the Dirac space $\Ker{L}$. To prove the converse we need to find subspaces $L_W\subset C$ and $\overline{W}\subset V^*$ such that $L_W$ is isotropic and they fit in the short exact sequence
\begin{equation*}
    \begin{tikzcd}
        0 \arrow[r] & \overline{W} \arrow[r, "j"] & L_W \arrow[r, "\rho"] & W \arrow[r]  & 0
    \end{tikzcd}.
\end{equation*}
This is achieved by setting
\begin{equation*}
    \overline{W}:=\{\alpha\in V^*|\quad \exists w\in W: \quad\alpha(v)=\omega(v,w) \quad\forall v\in V\}
\end{equation*}
and then defining $L_W^\nabla:=\nabla(W)+j(\overline{W})\subset C$, which is clearly isotropic.
\end{proof}

It trivially follows from this proposition that, given a 2-form $\omega\in\wedge^2 V^*$ and an isotropic splitting $\nabla:V\to C$, we have the Dirac space
\begin{equation*}
    \begin{tikzcd}
        0 \arrow[r] & \omega^\flat(V) \arrow[r, "j"] & L_\omega^\nabla \arrow[r, "\rho"] & V \arrow[r]  & 0
    \end{tikzcd}.
\end{equation*}
A change of isotropic splitting $\nabla'=\nabla +j\circ b^\flat$ then induces a new Dirac space $L_{\omega'}^{\nabla'}\subset C$ with associated 2-form $\omega'\in\wedge^2V^*$ satisfying
\begin{equation*}
    \omega'=\omega+b.
\end{equation*}
An important example of Dirac space is given by a bivector $\pi\in \wedge^2V$. We define $L_\pi^\nabla\subset C$, for a choice of isotropic splitting $\nabla:V\to C$, as: 
\begin{equation*}
\begin{tikzcd}
0 \arrow[r] & V^* \arrow[r, "j"] & L^\nabla_\pi \arrow[r, "\rho"] & \pi^\sharp(V^*) \arrow[r]  & 0.
\end{tikzcd}
\end{equation*}
The 2-form of these Dirac structures is given by $\omega_{\pi}^\nabla(v_1,v_2):=p(\alpha_1,\alpha_2)$ where $\nabla(v_i)+j(\alpha_i)\in L^\nabla_\pi$ $i=1,2$. Note that the datum of subspace $W\subset V$, subject to the choice of isotropic splitting, is enough to define Dirac structures of the form:
\begin{equation*}
\begin{tikzcd}
0 \arrow[r] & W^0 \arrow[r, "j"] & L_W \arrow[r, "\rho"] & W \arrow[r]  & 0,
\end{tikzcd}
\end{equation*}
in which case the associated 2-form vanishes.\newline

Let $(C,\langle, \rangle_C,\rho_C:C\to V)$ and $(D,\langle, \rangle_D,\rho_D:D\to W)$ be two exact Courant spaces, a \textbf{linear Lagrangian relation} between them is a Lagrangian subspace $\Lambda \subset D\oplus \overline{C}$ and we say that two elements $c\in C$ and $d\in D$ are \textbf{$\Lambda$-related} if $d\oplus c \in \Lambda$, we then write $d \sim_\Lambda c$. It is clear from the definition that, given two pairs of related elements $d \sim_\Lambda c$ and $d' \sim_\Lambda c'$, we have:
\begin{equation*}
\langle d ,d' \rangle_D = \langle c,c' \rangle_C .
\end{equation*}
Taking into account not only the bilinear form but also the extra structure given to Courant spaces by the anchor maps, we define a \textbf{morphism of Courant spaces} as a Lagrangian relation $\Gamma\subset D\oplus \overline{C}$ for which there exists a linear map $\psi : V \to W$ such that for all $d\oplus c \in \Gamma$ we have $\rho_D(d)=\psi(\rho_C(c))$. A Courant morphism then becomes a Dirac space $\Gamma_\psi\subset D\oplus \overline{C}$ that fits in the exact sequence:
\begin{equation*}
\begin{tikzcd}
0 \arrow[r] & \text{graph}(\psi^*) \arrow[r, "j_{D\oplus C}"] & \Gamma_\psi \arrow[r, "\rho_{D\oplus C}"] & \text{graph}(\psi) \arrow[r]  & 0
\end{tikzcd}
\end{equation*}
and we say that $\Gamma_\psi$ extends the linear map $\psi : V \to W$. A Courant morphism will be denoted by $\Gamma_\psi : C \dashrightarrow D$ to reflect the fact that it covers a linear map between the anchoring spaces, even if $\Gamma_\psi$ does not necessarily correspond to a map between the Courant spaces themselves. We define the \textbf{kernel} and \textbf{range} of a Courant morphism in a natural way:
\begin{align*}
\text{ker}(\Gamma_\psi) & :=\{c\in C|\quad 0\sim_{\Gamma_\psi}c\}=\Gamma_\psi\cap C\\
\text{ran}(\Gamma_\psi) & :=\{d\in D|\quad \exists c\in C:\quad d\sim_{\Gamma_\psi}c\}= \text{pr}_D(\Gamma_\psi).
\end{align*}
Given two Lagrangian relations $\Lambda \subset D\oplus \overline{C}$ and $\Sigma \subset E\oplus \overline{D}$ we can define the composition in a natural way:
\begin{equation*}
E\oplus \overline{C} \supset \Sigma \circ \Lambda := \{e\oplus c |\quad \exists d\in D,\quad e\oplus d\in \Sigma\quad \&\quad d\oplus c\in \Lambda \}.
\end{equation*}
This operation for Lagrangian relations is clearly associative by construction.
\begin{prop}[Composition of Linear Courant Morphisms, {\cite[Lemma 5.2]{bursztyn2013brief}}] \label{LinearCourantComposition}
Let two Courant morphisms $\Gamma_\psi : C \dashrightarrow D$ and $\Gamma_\phi : D \dashrightarrow E$ extending maps $\psi: V\to W$ and $\phi: W\to U$, then the composition $\Gamma_\phi \circ \Gamma_\psi:C \dashrightarrow E$ is a Courant morphism extending the map $\phi\circ \psi : V\to U$.
\end{prop}
\begin{proof}
First we show that $\Gamma_\phi \circ \Gamma_\psi$ is a Lagrangian subspace of $E\oplus \overline{C}$. To see this, consider the \emph{intermediate} Courant space $R:=E\oplus \overline{D}\oplus D \oplus \overline{C}$ and the subspace $I=E\oplus\text{Diag}(D)\oplus \overline{C}$. Note that $I^\perp = 0\oplus \text{Diag}(D)\oplus 0 \subset I$ so we can rewrite the relation space for the composition as the quotient $I/I^\perp=E\oplus \overline{C}$. It is easy to verify that the subspace  given by $\Gamma_\phi \circ \Gamma_\psi\subset E\oplus \overline{C}$ corresponds to the image of $\Gamma_\phi \oplus \Gamma_\psi \subset R$ in the quotient $I/I^\perp$ which is given by:
\begin{equation*}
\Gamma_\phi \circ \Gamma_\psi=\frac{(\Gamma_\phi \oplus \Gamma_\psi)\cap I + I^\perp}{I^\perp}.
\end{equation*}
The quotient map preserves the bilinear forms as $I/I^\perp$ is indeed just the coisotropic reduction of $R$ by $I$ which implies that $\Gamma_\phi \circ \Gamma_\psi$ will be Lagrangian if and only if $(\Gamma_\phi \oplus \Gamma_\psi)\cap I + I^\perp$ is Lagrangian. Using the properties of bilinear forms of split signature, this is now a simple check:
\begin{equation*}
((\Gamma_\phi \oplus \Gamma_\psi)\cap I + I^\perp)^\perp = (\Gamma_\phi \oplus \Gamma_\psi)\cap I + I\cap I^\perp=( \Gamma_\phi \oplus \Gamma_\psi)\cap I + I^\perp.
\end{equation*}
Then it only remains to show that the canonical relation $\Gamma_\phi \circ \Gamma_\psi$ extends the composition of maps, this is a direct consequence of the definition: elements $\rho_D(d)$ will be of the form $\psi(\rho_C(c))$ for $d\oplus c \in\Gamma_\psi$ and similarly for $e\oplus d \in\Gamma_\phi$, which implies that elements $\rho_E(e)$ will be of the form $\phi(\psi(\rho_C(c)))$.
\end{proof}

We then identify \textbf{the category of exact Courant spaces} $\Crnt$. The morphisms of this category will be morphisms between exact Courant spaces extending linear maps with the identity morphisms given by the diagonals Diag$(C)\subset C\oplus \overline{C}$. Courant morphisms may correspond to maps between the Courant spaces in special cases. Indeed, it is a simple calculation to show that a linear map $\Psi: C\to D$ satifies $\Psi^*\langle , \rangle_D =\langle , \rangle_C$ and $\rho_D \circ \Psi = \psi \circ \rho_C$, for some linear map $\psi: V\to W$, if and only if $\Grph{\psi}_\psi\subset D\oplus \overline{C}$ is a Courant morphism. Such a map is called a \textbf{linear Courant map}. It is also a quick check to show that the composition of two Courant maps $\Phi \circ \Psi$ gives a Courant map and that the corresponding Courant morphisms satisfy $\Grph{\Phi}\circ\Grph{\Psi}=\Grph{\Phi \circ \Psi}$. Since Courant morphisms are, in particular Dirac spaces, it follows from proposition \ref{LinearDirac2Form} that 2-forms characterize Courant morphisms extending linear maps up to isotropic splitting and then the set of Courant morphisms between two Courant spaces $C$ and $D$ sits in the following fibration sequence:
\begin{equation*}
\begin{tikzcd}
 \wedge^2V^* \arrow[r] & \text{Hom}_{\Crnt}(C,D) \arrow[r] & \text{Hom}_\Vect(V,W).
\end{tikzcd}
\end{equation*}
Invertible Courant morphisms indeed correspond to invertible Courant maps and are thus called \textbf{Courant isomorphisms}. We define the \textbf{Courant automorphisms} of $C\in\Crnt$ as the set of Courant isomorphisms of $C$ onto itself, $\text{Aut}_{\Crnt}(C)$. Then it is clear that the fibration sequence above becomes a group extension sequence for the automorphisms:
\begin{equation*}
\begin{tikzcd}
0 \arrow[r] & \wedge^2V^* \arrow[r] & \text{Aut}_{\Crnt}(C) \arrow[r] & \text{GL}(V) \arrow[r] & 0.
\end{tikzcd}
\end{equation*}

The following proposition shows that exact Courant spaces are isomorphic to the standard Courant space.
\begin{prop}[Exact Courant Spaces] \label{StandardCourantSpace}
Any exact Courant space $\begin{tikzcd}[cramped, sep=small]0 \arrow[r] & V^* \arrow[r, "j"] & C \arrow[r, "\rho"] & V \arrow[r]  & 0 \end{tikzcd}$ is (non-canonically) isomorphic to the standard Courant space $\mathbb{V}:= V\oplus V^*$. The choice of isomorphism is equivalent to a choice of isotropic splitting.
\end{prop}
\begin{proof}
To prove this statement, we simply construct a vector space isomorphism from a choice of isotropic splitting $\nabla: V\to C$ as follows:
\begin{align*}
\Phi^\nabla: V\oplus V^* & \to C\\
v\oplus \alpha & \mapsto \nabla(v) + j(\alpha)
\end{align*}
and directly check the Courant map conditions:
\begin{align*}
\rho(\Phi^\nabla(v\oplus \alpha)) & = \rho(\nabla(v) + j(\alpha))=v=\text{id}_V(\text{pr}_V(v\oplus \alpha)),\\
\langle \Phi^\nabla(v\oplus \alpha),\Phi^\nabla(w\oplus \beta) \rangle_C & = \langle\nabla(v) ,\nabla(w) \rangle_C + \langle\nabla(v) ,j(\beta) \rangle_C + \langle j(\alpha), \nabla(w)\rangle_C + \langle j(\alpha),j(\beta) \rangle_C=\\
& = \beta(v)+\alpha(w)\\
& = \langle v\oplus \alpha,w\oplus \beta \rangle_{\mathbb{V}},
\end{align*}
where we have used exactness, $\rho\circ \nabla =$id$_V$ and the fact that the splitting is isotropic. We then see that $\Phi^\nabla$ is a Courant isomorphism extending the identity id$_V:V\to V$. The converse is obvious.
\end{proof}

Standard Courant spaces can indeed be regarded as a subcategory of exact Courant spaces $\mathbb{V}\textsf{ect}:=\{\mathbb{V}:=V\oplus V^*, \quad V\in\Vect\}\subset\Crnt$. The general fibration sequence above becomes
\begin{equation*}
\text{Hom}_{\mathbb{V}\textsf{ect}}(\mathbb{V},\mathbb{W})=\wedge^2V^* \times \text{Hom}_\Vect(V,W),
\end{equation*}
with the automorphisms now forming a semi-direct product of groups
\begin{equation*}
    \text{Aut}_{\mathbb{V}\textsf{ect}}(\mathbb{V})=\wedge^2V^*\rtimes \GL{V}.
\end{equation*}
More explicitly, for a morphism of standard Courant spaces $(\psi,B):\mathbb{V}\dashrightarrow \mathbb{W}$, we find the condition of elements being related:
\begin{equation*}
w\oplus \beta \sim_{(\psi,B)} v\oplus \alpha \quad \Leftrightarrow \quad w=\psi(v) \quad \& \quad \psi^*\beta = \alpha + i_v B.
\end{equation*}
The composition of morphisms of standard Courant spaces is then given by
\begin{equation*}
(\phi,B')\circ(\psi, B)=(\phi\circ \psi,B+\psi^*B'),
\end{equation*}
and it is also easy to verify that
\begin{align*}
\Ker{(\psi,B)} & =\{v\oplus -i_vB\in \mathbb{V}:v\in \Ker{\psi}\}\\
\text{ran}((\psi,B)) & =\text{ran}(\psi)\oplus W^*.
\end{align*}

The next proposition ensures that there is a natural notion of reduction in the category of Courant spaces.

\begin{prop}[Reduction of Courant Spaces] \label{LinearCourantReduction}
Let $\begin{tikzcd}[cramped, sep=small]0 \arrow[r] & V^* \arrow[r, "j"] & C \arrow[r, "\rho"] & V \arrow[r]  & 0 \end{tikzcd}$ be an exact Courant space and consider a subspace $i:K\hookrightarrow C$, then there is a commutative diagram of vector spaces
\begin{equation*}
\begin{tikzcd}
 C \arrow[d,"\rho"] & K \arrow[d,"\rho|_K"]\arrow[l, "i"'] \arrow[r,  "q"] & C' \arrow[d,"\rho'"]  \\
 V  & \rho(K) \arrow[l, "i"']\arrow[r, "q'"] &  V'
\end{tikzcd}
\end{equation*}
defining a Courant space $\rho':C'\to V'$ by
\begin{equation*}
    C':=\frac{K}{K\cap K^\perp},\qquad V':=\frac{\rho(K)}{\rho(K\cap K^\perp)}.
\end{equation*}
The maps $q,q'$ denote the canonical projections of the vector space quotients. Furthermore, $C'$ is an exact Courant space iff
\begin{equation*}
    \rho(K)\cap\rho(K^\perp)=\rho(K\cap K^\perp).
\end{equation*}
\end{prop}
\begin{proof}
This result for Courant spaces is a direct consequence of the following simple fact about general bilinear forms: any subspace $B\subset A$ of a vector space with a non-degenerate bilinear form $(A,\alpha)$ induces a reduction to a new vector space with a non-degenerate bilinear form $(A',\alpha')$ by first restricting $(B,\alpha|_B$ and then factoring out the kernel of the possibly degenerate bilinear form $\alpha_B$. Noting that $\Ker{\alpha_B^\flat}=B\cap B^\perp$ we obtain the desired result. The compatibility of the anchors is induced by construction since they are linear maps and thus descend to linear maps of the quotient. The exactness condition for $C'$ then follows from a simple dimension count.
\end{proof}

Courant morphisms may provide, under suitable conditions, a notion of morphism between Dirac spaces. Let us first note that it is possible to regard Lagrangian subspaces of a Courant space $L\subset C$ as Courant morphisms in a trivial way: $0 \oplus C \supset \Lambda_L\cong L$ and similarly $C\oplus 0 \supset \Lambda_L\cong L$. Therefore, given a Courant morphism $\Gamma: C \dashrightarrow D$ and Dirac spaces $L\subset C$ and $N\subset D$, we can define the \textbf{backwards image} of $N$ and \textbf{forward image} of $L$ under $\Gamma$ as follows:
\begin{align*}
\mathcal{B}_\Gamma (N)& := N\circ \Gamma\subset C\\
\mathcal{F}_\Gamma (L)& := \Gamma \circ L \subset D.
\end{align*}
proposition \ref{LinearCourantComposition} ensures that these are Dirac spaces. These operations on Dirac spaces are not inverses of one another in general.It is easy to see that $\mathcal{B}_\Gamma\circ \mathcal{F}_\Gamma =\Id=\mathcal{F}_\Gamma\circ \mathcal{B}_\Gamma$ iff $\Gamma$ is a Courant isomorphism.\newline

Let two pairs of Courant spaces with a Dirac space, $L\subset C$ and $N\subset D$, a Courant morphism $\Gamma: C \dashrightarrow D$ is called a \textbf{Dirac morphism} between $L$ and $N$ if $N=\mathcal{F}_\Gamma(L)$. If, furthermore, ker$(\Gamma)\cap L=0$ then we say that $\Gamma$ is a \textbf{strong Dirac morphism} between $L$ and $N$. Composition of Dirac morphisms and strong Dirac morphisms clearly gives Dirac morphisms and strong Dirac morphisms by construction. The following proposition gives a characterization of Dirac morphisms.
\begin{prop}[Morphisms of Dirac Spaces, {\cite[Lemma 1.12]{alekseev2007pure}}] \label{LinearDiracMorphism}
A Courant morphism $\Gamma: C \dashrightarrow D$ gives a Dirac morphism between the Dirac spaces $L\subset C$ and $N\subset D$ iff there exists a map $a:N\to L$ such that
\begin{equation*}
n\sim_\Gamma a(n) \quad \forall n\in N.
\end{equation*}
The Dirac morphism $\Gamma: C \dashrightarrow D$ is strong iff $a$ is unique.
\end{prop}
\begin{proof}
It is clear that the condition of the target Dirac space being the forward image under $\Gamma$ of the source Dirac space is equivalent to the existence of a map $a:N\to L$. Furthermore, given the fact that $\Gamma$ is assumed to be Courant morphism, this map will be completely determined by its base component $\overline{a}=\rho_C \circ a: N \to V$. Different choices of these component correspond to maps that differ by an element in $\Ker{\Gamma\cap L}$, therefore the choice is unique precisely when the Dirac morphism is strong.
\end{proof}

This last result motivates the definition of \textbf{the category of Dirac spaces} $\Dir$ with objects the Lagrangian subspaces of Courant spaces and morphisms the strong Dirac morphisms between them. The \textbf{direct sum of Dirac spaces} is defined in the obvious way, as for general subspaces, from the direct sum of Courant spaces. \textbf{Reduction of Dirac spaces} is also induced from the general notion of reduction of Courant spaces: let $C$ be an exact Courant space with reduction $C'$ via subspace $K\subset C$, then if $L\subset C$ is a Dirac space the subspace
\begin{equation*}
L':=\frac{L\cap K_*+K_*^\perp}{K_*^\perp}\subset C'
\end{equation*}
where $K_*=K^\perp\cap (K+j(V^*))$ is also a Dirac space. When $K\subset C$ is coisotropic, this construction becomes the usual reduction of Lagrangian subspaces.

\subsection{Courant Algebras} \label{CourantAlgebras}

A \textbf{Leibniz algebra} is a vector space $\mathfrak{a}$ together with a bilinear bracket $[,]$ that satisfies:
\begin{equation*}
[a,[b,c]]=[[a,b],c]+[b,[a,c]] \quad \forall a,b,c\in \mathfrak{a}.
\end{equation*}
A \textbf{morphism of Leibniz algebras} is a linear map $\Psi:\mathfrak{a}\to \mathfrak{a}'$ that is an algebra morphism with respect to the Leibniz algebra structures. Note that a Leibniz algebra $(\mathfrak{a},[,])$ becomes a Lie algebra when $[,]$ is antisymmetric. A \textbf{Leibniz representation} of a Leibniz algebra $\mathfrak{a}$ on a vector space $\mathfrak{h}$ is a pair of maps
\begin{equation*}
[,]_L:\mathfrak{a}\times\mathfrak{h}\to \mathfrak{h}\quad [,]_R:\mathfrak{h}\times\mathfrak{a}\to \mathfrak{h}
\end{equation*}
such that all compositions that will generate triple products of an element $h\in\mathfrak{h}$ and two elements $a,b\in \mathfrak{a}$ satisfy the Leibniz product rule, reflecting the defining property of the Leibniz algebra, explicitly
\begin{align*}
 [h,[a,b]]_R &=[[h,a]_R,b]_R+[a,[h,b]_R]_L\\
 [a,[b,h]_L]_L &=[[a,b],h]_L+[b,[a,h]_L]_L \\
 [a,[h,b]_R]_L &=[[a,h]_L,b]_R+[h,[a,b]]_R .
\end{align*}
From a Leibniz representation we can define a differential complex using the tensor powers of $\mathfrak{a}$ instead of the wedge powers, in contrast with the case of Lie algebra cohomology, to account for the lack of antisymmetry of the Leibniz bracket. We define the following cochain complex
\begin{equation*}
CL^k(\mathfrak{a};\mathfrak{h}):= \text{Hom}_\Real(\mathfrak{a}^{\otimes k},\mathfrak{h})\cong (\mathfrak{a}^*)^{\otimes k} \otimes \mathfrak{h}
\end{equation*}
which together with the map $\delta : CL^k\to CL^{k+1}$ defined by
\begin{align*}
(\delta \eta)(a_1,\dots ,b_{k+1}):= & \sum_{1\leq i < j \leq k+1}(-)^{j+1}\eta(a_1,\dots,a_{i-1},[a_i,a_j],\dots,\hat{a_j},\dots,a_{k+1})\\
& + [a_1,\eta(a_2,\dots,a_{k+1})]_L + \sum_{i=2}^{k+1}(-)^i[\eta(a_1,\dots,\hat{a_i},\dots,a_{k+1}),a_1]_R
\end{align*}
forms a differential complex whose homology we define as the \textbf{Leibniz cohomology of $\mathfrak{a}$ with coefficients in $\mathfrak{h}$} denoted $HL^k(\mathfrak{a};\mathfrak{h})$. For example, for a 1-cochain $\eta:\mathfrak{a}\to \mathfrak{h}$ we obtain the 2-cochain:
\begin{align*}
(\delta \eta)(a,b)=  [a,\eta(b)]_L+[\eta(a),b]_R-\eta([a,b]).
\end{align*}

A \textbf{Courant algebra $\mathfrak{a}$ over a Lie algebra $\mathfrak{g}$} is a Leibniz algebra $(\mathfrak{a},[,])$ and a linear map $\rho : \mathfrak{a}\to \mathfrak{g}$ such that
\begin{equation*}
\rho([a,b]_{\mathfrak{a}})=[\rho(a),\rho(b)]_{\mathfrak{g}} \quad \forall a,b\in \mathfrak{a},
\end{equation*}
in other words, $\rho$ is a Leibniz algebra morphism. Without any further requirements, the map $\rho$ induces a two-sided Leibniz ideal, $\mathfrak{h}:=$ker$(\rho)\subset \mathfrak{a}$, and a Lie subalgebra, $\rho(\mathfrak{a})\subset \mathfrak{g}$. A Courant algebra $(\mathfrak{a},\rho,[,])$ is called an \textbf{exact Courant algebra} when
\begin{equation*}
\begin{tikzcd}
0 \arrow[r] & \mathfrak{h} \arrow[r, "i"] & \mathfrak{a} \arrow[r, "\rho"] & \mathfrak{g} \arrow[r]  & 0.
\end{tikzcd}
\end{equation*}
is a short exact sequence of Leibniz algebras and $\mathfrak{h}$ is abelian, $[\mathfrak{h},\mathfrak{h}]=0$. Given a Lie algebra $\mathfrak{g}$ acting on a vector space $\mathfrak{h}$ we find the first natural example of an exact Courant algebra in the so-called \textbf{hemisemidirect product of $\mathfrak{g}$ with $\mathfrak{h}$}: the vector space of the Courant algebra is simply the direct sum $\mathfrak{a}=\mathfrak{h} \oplus\mathfrak{g}$ and the Leibniz bracket is uniquely determined by the Lie algebra representation and the Lie algebra bracket:
\begin{equation*}
[h_1 \oplus g_1,h_2 \oplus g_2]:= g_1\cdot h_2\oplus [g_1,g_2].
\end{equation*}

Let two Courant algebras $(\mathfrak{a},[,],\rho)$ and $(\mathfrak{a}',[,]',\rho')$ over the Lie algebras $\mathfrak{g}$ and $\mathfrak{g}'$ respectively. A map $\Psi:\mathfrak{a}\to \mathfrak{a}'$ is said to be a \textbf{Courant algebra morphism} if it is a Leibniz algebra morphism and there exists a Lie algebra morphism $\psi:\mathfrak{g}\to \mathfrak{g}'$ such that $\rho'\circ \Psi = \rho \circ \psi$. When the map $\Psi$ is invertible we call it a \textbf{Courant algebra isomorphism}. In the case of exact Courant algebras, Courant morphisms equivalently correspond to morphisms of exact sequences, i.e. commutative diagrams:
\begin{equation*}
\begin{tikzcd}
0 \arrow[r] & \mathfrak{h}\arrow[d, "\eta"'] \arrow[r, "i"] & \mathfrak{a} \arrow[r, "\rho"] \arrow[d,"\Psi"] & \mathfrak{g}\arrow[d, "\psi"] \arrow[r] & 0 \\
0 \arrow[r] & \mathfrak{h}' \arrow[r, "i'"] & \mathfrak{a}' \arrow[r, "\rho'"] & \mathfrak{g}' \arrow[r] & 0 
\end{tikzcd}
\end{equation*}
where $\Psi$ is a Leibniz morphism, $\psi$ is a Lie algebra morphism and $\eta$ is induced by commutativity.\newline

An exact Courant algebra $\mathfrak{a}$ over $\mathfrak{g}$ always provides a Leibniz representation of $\mathfrak{g}$ on the kernel $\mathfrak{h}$ (note that we regard the three vector spaces as Leibniz algebras in general). The representation maps are
\begin{align*}
[,]_L:\mathfrak{g}\times\mathfrak{h} & \to \mathfrak{h}\\
(g,h) & \mapsto [a,h] \quad , \quad g=\rho(a)
\end{align*}
and similarly for $[,]_R$. These maps take values in the correct space because $\mathfrak{h}$ is an ideal and they are well-defined from the fact that if $\rho(a)=g=\rho(a')$ then $\rho(a-a')=0$ so $a-a'=h'$ for some $h\in \mathfrak{h}$, then
\begin{equation*}
[g,h]_L=[a',h]=[a+h',h]=[a,h]+[h',h]=[a,h]=[g,h]_L,
\end{equation*}
where we have used the fact that $\mathfrak{h}$ is abelian. To see a more explicit connection between an exact Courant algebra and Leibniz cohomology we note that we may always find a splitting of vector spaces for the exact sequence
\begin{equation*}
\begin{tikzcd}
0 \arrow[r] & \mathfrak{h} \arrow[r, "i"] & \mathfrak{a} \arrow[r, "\rho"] & \mathfrak{g} \arrow[r] \arrow[l, bend left, "n"]  & 0.
\end{tikzcd}
\end{equation*}
that in general will not be a splitting of Leibniz algebras. We can define the curvature that measures to what degree $n$ fails to be a Leibniz morphism by
\begin{equation*}
C_n(g,g'):= [n(g),n(g')]-n([g,g']).
\end{equation*}
This is clearly linear in both arguments and verifies $\rho(C_n(g,g'))=0$ by construction so we can see $C_n$ as a linear map of the form
\begin{equation*}
C_n : \mathfrak{g}\otimes \mathfrak{g}\to \mathfrak{h}
\end{equation*}
hence implying $C_n\in CL^2(\mathfrak{g};\mathfrak{h})$. Noting that that $n$ allows to write the Leibniz representation maps explicitly as $[g,h]_L=[n(g),h]$ and $[h,g]_R=[h,n(g)]$, a direct computation using the Leibniz property shows that $\delta C_n = 0$. Therefore, an exact Courant algebra defines a Leibniz 2-cohomology class $[C_n]\in HL^2(\mathfrak{g};\mathfrak{h})$. If we now consider a different splitting $n':\mathfrak{g}\to \mathfrak{a}$, which clearly relates to our previous choice by $n'=n+\eta$ with $\eta: \mathfrak{g}\to \mathfrak{h}$, it is a straightforward calculation, using the above expressions for the right and left actions, to show that:
\begin{equation*}
C_{n'}(g,g')=C_n(g,g')+\delta\eta (g,g')
\end{equation*}
and so we see that $[C_n]=[C_{n'}]\in HL^2(\mathfrak{g};\mathfrak{h})$. This indicates that an exact Courant algebra $\mathfrak{a}$ over the Lie algebra $\mathfrak{g}$ determines a class in the Leibniz 2-Cohomology, we call this the \textbf{characteristic class of the exact Courant algebra} and we denote it by $[C_\mathfrak{a}]\in HL^2(\mathfrak{g};\mathfrak{h})$. The following proposition shows that exact Courant algebras are specified by their characteristic classes up to isomorphism.
\begin{prop}[Isomorphic Exact Courant Algebras] \label{IsomorphicExactCourantAlgebras}
Two exact Courant algebras $\mathfrak{a}$ and $\mathfrak{a}'$ over the same Lie algebra $\mathfrak{g}$ are isomorphic iff their characteristic classes agree, $[C_\mathfrak{a}]=[C_{\mathfrak{a}'}]$.
\end{prop}
\begin{proof}
We first fix vector isomorphisms by choosing splittings $n:\mathfrak{g}\to \mathfrak{a}$ and $n':\mathfrak{g}\to \mathfrak{a}'$ so that we can identify $\mathfrak{a} = i(\mathfrak{h})\oplus n(\mathfrak{g})$ and $\mathfrak{a}= i'(\mathfrak{h})\oplus n'(\mathfrak{g})$. Then assume there exists a Courant algebra isomorphism so that elements from the two spaces above are in 1-to-1 correspondence via a Leibniz map that also gives an isomorphism of exact sequences. Commutativity of the diagrams allows us to relate the direct sum elements via:
\begin{equation*}
\tilde{\Psi}(h\oplus g)=h+\eta_\Psi(g)\oplus g
\end{equation*}
for some $\eta_\Psi: \mathfrak{g}\to \mathfrak{h}$. Writing the cocycles for each of the direct sums and using the fact that $\Psi$ is Leibniz and using commutativity again we can see:
\begin{equation*}
C_n(g,g')\oplus 0=C_{n'}(g,g')+\delta\eta_\Psi(g,g')\oplus 0
\end{equation*}
hence showing that the two cocycles are cohomologous and proving that the charateristic classes agree. Conversely, let us assume that the two exact Courant algebras share the same characteristic class. Again making use of the two splittings we define a linear isomorphism via $\mathfrak{a} = i(\mathfrak{h})\oplus n(\mathfrak{g}) \cong \mathfrak{h}\oplus \mathfrak{g} \cong i'(\mathfrak{h})\oplus n'(\mathfrak{g})=\mathfrak{a}'$ in the obvious way so that it fits in the commutative diagrams. From our assumption we know that the classes given by representatives $[C_n],[C_{n'}]\in HL^2(\mathfrak{g};\mathfrak{h})$ must agree, i.e. there exists a map $\eta_{n,n'}:\mathfrak{g}\to \mathfrak{h}$ so that $C_n=C_{n'}+\delta \eta_{n,n'}$; an elementary calculation shows that this additional condition allows to make the linear isomorphism above into a Leibniz algebra morphism hence proving the converse statement.
\end{proof}
To see the implications of the above proposition more explicitly, consider the isomorphism induced by a concrete splitting $n$:
\begin{align*}
\Phi_n: \mathfrak{h}\oplus \mathfrak{g} & \to \mathfrak{a}\\
h\oplus g & \mapsto i(h)+n(g)
\end{align*}
which, in general, defines a Leibniz bracket on $\mathfrak{h}\oplus \mathfrak{g}$ via:
\begin{equation*}
[h\oplus g,h'\oplus g']:= [h,g']_R+[g,h']_L+C_n(g,g')\oplus[g,g'].
\end{equation*}
proposition \ref{IsomorphicExactCourantAlgebras} ensures that if $[C_n]=0$, we can find a (possibly different) isomorphism of Leibniz algebras for which the last term in the first factor vanishes. If this is the case, the exact Courant algebra $\mathfrak{a}$ is isomorphic to the Leibniz algebra $\mathfrak{h}\oplus\mathfrak{g}$ with bracket given simply by a Leibniz representation of $\mathfrak{g}$ on $\mathfrak{h}$, this is called the \textbf{semidirect product Courant algebra}.

\subsection{Courant Algebroids and Dirac Structures} \label{CourantAlgebroids}

Note that both the linear Courant spaces of Section \ref{LinearCourant} and Courant algebras of Section \ref{CourantAlgebras} were introduced as extensions of the familiar notions of vector space and Lie algebra respectively. Recall that Lie algebroids were presented in Section \ref{LieAlgebroids} as fibrations of vector spaces over smooth manifolds whose modules of sections carried a Lie algebra structure, in the same spirit, Courant algebroids will be defined below as fibrations of Courant spaces over smooth manifolds whose modules of sections carry a Courant algebra structure.\newline

Before giving the general definition, let us discuss a motivating example. A \textbf{Lie bialgebroid} is defined as the natural algebroid generalization of Lie bialgebra: a pair of a vector bundle $\pi:A\to M$ and its dual, $(A,A^*)$, such that there are Lie algebroid structures on both, denoted by $(A,\rho,[,])$ and $(A^*,\rho_*,[,]_*)$, which are compatible in the following sense:
\begin{equation*}
    d_{A^*}[a,b]=[d_{A^*}a,b]+[a,d_{A^*}b]
\end{equation*}
for all $a,b\in\Sec{A}$ and where $[,]$ denotes the Schouten bracket of the Lie algebroid Gerstenhaber algebra $\Sec{\wedge^\bullet A}$ and $d_{A^*}$ is the graded differential of $A^*$ defined on the complex $\Sec{\wedge^\bullet (A^*)^*}\cong \Sec{\wedge^\bullet A}$. Forming the direct sum $E:=A\oplus A^*$ with anchor $\rho_E:=\rho + \rho_*:E\to \Tan M$ and defining the obvious fibre-wise symmetric bilinear form from the dual pairing
\begin{equation*}
    \langle a\oplus\alpha ,b\oplus \beta\rangle:=\tfrac{1}{2}(i_b\alpha+i_a\beta),
\end{equation*}
clearly makes $E$ into a vector bundle with each fibre being a Courant space. The identities of the Cartan calculus of the Lie algebroids together with the bialgebroid compatibility condition ensure that $\Sec{A}$ and $\Sec{A^*}$ have Leibniz representations on each other, thus inducing two semidirect product Courant algebras on $\Sec{A}\oplus\Sec{A^*}$. Putting these together we find a Courant algebra on $\Sec{E}$ with bracket explicitly given by
\begin{equation*}
    [a\oplus\alpha ,b\oplus \beta]:= [a,b] + \LDer_\alpha b -i_\beta d_{A^*}a \oplus [\alpha,\beta] + \LDer_a\beta -i_bd_{A}\alpha.
\end{equation*}
Furthermore, the fact that the Cartan derivations $i$, $\LDer$, $d_{A,A^*}$ were used to define both the bilinear form and the Leibniz bracket on sections will imply some compatibility conditions between them. These will be captured in the general definition of Courant algebroid below.\newline

A \textbf{Courant algebroid} structure on a vector bundle $\pi:E\to M$, denoted by $(E,\langle,\rangle,\rho:E\to \Tan M,[,])$, is a non-degenerate symmetric bilinear form $\langle,\rangle\in\Sec{\odot^2 E^*}$ and Courant algebra structure $\rho:(\Sec{E},[,])\to(\Sec{\Tan M},[,])$ such that
\begin{align*}
    [a,f\cdot b] &=f\cdot [a,b]+\rho(a)[f]\cdot b\\
    \rho(a)[\langle b , c \rangle] &=\langle [a,b] , c \rangle + \langle b , [a,c] \rangle\\
    [a,a] &=D\langle a , a \rangle
\end{align*}
for all $a,b\in\Sec{E}$, $f\in\Cin{M}$, where $D:=\sharp\circ \rho^* \circ d:\Cin{M}\to \Sec{E}$ is defined from the musical isomorphism $\sharp:E^*\to E$ induced by the non-degenerate bilinear form. Since the anchor map is $\Cin{M}$-linear on sections, the above identities imply that $\rho\circ \sharp \circ \rho^*=0$ as bundle maps, then, denoting $j:=\sharp \circ \rho^*$, we see that Courant algebroids always sit in a complex of vector bundles of the form
\begin{equation*}
    \begin{tikzcd} \Cot M \arrow[r, "j"] & E \arrow[r, "\rho"] & TM\end{tikzcd}.
\end{equation*}
This is the vector bundle analogue of the Leibniz ideal given by the kernel of the anchor in a general Courant algebra, which in this case further appears as the sections of an isotropic subbundle $j(\Cot M)\subset \Ker{\rho}\subset E$. When the complex above is a short exact sequence of vector bundle morphisms covering the identity, $(E,\langle,\rangle,\rho,[,])$ is called an \textbf{exact Courant algebroid}. Clearly, the Courant algebra of sections of an exact Courant algebroid is exact and the failure of the Leibniz bracket to be antisymmetric is precisely controlled by the symmetric bilinear form on sections since one readily checks:
\begin{equation*}
    [a,b]=-[b,a]+2D\langle a,b\rangle
\end{equation*}
for all $a,b\in\Sec{E}$. For the reminder of this thesis Courant algebroids will be assumend to be exact unless otherwise stated. We note in passing that the Leibniz identity of the anchor with the bracket in the definition of Courant algebroid above precisely corresponds to the condition of locality of a derivative bracket, in direct analogy with our definition of local Lie algebra of Section \ref{DerivativeAlgebras}. Although it will not be discussed further in this thesis, we see that Courant algebroids could be regarded as a subclass of local Leibniz algebras.\newline

Let $(E\to M,\langle,\rangle,\rho,[,])$ be a Courant algebroid, a \textbf{Dirac structure} supported on the submanifold $Q\subset M$ is a subbundle $L\subset E|_Q$ over $Q$ satisfying:
\begin{itemize}
    \item[1)] $L_q\subset E_q$ is maximally isotropic for all $q\in Q$,
    \item[2)] $\rho(L)\subset \Tan Q$,
    \item[3)] for all $a,b\in \Sec{E}$, $a|_Q,b|_Q \in \Sec{L} \Rightarrow [a,b]|_Q\in \Sec{L}$.
\end{itemize}
Note that the properties of the bracket in the ambient Courant algebroid ensure that 3) suffices to hold on sections that locally span $L\subset E|_Q$. From this definition and the axioms of Courant algebroid, it is clear that $(L,[,]|_L,\rho|_L)$ becomes a Lie algebroid. To see this more explicitly we can define the \textbf{Courant tensor} of a maximally isotropic subbundle $L\subset E$ (possibly with support) as
\begin{equation*}
\Upsilon_L(x,y,z):=\langle[x,y]|_L,z\rangle|_L
\end{equation*}
for sections $x,y,z\in \Sec{L}$. From isotropy of $L$ and the properties of the bracket on $E$, we can verify that this form is indeed tensorial and antisymmetric in all its entries, $\Upsilon_L\in \Sec{\wedge^3L^*}$. The involutivity condition for $L$ to be a Dirac structure, $[\Sec{L,\Sec{L}}]|_L\subset \Sec{L}$, can be easily checked to reduce to the condition $\Upsilon_L=0$.\newline

As mentioned above, Dirac structures are Lie algebroids and so they induce an involutive tangent distribution $R:=\rho(L)\subset \Tan M$ which will be integrable by a singular foliation: through any point $q\in Q$ there is a submanifold $O_q$ such that $R_q=T_qO_q$. We call $R$ the \textbf{range distribution} of the Dirac structure $L$. This gives the fibre-wise linear Dirac spaces:
\begin{equation*}
\begin{tikzcd}
0 \arrow[r] & \overline{R}_q \arrow[r, "j"] & L_q \arrow[r, "\rho"] & T_qO_q \arrow[r]  & 0.
\end{tikzcd}
\end{equation*}
Then, it follows from proposition \ref{LinearDirac2Form} that there are 2-forms at each tangent space of the distribution $T_qO_q$ which, from smoothness of all the maps involved, fit into a well-defined 2-form $\omega_L^q\in \Omega^2(O_q)$. Recalling the point-wise definition of the Courant tensor $\Upsilon_L$ above it follows by construction that
\begin{equation*}
d\omega_L^q(X,Y,Z)=\langle [x,y],z \rangle|_L=\Upsilon_L(x,y,z)
\end{equation*}
for all $X,Y,Z\in\Sec{\Tan O_q}$ and any $x,y,z\in\Sec{L}$ such that $\rho(x)=X$, $\rho(y)=Y$, $\rho(z)=Z$. Since $L$ is involutive, $\Upsilon_L=0$, which in turn makes the 2-form $\omega^q_L$ closed, and thus each integral manifold $(O_q,\omega^q_L)$ becomes a presymplectic manifold. The singular foliation integrating the range distribution is then called the \textbf{presymplectic foliation} of the Dirac structure $L$. The \textbf{characteristic distribution} $K\subset \Tan Q$ of a Dirac structure $L$ is as the leaf-wise kernel of $\omega_L^q$.\newline

Associated to Dirac structures we also find a bracket on functions of the base manifold. Let us define the subring of \textbf{admissible functions} on $Q$ associated to the Dirac structure $L$ as follows:
\begin{equation*}
\Cin{Q}_L=:\{f\in \Cin{Q}|\quad \exists x\in \Sec{E|_Q}, \quad x+j(df)\in L\}.
\end{equation*}
Any vector field $\rho(x)\in \Sec{\Tan Q}$, with $x\in \Sec{E|_Q}$ any section that makes $f\in \Cin{Q}$ admissible, is called a \textbf{Hamiltonian vector field} relative to the Dirac structure $L$. Clearly, Hamiltonian vector fields associated with a Dirac structure are unique, if they exist, when the characteristic distribution vanishes. The \textbf{Dirac bracket} of a pair of admissible functions $f,g\in \Cin{Q}_L$ can be explicitly given as:
\begin{equation*}
\{f,g\}_L:=\rho(x)[g]=-\rho(y)[f],
\end{equation*}
where $\rho(x),\rho(y)\in \Sec{\Tan Q}$ are Hamiltonian vector fields for $f$ and $g$ respectively. This bracket fails to satisfy the Jacobi identity generically since it can be shown to be controlled by the Severa class (identified in proposition \ref{SeveraClass} below) of the underlying exact Courant algebroid.\newline

Let $(E\to M,\langle,\rangle_E,\rho_E,[,]_E)$ and $(D\to N,\langle,\rangle_D,\rho_D,[,]_D)$ be two Courant algebroids, we define the \textbf{opposite Courant algebroid} of $E$ as $\overline{E}:=(E,-\langle,\rangle_E,\rho_E,[,]_E)$ and the \textbf{product Courant algebroid} as $E \boxplus D$ with bilinear form $\langle,\rangle_{E \boxplus D}$ defined fibre-wise in the obvious way, anchor map $\rho_E\boxplus\rho_D:E \boxplus D\to \Tan M\boxplus\Tan N$ and Leibniz bracket $[,]_{E \boxplus D}$ defined from $[,]_E\oplus [,]_D$ and extending by the Leibniz identity with the anchor. A \textbf{morphism of Courant alegrbroids} is a Dirac structure $\Gamma_\phi\subset D\boxplus \overline{E}$ supported on $\Grph{\varphi}\subset N\times M$ for some smooth map $\varphi: M\to N$. We say that $\Gamma_\phi$ is a Courant morphism extending the smooth map $\varphi$ and denote $\Gamma_\varphi : E \dashrightarrow D$. More explicitly, a Courant morphism $\Gamma_\varphi$ sits in a short exact sequence of vector bundles
\begin{equation*}
\begin{tikzcd}
0 \arrow[r] & \Cot\varphi \arrow[r, "j_{D\boxplus E}"] & \Gamma_\phi \arrow[r, "\rho_{D\boxplus E}"] & \Grph{\Tan \varphi} \arrow[r]  & 0
\end{tikzcd}
\end{equation*}
where $\Cot \varphi\subset \Cot(N\times M)$ is the Lagrangian relation induced by cotangent lift and $\Grph{\Tan\varphi}\subset \Tan (M\times N)$ is simply the graph of the tangent map. A vector bundle morphism $\Phi:E\to D$ covering a smooth map $\phi:M\to N$ is called a \textbf{Courant map} when $\Grph{\Phi}\subset D\oplus \bar{E}$ is a morphism of Courant algebroids. A Courant map $\Phi:E\to D$ is shown to preserve the Courant algebroid structures explicitly in the following sense:
\begin{align*}
\rho_D\circ \Phi & = T\phi\circ \rho_E\\
\Phi^*\langle,\rangle_D & = \langle,\rangle_E\\
[\Phi_*\cdot,\Phi_*\cdot]_D &= \Phi_*[\cdot,\cdot]_E,
\end{align*}
where in the second identity the pull-back of covariant tensors of vector bundles has been used and noting that the third equation must be read in terms of $\Phi$-relatedness in general. When $\Phi:E\to D$ covers a diffeomorphism $\varphi:M\to N$ the third equation indeed becomes the condition that the push-forward $\Phi_*:\Sec{E}\to\Sec{D}$ is a morphism of Leibniz algebras. Comparing this notion of morphism of Courant algebroids with the notion of linear Courant morphism of Section \ref{LinearCourant} we see that it corresponds to the natural vector bundle generalization. In particular, it follows from our construction above that a morphism of Courant algebroids $\Gamma_\varphi : E \dashrightarrow D$ induces morphisms of linear Courant spaces $\Gamma_{\Tan_x\varphi}:E_x\dashrightarrow D_{\varphi(x)}$ for all $x\in M$. The \textbf{kernel} $\Ker{\Gamma_\varphi}\subset E$ and \textbf{range} $\text{ran}(\Gamma_\varphi)\subset D$ of a morphism of Courant algebroids $\Gamma_\varphi : E \dashrightarrow D$ are defined in an analogous way to the linear case.\newline

Given two Courant morphisms $\Gamma_\phi : E \dashrightarrow D$ and $\Gamma_\varphi : D \dashrightarrow C$ extending the smooth maps $\phi: M\to N$ and $\varphi: N\to S$, respectively, the construction of the \textbf{composition} $\Gamma_\varphi\circ \Gamma_\phi: E \dashrightarrow C$ as a Dirac structure on $C\boxplus \overline{E}$ supported on graph$(\varphi\circ \phi)\subset S\times M$ can always be attempted and, at least fibre-wise, it will give a collection of linear Courant morphisms. However, much like in the Weinstein category of symplectic manifolds, this construction involves intersections of general submanifolds in the total spaces of the vector bundles involved and thus they may fail to be smooth. Note, however, that this technical complication disappears when the Courant morphisms are graphs of Courant maps, since for any two vector bundle maps we have $\Grph{\Psi}_\varphi\circ\Grph{\Phi}_\phi=\Grph{\Psi \circ \Phi}_{\varphi\circ \phi}$. We then identify \textbf{the category of Courant algebroids} over smooth manifolds $\Crnt_\Man$ whose objects are Courant algebroids and morphisms are Courant maps. Following the analogy with the linear case we define \textbf{morphisms of Dirac structures} and \textbf{the category of Dirac manifolds} $\Dir_\Man$ whose objects are Dirac structures of Courant algebroids over smooth manifolds.\newline

Any smooth manifold $M$ carries a Courant algebroid structure on the vector bundle $\mathbb{T}M:=\Tan M\oplus \Cot M$, called the \textbf{standard Courant algebroid}, and constructed by regarding $(\Tan M,\Cot M)$ as a Lie bialgebroid with zero anchor and zero Lie bracket on $\Cot M$. In particular, this gives the Leibniz bracket
\begin{equation*}
[X\oplus \alpha,Y\oplus \beta]=[X,Y]\oplus\LDer_X\beta -i_yd\alpha
\end{equation*}
referred to as the \textbf{Dorfman bracket} in the standard literature. Recall from our comments at the end of Section \ref{CourantAlgebras} that any Leibniz representation can be twisted by the addition a 2-cocycle in the Leibniz cohomology. In the case of the Lie bialgebroid $(\Tan M,\Cot M)$, it is clear that the Leibniz cohomology corresponds to the usual de Rham cohomology $HL^2(\Sec{\Tan M};\Sec{\Cot M})\cong H_{\text{dR}}^3(M)$. Then, for any closed 3-form $H\in\Omega_{\text{cl}}^3(M)$, we define the \textbf{$H$-twisted standard Courant algebroid} $\mathbb{T}M_H$ with the same underlying vector bundle, bilinear form and anchor map but altering the Dorfman bracket by the addition of the 3-form:
\begin{equation*}
[X\oplus \alpha,Y\oplus \beta]_H=[X,Y]\oplus\LDer_X\beta -i_yd\alpha +i_Yi_XH.
\end{equation*}
The next proposition establishes a close connection between the Leibniz cohomology of the exact Courant algebras arising as modules of sections of Courant algebroids and the usual de Rham cohomology of the base manifolds of the vector bundle underlying the Courant algebroid structure.

\begin{prop}[The Severa Class of a Courant Algebroid, {see \cite{vsevera2001poisson}}] \label{SeveraClass}
A Courant algebroid structure on the vector bundle $\pi:E\to M$ uniquely determines a cohomology class $\normalfont[H_E]\in H^3_{\text{dR}}(M)$, called the \textbf{Severa class} of the Courant algebroid $E$, corresponding to the characteristic class of the exact Courant algebra structure of the module of sections $(\Sec{E},[,])$. Furthermore, two Courant algebroids are isomorphic iff their Severa classes agree.
\end{prop}
\begin{proof}
Under the mild assumption that $M$ is a paracompact manifold, splittings of exact sequences of vector bundles over $M$ always exist, then, by a fibre-wise argument, proposition \ref{IsotropicSplittingLinearCourant} implies that isotropic splittings $\nabla: TM\to E$ always exist for Courant algebroids and that, for any two such splittings, we have $\nabla - \nabla' = j\circ b^\flat$ with $b\in \Omega^2(M)$. Applying proposition \ref{StandardCourantSpace} fibre-wise, we see that a choice of isotropic splitting allows us to construct a vector bundle isomorphism explicitly as:
\begin{align*}
\Phi_\nabla: TM\oplus \Cot M  &\to E\\
X + \alpha  &\mapsto \nabla X + j \alpha
\end{align*}
Recalling that $j=\sharp\rho^*$, we verify
\begin{align*}
\rho_{\mathbb{T}M}(X+\alpha) &= \rho (\Phi(X+\alpha))=\rho\nabla X = X\\
\langle X+\alpha ,Y+\beta \rangle_{\mathbb{T}M} &= \langle\nabla X+\sharp\rho^*\alpha ,\nabla Y+\sharp\rho^* \beta \rangle=\langle\sharp\rho^*\alpha ,\nabla Y \rangle+\langle\nabla X, \sharp\rho^*\beta \rangle = i_X\beta + i_Y\alpha
\end{align*}
so the fibre-wise linear Courant structures are preserved under the isomorphism. The Leibniz bracket on $E$ is pulled-back to the bracket:
\begin{equation*}
\Phi^{-1}([\Phi(X+\alpha),\Phi(Y+\beta)]) = \rho([\Phi(X+\alpha),\Phi(Y+\beta)]) + \nabla^*\flat([\Phi(X+\alpha),\Phi(Y+\beta)])
\end{equation*}
on sections $X+\alpha,Y+\beta\in \Sec{\Tan M}$. The $\Tan M$ component is easily checked to be simply $[X,Y]$ and using the identities for general Courant algebroids we find the $\Cot M$ component to be
\begin{align*}
\nabla^*\flat([\Phi(X+\alpha),\Phi(Y+\beta)]) & = \nabla^*\flat([\nabla X,\nabla Y]+[\nabla X,+\sharp\rho^*\beta)]+[\sharp\rho^*\alpha,\nabla Y]+[\sharp\rho^*\alpha,\sharp\rho^*\beta)])\\
& = \nabla^*\flat(\sharp\rho^*(i_Yi_XH_\nabla)+\sharp\rho^*(\mathcal{L}_X\beta)+\sharp\rho^*(-i_Yd\alpha)) \\
& = \mathcal{L}_X\beta -i_Yd\alpha + i_Yi_XH_\nabla,
\end{align*}
where we have identified 
\begin{equation*}
    H_\nabla(X,Y,Z):=\langle[\nabla(X),\nabla(Y)],\nabla(Z)\rangle
\end{equation*}
to be tensorial and totally antisymmetric, hence a 3-form on $M$, from the fact that $\nabla$ is isotropic. This shows that, after choosing an isotropic splitting $\nabla$, the exact Courant algebroid $E\to M$ is isomorphic to the $H_\nabla$-twisted standard Courant algebroid over its base, $E\cong \mathbb{T}M_{H_\nabla}$. Let two isotropic splittings $\nabla,\nabla'$ of the same exact Courant algebroid $E$, as we saw above, they are related by $\nabla' = \nabla+\lambda$ with $\lambda = j \circ b^\flat$ for some $b\in \Omega^2(M)$. The brackets induced by the two splittings only depend on the choice of splitting via the twisting 3-form so we compare $H_\nabla$ with $H_{\nabla'}$ as follows
\begin{align*}
i_Zi_Yi_XH_{\nabla'} =i_Zi_Yi_XH_{\nabla+\lambda} & = \langle [\nabla X+\lambda X,\nabla Y+\lambda Y ],\nabla Z+\lambda Z\rangle \\
& = \langle [\nabla X,\nabla Y],\nabla Z\rangle + \langle [\lambda X,\nabla Y],\nabla Z\rangle + \langle [\nabla X,\lambda Y ],\nabla Z\rangle + \langle [\nabla X,\nabla Y],\lambda Z\rangle
\end{align*}
where we have cancelled all the terms that contain more than one $\lambda$. We are now able to rewrite the remaining terms as
\begin{align*}
i_Zi_Yi_XH_{\nabla'} & = i_Zi_Yi_XH_\nabla + \mathcal{L}_Xi_Yi_Zb-i_Y\mathcal{L}_Xi_Zb-i_Zi_Y\mathcal{L}_Xb+i_Z\mathcal{L}_Xi_Yb+i_Zi_Yi_Xdb\\
& = i_Zi_Yi_X(H_\nabla+db) + [\mathcal{L}_X,i_Y]i_Zb + i_Z[\mathcal{L}_X,i_Y]b \\
& = i_Zi_Yi_X(H_\nabla+db)
\end{align*}
hence showing that the relation between twistings of two different isotropic splittings is given by:
\begin{equation*}
H_{\nabla'}=H_\nabla + db.
\end{equation*}
Then, $[H_\nabla]=[H_{\nabla'}]$, which proves that the cohomology class defined by any particular isotropic splittings is, in fact, splitting independent. This defines the \textbf{Severa class} of the Courant algebroid $[H_E]:=[H_\nabla]\in H^3_{\text{dR}}(M)$. Note that an isotropic splitting $\nabla$ for $E$ gives a splitting for the sequence of sections as well thus inducing an splitting of exact Courant algebras with explicit Leibniz representation of $\Sec{\Tan M}$ on $\Sec{\Cot M}$ given by:
\begin{equation*}
[\alpha,X]_R=[j\alpha,\nabla X]=j(-i_Xd\alpha) \quad\quad [X,\alpha]_L=[\nabla X,j\alpha]= j(\mathcal{L}_X\alpha).
\end{equation*}
Take the cocycle of Leibniz cohomology for this particular splitting $C_\nabla(X,Y):= [\nabla X,\nabla Y]-\nabla([X,Y])$, which by construction satisfies
\begin{align*}
C_\nabla(X,Y) &=-C_\nabla(Y,X)\\
C_\nabla(X,f\cdot Y) &= f [\nabla X,\nabla Y] + \rho(\nabla X)[f]\nabla Y -f \nabla [X,Y]-\nabla(X[f]Y)\\
 &= fC_\nabla(X,Y) +X[f]\nabla Y - X[f]\nabla Y\\
 &= f\cdot C_\nabla(X,Y)
\end{align*}
hence $C_\nabla : \Sec{\Tan M}\otimes \Sec{\Tan M}) \to \Sec{E}$. But also, trivially $\rho(C_\nabla(X,Y))=0$ so $C_\nabla = j(\eta_\nabla)$ for some $\Cin{M}$-linear map $\eta_\nabla:\Sec{\Tan M}\otimes \Sec{\Tan M} \to \Sec{\Cot M}$. We can write the action of $\eta_\nabla$ on arbitrary arguments as
\begin{align*}
 \eta_\nabla(X,Y)(Z) & =\langle j(\eta_\nabla(X,Y)),\nabla Z \rangle = \langle C_\nabla(X,Y),\nabla Z\rangle \\
& = \langle [\nabla X,\nabla Y],\nabla Z\rangle - \langle \nabla [X,Y],\nabla Z \rangle \\
& = \langle [\nabla X,\nabla Y],\nabla Z\rangle = H_\nabla(X,Y,Z)
\end{align*}
thus giving
\begin{equation*}
C_\nabla(X,Y)=j(i_Yi_XH_\nabla).
\end{equation*}
Which shows the desired relation between the Severa class of an exact Courant algebroid $[H_\nabla]\in H^3_{\text{dR}}(M)$ and the Leibniz cohomology class $[C_\nabla]\in HL^2(\Gamma(TM);\Gamma(\Cot M))$ of its associated exact Courant algebra. This identification then allows us to use proposition \ref{IsomorphicExactCourantAlgebras} for the case of $\Cin{M}$-linear isomorphisms of exact Courant algebras to complete the proof.
\end{proof}

It follows from the previous proposition that, under the choice of isotropic splittings, an arbitrary Courant morphism between two exact Courant algebroids $\Gamma_\varphi:E \dashrightarrow D$ is characterised by a smooth map between the base spaces $\phi:M\to N$ and a 2-form $B\in \Omega^2(M)$. As for general Dirac structures, the involutivity condition for $\Gamma_\varphi$ can be expressed as the vanishing of the Courant tensor on $D\oplus \bar{E}$. Note that, by construction, we have $\Upsilon^{D\oplus \bar{E}}=\Upsilon^D\oplus -\Upsilon^E$, then choosing a pair of splittings $\nabla:\Tan M\to E$ and $\Delta:\Tan N\to D$, so that the Leibniz brackets on each factor are given by the twisted Dorfman brackets $[,]_{H_\nabla}$ and $[,]_{H_\Delta}$, it is a simple calculation to show
\begin{equation*}
\Upsilon^{D\oplus \bar{E}}|_{\Gamma_\varphi}= 0 \quad \Leftrightarrow \quad \phi^*H_{\Delta}=H_\nabla+dB.
\end{equation*}
In other words, $\phi^*[H_D]=[H_E]$ as the Severa classes are splitting independent. This implies that Courant morphisms between exact Courant algebroids extend smooth maps that preserve the Severa classes and shows that morphisms between to fixed Courant algebroids fit into the following fibration sequence:
\begin{equation*}
\begin{tikzcd}
 \Omega^2_{\text{cl}}(M) \arrow[r] & \text{Hom}_{\Crnt_\Man}(E,D) \arrow[r] & \Cin{M,N}_{H_{\text{dR}}^3}
\end{tikzcd}
\end{equation*}
where $\Cin{M,N}_{H_{\text{dR}}^3}$ is the set of all smooth maps from $M$ to $N$ preserving 3-cohomology classes under pull-back. This fibration sequence gives the following group extension sequence for the group of \textbf{Courant automorphisms} of a given Courant algebroid:
\begin{equation*}
\begin{tikzcd}
0 \arrow[r] & \Omega^2_{\text{cl}}(M) \arrow[r] & \text{Aut}_{\Crnt_\Man}(E) \arrow[r] & \text{Diff}(M)_{[H_E]} \arrow[r] & 0
\end{tikzcd}
\end{equation*}
where $\text{Diff}(M)_{[H_E]}$ is the group of diffeomorphisms that preserve the Severa class of $E$ under pull-back. In the particular case of a $H$-twisted standard Courant algebroid $\mathbb{T}M_H$, the extension sequence above becomes the semi-direct product of groups
\begin{equation*}
\text{Aut}(\mathbb{T}M_H) \cong \Omega^2(M)\rtimes \text{Diff}(M) |_{\varphi^*H=H-dB}.
\end{equation*}
This shows that the automorphisms of a $H$-twisted standard Courant algebroid are given by pairs $(B,\varphi)$ with $B$ a 2-form acting by addition, commonly referred to as the $B$\textbf{-field transformation} in the literature, and $\varphi$ a diffeomorphism acting by pull-back and preserving the cohomology class of $H$ via the relation $\varphi^* H=H-dB$. Infinitesimally, the group extension sequence for the Courant automorphisms above gives the \textbf{Courant derivations} as a Lie algebra extension sitting in the following short exact sequence:
\begin{equation*}
\begin{tikzcd}
0 \arrow[r] & \Omega_{\text{cl}}^2(M) \arrow[r] & \Dr{E} \arrow[r] & \Sec{\Tan M}_{[H_E]} \arrow[r]  & 0
\end{tikzcd}
\end{equation*}
where $\Sec{\Tan M}_{[H_E]}$ denotes the vector fields preserving the Severa class, $[\LDer_X H_E]=[0]\in H_{\text{dR}}^3(M)$, and $\Omega_{\text{cl}}^2(M)$ is regarded as an abelian Lie algebra. The extension class of $\Dr{E}$ is given by the cocycle in Lie algebra cohomology $c$ defined from any representative of the Severa class $H$ as follows:
\begin{equation*}
c(X,Y)=di_Xi_YH.
\end{equation*}
In the case of a $H$-twisted standard Courant algebroid, the extension sequence above becomes a semi-direct product of Lie algebras:
\begin{equation*}
    \Dr{\mathbb{T}M_H}\cong \Omega_{\text{cl}}^2(M) \rtimes \Sec{\Tan M}|_{\LDer_XH=dB}.
\end{equation*}
Sections of a Courant algebroid $E$ naturally act as Courant derivations on $E$ itself via the usual definition of an \textbf{adjoint map}, $\text{ad}:a\in\Sec{E}\mapsto [a,-]\in\Dr{E}$. The next proposition gives a characterization of $\text{ad}$ as a map of Leibniz algebras.
\begin{prop}[The Adjoint Map of a Courant Algebroid, {\cite[Proposition 2.5]{bursztyn2007reduction}}] \label{AdjointMapCourant}
The adjoint map $\normalfont\text{ad}$ sits in the following exact sequence of Leibniz algebras
\begin{equation*}
\begin{tikzcd}
0 \arrow[r] & \normalfont\Omega_{\text{cl}}^1(M) \arrow[r,"j"] & \normalfont\Sec{E} \arrow[r, "\text{ad}"]& \normalfont\Dr{E} \arrow[r,"\sigma"] & \normalfont H^2_{\text{dR}}(M) \arrow[r] & 0
\end{tikzcd}
\end{equation*}
where $\normalfont\Omega_{\text{cl}}^1(M)$ and $\normalfont H^2_{\text{dR}}(M)$ are regarded as abelian.
\end{prop}
\begin{proof}
Recall that $j=\sharp\rho^*$, then we show that $j(\Omega^1_{\text{cl}}(M))=\Ker{\text{ad}}$ using exactness of the Courant algebroid and the fact that the condition $[k,e]=0$ implies $i_{\rho(e)}d\kappa=0$ for all $e\in \Sec{E}$ but $\rho$ is surjective so $d\kappa =0$. Let us choose an isotropic splitting so that $E\cong \mathbb{T}M_H$, then the Courant derivations become isomorphic to $\Omega_{\text{cl}}^2(M) \rtimes \Sec{\Tan M}|_{\LDer_XH=dB}$ and the image of the adjoint map can be written explicitly as:
\begin{equation*}
\text{ad}(\Sec{E})\cong \{(X,i_XH-d\alpha) \quad X\in \Sec{\Tan M}, \alpha\in \Omega^2(M) \}.
\end{equation*}
Let us define
\begin{align*}
\sigma: \text{Der}(E) & \to H^2_{\text{dR}}(M)\\
(X,b) & \mapsto [i_XH-b]
\end{align*}
which is well-defined, since $d(i_xH-b)=\mathcal{L}_XH+i_XdH-db=db-db=0$, and surjective, since $b$ takes all possible values in $\Omega^2(M)$. Finally, for the kernel of this map we have:
\begin{align*}
\text{ker}(\sigma) & = \{(X,b)\in \text{Der}(E) : [i_XH-b]=[0]\}\\
& = \{(X,b)\in \text{Der}(E) : \exists \alpha\in \Omega^1(M) : i_XH-b-d\alpha=0\}\\
& = \{(X,i_XH-d\alpha) \quad X\in \Sec{\Tan M}, \alpha\in \Omega^2(M) \}=\text{ad}(\Sec{E}).
\end{align*}
\end{proof}

\subsection{Reduction of Courant Algebroids} \label{CourantReduction}

In this section we will treat the problem of reduction in the category of Courant algebroids via Lie group actions. In order to do so, we will generalize the notion of an infinitesimal action of a Lie algebra in the context of Courant algebras. Let $(E,\langle,\rangle,\rho:E\to \Tan M,[,])$ be a Courant algebroid, $\rho':\mathfrak{a}\to \mathfrak{g}$ an exact Courant algebra and recall that a morphism between the Courant algebras $\Sec{E}$ and $\mathfrak{a}$ is given by a commutative diagram
\begin{equation*}
\begin{tikzcd}
0 \arrow[r] & \mathfrak{h} \arrow[r, "i"] \arrow[d, "\nu"'] & \mathfrak{a} \arrow[r, "\rho'"] \arrow[d, "\lambda"'] & \mathfrak{g} \arrow[r] \arrow[d, "\psi"] & 0 \\
0 \arrow[r] & \Sec{\Cot M} \arrow[r, "j"] & \Sec{E} \arrow[r, "\rho"] & \Sec{\Tan M} \arrow[r] & 0 
\end{tikzcd}
\end{equation*}
where $\lambda:\mathfrak{a}\to \Sec{E}$ is a Leibniz algebra morphism and $\psi:\mathfrak{g}\to\Sec{\Tan M}$ is a Lie algebra morphism. We then see how this notion is naturally generalizing that of an infinitesimal action on the base manifold, given by the map $\psi$. The following proposition further connects morphisms of Courant algebras with a generalized notion of actions on Courant algebroids.
\begin{prop}[Extended Actions on a Courant Algebroid] \label{ExtendedActionsCourantAlgebroids}
Let $E$ be a Courant algebroid, $\rho':\mathfrak{a}\to \mathfrak{g}$ an exact Courant algebra and $\lambda:\mathfrak{a}\to \Sec{E}$ a morphism of Courant algebras as above, then, if $\mathfrak{h}:=\Ker{\rho'}$ is mapped trivially into $\Sec{E}$, there is an infinitesimal vector bundle action of $\mathfrak{g}$ on $E$ via Courant derivations $\normalfont\Psi:\mathfrak{g}\to \Dr{E}$. Furthermore, if there is a group action $G\Acts M$ integrating the infinitesimal action $\psi:\mathfrak{g}\to \Sec{\Tan M}$, then $G$ acts on $E$ via Courant automorphisms. We say that $\lambda$ is an \textbf{infinitesimal Courant algebra action extending} the infinitesimal Lie algebra action $\psi$.
\end{prop}
\begin{proof}
Note that the condition that $\mathfrak{h}$ is mapped trivially into $\Sec{E}$, written explicitly $\text{ad}\circ \lambda \circ i (\mathfrak{h)}=0$, is equivalent to $\nu (\mathfrak{h})\subset \Omega^1_{\text{cl}}(M)\subset\Sec{\Cot M}$ as a direct consequence of proposition \ref{AdjointMapCourant}. Also, from the fact that $[\mathfrak{a},\mathfrak{h}]\subset \mathfrak{h} \supset [\mathfrak{h},\mathfrak{a}]$ and $[\mathfrak{h},\mathfrak{h}]=0$, we see that $\mathfrak{a}/\mathfrak{\text{ker}(\rho')}=\mathfrak{a}/\mathfrak{h}$ inherits a well defined Leibniz bracket:
\begin{equation*}
[a+\mathfrak{h},b+\mathfrak{h}]_{\mathfrak{a}/\mathfrak{h}} := [a,b]+\mathfrak{h}.
\end{equation*}
Now, we readily check
\begin{equation*}
[a+\mathfrak{h},b+\mathfrak{h}]_{\mathfrak{a}/\mathfrak{h}} + [b+\mathfrak{h},a+\mathfrak{h}]_{\mathfrak{a}/\mathfrak{h}} = [a,b]+[b,a]+\mathfrak{h}
\end{equation*}
and
\begin{equation*}
\rho'([a,b]+[b,a])=[\rho'(a),\rho'(b)]_\mathfrak{g}+[\rho'(b),\rho'(a)]_\mathfrak{g}=[\rho'(a),\rho'(b)]_\mathfrak{g}-[\rho'(a),\rho'(b)]_\mathfrak{g}=0,
\end{equation*}
which imply $[a,b]+[b,a]\in \mathfrak{h}$, thus showing that the bracket is antisymmetric. This indeed makes the quotient $\mathfrak{a}/\mathfrak{h}$ together with this bracket into a Lie algebra which, by exactness of $\rho'$, is isomorphic to $(\mathfrak{g},[,]_\mathfrak{g})$. From now on we will drop all subindices in the brackets. This further allows to define a map $\tilde{\lambda}:\mathfrak{a}/\mathfrak{h} \to \Gamma(E)$ simply by $\tilde{\lambda}(a+\mathfrak{h}):=\lambda(a)$ which, because $\mathfrak{h}$ acts trivially, makes $\text{ad}\circ \tilde{\lambda}:\mathfrak{a}/\mathfrak{h}\to$Der$(E)$ into a well-defined Lie algebra morphism. We can write the explicit lifted action of $\mathfrak{g}$ on $\Gamma(E)$ as:
\begin{align*}
\Psi: \mathfrak{g}& \to \Gamma(E)\\
g & \mapsto \lambda(a) \quad , \quad \rho'(a)=g.
\end{align*}
This is now a vector bundle infinitesimal $\mathfrak{g}$-action in the usual sense since, from the axioms of Courant algebroid and the defining properties of $\lambda$, we can check that $\text{ad}\circ \Psi:\mathfrak{g} \to \text{Der}(E)$ is a well-defined Lie algebra morphism that further satisfies
\begin{equation*}
    \Psi(g)[f\cdot s]=s\cdot \Psi(g)[s]+\psi(g)[f]\cdot s
\end{equation*}
for all $s\in\Sec{E}$ and $f\in\Cin{M}$. Adjoint derivations are indeed the internal Courant derivations and thus they are integrated by Courant automorphisms, then when the infinitesimal action $\psi$ is integrated by a $G$-action on the base manifold $M$, the infinitesimal action $\Psi$ will be integrated by a vector bundle $G$-action via Courant automorphisms.
\end{proof}

The next proposition shows that there is a natural notion of equivalence between extended actions and that, in general, there are many Courant algebra actions on a fixed Courant algebroid extending a given Lie algebra action. 

\begin{prop}[Equivalence of Extended Actions] \label{EquivalenceExtendedActionsCourant}
Let $\lambda:\mathfrak{a}\to \Sec{E}$ be a Courant algebra action extending an integrable Lie algebra action $\psi:\mathfrak{g}\to \Sec{\Tan M}$ as in proposition \ref{ExtendedActionsCourantAlgebroids} and $f\in \Cin{M,\mathfrak{g}^*}^G$ a $G$-equivariant function, then the map
\begin{equation*}
\lambda_f(a):=\lambda(a)+D\langle f,\rho'(a)\rangle
\end{equation*}
also defines a action $\lambda_f:\mathfrak{a}\to \Sec{E}$ extending $\psi:\mathfrak{g}\to \Sec{\Tan M}$. Two Courant algebra actions extending $\lambda,\lambda':\mathfrak{a}\to \Sec{E}$ both extending $\psi:\mathfrak{g}\to \Sec{\Tan M}$ are called \textbf{equivalent} if there exists $f\in \Cin{M,\mathfrak{g}^*}^G$ such that $\lambda-\lambda=D\langle f,\rho'\rangle$.
\end{prop}
\begin{proof}
By $\langle,\rangle:\mathfrak{g}^*\times \mathfrak{g}\to \mathbb{R}$ we mean the usual dual pairing. For this to be a well-defined equivalence we need to check that, given an extended action $\lambda$ of $\psi$, the map $\lambda_f(-):=\lambda(-)+D\langle f,\rho'(-)\rangle$ is again an extended action of $\psi$. The map $\lambda_f$ clearly fits in the commutative diagram because $\rho(\lambda + D\langle f,\rho'\rangle)=\rho\lambda=\psi\rho'$ as $\rho D = 0$ and $\text{ad}\circ\lambda_f(\mathfrak{h})=\text{ad}\circ(\lambda(\mathfrak{h}))+D\langle f,\rho'(\mathfrak{h})\rangle)= \text{ad} \circ\lambda(\mathfrak{h})=0$. This last remark implies that the induced $\mathfrak{g}$-action will be defined via the same map and hence it will also be integrable. Then, it only remains to verify that $\lambda_f$ is a Leibniz morphism. To see this, let us first derive the infinitesimal counterpart of $G$-equivariance for $f\in \text{C}^\infty(M,\mathfrak{g}^*)$ using the usual identifications $\Tan_{\Id_G}\text{Ad}=\text{ad}$ and $\Tan_{f(x)}\mathfrak{g}^*\cong \Tan_0\mathfrak{g}^*\cong \mathfrak{g}^*$ to find 
\begin{equation*}
f \text{ $G$-equivariant } \quad \Leftrightarrow \quad \langle f, [g,g']\rangle = \langle \Tan f\circ \psi(g), g' \rangle \quad \forall g,g'\in\mathfrak{g}
\end{equation*}
By linearity of the brackets, the new action $\lambda_f$ will be a Leibniz algebra morphism, i.e. $\lambda_f([a,b])=[\lambda_f(a),\lambda_f(b)]$ for all $a,b\in \mathfrak{a}$, when the following equation holds
\begin{equation*}
D\langle f, [\rho'(a),\rho'(b)]\rangle = [\lambda(a),D\langle f,\rho'(b)\rangle]+[D\langle f,\rho'(a)\rangle,\lambda(b)]+[D\langle f,\rho'(a)\rangle,D\langle f,\rho'(b)\rangle].
\end{equation*}
Note that by the properties of the Leibniz bracket of the Courant algebroid the last term vanishes and the two other terms can be rewritten as
\begin{align*}
[\lambda(a),D\langle f,\rho'(b)\rangle]+[D\langle f,\rho'(a)\rangle,\lambda(b)]&=j(\rho(\lambda(a))[\langle f,\rho'(b)]\rangle-i_{\rho(\lambda(b))}dd\langle f,\rho'(a)\rangle)\\
& =jd\rho(\lambda(a))[\langle f,\rho'(b)\rangle]\\
& =D\rho(\lambda(a))[\langle f,\rho'(b)\rangle].
\end{align*}
Note that the Lie derivative of the pairing satisfies $X[\langle f,g\rangle]=\langle \Tan f\circ X,g\rangle$, and then the condition above now becomes simply:
\begin{equation*}
D\langle f, [\rho'(a),\rho'(b)]\rangle = D\langle \Tan f\circ \rho\lambda(a),\rho'(b)\rangle
\end{equation*}
but $\rho \circ \lambda = \psi \circ \rho'$ so this is directly implied by infinitesimal $G$-equivariance of $f$.
\end{proof}

Recall that in  the context of linear Courant spaces proposition \ref{LinearCourantReduction} showed how a subspace in a linear Courant space induced a reduction scheme. The two propositions below show how extended Courant actions integrating to Lie group actions induce an analogous reduction scheme in the context of Courant algebroids. Consider an extended Courant algebra action $\lambda:\mathfrak{a}\to \Sec{E}$, let us identify the following distributions on the Courant algebroid $K:=\lambda(\mathfrak{a})$ and $K^\perp$, and on the tangent bundle
\begin{equation*}
\Delta_s = \rho(K^\perp) \subset \Tan M, \quad \quad \Delta_b = \rho(K+K^\perp) \subset \Tan M,
\end{equation*}
called the \textbf{small distribution} and the \textbf{big distribution} respectively. These are clearly smooth distributions as they are generated by sections.  
\begin{prop}[Distributions of an Extended Actions, {\cite[Lemma 3.2]{bursztyn2007reduction}}] \label{DistributionsExtendedActions}
Let $G$ be a connected Lie group acting on $M$ with infinitesimal action $\psi:\mathfrak{g}\to \Sec{\Tan M}$ and $(E,\langle,\rangle,\rho:E\to \Tan M,[,])$ be a Courant algebroid with infinitesimal Courant algebra action $\lambda:\mathfrak{a}\to \Sec{E}$ extending $\psi$. Assume there is an integral manifold $P\subset M$ of the big distribution,  $TP=\Delta_b|_P$, on which $G$ acts freely and properly and along which $\lambda(\mathfrak{h})$ has constant rank, then $\Delta_s|_P\subset TP$ is an integrable tangent distribution on $P$ and both $K$ and $K\cap K^\perp$ have constant rank along $P$.
\end{prop}

This proposition implies, in particular, that $K|_P$, $K^\perp|_P$ and $K\cap K^\perp|_P$ are all vector bundles over $P$ equipped with a vector bundle $G$-action. This allows us to use proposition \ref{LinearCourantReduction} fibre-wise along $P$ to construct the quotient vector bundle and use the vector bundle $G$-action to form the orbit vector bundle thus giving:
\begin{equation*}
E_r:=\frac{K^\perp|_P}{K\cap K^\perp|_P}\big/ G.
\end{equation*}
which is a vector bundle over the base space $P/G$. Since $G$ acts via Courant automorphisms, we can define an anchor map $\rho_r:E_r \to \Tan (P/G)$ and a bilinear form $\langle,\rangle_r\in\Sec{\odot^2E_r}$ simply from the fibre-wise construction of the quotient. It follows from the general results about vector bundle actions presented at the end of Section \ref{VBundleAutomorphisms} that $\Sec{E_r}\cong \Sec{K^\perp|_P}^G/\Sec{K\cap K^\perp|_P}^G$, then we can give the following explicit definition of a Leibniz bracket on the sections of $E_r$:
\begin{align*}
[ ,]_r : \Sec{E_r} \times \Sec{E_r} &\to \Sec{E_r}\\
(k + \Sec{K\cap K^\perp|_P}^G,k' + \Sec{K\cap K^\perp|_P}^G)& \mapsto [e,e']|_P + \Sec{K\cap K^\perp|_P}^G
\end{align*}
for $e,e'\in \Sec{E}$ any extensions such that $e|_P=k$ and $e'|_P=k'$. We are now in the position to formulate the main reduction theorem for Courant algebroids.
\begin{prop}[Reduction of Courant Algebroids, {\cite[Theorem 3.3]{bursztyn2007reduction}}] \label{ReductionExactCourantAlgebroids}
Let $E$ be a Courant algebroid and $\lambda$ an extended action as in proposition \ref{DistributionsExtendedActions}, then $(E_r,\rho_r,\langle,\rangle_r,[,]_r)$ is a (non-necessarily exact) Courant algebroid over $P/G$ with surjective anchor. $E_r$ is an exact Courant algebroid iff
\begin{equation*}
\rho(K)\cap\rho(K)^\perp=\rho(K\cap K^\perp)
\end{equation*}
is satisfied along $P$. We call $E_r$ the \textbf{reduced Courant algebroid}.
\end{prop}

Recall from our discussion at the end of Section \ref{LinearCourant} that Dirac spaces in linear Courant spaces always reduced to the quotient Courant space. Consider now a Dirac structure of a Courant algebroid $L\subset E$ for which there is a reduction $E_r$ via the extended $G$-action $\lambda:\mathfrak{a}\to \Sec{E}$ as in proposition \ref{ReductionExactCourantAlgebroids}. We say that a Dirac structure $L\subset E$ is \textbf{preserved} by the action $\lambda$ if
\begin{equation*}
[\lambda(\mathfrak{a}),\Sec{L}]\subset \Sec{L}.
\end{equation*}
As a direct consequence of proposition \ref{ExtendedActionsCourantAlgebroids} we see that $L$ becomes a $G$-invariant subbundle and thus one may define a subset $L_r\subset E_r$ simply by applying the reduction of linear Dirac spaces fibre-wise. When $L_r$ is smooth it is clearly a subbundle and if $E_r$ is an exact Courant algebroid it follows from dimension counting that $L_r$ is Lagrangian. Under these assumptions we can see that the invariance condition for $L\subset E$ above implies that $L_r\subset E_r$ is an involutive distribution and thus becomes a Dirac structure of the reduced Courant algebroid. We call $L_r$ the \textbf{reduced Dirac structure}.

\section{Line Bundles} \label{LineBundles}

\subsection{The Category of Line Bundles} \label{CategoryOfLineBundles}

We define the \textbf{category of line bundles} $\Line_\Man$ as the subcategory of $\Vect_\Man$ whose objects are rank $1$ vector bundles $\lambda: L \to M$ and whose morphisms are regular, i.e. fibre-wise invertible, bundle morphisms covering general smooth maps
\begin{equation*}
\begin{tikzcd}
L_1 \arrow[r, "B"] \arrow[d, "\lambda_1"'] & L_2 \arrow[d, "\lambda_2"] \\
M_1 \arrow[r, "b"'] & M_2
\end{tikzcd}
\end{equation*}
In the interest of brevity, we may refer to line bundles $L\in \Line_\Man$ as \textbf{lines} and regular line bundle morphisms $B\in\text{Hom}_{\Line_\Man}(L_1,L_2)$ as \textbf{factors}. This nomenclature aligns with that of Section \ref{CategoryOfLines} as we can indeed regard the category of vector lines $\Line$ as the subcategory of $\Line_\Man$ of line bundles over a single point. A factor covering a diffeomorphism, i.e. a line bundle isomorphism, is called a \textbf{diffeomorphic factor}. Similarly, a factor covering an embedding or submersion is called an \textbf{embedding factor} or \textbf{submersion factor}, respectively. By fixing the base manifold $M$ we can use the tensor product and duality functor of $\Line$ to give fibre-wise analogues for general line bundles over $M$, denoted by $\Line_M$. This makes $(\Line_M,\otimes)$ into a symmetric monoidal category with tensor unit $\Real_M$. The duality assignment $L\in\Line_\Man\mapsto L^*\in\Line_\Man$ is well defined for objects but, following our discussion for general vector bundles at the end of Section \ref{CategoryVectorBundles}, it fails to give a functor of factors in general.\newline

Let a $B:L_{M_1}\to L_{M_2}$ be a factor covering the smooth map $b:M_1\to M_2$, we can define the \textbf{factor pull-back} on sections as:
\begin{align*}
B^*: \Sec{L_{M_2}} & \to \Sec{L_{M_1}}\\
s_2 & \mapsto B^*s_2
\end{align*}
where
\begin{equation*}
    B^*s_2(x_1):=B^{-1}_{x_1}s_2(b(x_1))
\end{equation*}
for all $x_1\in M_1$. Note that the fibre-wise invertibility of $B$ has critically been used for this to be a well-defined map of sections. Linearity of $B$ implies the following interaction of factor pull-backs with the module structures:
\begin{equation*}
    B^*(f\cdot s)=b^*f\cdot B^*s \quad \forall f\in\Cin{M_2}, s\in \Sec{L_{M_2}}.
\end{equation*}
Then, if we consider $\textsf{RMod}_1$, the category of $1$-dimensional projective modules over rings with module morphisms covering ring morphisms and satisfying the pull-back condition above, the assignment of sections becomes a contravariant functor
\begin{equation*}
    \Gamma:\Line_\Man \to \textsf{RMod}_1.
\end{equation*}

Since factors restrict to linear fibre isomorphisms, there is a one-to-one correspondence between monomorphisms in the category $\Line_\Man$ and embeddings of smooth manifolds in the category $\Man$. Indeed, given a submanifold $i:S\hookrightarrow M$ and a line bundle $\lambda:L\to M$, restriction to $S$ gives a factor
\begin{equation*}
\begin{tikzcd}
i^*L \arrow[r, "\iota"] \arrow[d] & L \arrow[d] \\
S \arrow[r,hook, "i"'] & M
\end{tikzcd}
\end{equation*}
where $\iota$ is an embedding factor. The pull-back bundle will often be denoted simply by $L_S$. We may give an algebraic characterization of submanifolds of the base space of a line bundle as follows: for an embedded submanifold $i:S\hookrightarrow M$ we define its \textbf{vanishing submodule} as
\begin{equation*}
    \Gamma_S:=\Ker{\iota^*}=\{s\in\Sec{L_M}|\quad s(x)=0\in L_x \quad \forall x\in S\}.
\end{equation*}
This is the line bundle analogue to the characterization of submanifolds with their vanishing ideals. In fact, the two notions are closely related since (depending on the embedding $i$, perhaps only locally) we have
\begin{equation*}
    \Gamma_S=I_S\cdot \Sec{L_M},
\end{equation*}
where $I_S:=\Ker{i^*}$ is the multiplicative ideal of functions vanishing on the submanifold. This ideal gives a natural isomorphism of rings $\Cin{S}\cong \Cin{M}/I_S$ which, in turn, gives the natural isomorphism of $\Cin{S}$-modules
\begin{equation*}
    \Sec{L_S}\cong \Sec{L}/\Gamma_S.
\end{equation*}

Let $G$ be a Lie group with lie algebra $\mathfrak{g}$ and $\lambda:L\to M$ a line bundle $\lambda:L\to M$, we say that \textbf{$G$ acts on $L$} and denote $G\Acts L$ when there is a smooth map $\Phi:G\times L\to L$ such that $\Phi_g:L\to L$ is a diffeomorphic factor for all $g\in G$ and the usual axioms of a group action are satisfied
\begin{equation*}
    \Phi_g\circ \Phi_h = \Phi_{gh}\qquad \Phi_{e}=\Id_L \qquad \forall g,h\in G.
\end{equation*}
We call the map $\Phi$ a \textbf{line bundle $G$-action}. It follows by construction that any such action induces a standard group action of $G$ on the base space $\phi:G\times M\to M$, which gives the base smooth maps of the line bundle action factors and will be thus called the base action. The infinitesimal counterpart of a line bundle $G$-action is a morphism of Lie algebras
\begin{equation*}
    \Psi:\mathfrak{g}\to \Dr{L}.
\end{equation*}
Note that $G$ acts on both the domain and codomain of the infinitesimal action by the adjoint action and by push-forward respectively; the fact that $\Psi$ is defined as the infinitesimal counterpart of the action $\Phi$ manifest as $G$-equivariance in the following sense
\begin{equation*}
    \Psi\circ \text{Ad}_g = (\Phi_g)_*\circ \Psi \qquad \forall g\in G.
\end{equation*}
Denoting by $\psi:\mathfrak{g}\to \Sec{\Tan M}$ the infinitesimal counterpart of $\phi$, the infinitesimal line bundle action satisfies 
\begin{equation*}
    \Psi(\xi)[f\cdot s]=f\cdot \Psi(\xi)[s]+\psi(\xi)[f]\cdot s \qquad\forall f\in\Cin{M},s\in\Sec{L},\xi\in\mathfrak{g}.
\end{equation*}
The \textbf{orbits} of a line bundle $G$-action $\Phi$ can be simply defined as the images of all group elements acting on a single fibre, and thus they are naturally regarded as the line bundle restricted to the orbits of the base action $\phi$. In analogy with the case of smooth actions, we denote the set of orbits by $L/G$. Any notion defined for usual group actions on smooth manifolds extends to a corresponding notion for line bundle actions simply requiring the base action to satisfy the corresponding conditions. In particular, a \textbf{free and proper} line bundle $G$-action gives a well-defined submersion factor
\begin{equation*}
\begin{tikzcd}
L \arrow[r, "\zeta"] \arrow[d] & L/G \arrow[d] \\
M \arrow[r,two heads, "z"'] & M/G
\end{tikzcd}
\end{equation*}
where the line bundle structure on the set of orbits $L/G$ is induced from the fact that the orbit space of the base action $M/G$ is a smooth manifold and from the fact that all pairs of fibres over the same orbit are mapped isomorphically by some group element. Under these conditions, analogously to the ring isomorphism between functions on the orbit space and invariant functions $\Cin{M/G}\cong \Cin{M}^G$, we find the following isomorphism of modules
\begin{equation*}
    \Sec{L/G}\cong \Sec{L}^G:=\{s\in\Sec{L}| \quad \Phi_g^*s=s \quad \forall g\in G\}.
\end{equation*}

Let us now define the analogue of the Cartesian product of manifolds in the category $\Line_\Man$. Consider two line bundles $\lambda_i:L_i\to M_i$, $i=1,2$, we will use the notations $L_i$ and $L_{M_i}$ indistinctly. We begin by defining the set of all linear invertible maps between fibres:
\begin{equation*}
    M_1 \dtimes M_2:=\{B_{x_1x_2}:L_{x_1}\to L_{x_2}, \text{ linear isomorphism }, (x_1,x_2)\in M_1\times M_2\},
\end{equation*}
we call this set the \textbf{base product} of the line bundles. Let us denote by $p_i:M_1\dtimes M_2 \to M_i$ the obvious projections. This set is easily shown to be a smooth manifold by taking trivializations $L|_{U_i}\cong U_i\times \Real$, $i=1,2$ that give charts of the form $U_1\times U_2\times \Real^\times$ for the open neighbourhoods $(p_1\times p_2)^{-1}(U_1 \times U_2)\subset M_1\dtimes M_2$. In fact, the natural $\Real^\times$-action given by fibre-wise multiplication makes $M_1\dtimes M_2$ into a principal bundle
\begin{equation*}
    \begin{tikzcd}[column sep=0.1em]
        \Real^\times & \Acts & M_1\dtimes M_2  \arrow[d,"p_1\times p_2"]\\
         & & M_1\times M_2.
    \end{tikzcd}
\end{equation*}
Note that the construction of the base product of two line bundles $M_1\dtimes M_2$ as the space of fibre isomorphisms allows to identify the factor $B$ with a submanifold that we may regard as the line bundle analogue of a graph:
\begin{equation*}
    \LGrph{B}:=\{C_{x_1x_2}\in M_1\dtimes M_2|\quad x_2=b(x_1), \quad C_{x_1x_2}=B_{x_1}\},
\end{equation*}
we call this submanifold the \textbf{L-graph} of the factor $B$. We define the \textbf{line product} of the line bundles as $L_1\utimes L_2 :=p_1^*L_1$, which is a line bundle over the base product 
\begin{equation*}
    \lambda_{12}:L_1\utimes L_2\to M_1\dtimes M_2.
\end{equation*}
In the interest of simplicity of notation, the base product $M_1 \dtimes M_2$ does not explicitly keep track of the line bundles $\lambda_i:L_i\to M_i$, however, care has been taken throughout this thesis to always provide sufficient context for an unambiguous interpretation.\newline

Despite the apparent asymmetry of the definition, note that we can define the following factors
\begin{align*}
    P_1(B_{x_1x_2},l_{x_1})&:= l_1\in L_{x_1}\\
    P_2(B_{x_1x_2},l_{x_1})&:= B_{x_1x_2}(l_{x_1})\in L_{x_2}
\end{align*}
where $(B_{x_1x_2},l_{x_1})\in p_1^*L_1$, thus giving the following commutative diagram
\begin{equation}\label{LineProductCommutativeDiagram}
\begin{tikzcd}
L_1 \arrow[d, "\lambda_1"'] & L_1\utimes L_2 \arrow[l,"P_1"']\arrow[d, "\lambda_{12}"]\arrow[r,"P_2"] & L_2\arrow[d,"\lambda_2"] \\
M_1 & M_1 \dtimes M_2\arrow[l,"p_1"]\arrow[r,"p_2"'] & M_2
\end{tikzcd}
\end{equation}
where $P_1$ and $P_2$ are submersion factors. It is clear from this definition that taking the pull-back bundle $p_2^*L_2$ instead of $p_1^*L_1$ as our definition for line product will give the line bundle that we have denoted by $L_2\utimes L_1$. By construction of the base product $M_1\dtimes M_2$ we can construct the following smooth map:
\begin{align*}
c_{12}: M_1\dtimes M_2 & \to M_2\dtimes M_1\\
B_{x_1x_2} & \mapsto B_{x_1x_2}^{-1}
\end{align*}
which is clearly invertible, $c_{12}^{-1}=c_{21}$, and induces a factor on the line products by setting
\begin{equation*}
    C_{12}(B_{x_1x_2},l_{x_1})=(B_{x_1x_2}^{-1},B_{x_1x_2}(l_{x_1})).
\end{equation*}
We have thus found a canonical factor covering a diffeomorphism
\begin{equation*}
\begin{tikzcd}
L_1\utimes L_2 \arrow[r,"C_{12}"]\arrow[d, "\lambda_{12}"] & L_2\utimes L_1\arrow[d, "\lambda_{21}"] \\
M_1 \dtimes M_2\arrow[r,"c_{12}"'] & M_2 \dtimes M_1
\end{tikzcd}
\end{equation*}
hence proving that the definition of line product is indeed symmetrical, i.e the two products are canonically isomorphic as line bundles $L_1\utimes L_2\cong L_2\utimes L_1$. The following proposition establishes the line product as the line bundle version of the Cartesian product.
\begin{prop}[Line Product]\label{LineProduct}
The line product construction 
\begin{equation*}
    \normalfont\utimes:\Line_\Man \times \Line_\Man \to \Line_\Man
\end{equation*}
is a categorical product.
\end{prop}
\begin{proof}
For $\utimes$ to be a categorical product in $\Line_\Man$ first we need the two canonical projection morphisms, these are given by the factors $P_1$ and $P_2$ defined above
\begin{equation*}
\begin{tikzcd}
L_1 & L_1\utimes L_2 \arrow[l,"P_1"']\arrow[r,"P_2"] & L_2.
\end{tikzcd}
\end{equation*}
Now we must check the universal property that, given two morphisms $B_i:L\to L_i$, $i=1,2$, we get a unique morphism $B:L\to L_1 \utimes L_2$ such that the following diagram commutes
\begin{equation*}
\begin{tikzcd}
 & L\arrow[dl,"B_1"'] \arrow[dr,"B_2"]\arrow[d,"B"] & \\
L_1 & L_1\utimes L_2 \arrow[l,"P_1"]\arrow[r,"P_2"'] & L_2.
\end{tikzcd}
\end{equation*}
We can indeed define the \textbf{line product of factors} as
\begin{align*}
B_1\utimes B_2: L & \to L_1\utimes L_2\\
l_x & \mapsto (B_2|_x\circ B_1|_{x}^{-1},B_1(l_x))
\end{align*}
noting that, by assumption, $B_2|_x\circ B_1|_{x}^{-1}:L_{b_1(x)}\to L_{b_2(x)}\in M_1\dtimes M_2$ since factors are fibre-wise invertible. We readily check that setting $B=B_1\utimes B_2$ makes the above diagram commutative, thus completing the proof.
\end{proof}
The construction of the line product as a pull-back bundle over the base product indicates that the sections $\Sec{L_1\utimes L_2}$ are spanned by the sections of each factor, $\Sec{L_1}$, $\Sec{L_2}$, and the functions on the base product. More precisely, we have the following isomorphisms of $\Cin{M_1\dtimes M_2}$-modules:
\begin{equation*}
    \Cin{M_1\dtimes M_2}\cdot P_1^*\Sec{L_1}\cong \Sec{L_1 \utimes L_2} \cong \Cin{M_1\dtimes M_2}\cdot P_2^*\Sec{L_2},
\end{equation*}
where the second isomorphism is induced by the canonical factor $C_{12}$. We can summarize this by writing the image of the line product commutative diagram (\ref{LineProductCommutativeDiagram}) under the section functor
\begin{equation*}
\begin{tikzcd}[row sep=tiny]
\Sec{L_1} \arrow[r,"P_1^*"] & \Sec{L_1\utimes L_2}  & \Sec{L_2}\arrow[l,"P_2^*"'] \\
\bullet & \bullet & \bullet \\
\Cin{M_1}\arrow[r,hook,"p_1^*"'] & \Cin{M_1 \dtimes M_2} & \Cin{M_2}\arrow[l,hook',"p_2^*"]
\end{tikzcd}
\end{equation*}
where $\bullet$ denotes ring-module structure. We clearly see that pull-backs of functions on the Cartesian product of base manifolds $M_1\times M_2$ form a subring of $\Cin{M_1\dtimes M_2}$ but a quick computation introducing trivializations shows that differentials of these do not span the cotangent bundle of the base product everywhere. However, the proposition below shows that it is possible to identify a subspace of spanning functions on $M_1 \dtimes M_2$ defined from local sections of $L_1$ and $L_2$. Consider two trivializations $L|_{U_i}\cong U_i\times \Real$, $i=1,2$ so that the spaces of non-vanishing sections over the trivializing neighbourhoods, denoted by $\Sec{L^\times_i}$, are non-empty. Recall that the ratio maps introduced in Section \ref{CategoryOfLines} allowed us to define functions over the sets of invertible maps between lines. Working fibre-wise, we can use these together with the choice of two local sections to define the \textbf{ratio functions} on $M_1\dtimes M_2$ as follows: let the local sections $s_1\in\Sec{L_1}$ and $s_2\in\Sec{L_2^\times}$, then the function $\tfrac{s_1}{s_2}\in \Cin{M_1\dtimes M_2}$ is defined by
\begin{equation*}
    \frac{s_1}{s_2}(B_{x_1x_2}):=\frac{B_{x_1x_2}(s_1(x_1))}{s_2(x_2)}
\end{equation*}
where we have made use of the map $l^{12}$ in (\ref{ratio2}) in fraction notation. Similarly, using the map $r^{12}$ in (\ref{ratio2}) we can define $\tfrac{s_2}{s_1}\in \Cin{M_1\dtimes M_2}$ for some local sections $s_1\in\Sec{L_1^\times}$ and $s_2\in\Sec{L_2}$. If we consider two non-vanishing local sections $c_1\in\Sec{L_1^\times}$ and $c_2\in\Sec{L_2^\times}$, it is a direct consequence of (\ref{ratio1}) that the relation
\begin{equation*}
    \frac{c_1}{c_2}\cdot \frac{c_2}{c_1}=1.
\end{equation*}
holds on the open set $(p_1\times p_2)^{-1}(U_1 \times U_2)$ where $U_i\subset M_i$ are the neighbourhoods where the local sections are defined.
\begin{prop}[Spanning Functions of $M_1\dtimes M_2$]\label{SpanningFunctionsLineProduct}
Let $B_{x_1x_2}\in M_1\dtimes M_2$ any point in the base product of two line bundles $L_{M_1}$, $L_{M_2}$, consider trivializing neighbourhoods $x_i\in U_i\subset M_1$ so that $\Sec{L_i^\times}$ are non-empty; then the cotangent bundle in a neighbourhood of $B_{x_1x_2}$ is spanned by the differentials of all possible ratio functions, i.e.
\begin{equation*}
\normalfont
    \Cot(M_1\dtimes M_2)=d\frac{\Sec{L_1}}{\Sec{L_2^\times}}
\end{equation*}
when restricted to the open set $(p_1\times p_2)^{-1}(U_1 \times U_2)$.
\end{prop}
\begin{proof}
First note that it is a general fact from basic vector bundle geometry that such open neighbourhoods $U_i$ can be chosen for any points $x_i\in M_i$. The trivializations induce the following diffeomorphism
\begin{equation*}
    (p_1\times p_2)^{-1}(U_1 \times U_2) \cong U_1\times U_2\times \Real^\times
\end{equation*}
and thus we can assign coordinates $(y_1,y_2,b)$, where $y_i\in\Real^{\text{dim}M_i}$ and $b>0$, to any element $B_{x_1x_2}\in (p_1\times p_2)^{-1}(U_1 \times U_2)$. Local sections $s_i\in\Sec{L_i}$ are given by  functions $f_i:U_i\to \Real$ under the trivialization. It follows directly from the definition that ratio functions have the following coordinate expressions
\begin{equation*}
    \frac{s_1}{s_2}(y_1,y_2,b)=f_1(y_1) f_2^{-1}(y_2)b
\end{equation*}
where $f_1$, $f_2$ are the trivialized expressions for the sections $s_1\in\Sec{L_1}$, $s_2\in\Sec{L_2^\times}$. A simple computation shows that
\begin{equation*}
    d_{(y_1,y_2,b)}\frac{s_1}{s_2}=f^{-1}_2b\cdot d_{y_1}f_1-f_1f_2^{-2}b\cdot d_{y_2}f_2+f_1f_2^{-1}\cdot d_bb
\end{equation*}
Since $f_1(y_1)$ is allowed to vanish in general, these differentials span all the cotangent vectors at any coordinate point $(y_1,y_2,b)$.
\end{proof}

In analogy with the constant function $1\in\Cin{M}$ seen as the unit section of the trivial line bundle $\Real_M$, a local non-vanishing section $u\in \Sec{L|_U^\times}$ will be called a \textbf{local unit} of the generically-non-trivial line bundle $L$.

\subsection{The Category of L-Vector Bundles} \label{CategoryOfLVectorBundles}

The notion of L-vector space was introduced in Section \ref{CategoryOfLVectorSpaces} as the natural analogue of vector spaces within the context of the category of lines. The standard definition of the category of vector bundles, as in Section \ref{CategoryVectorBundles}, has objects the fibrations of vector spaces over smooth manifolds and morphisms the smooth maps between the base manifolds with fibre-wise vector space morphisms. In the same spirit, we define a \textbf{L-vector bundle} $\epsilon:E^L\to M$ as a pair $(E,L)$ with $E$ a vector bundle and $L$ a line bundle over the same manifold $M$. The fibres are indeed identified with objects in the category of L-vector spaces $E_x^{L_x}\in\LVect$ for all $x\in M$. A \textbf{morphism of L-vector bundles} $F^B:E_1^{L_1}\to E_2^{L_2}$ is defined as a pair $(F,B)$ with $F:E_1\to E_2$ a vector bundle morphism and $B:L_1\to L_2$ a factor of lines. Naturally, the fibre-wise components of a morphism of L-vector bundles are morphisms of L-vector spaces. We have thus identified the \textbf{the category of L-vector bundles}, which will be denoted by $\LVect_\Man$. The notations $E^L$ and $E_M^L$ will be used indistinctly to denote a L-vector bundle and $L$ may be omitted altogether as superscript if it is understood from the context. Similar notational conventions will apply for L-vector bundle morphisms.\newline

Fixing a base manifold $M$, we find the subcategory of L-vector bundles over a single base $\LVect_M$ with L-vector bundle morphisms covering the identity on $M$, this will be referred to as the restricted category of L-vector bundles. The category $\LVect_M$ admits \textbf{L-direct sums}, \textbf{L-tensor products}, and \textbf{L-duals} which are defined using the fibre-wise L-linear structure and replicating the corresponding constructions of L-vector spaces found in Section \ref{CategoryOfLVectorSpaces}. Equipped with these operations, $\LVect_M$ displays an essentially identical categorical structure to that of $\LVect$. We emphasize at this point that the feature that sets L-vector bundles apart from regular vector bundles is the fact that direct sums can only be taken between L-vector bundles sharing the same line bundle component; tensor products and duality, however, operate in an entirely analogous way.\newline

Let $\epsilon:E^L\to M$ be a L-vector bundle and $\phi:N\to M$ a smooth map, the \textbf{pull-back L-vector bundle} is defined in the obvious way by taking the simultaneous pull-back of both the vector bundle and the line bundle over the same base $\phi^*(E^L):=(\phi^*E)^{(\phi^*L)}$. From the peculiarity of the direct sum of L-vector bundles discussed above one may be discouraged to expect L-vector bundles to admit a generalization of the ordinary vector bundle product, since in the construction of the latter the direct sum plays a crucial role. The next proposition shows that this is, fortunately, not the case, and we can indeed find a categorical product in the category of L-vector bundles.

\begin{prop}[Product of L-Vector Bundles] \label{ProductLVectorBundles}
Let $\epsilon_i: E_i^{L_i}\to M_i$, $i=1,2$, be two L-vector bundles, then the \textbf{product of L-vector bundles} defined by means of the line product as
\begin{equation*}
    E_1^{L_1}\boxplus E_2^{L_2}:=(p_1^*E_1\oplus p_2^*E_2)^{L_1\utimes L_2},
\end{equation*}
is a well-defined a L-vector bundle over the base product $M_1\dtimes M_2$. This construction makes 
\begin{equation*}
    \normalfont\boxplus: \LVect \times \LVect \to \LVect
\end{equation*}
into a categorical product.
\end{prop}
\begin{proof}
To see that $E_1^{L_1}\boxplus E_2^{L_2}$ is well-defined we first note that the commutative diagram \ref{LineProductCommutativeDiagram} allows for L-vector bundle pull-backs $p_i^*E_i^{L_i}\to M_1\dtimes M_2$ which, by construction of the line product $L_1 \utimes L_2$, share the same line bundle component. It is then possible to take their direct sum in the restricted category $\LVect_{M_1\dtimes M_2}$. In order to show that $\boxplus$ is a categorical product it will suffice to construct the product of two morphisms of L-vector bundles $ F_i^{B_i}: E^L \to E_i^{L_i}$, $i=1,2$. The line product of factors of proposition \ref{LineProduct} allows for an explicit point-wise definition:
\begin{equation*}
    F_1^{B_1}\boxplus F_2^{B_2}(a^l):=(B_2|_x\circ B_1|_x^{-1},F_1|_x(a)\oplus F_2|_x(a))^{B_1\utimes B_2(l)},
\end{equation*}
for an arbitrary element $a^l\in E^L|_x$.
\end{proof}

\subsection{Unit-Free Manifolds} \label{UnitFreeManifolds}

Recall that in Section \ref{CategoryOfLines} we argued that objects of $\Line$, i.e. $1$-dimensional vector spaces, were to be thought of as ``unit-free'' fields of real numbers. Similarly, the identification of the category $\Line_\Man$ can be reinterpreted as an attempt to find a notion of ``unit-free'' manifold. This section is aimed at making this interpretation clear.\newline 

One thinks of a manifold $M$ as a smooth space which carries a ring of functions that has a unit, the global constant function $1_M\in\Cin{M}$, with the assignment $\text{C}^\infty:\Man\to \Ring$ being a contravariant functor sending smooth maps to pull-backs of functions. By taking an object in the category of line bundles $L_M\in\Line_\Man$ as our replacement for smooth space and its sections $\Sec{L_M}\in\textsf{RMod}_1$, a $1$-dimensional projective $\Cin{M}$-module, as our replacement for functions, we find $L_M$ to be our candidate for the notion of \textbf{unit-free smooth manifold}. Note that there is no canonical choice of global unit section in the module of sections of a generic line bundle $\Sec{L_M}$ so, indeed, our choice of line bundles as unit-free manifolds seems, at least, reasonable. In Section \ref{CategoryOfLineBundles} it was argued that line bundles exhibit notions of subobjects, quotients by group actions and products that are entirely analogous to the ones found in the category of smooth manifolds; which further support our proposed definition. L-vector bundles were introduced in Section \ref{CategoryOfLVectorBundles} as the direct analogues of vector bundles within the category of line bundles. The natural constructions of tangent and cotangent bundles of conventional manifolds establish a strong connection between the category of smooth manifolds and the category of vector bundles - a summary of these can be found in Section \ref{TangentFunctor}. Then, the natural question arises of whether unit-free manifolds (line bundles) exhibit an analogous connection with L-vector bundles. The reminder of this section is devoted to proving the necessary results to positively support the conclusion that such a connection between line bundles and L-vector bundles does indeed exist. Results of this section should be compared with those of Section \ref{TangentFunctor}.\newline

Following the construction of the bundle of derivations of a vector bundle given in Section \ref{DifferentialOperators}, one immediately sees that it can be explicitly defined as the \textbf{der bundle}:
\begin{equation*}
    \Der L:=\{a_x:\Sec{L}\to L_x|\quad a_x(f\cdot s)=f(x)a_x(s)+\delta(a_x)(f)(x)s(x)\quad \forall f\in \Cin{M},s\in\Sec{L}\}
\end{equation*}
where $\delta:\Der A\to \Tan M$ is a vector bundle map covering the identity, defined from the symbol of derivations $\sigma$. This is clearly a vector bundle whose fibres realize the unit-free version of local directional derivatives. In analogy with the usual notation of the action a vector field $X$ on a function $f$ as a directional derivative $X[f]$, given a line bundle $L_M$ we define the action of a section of the der bundle, $a\in \Sec{\Der L}$, as a derivation on sections by setting
\begin{equation*}
    a[s](x):=a(x)(s)\in L_x, \quad \forall x\in M.
\end{equation*}
We then write the Lie bracket of derivations as a bracket on sections of $\Der L$ explicitly
\begin{equation*}
    [a,a'][s]:=a[a'[s]]-a'[a[s]]
\end{equation*}
for all $a,a'\in \Sec{\Der L}$ and $s\in\Sec{L}$. Since sections of $L$ are interpreted as the line bundle version of functions, the Lie algebroid of derivations $(\Der L,\delta,[,] )$ appears as an obvious candidate for the generalization of the tangent bundle for a line bundle.\newline

Let a factor $B:L_1\to L_2$ covering a smooth map $b:M_1\to M_2$, we define its \textbf{der map} at a point $x_1\in M_1$ as follows:
\begin{align*}
\Der_{x_1} B: \Der_{x_1} L_1 & \to \Der_{b(x_1)} L_2\\
a_{x_1} & \mapsto \Der_{x_1}a_{x_1} 
\end{align*}
where
\begin{equation*}
    \Der_{x_1}a_{x_1}(s_2):=B_{x_1}(a_{x_1}(B^*s_2)) \quad \forall s_2\in\Sec{L_2}.
\end{equation*}
Given two derivations, $a_1\in\Sec{\Der L_1}$ and $a_2\in\Sec{\Der L_2}$ we say that they are $B$\textbf{-related}, and denote $a_1\sim_B a_2$, if the following diagram commutes
\begin{equation*}
    \begin{tikzcd}
    \Der L_1 \arrow[r, "\Der B"] & \Der L_2 \\
    M_1 \arrow[r, "b"']\arrow[u,"a_1"] & M_2.\arrow[u,"a_2"]
\end{tikzcd}
\end{equation*}
When the base map of the factor $b$ is a diffeomorphism we can define the \textbf{der push-forward} of derivations by
\begin{align*}
B_*: \Sec{\Der L_1} & \to \Sec{\Der L_2}\\
a & \mapsto \Der B \circ a \circ b^{-1}.
\end{align*}
With our identification of derivations with sections of the der bundle and the definition of pull-back of sections, we can readily check the following identity:
\begin{equation*}
    B_*a[s]=(B^{-1})^*a[B^*s] \qquad \forall s\in \Sec{L_2},
\end{equation*}
which gives an alternative definition of der push-forward only in terms of pull-backs of diffeomorphic factors.\newline

We see that the der bundle indeed recovers the two usual interpretations of the tangent bundle: first as the bundle realizing infinitesimal symmetries but also the bundle whose fibres are isomorphic to the spaces local directional derivatives. The propositions that follow below will show that the assignment of the der bundle to a line bundle has analogous functorial properties to those of the tangent functor, thus establishing the choice of the der bundle of a line bundle as the correct generalization of the tangent bundle.

\begin{prop}[The Der Functor]\label{DerFunctor}
The der map of a factor $B:L_1 \to L_2$ gives a well-defined L-vector bundle morphism
\begin{equation*}
\normalfont
\begin{tikzcd}
\Der L_1 \arrow[r, "\Der B"] \arrow[d, "\lambda_1"'] & \Der L_2 \arrow[d, "\lambda_2"] \\
M_1 \arrow[r, "b"'] & M_2
\end{tikzcd}
\end{equation*}
which is a Lie algebroid morphism, i.e.
\begin{equation*}
\normalfont
    \Tan b\circ \delta_1 = \delta_2\circ \Der B \qquad  \qquad a_1\sim_B a_2 ,\quad a'_1\sim_B a'_2 \quad \Rightarrow \quad [a_1,a'_1]_1\sim_B [a_2,a'_2]_2
\end{equation*}
for $\normalfont a_i,a_i'\in\Sec{\Der L_i}$. Furthermore, for any other factor $F:L_2 \to L_3$ and $\normalfont \Id_L:L\to L$ the identity factor, we have
\normalfont
\begin{align*}
    \Der (F\circ B) & = \Der F \circ \Der B\\
    \Der (\Id_L) & = \Id_{\Der L}.
\end{align*}
\end{prop}
\begin{proof}
Note first that  $B^*(f\cdot s)=b^*f\cdot B^*s$ for $f\in\Cin{M_2}$ and $s\in\Sec{L_2}$, then applying the definition of point-wise der map one obtains
\begin{equation*}
    \Der_xBa_x(f\cdot s)=f(b(x))\cdot \Der_x B a_x(s)+\delta_1(a_x)(b^*f)(x)\cdot s(b(x)).
\end{equation*}
Simply using the definition of tangent map $\Tan b$ we can rewrite the second term to find
\begin{equation*}
    \Der_xBa_x(f\cdot s)=f(b(x))\cdot \Der_x B a_x(s)+\Tan_xb\circ \delta_1(a_x)(f)(b(x))\cdot s(b(x))
\end{equation*}
which shows that indeed $\Der_xBa_x\in \Der_{b(x)}L_2$, making $\Der B$ a well-defined vector bundle morphism. We can then use the anchor $\delta_2$ to rewrite the LHS
\begin{equation*}
    \Der_xBa_x(f\cdot s)=f(b(x))\cdot \Der_x B a_x(s)+\delta_2(\Der_xB a_x)(f)(b(x))\cdot s(b(x))
\end{equation*}
and thus we obtain
\begin{equation*}
    \Tan_xb\circ \delta_1(a_x)=\delta_2(\Der_xB a_x),
\end{equation*}
which must hold for any $a_x\in\Der_xL_1$, thus giving the desired compatibility condition with the anchors. The $B$-relatedness condition follows directly from the definition of der map by noting that
\begin{equation*}
    a_2(b(x))(a'_2[s])=(\Der_x B a_1(x))(a'_2[s])=B_x(a_1(x)(B^*a_2[s]))
\end{equation*}
for all $s\in\Sec{L_2}$. The two functorial identities follow from contravariance of pull-backs $(F\circ B)^*=B^*\circ F^*$ and the trivial fact $\Id_L^*=\Id_{\Sec{L}}$.
\end{proof}

\begin{prop}[Der Bundle of the Dual] \label{DerBundleDual}
Let $E$ be a vector bundle and its dual $E$, then there is a canonical Lie algebroid isomorphism covering the identity
\begin{equation*}
\normalfont
    R:\Der E \to \Der E^*.
\end{equation*}
Furthermore, when $E$ is a line bundle, $R$ is also an isomorphism of L-vector bundles
\end{prop}
\begin{proof}
Consider the vector bundle $E\to M$ and let us denote the der bundle Lie algebroid structures by $(\Der E,\delta,[,])$ and $(\Der E,\delta_*,[,]_*)$, then we can define the Lie algebroid isomorphism explicitly on sections by setting
\begin{equation*}
    R(a)[\epsilon](e):=\delta(a)[\epsilon(e)]-\epsilon(a[e])
\end{equation*}
for $a\in \Dr{E}$ and $e\in\Sec{E}$, $\epsilon\in\Sec{E^*}$. This is clearly $\Cin{M}$-linear by construction and thus we simply have to check compatibility with the anchors and the Lie brackets. Compatibility with the anchors is evident by construction since one readily checks
\begin{equation*}
    \delta_*(R(a))[f]=\delta(a)[f]
\end{equation*}
for an arbitrary function $f$ that can be given as a pairing of sections $f=\epsilon(e)\in\Cin{M}$ without loss of generality. Compatibility with the Lie brackets follows by direct computation: first, using the fact that the anchor $\delta$ is a Lie algebra morphism write
\begin{equation*}
    R([a,b])[\epsilon](e)=[\delta(a),\delta(b)][\epsilon(e)]-\epsilon([a,b][e]),
\end{equation*}
then simply by definition
\begin{equation*}
    R(a)\circ R(b)[\epsilon](e)=\delta(a)[\delta(b)[\epsilon(e)]-\delta(a)[\epsilon( b[e])]-\delta(b)[\epsilon( a[e])]+\epsilon(a\circ b[e])
\end{equation*}
which, after cancellations, gives the desired result
\begin{equation*}
    R([a,b])[\epsilon](e)=R(a)\circ R(b)[\epsilon](e)-R(b)\circ R(a)[\epsilon](e)=[R(a),R(b)][\epsilon](e).
\end{equation*}
\end{proof}

\begin{prop}[Der Bundle of a Submanifold]\label{DerBundleLineSubmanifold}
Let $\lambda:L\to M$ be a line bundle and $i:S\hookrightarrow M$ an embedded submanifold of the base. The canonical embedding factor $\iota : L_S\hookrightarrow L$ gives an injective L-vector bundle morphism
\begin{equation*}
\normalfont
    \Der \iota : \Der L_S \hookrightarrow \Der L.
\end{equation*}
By a slight abuse of notation we will write $\normalfont\Der L_S\subset \Der L$ in the same way that we write $\normalfont\Tan S\subset \Tan M$, then we have
\begin{equation*}
\normalfont
    \delta_L(\Der L_S)=\Tan S,
\end{equation*}
where $\normalfont\delta_L:\Der L\to \Tan M$ is the anchor of the der bundle.
\end{prop}
\begin{proof}
That $\Der \iota$ is a L-vector bundle morphism follows by construction since $L_S:=i^*L$, so the line bundle morphism is simply the fibre-wise identity map. Then, injectivity of $\Der \iota$ follows simply from injectivity of $i:S\hookrightarrow M$. Note that using the full notation the second identity in the proposition reads
\begin{equation*}
    \delta_L(\Der \iota (\Der L_S))=\Tan i(\Tan S)
\end{equation*}
which is clearly a direct consequence of $\Der \iota$ being a Lie algebroid morphism, in particular compatible with the anchors, and the anchors being surjective so that $\delta_L(\Der L)=\Tan M$ and $\delta_{L_S}(\Der L_S)=\Tan S$.
\end{proof}

\begin{prop}[Derivations of a Submanifold]\label{DerivationsLineSubmanifold}
Let $\lambda:L\to M$ be a line bundle, $i:S\hookrightarrow M$ an embedded submanifold of the base and denote by $\Gamma_S\subset\Sec{L}$ its vanishing submodule. We define the derivations that tangentially restrict to $S$ as
\begin{equation*}
\normalfont
    \text{Der}_S(L):=\{D\in\Dr{L}|\quad D[\Gamma_S]\subset \Gamma_S\}
\end{equation*}
and the derivations that vanish on $S$ as
\begin{equation*}
\normalfont
    \text{Der}_{S_0}(L):=\{D\in\Dr{L}|\quad D[\Sec{L}]\subset \Gamma_S\},
\end{equation*}
then there is a natural isomorphism of $\normalfont\Cin{M}$-modules and Lie algebras
\begin{equation*}
\normalfont
    \Dr{L_S}\cong \text{Der}_S(L)/\text{Der}_{S_0}(L). 
\end{equation*}
\end{prop}
\begin{proof}
The isomorphism as modules follows directly from proposition \ref{DerBundleLineSubmanifold} using the correspondence between sections of the der bundle and derivations. The Lie algebra isomorphism is then a consequence of the following simple facts
\begin{align*}
    [\text{Der}_S(L),\text{Der}_S(L)] & \subset  \text{Der}_S(L)\\
    [\text{Der}_S(L),\text{Der}_{S_0}(L)] & \subset  \text{Der}_{0S}(L)\\
    [\text{Der}_{0S}(L),\text{Der}_{S_0}(L)] & \subset  \text{Der}_{0S}(L),
\end{align*}
easily derived from the definitions above, thus showing that $\text{Der}_S(L)\subset \Dr{L}$ is a Lie subalgebra and $\text{Der}_{0S}(L)\subset\text{Der}_S(L)$ is a Lie ideal making the subquotient $\text{Der}_S(L)/\text{Der}_{0S}(L)$ into a Lie algebra reduction.
\end{proof}

\begin{prop}[Der Bundle of a Group Action Quotient]\label{DerBundleGroupAction}
Let $\Phi:G\times L\to L$ be a free and proper line bundle $G$-action with infinitesimal counter part $\normalfont \Psi:\mathfrak{g}\to \Dr{L}$ and let us denote by $\zeta: L\to L/G$ the submersion factor given by taking the quotient onto the space of orbit line bundles. Defining the $G$-invariant derivations as
\begin{equation*}
\normalfont
    \Dr{L}^G:=\{D\in\Dr{L}|\quad (\Phi_g)_*D=D \quad \forall g\in G\}=\{D\in\Dr{L}|\quad D[\Sec{L}^G]\subset \Sec{L}^G\}
\end{equation*}
and the derivations that tangentially restrict to the orbits as
\begin{equation*}
\normalfont
     \Dr{L}^{G_0}:=\{D\in\Dr{L}|\quad D[\Sec{L}^G]=0\}
\end{equation*}
we find the following isomorphism of modules and Lie algebras
\begin{equation*}
\normalfont
    \Dr{L/G}\cong \Dr{L}^G/\Dr{L}^{G_0}.
\end{equation*}
When $G$ is connected, this isomorphism becomes
\begin{equation*}
\normalfont
    \Dr{L/G}\cong \Dr{L}^{\mathfrak{g}}/\Psi(\mathfrak{g}).
\end{equation*}
where the $\mathfrak{g}$-invariant derivations are defined as
\begin{equation*}
\normalfont
    \Dr{L}^{\mathfrak{g}}:=\{D\in\Dr{L}|\quad [\Psi(\xi),D]=0 \quad \forall \xi\in \mathfrak{g}\}.
\end{equation*}
When applied point-wise, this last isomorphism gives a fibre-wise isomorphism of der spaces over the base orbit space $M/G$:
\begin{equation*}
\normalfont
    \Der_{[x]}(L/G)\cong \Der_xL/\mathfrak{g}.
\end{equation*}
\end{prop}
\begin{proof}
The module isomorphism $\Dr{L/G}\cong \Dr{L}^G/\Dr{L}^{G_0}$ is a direct consequence of the ring isomorphism $\Cin{M/G}\cong \Cin{M}^G$ and the module isomorphism $\Sec{L/G}\cong\Sec{L}^G$. The isomorphism as Lie algebras follows from the fact that, by construction, we have
\begin{align*}
    [\Dr{L}^G,\Dr{L}^G] & \subset  \Dr{L}^G\\
    [\Dr{L}^G, \Dr{L}^{0G}] & \subset   \Dr{L}^{G_0}\\
    [ \Dr{L}^{0G}, \Dr{L}^{0G}] & \subset   \Dr{L}^{G_0},
\end{align*}
thus showing that $\Dr{L}^G\subset \Dr{L}$ is a Lie subalgebra and $\Dr{L}^{0G}\subset\Dr{L}^G$ is a Lie ideal making the subquotient $\Dr{L}^G/\Dr{L}^{0G}$ into a Lie algebra reduction. When $G$ is connected, its action is uniquely specified by the infinitesimal counterpart, which can be equivalently regarded as a Lie algebroid morphism $\Psi:(\mathfrak{g}_M,\psi,[,])\to (\Der L,\delta,[,])$, where $\mathfrak{g}_M=M\times \mathfrak{g}$ here denotes the action Lie algebroid with anchor given by the infinitesimal action on the base manifold $\psi:\mathfrak{g}\to \Sec{\Tan M}$. It is clear by construction that for the action of a connected $G$ we have $\Dr{L}^{0G}=\Psi(\mathfrak{g})$ and that $G$-invariance under push-forward, $(\Phi_g)_*D=D$, becomes vanishing commutator with the infinitesimal action, $[\Psi(\xi),D]=0$. Then the second isomorphism follows. Since the action is free and proper, the map $\Psi$ will be injective as a vector bundle morphism covering the identity map on $M$, then applying the second isomorphism point-wise and injectivity of the infinitesimal action, so that $\Psi(\mathfrak{g})_x\cong \mathfrak{g}$, we find the last desired result.
\end{proof}

\begin{prop}[Der Bundle of the Line Product]\label{DerLineProduct}
Let two line bundles $L_1, L_2\in \Line_\Man$ and take their line product $L_1\utimes L_2$, then there is a canonical isomorphism of L-vector bundles covering the identity map on the base product:
\begin{equation*}
\normalfont
    \Der (L_1\utimes L_2)\cong \Der L_1 \boxplus \Der L_2.
\end{equation*}
\end{prop}
\begin{proof}
Note first that both vector bundles are, by construction, vector bundles over the same base manifold $M_1\dtimes M_2$. The first term $p_1^*\Der L_1$ is clearly an L-vector bundle with line $p_1^*L_1$ and although the second term is an L-vector bundle with line $p_2^*L_2$ we can use the canonical factor $C_{12}$ to isomorphically regard $p_2^*\Der L_2$ as an L-vector bundle with line $p_1^*L_1$. Then the direct sum as L-vector bundles is well-defined and giving an L-vector bundle isomorphism now reduces to finding a vector bundle isomorphism. A quick check introducing trivializations (see the end of this section for details) shows that both vector bundles have the same rank and so it is enough to find a fibre-wise surjective map $\Phi$ between the vector bundles. We can write this map explicitly as
\begin{align*}
\Phi_{B_{x_1x_2}}:\Der_{B_{x_1x_2}}(L_1\utimes L_2) & \to \Der_{x_1}L_1\oplus \Der_{x_2}L_2\\
a & \mapsto \Der_{B_{x_1x_2}}P_1(a) \oplus \Der_{B_{x_1x_2}}P_2(a)
\end{align*}
where $P_1$ and $P_2$ are the line product projection factors as in (\ref{LineProductCommutativeDiagram}). Surjectivity of this map follows directly from the definition of the der map and the fact that projection factors are surjective. The line bundle component of the isomorphism of L-vector bundles is given simply by the fact that the L-vector bundle product $\boxplus$ is defined as a L-vector bundle with line bundle component given by the line product $\utimes$.
\end{proof}
\begin{prop}[Derivations of the Line Product]\label{DerivationsLineProduct}
Let two line bundles $L_1, L_2\in \Line_\Man$ and take their line product $L_1\utimes L_2$, then we find the derivations of each factor as submodules of the derivations of the product
\begin{equation*}
\normalfont
\begin{tikzcd}
\Dr{L_1}\arrow[r,hook, "k_1"] & \Dr{L_1\utimes L_2} & \Dr{L_2}\arrow[l,hook',"k_2"'].
\end{tikzcd}
\end{equation*}
The maps $k_i$ are Lie algebra morphisms making $\normalfont \Dr{L_i}\subset\Dr{L_1\utimes L_2}$ into Lie subalgebras which, furthermore, satisfy
\begin{equation*}
\normalfont
    [\Dr{L_1},\Dr{L_1}]\subset \Dr{L_1} \qquad [\Dr{L_1},\Dr{L_2}]=0 \qquad [\Dr{L_2},\Dr{L_2}]\subset\Dr{L_2}.
\end{equation*}
\end{prop}
\begin{proof}
This can be proved using the isomorphism in proposition \ref{DerLineProduct} above and the fact that sections of the der bundle are identified with derivations. However, we give an independent, more algebraic proof that explicitly involves the module structure of the sections of the line product. Recall that sections of the line product $\Sec{L_1\utimes L_2}$ are spanned by the pull-backs $P^*_i\Sec{L_i}$ over the functions on the base product $\Cin{M_1\dtimes M_2}$. This means that a derivation $D$ is characterised by its action on projection pull-backs and the action of its symbol $X$ on spanning functions of the base product, which are the ratio functions defined in proposition \ref{SpanningFunctionsLineProduct}. With this in mind, given derivations $D_i\in\Dr{L_i}$ we give the derivations on the line product $k_i(D_i)=\overline{D}_i\in\Dr{L_1\utimes L_2}$ determined uniquely by the conditions
\begin{align*}
    \overline{D}_1(P^*_1s_1)&=P^*D_1(s_1) &\overline{D}_1(P^*_2s_2)&=0\\
    \overline{D}_2(P^*_1s_1)&=0 &\overline{D}_2(P^*_2s_2)&=P^*_2D_2(s_2)
\end{align*}
and with symbols $\overline{X}_i\in\Sec{T(M_1\dtimes M_2)}$ defined on ratio functions by
\begin{align*}
    \overline{X}_1[\tfrac{s_1}{u_2}]=\tfrac{D_1(s_1)}{u_2}\\
    \overline{X}_2[\tfrac{s_2}{u_1}]=\tfrac{D_2(s_2)}{u_1}
\end{align*}
for all sections $s_i\in\Sec{L_i}$ and local non-vanishing sections $u_i\in\Sec{L_i^\times}$. The fact that $k_i$ are Lie algebra morphisms and that $[\overline{D}_1,\overline{D}_2]=0$ follows directly from the defining conditions above.
\end{proof}

In summary, we find the \textbf{der functor} for line bundles
\begin{equation*}
    \Der : \Line_\Man \to \Lie_\Man\subset \LVect_\Man,
\end{equation*}
which plays a categorical role entirely analogous to that of the tangent functor for smooth manifolds as described Section \ref{TangentFunctor}.\newline

Recall that, for a conventional manifold, one defines the cotangent bundle as the vector bundle of point-wise differentials of functions, which is indeed dual to the tangent bundle $\Cot M\cong (\Tan M)^*$. Differentials of functions give sections of this bundle and one defines the exterior algebra from the exterior derivative $d:\Cin{M}\to \Sec{\Cot M}$, which also induces the seminal identity of Cartan calculus on $M$: $X[f]=i_Xdf=\LDer_Xf$. When we consider a section of a line bundle $s\in\Sec{L}$ as the generalization of a unit-free function, the analogue of a point-wise differential is the 1-jet value $j_xs\in \Jet^1L$, thus hinting to the bundle of 1-jets as the analogue of the cotangent bundle for unit-free manifolds. It follows from the definitions presented in Section \ref{DifferentialOperators} and proposition \ref{LocalLieRk1} that the \textbf{jet bundle} $\Jet^1 L$ indeed corresponds to the L-dual to the der bundle
\begin{equation*}
    \Jet^1 L:=(\Der L)^*\otimes L = (\Der L)^{*L}.
\end{equation*}
The action of a section of the der bundle $a\in\Sec{\Der L}$ on a section $s\in\Sec{L}$ as a local derivation can now be rewritten as
\begin{equation*}
    a[s]=j^1s(a)
\end{equation*}
where the jet prolongation $j^1:\Sec{L}\to \Sec{\Jet^1 L}$ gives the unit-free generalization of the exterior derivative. In fact, the Cartan calculus of the Lie algebroid $\Der L$ with the natural representation on $L$ shows that this map extends to a graded differential on the \textbf{der complex} of the line bundle $(\Sec{\wedge^\bullet(\Der L)^*\otimes L},d_L)$. These remarks confirm that the jet bundle of a line bundle is indeed the unit-free analogue of the cotangent bundle.\newline

The symbol and Spencer short exact sequences, defined in general for differential operators between vector bundles in Section \ref{DifferentialOperators}, find a straightforward reformulation in the category of L-vector bundles. Note that, for any line bundle $\lambda:L\to M$, the anchor map $\delta:\Der L\to \Tan M$ is surjective by construction and elements of its kernel correspond to bundle endomorphisms induced from fibre-wise scalar multiplication, $\Ker{\delta}\cong \Real_M$. Then, regarding $\Der L$ and $\Tan M$ as L-vector bundles and $\delta$ as a L-vector bundle morphism in a trivial way, the symbol sequence is simply the short exact sequence of L-vector bundles induced by the fact that $\delta$ is surjective
\begin{equation*}
\begin{tikzcd}
0 \arrow[r] & \Real_M \arrow[r] & \Der L \arrow[r, "\delta"] & \Tan M \arrow[r] & 0.
\end{tikzcd}
\end{equation*}
This will be sometimes called the \textbf{der sequence} of the line bundle $L$. The Spencer sequence then corresponds precisely to the L-dual of the short exact sequence above
\begin{equation*}
\begin{tikzcd}
0 &  L \arrow[l] & \Jet^1L \arrow[l]& (\Tan M)^{*L} \arrow[l,"i"'] & 0 \arrow[l]
\end{tikzcd}
\end{equation*}
where $i=\delta^{L*L}$ is injective. We will refer to this sequence as the \textbf{jet sequence} of the line bundle $L$. We will use the notation $\Tan^{*L} M:= (\Tan M)^{*L}$ in analogy with the usual notation for cotangent bundles.\newline

We end this section by giving a summary of the results proved above for the particular case of trivial line bundles. Note that this discussion will also apply to general line bundles when restricted to trivializing neighbourhoods. Sections of a trivial line bundle $\Real_M$ are isomorphic to the smooth functions of the base $\Sec{\Real_M}\cong \Cin{M}$ with the module structure being simply point-wise multiplication, this can be seen explicitly on generic line bundles by the choice of a local unit. This implies that there is now a natural inclusion $\Dr{\Cin{M}}\subset\Dr{\Real_M}$ making the der short exact sequence split and thus giving an isomorphism of modules
\begin{equation*}
    \Dr{\Real_M} \cong \Sec{\Tan M}\oplus \Cin{M}.
\end{equation*}
The action of a derivation on a function $s\in \Sec{\Real_M}$, regarded as a section of the trivial bundle $\Real_M$, is given by
\begin{equation*}
    (X\oplus f)[s]=X[s]+fs.
\end{equation*}
It then trivially follows that the der bundle is
\begin{equation*}
    \Der \Real_M \cong \Tan M \oplus \Real_M
\end{equation*}
with anchor $\delta=\Proj_1$ and Lie bracket bracket
\begin{equation*}
    [X\oplus f,Y\oplus g]=[X,Y]\oplus X[g]-Y[f].
\end{equation*}
Note how this Lie bracket is entirely determined by the fact that vector fields are the Lie algebra of derivations on functions. Taking $\Real_M$-duals corresponds to taking ordinary duals, therefore the jet bundle is
\begin{equation*}
    \Jet^1 \Real_M\cong \Cot M\oplus \Real_M.
\end{equation*}
The jet prolongation map then becomes
\begin{align*}
j^1: \Cin{M} & \to \Sec{\Cot M}\oplus \Cin{M}\\
s & \mapsto ds\oplus s,
\end{align*}
which indeed only carries the information of the ordinary exterior differential. The base product of two trivial line bundles $\Real_M$ and $\Real_N$ is
\begin{equation*}
    M\dtimes N\cong M\times N\times \Real^\times
\end{equation*}
and the line product is again a trivial line bundle
\begin{equation*}
    \Real_M\utimes \Real_N \cong \Real_{M\dtimes N}.
\end{equation*}
A factor between trivial line bundles $B:\Real_M\to \Real_N$ is given by a pair $B=(b,\beta)$ with $b:M\to N$ a smooth map and $\beta\in\Cin{M}$ a nowhere-vanishing function. We have explicitly
\begin{align*}
(b,\beta): \Real_M & \to \Real_N\\
(x,l) & \mapsto (b(x),\beta(x)l).
\end{align*}
Pull-backs then become
\begin{equation*}
    (b,\beta)^*s=\tfrac{1}{\beta}\cdot b^*s
\end{equation*}
for all $s\in\Cin{N}$. A simple computation shows that the der map of a factor $(b,\beta)$ gives a map of the form
\begin{align*}
\Der (b,\beta): \Tan M\oplus \Real_M & \to \Tan N\oplus \Real_N\\
v_x\oplus a_x & \mapsto \Tan_x b(v_x)\oplus a_x-\tfrac{1}{\beta(x)}d_x\beta(v_x).
\end{align*}
For a diffeomorphic factor $(b,\beta):\Real_M\to \Real_M$, i.e. when $b$ is a diffeomorphism, the der push-forward of a derivation is given by
\begin{equation*}
    (b,\beta)_*(X\oplus f)=b_* X \oplus b_*f +\beta\cdot X[\tfrac{1}{\beta}].
\end{equation*}

\section{Jacobi, Precontact and Contact Geometry} \label{ContactGeometry}

One of the strengths of the notion of unit-free manifold introduced in Section \ref{UnitFreeManifolds} is that it allows to identify Jacobi, precontact and contact manifolds as the analogues of Poisson, presymplectic and symplectic manifolds, respectively. In order to make this claim precise, this section follows a structure that mirrors our presentation of Poisson and presymplectic geometry in Section \ref{SymplecticGeometry}.\newline

A \textbf{Jacobi manifold} or \textbf{Jacobi structure} is defined as a line bundle $\lambda:L\to M$, i.e. an object the category $\Line_\Man$, whose module of sections carries a local Lie algebra structure $(\Sec{L},\{,\})$. In virtue of proposition \ref{LocalLieRk1}, the locality condition is tantamount to the adjoint map of the Lie bracket mapping into derivations
\begin{equation*}
    \text{ad}_{\{,\}}:\Sec{L}\to \Dr{L}.
\end{equation*}
Note that this is entirely analogous to the Leibniz condition of a Poisson bracket where the smooth functions are replaced by the sections of the line bundle and vector fields are replaced by derivations. Since the symbol of any derivation defines a vector field, the analogue of the Hamiltonian map for a Jacobi manifold becomes a pair of $\Real$-linear maps
\begin{align*}
D: \Sec{L} & \to \Dr{L}\\
u & \mapsto D_u:=\{u,-\},
\end{align*}
\begin{align*}
X: \Sec{L} & \to \Sec{\Tan M}\\
u & \mapsto X_u:=\delta_*(D_u).
\end{align*}
For a section $u\in\Sec{L}$, $D_u$ is called its \textbf{Hamiltonian derivation} and $X_u$ its \textbf{Hamiltonian vector field}. As for any derivation of a line bundle, the Hamiltonian vector field map is uniquely determined by the following Leibniz identity
\begin{equation*}
    D_u(f\cdot v)=f\cdot D_u(v)+X_u[f]\cdot v,
\end{equation*}
which must hold for all sections $u,v\in\Sec{L}$ and smooth functions $f\in\Cin{M}$. The following proposition establishes these as the Jacobi analogue of the Hamiltonian map.
\begin{prop}[Hamiltonian Maps of a Jacobi Structure]\label{HamiltonianMapsJacobi}
The Hamiltonian maps $D$ and $X$ are Lie algebra morphisms, i.e.
\begin{equation*}
    D_{\{u,v\}}=[D_u,D_v] \qquad X_{\{u,v\}}=[X_u,X_v]
\end{equation*}
for all $u,v\in\Sec{L}$.
\end{prop}
\begin{proof}
The Lie algebra morphism condition for $D$ follows directly from the Jacobi identity of the bracket. To show that $X$ is a Lie algebra morphism we consider the triple bracket $\{u,\{v,f\cdot w\}\}$. We can expand the expression by using the Jacobi identity of the Lie bracket first and then applying the Leibniz property repeatedly to obtain
\begin{equation*}
   \{u,\{v,f\cdot w\}\}=f\cdot\{u,\{v,w\}\}+X_{u}[f]\cdot \{v,w\} + X_{v}[f]\cdot \{u,w\} + X_{\{u,v\}}[f]\cdot w + X_{v}X_{u}[f]\cdot w,
\end{equation*}
or we can expand only using the Leibniz property obtaining
\begin{equation*}
   \{u,\{v,f\cdot w\}\}=f\cdot\{u,\{v,w\}\}+X_{u}[f]\cdot \{v,w\} + X_{v}[f]\cdot \{u,w\} + X_{u}X_{v}[f]\cdot w.
\end{equation*}
Since both expressions must agree, it follows that
\begin{equation*}
   X_{\{u,v\}}[f]\cdot w + X_{v}X_{u}[f]\cdot w=X_{u}X_{v}[f]\cdot w
\end{equation*}
which must hold for all sections $u,v,w\in\Sec{A}$ and functions $f\in\Cin{M}$, thus implying the desired Lie algebra morphism condition.
\end{proof}
The fact that the Lie bracket $\{,\}$ is a derivation on each argument allows us to write
\begin{equation*}
    \{u,v\}=\Tilde{\Lambda}(j^1u,j^1v) \quad \text{ with }\quad  \Tilde{\Lambda}\in \Sec{\wedge^2(\Jet^1L)^*\otimes L}.
\end{equation*}
Note that $\tilde{\Lambda}$ is simply the squiggle of the Lie bracket $\{,\}$, as identified in proposition \ref{SymbolLieAlgebra}. We thus see that the bilinear form $\Tilde{\Lambda}$ is the analogue of the Poisson bivector and we appropriately call it the \textbf{Jacobi biderivation}. Noting that $(\Jet^1L)^*\otimes L=((\Der L)^{*L})^{*L}\cong \Der L$, the biderivation induces a musical map
\begin{equation*}
    \Tilde{\Lambda}^\sharp:\Jet^1 L \to \Der L
\end{equation*}
which allows us to rewrite the Hamiltonian derivation of a section $u\in\Sec{L}$ as
\begin{equation*}
    D_u = \Tilde{\Lambda}^\sharp (j^1u).
\end{equation*}
This musical map connects the, otherwise disconnected, der and jet sequences of the underlying line bundle in the following way
\begin{equation*}
\begin{tikzcd}
0 \arrow[r] &  \Real_M \arrow[r] & \Der L \arrow[r, "\delta"]  & \Tan M \arrow[r]  & 0 \\
0 &  L \arrow[l] & \Jet^1L  \arrow[u, "\Tilde{\Lambda}^\sharp"] \arrow[l]& \Tan^{*L} M \arrow[u,"\Lambda^\sharp"] \arrow[l,"i"] & 0 \arrow[l] 
\end{tikzcd}
\end{equation*}
where we have defined the bundle map $\Lambda^\sharp$ so that the diagram commutes. By construction, we see that this bundle map is indeed a musical map for the bilinear form defined as the fibre-wise pull-back of the Jacobi biderivation via the injective map of the jet sequence: $\Lambda=i^*\Tilde{\Lambda}\in\Sec{\wedge^2(\Tan^{*L} M)\otimes L}$. We call this form the \textbf{Jacobi L-bivector}. It follows from proposition \ref{SymbolLieAlgebra} that these bilinear forms should be interpreted as the symbols of the Hamiltonian maps when regarded as differential operators themselves; these indeed control the interaction of the Lie bracket with the $\Cin{M}$-module structure via the following identities
\begin{align*}
    D_{f\cdot u}&=f\cdot D_u+\Tilde{\Lambda}^\sharp(i(df\otimes u))\\
    X_{f\cdot u}&=f\cdot X_u+\Lambda^\sharp(df\otimes u).
\end{align*}
for all $f\in\Cin{M}$ and $u,v\in\Sec{L}$. Note that $X$ and $\Lambda$ are nothing but the symbol and squiggle of the Jacobi bracket seen as a general derivative Lie algebra bracket.\newline

In some situations it will be convenient to have a characterization of a Jacobi structure in terms of the Lie bracket and the Hamiltonian maps seen as $\Real$-linear maps satisfying some compatibility conditions with the module structure on sections. The next proposition gives this characterization as a form of converse statement of proposition \ref{SymbolLieAlgebra}, where the symbol and squiggle are defined for general derivative Lie algebras.
\begin{prop}[Extension by Symbol]\label{ExtensionBySymbol}
Let $\lambda:L\to M$ be a line bundle and let $\Sigma\subset \Sec{L}$ be a subspace of spanning sections such that $\Sec{L}\cong \Cin{M}\cdot \Sigma$. The datum of a Jacobi structure on $L$, i.e. a local Lie algebra $(\Sec{L},\{,\})$, is equivalent to a triple $((\Sigma,[,]),X,\Lambda)$ where
\begin{itemize}
\normalfont
    \item $(\Sigma,[,])$ is a $\Real$-linear Lie bracket,
    \item $X:\Sec{L}\to \Sec{\Tan M}$ is a $\Real$-linear map and,
    \item $\Lambda\in\Sec{\wedge^2\Tan^{L} M\otimes L}$ inducing a vector bundle morphism $\Lambda^\sharp:\Tan^{*L} M\to \Tan M$
\end{itemize}
satisfying the compatibility conditions
\begin{align*}
    1.\quad & [X_s,X_{r}]=X_{[s,r]},\\
    2.\quad & X_{f\cdot s}=f\cdot X_s+\Lambda^\sharp(df\otimes s), \\
    3.\quad & [X_s,\Lambda^\sharp(df\otimes r)]=\Lambda^\sharp(dX_{s}[f]\otimes r + df\otimes [s,r]),\\
    4.\quad & \Lambda(df\otimes s, \Lambda(dg \otimes r , dh\otimes t))+\Lambda(dg\otimes r, \Lambda(dh \otimes t , df\otimes s))+\Lambda(dh\otimes t, \Lambda(df \otimes s , dg\otimes r)) =\\
    =&X_r[f]\cdot \Lambda( dh\otimes t, dg\otimes s)+X_t[g]\cdot \Lambda(df\otimes s, dh\otimes r)+X_s[h]\cdot \Lambda(dg\otimes r,df\otimes t),
\end{align*}
for all $f,g,h\in\Cin{M}$ and $s,r,t\in\Sigma$. The equivalence is realized by setting
\begin{align*}
    \{s,s'\}&:=[s,s'] \\
    \{s,f\cdot s'\}&:=f\cdot \{s,s'\} + X_{s}[f]\cdot s'
\end{align*}
which uniquely determine a local Lie bracket $(\Sec{L},\{,\})$.
\end{prop}
\begin{proof}
Since $\Sigma$ is a subspace of spanning sections, a Jacobi bracket will be completely determined by specifying brackets of the form $\{f\cdot s,g\cdot s'\}$ for all $f,g\in\Cin{M}$ and $s,s'\in\Sigma$., then our aim will be to reconstruct a Jacobi bracket from its symbol-squiggle identity
\begin{equation*}
    \{f\cdot s,g\cdot s'\}:=fg\cdot\{s,s'\}+fX_s[s]\cdot s' - gX_{s'}[f]\cdot s +\Lambda^\sharp(df\otimes s)[g]\cdot s'.
\end{equation*}
By setting $\{s,s'\}=[s,s']$ it is clear that the triple $((\Sigma,[,]),X,\Lambda)$ determines an antisymmetric bracket $(\Sec{L},\{,\})$, thus we are only left with showing that it does indeed satisfy the Jacobi identity. Since the bracket $[,]$ is assumed to be Lie, we have to verify the following general Jacobi identity:
\begin{equation*}
    \{f\cdot u,\{g\cdot v,h\cdot w\}\}=\{\{f\cdot u,g\cdot v\},h\cdot w\}+\{g\cdot v,\{f\cdot u,h\cdot w\}\}.
\end{equation*}
This identity follows directly from the compatibility conditions of the symbol and squiggle 1. - 4. assumed in the statement of the proposition but checking it is quite a formidable computation. The length of the task is somewhat reduced if we prove the Jacobi identity for some of the functions $f,g,h$ being $1$. We first check the Jacobi identity for $\{u,\{v,h\cdot w\}\}$, by expanding both sides of the equation above for the case $f=1,g=1$ with the defining identities of the bracket $\{,\}$ and simplifying whenever possible, we easily find
\begin{equation*}
    X_u[X_v[h]]\cdot w=X_v[X_u[h]\cdot w +X_{\{u,v\}}[h]\cdot w
\end{equation*}
but this is precisely the compatibility condition 1, then the identity holds. For the bracket $\{w,\{f\cdot u,g\cdot v\}\}$ we expand both sides of the Jacobi identity, again using the defining identities of $X$ and $\Lambda$ and simplifying whenever possible with the Jacobi identity of the bracket $[,]$, to obtain
\begin{align*}
& X_w[g]X_v[f]\cdot u - gX_{\{v,w\}}[f]\cdot u -X_w[f]X_u[g]\cdot v - fX_{\{w,u\}}[g]\cdot v = \\
& \Lambda^\sharp(df\otimes u)[g]\cdot \{w,v\} +X_w[\Lambda^\sharp(df\otimes u)[g]]\cdot v -[X_w,X_{f\cdot u}][g]\cdot v + [X_w,X_{g\cdot v}][f]\cdot u.
\end{align*}
Then, further expanding by the defining identities of $X$ and crucially using the $\Cin{M}$-linearity of $\Lambda^\sharp$, after simplifications we get
\begin{equation*}
    X_w[\Lambda^\sharp(df\otimes u)[g]]\cdot v -\Lambda^\sharp(df\otimes \{w,u\})[g]\cdot v = \Lambda^\sharp(dX_w[f]\otimes u)[g]\cdot v - \Lambda^\sharp(dX_w[g]\otimes v)[f]\cdot u.
\end{equation*}
We can use the antisymmetry of $\Lambda$ to rewrite the last term and make the identity into $\Cin{M}$-linear combination of terms in $v$, thus obtaining
\begin{equation*}
    [X_w,\Lambda^\sharp(df\otimes u)][g]\cdot v = \Lambda^\sharp(dX_w[f]\otimes u + df\otimes \{w,u\})[g]\cdot v,
\end{equation*}
which is indeed the compatibility condition 3. Lastly, for the bracket $\{f\cdot u,\{g\cdot v,h\cdot w\}\}$ it will suffice to use the bilinear form expression of the squiggle to write
\begin{equation*}
    \Lambda^\sharp(df\otimes u)[g]\cdot v=\Lambda(df\otimes u, dg\otimes v)=-\Lambda^\sharp(dg\otimes v)[f]\cdot u
\end{equation*}
in order to directly apply the compatibility conditions 1. and 3. to the Leibniz expansion of the Jacobi identity. Then after a long, but routine, calculation cancelling whenever possible, it follows that the Jacobi identity reduces to the compatibility condition 4.
\end{proof}

Analogously to the linear Poisson structures defined on vector bundles, we define the corresponding notion for Jacobi structures on L-vector bundles. Let $\epsilon:E^L\to M$ be a L-vector bundle and consider the pull-back line bundle over the total space $L_E:=\epsilon^*L$. We identify the fibre-wise constant sections as $\epsilon^*:\Sec{L}\hookrightarrow \Sec{L_E}$ and the fibre-wise linear sections as $l:\Sec{E^{*L}}\hookrightarrow \Sec{L_E}$, which form a submodule of spanning sections for $\Sec{L_E}$. A \textbf{fibre-wise linear Jacobi structure} is defined as a Jacobi structure $(\Sec{L_E},\{,\})$ that is completely determined by its restriction to spanning sections in the following way
\begin{align*}
    \{l(\Sec{E^{*L}}),l(\Sec{E^{*L}})\} & \subset  l(\Sec{E^{*L}})\\
    \{l(\Sec{E^{*L}}),\epsilon^*(\Sec{L}\} & \subset  \epsilon^*(\Sec{L})\\
    \{\epsilon^*(\Sec{L}),\epsilon^*(\Sec{L})\} & = 0.
\end{align*}

Let $i:S\hookrightarrow M$ be an embedded submanifold of a line bundle $\lambda:L\to M$ and let us denote by $L_S:=i^*L$ the induced line bundle over $S$ and by $\Gamma_S\subset \Sec{L}$ the submodule of vanishing sections. Note that at each point $x\in M$ the jet space $\Jet^1_x L$ carries a $L_x$-valued bilinear form given by the Jacobi biderivation $\Tilde{\Lambda}_x$ and so the notion of isotropic subspace of a jet space (one where the bilinear form restricts to zero) is well-defined point-wise. We can then define coisotropic submanifolds of Jacobi manifolds in an entirely analogous way to the coisotropic submanifolds of Poisson manifolds. We say that $S$ is \textbf{coisotropic} if the L-annihilator of its der bundle $(\Der L_S)^{0L}\subset \Jet^1L$ is an isotropic subbundle with respect to the bilinear form $\Tilde{\Lambda}$. Equivalently, $S$ is coisotropic if
\begin{equation*}
    \Tilde{\Lambda}^\sharp_x((\Der_xL_S)^{0L})\subset \Der_xL_S \qquad \forall x\in S.
\end{equation*}
The following proposition gives several equivalent characterizations of coisotropic submanifolds.
\begin{prop}[Coisotropic Submanifolds of a Jacobi Manifold, cf. {\cite[Lemma 2.44]{tortorella2017deformations}}]\label{CoisotropicSubmanifoldsJacobi}
Let $i:S\hookrightarrow M$ be an embedded submanifold of a Jacobi structure $(\Sec{L},\{,\})$, then the following conditions are equivalent
\begin{enumerate}
\normalfont
    \item $S$ is a coisotropic submanifold;
    \item $\Lambda^\sharp_x(\Tan_x S^0\otimes L_x)\subset \Tan_x S \quad \forall x\in S$;
    \item the vanishing sections $\Gamma_S$ form a Lie subalgebra, $\{\Gamma_S,\Gamma_S\}\subset \Gamma_S$;
    \item the Hamiltonian vector fields of vanishing sections are tangent to $S$, $X_s\in\Sec{\Tan S}$ for $s\in\Gamma_S$.
\end{enumerate}
\end{prop}
\begin{proof}
Equivalence between 1. and 2. follows from the definition of the bilinear forms, $\Lambda=i^*\Tilde{\Lambda}$, and the fact that
\begin{equation*}
    i(\alpha\otimes l)\in (\Der_xL_S)^{0L} \quad \Leftrightarrow \quad i(\alpha\otimes l)(a)=0 \quad \forall a\in \Der_xL_S,
\end{equation*}
which, since $i$ is the $L$-dual to the anchor $\delta$, is tantamount to demanding
\begin{equation*}
    \alpha(\delta(a))=0 \quad \forall a\in \Der_xL_S.
\end{equation*}
Recall that $\delta(\Der L_S)=\Tan S$, then 
\begin{equation*}
    i(\alpha\otimes l)\in (\Der_xL_S)^{0L} \quad \Leftrightarrow \quad \alpha\otimes l \in \Tan_x S^0\otimes L_x.
\end{equation*}
To show equivalence between 2. and 3. we can use a local trivialization to write $\Ker{i^*}=\Gamma_S=I_S\cdot \Sec{L}$, then it follows from the definition of the Jacobi L-bivector that
\begin{equation*}
    i^*\{f\cdot s,g\cdot s\}=\Lambda^\sharp(df\otimes s)[g]\cdot i^*s=0
\end{equation*}
for $f,g\in I_S$ and any section $s\in\Sec{L}$. Finally, equivalence between 3. and 4. follows by considering $f\in I_S$, $s\in\Sec{L}$, $s'\in\Gamma_S$ and writing
\begin{equation*}
    i^*\{s',f\cdot s\}=i^*f\cdot i^*\{s',s\} + i^*(X_{s'}[f])\cdot i^*s= i^*(X_{s'}[f])\cdot i^*s
\end{equation*}
thus implying $X_{s'}[f]\in I_S$, completing the equivalence.
\end{proof}

Let $(\Sec{L_1},\{,\}_1)$ and $(\Sec{L_2},\{,\}_2)$ be two Jacobi structures, a \textbf{Jacobi map} is a factor, i.e. morphism in the category of line bundles $\Line_\Man$, $B:L_1\to L_2$ such that its pull-back on sections
\begin{equation*}
    B^*:(\Sec{L_2},\{,\}_2)\to (\Sec{L_1},\{,\}_1)
\end{equation*}
is a Lie algebra morphism. The next proposition shows that the categorical product in $\Line_\Man$ allows us to define a product of Jacobi manifolds.
\begin{prop}[Product of Jacobi Manifolds, cf. {\cite[Proposition 5.1]{ibanez1997coisotropic}}]\label{ProductJacobi}
Let Jacobi structures $(\Sec{L_1},\{,\}_1)$ and $(\Sec{L_2},\{,\}_2)$, then there exists a unique Jacobi structure in the line product $(\Sec{L_1\utimes L_2},\{,\}_{12})$ such that the canonical projection factors
\begin{equation*}
\begin{tikzcd}
L_1 & L_1\utimes L_2 \arrow[l,"P_1"']\arrow[r,"P_2"] & L_2.
\end{tikzcd}
\end{equation*}
are Jacobi maps. The Jacobi structure $(\Sec{L_1\utimes L_2},\{,\}_{12})$ is called the \textbf{product Jacobi structure}.
\end{prop}
\begin{proof}
Recall from the definition of line product as a pull-back bundle that $\Sec{L_1\utimes L_2}=\Cin{M_1\dtimes M_2} \cdot p_1^*\Sec{L_1}\cong \Cin{M_1\dtimes M_2} \cdot p_2^*\Sec{L_2}$, therefore it will suffice to determine the bracket on pull-backs of sections of each factor and then use the extension by symbol formula of proposition \ref{ExtensionBySymbol} for some appropriately chosen Hamiltonian vector field $X^{12}:\Sec{L_1\utimes L_2}\to \Sec{\Tan M_1\dtimes M_2}$ and L-bivector musical map $\Lambda^{12}:\Sec{\Tan ^{*L} (M_1\dtimes M_2)\otimes L_1\utimes L_2}\to \Sec{\Tan (M_1\dtimes M_2)}$. The bracket on the line product is then defined as a $\Real$-linear bracket given by
\begin{equation*}
    \{P_1^*s_1,P_1^*s_1'\}_{12}:=P_1^*\{s_1,s_1'\}_1 \qquad \{P_2^*s_2,P_2^*s_2'\}_{12}:=P_2^*\{s_2,s_2'\}_2\qquad \{P_1^*s_1,P_2^*s_2\}_{12}:=0
\end{equation*}
for all $s_i,s_i'\in\Sec{L_i}$, $i=1,2$. Proposition \ref{SpanningFunctionsLineProduct} ensures that it suffices to define the action of the Hamiltonian vector field and the L-bivector musical map locally on spanning functions of the base product, then we set
\begin{align*}
    X^{12}_{P_1^*s_1}[\tfrac{a}{b}]&=\tfrac{\{s_1,a\}_1}{b}   & X^{12}_{P_2^*s_2}[\tfrac{b}{a}]&=\tfrac{\{s_2,b\}_2}{a}\\
    \Lambda^{12}(d\tfrac{a}{b}\otimes P_1^*s_1)[\tfrac{a'}{b'}]&=\tfrac{\{a,a'\}_1}{b}\tfrac{s_1}{b'} &\Lambda^{12}(d\tfrac{b}{a}\otimes P_2^*s_2)[\tfrac{b'}{a'}]&=\tfrac{\{b,b'\}_2}{a}\tfrac{s_2}{a'}
\end{align*}
for all $s_i\in\Sec{L_i}$, $a,a'\in\Sec{L_1^\times}$ and $b,b'\in\Sec{L_2^\times}$. Checking that $X^{12}$ and $\Lambda^{12}$ so defined satisfy the extension-by-symbol compatibility conditions of proposition \ref{ExtensionBySymbol} follows by a long direct computation that becomes routine once formulas for the action of $X^{12}$ and $\Lambda^{12}$ on all the possible combinations of spanning functions are derived from the symbol-squiggle identity. Appendix \ref{SymbolSquiggleProductJacobi} contains all these formulas. We also point out that, although the above formulas may appear asymmetric for the local non-vanishing sections $a\in\Sec{L_1^\times}$ and $b\in\Sec{L_2^\times}$, the definition of local ratio functions is such that $\tfrac{a}{b}\tfrac{b}{a}=1$ and so, by construction, we can relate both sides of the definition by the identities
\begin{equation*}
    d\tfrac{a}{b}=-(\tfrac{a}{b})^2d\tfrac{b}{a} \qquad P_1^*a=\tfrac{a}{b}P_2^*b.
\end{equation*}
\end{proof}

Given a Jacobi structure $(\Sec{L},\{,\})$ we define the \textbf{opposite Jacobi structure} as $\overline{L}:=(\Sec{L},-\{,\})$. The proposition below shows that we can regard Jacobi maps as a particular case of \textbf{coistropic relations}, which are defined, in general, as coisotropic submanifolds of the product Jacobi manifold $L_1\utimes \overline{L}_2$.
\begin{prop}[Jacobi Maps as Coisotropic Relations, cf. {\cite[Theorem 5.3]{ibanez1997coisotropic}}]\label{JacobiMapCoisotropicRelation}
Let two Jacobi structures $(\Sec{L_{M_1}},\{,\}_1)$ and $(\Sec{L_{M_2}},\{,\}_2)$ and a factor $B:L_{M_1}\to L_{M_2}$, then $B$ is a Jacobi map iff its L-graph
\begin{equation*}
\normalfont
    \LGrph{B}\subset M_1\dtimes M_2
\end{equation*}
is a coisotropic submanifold of the product Jacobi structure $L_{M_1}\utimes \overline{L}_{M_2}$.
\end{prop}
\begin{proof}
We will use characterization 3 in proposition \ref{CoisotropicSubmanifoldsJacobi} to identify a coisotropic submanifold with a Lie subalgebra of the Jacobi bracket. Note that, by construction of $\LGrph{B}$, the vanishing sections $\Gamma_{\LGrph{B}}\subset \Sec{L_1\utimes \overline{L}_2}$ are generated by those of the form
\begin{equation*}
    P_1^*B^*s-P_2^*s\qquad \forall s\in\Sec{L_2}.
\end{equation*}
We can directly compute the brackets using the defining relations of the product Jacobi structure
\begin{equation*}
    \{P_1^*B^*s-P_2^*s,P_1^*B^*s'-P_2^*s'\}_{12}=P_1^*\{B^*s,B^*s'\}_1-P_2^*\{s,s'\}_2
\end{equation*}
which must hold for all $s.s'\in\Sec{L_2}$. Then it follows that
\begin{equation*}
    \{P_1^*B^*s-P_2^*s,P_1^*B^*s'-P_2^*s'\}_{12}\in \Gamma_{\LGrph{B}}\quad \Leftrightarrow \quad \{B^*s,B^*s'\}_1 = B^*\{s,s'\}_2
\end{equation*}
as desired.
\end{proof}

Let a line bundle $\lambda:L\to M$ with a Jacobi structure $(\Sec{L},\{,\})$ and consider a submanifold $i:C\hookrightarrow M$ with its corresponding embedding factor $\iota: L_C\to L$. For another line bundle $\lambda':L'\to M'$, assume there exists a submersion factor $\pi:L_C\to L'$ covering a surjective submersion $p:C\twoheadrightarrow M'$. Then, we say that a Jacobi structure $(\Sec{L},\{,\})$ \textbf{reduces} to the Jacobi structure $(\Sec{L'},\{\,\}')$ via $\pi:L_C\to L'$ when for all pairs of sections $s_1,s_2\in\Sec{L'}$ the identity
\begin{equation*}
    \pi^*\{s_1,s_2\}'=\iota^*\{S_1,S_2\}
\end{equation*}
holds for all choices of extensions $S_1,S_2$, i.e all choices of sections $S_1,S_2\in\Sec{L}$ satisfying 
\begin{equation*}
    \pi^*s_1=\iota^*S_1 \qquad \pi^*s_2=\iota^*S_2.
\end{equation*}
The following proposition shows that, in analogy with the case of reduction of Poisson manifolds, coistropic play a distinguished role in the reduction of Jacobi manifolds.

\begin{prop}[Coisotropic Reduction of Jacobi Manifolds, cf. {\cite[Proposition 2.56]{tortorella2017deformations}}] \label{CoisotropicReductionJacobi}
Let $\lambda:L\to M$ and $\lambda':L'\to M'$ be line bundles,  $(\Sec{L},\{,\})$ be a Jacobi structure, $i:C\hookrightarrow M$ a closed coisotropic submanifold and $\pi:L_C\to L'$ a submersion factor covering a surjective submersion $p:C\twoheadrightarrow M'$ so that we have the reduction diagram:
\begin{equation*}
    \begin{tikzcd}[sep=tiny]
    L_C \arrow[rr,"\iota"] \arrow[dd,"\pi"'] \arrow[dr]& & L \arrow[dr] & \\
    & C \arrow[rr,"i"', hook] \arrow[dd, "p",twoheadrightarrow] & & M \\
    L' \arrow[dr] & & & \\
    & M' &  & 
    \end{tikzcd}
\end{equation*}
Assume the following compatibility condition between the coisotropic submanifold and the submersion factor:
\begin{equation*}
\normalfont
    \delta(\Ker{\Der \pi})=\Lambda^\sharp((\Tan C)^{\text{0}L}),
\end{equation*}
where $\delta:\Der L \to \Tan M$ is the anchor of the der bundle, $\normalfont(\Tan C)^{0L}\subset \Tan^{*L}M$ is the $L$-annihilator and $\normalfont\Lambda^\sharp: \Tan^{*L}M\to \Tan M$ is the musical map, then there exists a unique Jacobi structure $(\Sec{L'},\{,\}')$ such that $(\Sec{L},\{,\})$ reduces to it via the submersion factor $\pi:L_C\to L'$.
\end{prop}
\begin{proof}
Recall the proposition \ref{CoisotropicSubmanifoldsJacobi} identified a coisotropic submanifold $C$ with the Lie subalgebra of its vanishing submodule $\Gamma_C\subset \Sec{L}$. Since, at least locally, this characterization is accompanied by the vanishing ideal $I_C\subset \Cin{M}$ so that $\Gamma_C=I_C\cdot \Sec{L}$, we can use the constructions for submersion and embedding factors introduced in Section \ref{CategoryOfLineBundles} to regard the functions and sections as quotients of functions and sections that restrict to fibres of the surjective submersion: $\Cin{M'}\cong (\Cin{M}/I_C)_p$ and $\Sec{L'}\cong (\Sec{L}/\Gamma_C)_\pi$. Then we can construct a Lie algebra $(\Sec{L'},\{,\})$ by noting that the compatibility condition ensures that sections restricting to $p$-fibres correspond to the elements in the Lie idealizer $N(\Gamma_C)$; hence $\Sec{L'}\cong N(\Gamma_C)/\Gamma_C$ as Lie algebras. To complete the proof we need to check that this Lie algebra indeed corresponds to a local Lie algebra structure on the line bundle $\lambda':L'\to M'$. This follows from the fact that the compatibility condition can be read as the point-wise requirement of Hamiltonian vector fields of the vanishing sections to be tangent to the $p$-fibres, then a few routine calculations show that the local Lie algebra properties of the bracket $\{,\}$ carry over to the bracket $\{,\}'$ in a natural way.
\end{proof}
As a particular case of coisotropic reduction we find \textbf{Jacobi submanifolds}, i.e. submanifolds $i:C\hookrightarrow M$ whose vanishing submodule $\Gamma_C$ is not only a Lie subalgebra but a Lie ideal, $\{\Gamma_C,\Sec{L}\}\subset \Gamma_C$. In this case, the distribution of Hamiltonian vector fields clearly vanishes, thus satisfying all the requirements of the proposition above trivially, so that the restricted line bundle itself inherits a Jacobi structure $(L_C,\{,\}_C)$. Another example of coistropic reduction is given by the presence of a \textbf{Hamiltonian group action}: a line bundle Lie group action on a Jacobi manifold $G \Acts L$ via Jacobi maps with infinitesimal action $\Psi:\mathfrak{g}\to \Dr{L}$, and a (co)moment map $\overline{\mu}:\mathfrak{g}\to \Sec{L}$ satisfying the defining conditions
\begin{equation*}
    \Psi(\xi)=D_{\overline{\mu}(\xi)},\qquad \{\overline{\mu}(\xi),\overline{\mu}(\zeta)\}=\overline{\mu}([\xi,\zeta]) \qquad \forall \xi,\zeta\in\mathfrak{g}.
\end{equation*}

Comparing the results presented in this section thus far with their analogues for Poisson manifolds in Section \ref{SymplecticGeometry} confirms that there are obvious structural similarities between Jacobi and Poisson geometry. In Section \ref{UnitFreeManifolds} we argued that line bundles conceptually represent unit-free manifolds from the fact that there is no canonical choice of global unit section, in contrast with the global constant function $1$ in ordinary manifolds. Jacobi structures are indeed Lie brackets for the unit-free generalization of functions, again in direct analogy with Poisson structures being Lie brackets on ordinary functions. Since local trivializations ensure the existence of local non-vanishing sections, one should expect to somehow recover Poisson algebras in the context of a Jacobi manifold, at least locally. The proposition below shows that this is indeed the case when we consider local units.
\begin{prop}[Local Units of a Jacobi Manifold] \label{LocalUnitsJacobi}
Let $\lambda:L\to M$ be a line bundle and $(\Sec{L},\{,\})$ a Jacobi structure on it, then a local unit $u\in\Sec{L^\times}$ on an open neighbourhood $U\subset M$ induces a Poisson structure on the $X_u$-invariant local functions $\normalfont (\text{C}_u^\infty(U),\{,\}_u)$. 
\end{prop}
\begin{proof}
Since $u$ is non-vanishing, any local section $s\in\Sec{L|_U}$ will be uniquely identified with a local function via $s=f\cdot u$. Then, considering the Jacobi bracket of any two such sections with $X_u$-invariant factors, we find
\begin{equation*}
    \{f\cdot u,g\cdot u\}=fg\cdot \{u,u\}+fX_u[g]\cdot u -gX_u[f]\cdot u + \Lambda^\sharp(df\otimes u)[g]\cdot u =\Lambda^\sharp(df\otimes u)[g]\cdot u
\end{equation*}
where we have made use of the properties of the symbol and squiggle of the Jacobi bracket and the invariance condition $X_u[f]=0$, $X_u[g]=0$. This clearly specifies a bracket on $X_u$-invariant functions that, by the definition of the musical map $\Lambda^\sharp$, can be written explicitly as
\begin{equation*}
    \{f,g\}_u\cdot u:=\Lambda(df\otimes u, dg\otimes u)
\end{equation*}
for $f,g\in \text{C}_u^\infty(U)$. This bracket is clearly antisymmetric and Leibniz by construction, and its Jacobi identity can be easily checked by considering a bracket of the form $\{f\cdot u,\{g\cdot u,h\cdot u\}\}$ and then using the Jacobi identity of $\{,\}$ and the correspondence $\Cin{U}\cong \Sec{L|_U}$.
\end{proof}

A local unit $u\in\Sec{L|_U^\times}$ whose Hamiltonian vector field vanishes $X_u=0$ induces, from the construction above, a Poisson structure on all the local functions $\Cin{U}$ and it is thus called a \textbf{Poisson local unit}. As it will be apparent from our discussion on contact and Lichnerowicz manifolds later in this section, Jacobi manifolds do not generically admit Poisson local units. A Jacobi structure $(\Sec{L},\{,\})$ that admits local trivializations via Poisson local units is called a \textbf{conformal Poisson structure}. If a conformal Poisson structure on $L$ admits a global unit, thus trivializing the line bundle $L\cong \Real_M$, we recover the definition of \textbf{Poisson structure}. We then conclude that Poisson manifolds can be regarded as a special subclass of Jacobi manifolds: the ones defined on trivial line bundles and characterised by the fact that the Hamiltonian vector field of a function only depends on the differential of the function but not on its point-wise values. It is a routine check to show that setting these conditions in all the definitions and propositions about Jacobi manifolds of this section recovers the definitions and results given in Section \ref{SymplecticGeometry} for Poisson manifolds. This now retroactively justifies our choice of common terminology in both Poisson and Jacobi geometry.\newline

A \textbf{precontact manifold} is a pair $(M,\text{H}\subset \Tan M)$ with $\text{H}$ a hyperplane distribution, i.e. $\dimm(\text{H}_x)=\dimm (M)-1$ for all $x\in M$. Note that a hyperplane distribution on the tangent bundle is equivalent to the datum of a (generically non-trivial) line bundle $\lambda:L\to M$ and a non-vanishing $L$-valued 1-form $\theta:\Tan M\to L$. The equivalence is realized by setting
\begin{equation*}
    \text{H}=\Ker{\theta}, \quad \text{ which then gives } \quad \Tan M/\text{H}\cong L.
\end{equation*}
Let us denote the $\Cin{M}$-submodule of vector fields tangent to the hyperplane distribution by $\Sec{\text{H}}$. We can define the following antisymmetric map for vector fields tangent to the hyperplane distribution
\begin{align*}
\omega: \Sec{\wedge^2 \text{H}} & \to \Sec{L}\\
(X,Y) & \mapsto \theta([X,Y]).
\end{align*}
The kernel of this map clearly measures the degree to which $\text{H}$ is integrable as a tangent distribution noting that, in particular, when $[\Sec{\text{H}},\Sec{\text{H}}]\subset\Sec{\text{H}}$ the map $\omega$ is identically zero. It follows by construction that $\omega$ is in fact $\Cin{M}$-bilinear and thus defines a bilinear form $\omega_\text{H}:\wedge^2 \text{H}\to L$ called the \textbf{curvature form} of the hyperplane distribution $\text{H}$. A hyperplane distribution $\text{H}$ is called \textbf{maximally non-integrable} when its curvature form $\omega_\text{H}$ is non-degenerate, i.e. when the musical map $\omega^\flat_\text{H}:\text{H}\to \text{H}^{*L}$ has vanishing kernel. Such a hyperplane distribution $\text{H}\subset \Tan M$ is called a \textbf{contact structure} on $M$ and we refer to the pair $(M,\text{H})$ as a \textbf{contact manifold}. Simple dimension counting applied to any tangent space of a contact manifold reveals that a manifold supporting a contact hyperplane is necessarily odd-dimensional. A \textbf{contact map}, defined as a smooth map $\varphi:(M_1,\text{H}_1)\to(M_2,\text{H}_2)$ whose tangent assigns the hyperplane distributions isomorphically $\Tan \varphi (\text{H}_1)=\text{H}_2$, is necessarily a local diffeomorphism from the fact that hyperplane distributions have codimension $1$ everywhere. This condition is equivalent to the tangent map $\Tan \varphi$ inducing a well-defined morphism of line bundles $\Phi:L_1\cong \Tan M_1/\text{H}_1\to L_2\cong \Tan M_2/\text{H}_2$. The presence of a hyperplane distribution on a precontact manifold allows for the identification of \textbf{isotropic submanifolds} as integral submanifolds of the tangent distribution of hyperplanes, i.e. submanifolds $S\subset M$ with $\Tan S\subset \text{H}$. \newline

A Jacobi structure $(\Sec{L}, \{,\})$ is called \textbf{non-degenerate} when its Jacobi biderivation $\Tilde{\Lambda}$ induces a musical isomorphism of L-vector bundles
\begin{equation*}
        \begin{tikzcd}
        \Jet^1 L\arrow[r, "\Tilde{\Lambda}^\sharp",yshift=0.7ex] & \Der L \arrow[l,"\Tilde{\Lambda}^\flat",yshift=-0.7ex]
        \end{tikzcd}
\end{equation*}
In Section \ref{SymplecticGeometry} we proved that symplectic manifolds, conventionally regarded as non-degenerate presymplectic manifolds, are equivalent to non-degenerate Poisson manifolds. It turns out that there is an entirely analogous connection between contact manifolds and non-degenerate Jacobi manifolds.

\begin{prop}[Non-Degenerate Jacobi Manifolds are Contact Manifolds, {\cite[Proposition 2.31]{tortorella2017deformations}}]\label{NonDegenerateJacobiContact}
Let $M$ be any smooth manifold, then the datum of a contact structure $(M,\text{H}\subset \Tan M)$ is equivalent to the datum of a non-degenerate Jacobi structure $(\Sec{L_M},\{,\})$.
\end{prop}
\begin{proof}
Let us assume $(M,\text{H}\subset \Tan M)$ is a contact structure first and show that it defines a Hamiltonian vector field and L-bivector, which, in virtue of proposition \ref{ExtensionBySymbol}, give a Jacobi structure that is non-degenerate. The curvature form is non-degenerate and thus has an inverse $\omega^\sharp:\text{H}^{*L}\to \text{H}\subset \Tan M$. Define $L:=\Tan M/\text{H}$ with canonical projection defining the contact form $\theta\in\Sec{\Tan^{*L} M}$, which in turn gives an $\Real$-linear map $X:\Sec{L}\to\Sec{\Tan M}$ defined by the condition $\theta(X_u)=u$. Let us define the L-bivector musical map
\begin{equation*}
    \Lambda^\sharp(df\otimes u):=\omega^\sharp(df|_\text{H}\otimes u).
\end{equation*}
It follows from the definition of curvature form and the Jacobi identity of Lie bracket of vector fields that $X$ and $\Lambda$ so defined satisfy the extension by symbol conditions. The point-wise decomposition of tangent spaces as $T_xM\cong \text{H}\oplus \Real\cdot X_u(x)$ for some locally non-vanishing section $u$ then shows that the Jacobi biderivation associated to this Jacobi structure is indeed non-degenerate. Conversely, let us assume $\lambda:L\to M$ is a line bundle whose sections carry a non-degenerate Jacobi structure $(\Sec{L},\{,\})$ with L-bivector musical map $\Lambda^\sharp:\Tan^{*L}M\to \Tan M$. A point-wise dimension count shows that $\text{H}=\Lambda^\sharp(\Tan^{*L}M)\subset \Tan M$ is a hyperplane distribution and by direct computation we show that it is, in fact, contact and precisely the converse construction to the previous case.
\end{proof}
A submanifold of a contact manifold $S\subset M$ that is isotropic with respect to the hyperplane distribution and that is coisotropic with respect to the associated non-degenerate Jacobi structure is called \textbf{Legendrian}. The non-degeneracy of the Jacobi structure forces Legendrian submanifolds to be maximally isotropic, then, if the odd dimension of the contact manifold is $\dimm M=2n+1$, the dimension of a Legendrian submanifold must be $\dimm S=n$. Similarly to the definition of the Weinstein category of symplectic manifolds in Section \ref{SymplecticGeometry}, we can use Legendrian submanifolds to define \textbf{the category of contact manifolds} $\Cont_\Man$ with objects smooth manifolds with a contact hyperplane distribution $(M,\text{H})$, or equivalently non-degenerate Jacobi structures on line bundles $(\Sec{L_M},\{,\})$, with morphisms given by Legendrian relations $R:L_1\dashrightarrow L_2$, i.e. Legendrian submanifolds $R\subset M_1\dtimes M_2$ of the line product of non-degenerate Jacobi structures $L_1\utimes \overline{L_2}$. In particular, note that the L-graph of a contact map is clearly maximally coisotropic, thus we see that contact maps are recovered as isomorphisms in this category. As in the case of composition of Lagrangian relations, the composition of Legendrian relations is subject to the same issues of possible non-smoothness. These shall be addressed in the particular subclasses of contact manifolds that will be considered in this thesis.\newline

We close this section by giving a brief description of Jacobi manifolds with trivial line bundle. It follows from the results presented at the end of Section \ref{CategoryOfLineBundles} for a trivial line bundle $\Real_M$, that a Jacobi structure corresponds to a Lie bracket on functions $(\Cin{M},\{,\})$ defined by a bivector field $\pi\in\Sec{\wedge^2 \Tan M}$ and a vector field $R\in\Sec{\Tan M}$ in the following way:
\begin{equation*}
    \{f,g\}=\pi(df,dg)+fR[g]-gR[f]
\end{equation*}
for all $f,g\in\Cin{M}$. The Jacobi identity of this bracket is tantamount to the following Gerstenhaber algebra conditions
\begin{equation*}
    [R,\pi]=0\qquad [\pi,\pi]+2R\wedge\pi =0.
\end{equation*}
When these identities hold, $(\pi,R)$ is called a \textbf{Jacobi pair} and $(M,\pi,R)$ is known as a \textbf{Lichnerowicz manifold}, which is indeed simply a Jacobi manifold with globally trivialized line bundle. Lichnerowicz structures are found in generically-non-trivial line bundles locally, we can see this explicitly by considering a local unit $u\in\Sec{L|_U^\times}$ and setting
\begin{equation*}
    \pi^u(df,dg):=\Lambda(df\otimes u, dg \otimes u), \qquad R^u[f]:=X_u[f]
\end{equation*}
for $f,g\in \Cin{U}$, which are easily shown to form a Jacobi pair $(\pi^u,R^u)$ from the basic symbol-squiggle identities of a Jacobi bracket and the identification $\Sec{L|_U}\cong \Cin{U}$. Note that a Lichnerowicz manifold $(M,\pi,R)$ with $R=0$ makes $(M,\pi)$ precisely into a Poisson manifold. Recall that for a trivial line bundle we have the isomorphism $\Der \Real_M\cong \Tan M\oplus \Real_M$, then the Hamiltonian derivation and Hamiltonian vector field of a function $f\in\Cin{M}$ take the following concrete form
\begin{equation*}
    D_f= \pi^\sharp(df)+fR\oplus -R[f] \qquad X_f=\pi^\sharp(df)+fR.
\end{equation*}
A \textbf{Lichnerowicz map} is a smooth map between Lichnerowicz manifolds whose pull-back on functions is a Lie algebra morphisms for the Jacobi brackets, a condition, in fact, equivalent to the bivector and vector fields being related by the tangent map of the smooth map. Lichnerowicz maps are special cases of Jacobi maps of trivial line bundles where the conformal factor is set to the constant function $1$.\newline

Suppose that the Jacobi structure of a Lichnerowicz manifold $(M,\pi,R)$ is non-degenerate so that the manifold is endowed with a contact structure $(M,\text{H})$. Since the line bundle is now globally trivial, the hyperplane distribution will be given by the kernel of some non-vanishing ordinary $1$-form 
\begin{equation*}
    \text{H}=\Ker{\theta},\qquad \theta\in\Omega^1(M)
\end{equation*}
called the \textbf{contact form}. In this case $(M,\theta)$ is called an \textbf{exact contact manifold}. It follows by construction that the curvature form of the hyperplane distribution is given by the restriction of the exterior derivative of the contact form
\begin{equation*}
    \omega_\text{H}=-d\theta|_\text{H}=:\omega|_\text{H},
\end{equation*}
and the fact that it is non-degenerate on $\text{H}$ is tantamount to the bilinear form
\begin{equation*}
    \eta:=\theta \otimes \theta + \omega \in \Sec{\otimes^2\Cot M}
\end{equation*}
being non-degenerate, thus giving a musical isomorphism
\begin{equation*}
        \begin{tikzcd}
        \Tan M \arrow[r, "\eta^\flat",yshift=0.7ex] & \Cot M \arrow[l,"\eta^\sharp",yshift=-0.7ex]
        \end{tikzcd}.
\end{equation*}
It then follows that the datum of an exact contact structure $(M,\theta)$ is equivalent to the datum of a non-degenerate Lichnerowicz manifold $(M,\pi,R)$. This is easily seen from the fact that the presence of a musical isomorphism allows for the natural assignment:
\begin{equation*}
    \theta = \eta^\flat R \qquad \omega=\eta^{\flat\flat}\pi
\end{equation*}
where
\begin{equation*}
    \eta^{\flat\flat}\pi(v,w):=\pi(\eta^\flat v, \eta^\flat w) \qquad \forall v,w\in\Tan M.
\end{equation*}

\section{L-Dirac Geometry} \label{LDiracGeometry}

In Section \ref{DiracGeometry} Dirac structures were shown to encompass Poisson and presymplectic manifolds unified in a formalism of involutive Lagrangian subbundles of Courant algebroids. With the introduction of the category of line bundles in Section \ref{CategoryOfLineBundles} and their interpretation as unit-free manifolds, Jacobi and precontact manifolds were presented as the unit-free analogues of Poisson and presymplectic manifolds. From this perspective, one expects that unit-free analogues of Courant algebroids and Dirac structure exist in the context of line bundles and L-vector bundles. This is indeed the case and this section is dedicated to the presentation of the basic definitions of such formalism.\newline

Recall that the basic example of a Courant space is the direct sum of a generic vector space with its dual, in the same spirit, we can consider a L-vector space $V^L\in\LVect$ and construct the L-direct sum
\begin{equation*}
    \mathbb{V}^L:=V^L\oplus V^{*L}
\end{equation*}
where we note that the direct sum is possible since the L-dual of a L-vector space is defined to have the same line component. It is clear that, in analogy with Courant spaces, there is a L-vector space morphism $\Proj_1:\mathbb{V}^L\to V^L$ and a non-degenerate symmetric $L$-valued bilinear form defined by
\begin{equation*}
    \langle v+\alpha ,w+\beta\rangle_{\mathbb{V}}:=\tfrac{1}{2}(\alpha(v)+\beta(w))
\end{equation*}
where elements of the L-dual are regarded as linear maps $\alpha,\beta : V\to L$, the L-vector space $\mathbb{V}^L$ is called the \textbf{standard L-Courant space}. This motivates the definition of a general \textbf{L-Courant space} as a L-vector space $C^L$ together with a L-vector space morphism to another L-vector space with the same line component $\rho:C^L\to V^L$ called the \textbf{anchor} and a non-degenerate symmetric $L$-valued bilinear form $\langle,\rangle:\odot^2 C\to L$. Note that the presence of a non-degenerate $L$-valued bilinear form, for some fixed line $L$, is equivalent to the presence of a conventional non-degenerate bilinear form. In fact, one can regard the datum of a $L$-valued bilinear form as coherent family of ordinary bilinear forms indexed by the choices of basis in $L$. Following from this remark, we see that the usual notions of musical isomorphism, orthogonal complements, isotropic subspaces and coisotropic subspaces appear naturally in L-Courant spaces in direct analogy to how they appear in ordinary Courant spaces, as presented in Section \ref{LinearCourant}. A L-Courant is called \textbf{exact} when it sits in the following short exact sequence of L-vector spaces
\begin{equation*}
\begin{tikzcd}
0 \arrow[r] & V^{*L} \arrow[r, "j"] & C^L \arrow[r, "\rho"] & V^L \arrow[r]  & 0
\end{tikzcd}
\end{equation*}
where $j:=\sharp \circ \rho^{*L}$ and $\sharp:C^{*L}\to C^L$ is the musical isomorphism induced by $\langle,\rangle$. Proposition \ref{IsotropicSplittingLinearCourant} can be generalized to L-Courant spaces in an obvious way thus ensuring that isotropic splittings of the exact sequence above always exist and then giving isomorphisms of exact and standard L-Courant spaces. For the reminder of this section, all L-Courant spaces will be assumed to be exact unless otherwise stated.\newline

Dirac spaces of ordinary Courant spaces find their obvious generalization in L-Courant spaces as \textbf{L-Dirac spaces}: we say that $D\subset C^L$ is a L-Dirac space in the L-Courant space $(C^L, \rho:C^L\to V^L,\langle,\rangle)$ when it is a Lagrangian subspace that sits on a short exact sequence of L-vector spaces
\begin{equation*}
\begin{tikzcd}
0 \arrow[r] & \overline{W} \arrow[r, "j"] & D \arrow[r, "\rho"] & W \arrow[r]  & 0
\end{tikzcd}
\end{equation*}
for some $W\subset V^L$ and $\overline{W}\subset V^{*L}$. Once more, the results of ordinary Dirac spaces presented in section naturally generalize to the L-Courant setting. We note that the 2-forms that naturally appear in L-Courant spaces, namely, differences of isotropic splittings and the 2-forms associated to L-Dirac spaces, are now $L$-valued. It is then clear that a discussion parallel to that of Section \ref{LinearCourant} regarding Courant morphisms and the categories $\Crnt$ and $\Dir$ follows from the proposed definitions in the context of L-Courant spaces. We find \textbf{L-Courant morphisms} as L-Dirac spaces covering graphs of L-vector space morphisms that compose as Lagrangian relations to form \textbf{the category of L-Courant spaces} $\LCrnt$. Following from this it is immediate to define backward and forward images of L-Dirac spaces then naturally identifying \textbf{L-Dirac morphisms} and \textbf{the category of L-Dirac spaces} $\LDir$.\newline

As it will become apparent from the discussion below, the notions of Courant algebras and Leibniz cohomology are broad enough to capture the algebraic aspects of the objects that generalize Lie and Courant algebroids in the context of line bundles. It is, therefore, not necessary to identify L-Leibniz algebras or L-Courant algebras as new mathematical objects different from Leibniz and Courant algebras introduced in Section \ref{CourantAlgebras}.\newline 

Let us now motivate the natural generalization of Lie algebroids in the context of line bundles. First note that the general definition of Lie algebroid of Section \ref{LieAlgebroids} can be summarized as the datum of a vector bundle $\alpha:A \to M$ with a vector bundle morphism to the tangent bundle of the base $\rho: A\to \Tan M$ and a Lie bracket on sections $(\Sec{A},[,])$ satisfying the Leibniz compatibility condition. Inspired by the interpretation of line bundles as unit-free manifolds of Section \ref{UnitFreeManifolds} where the der bundle was argued to be the line bundle analogue of the tangent bundle, we could simply generalize the definition of Lie algebroid by replacing $M$ with a line bundle $\lambda:L\to M$ and $\Tan M$ with the der bundle of the line bundle $\Der L$. This is indeed gives the desired generalization, as it will be clear from all the constructions follow below, and so a L-vector bundle $\alpha:A^L\to M$ is called a \textbf{Jacobi algebroid}\footnote{We could have coined the term \textbf{L-Lie algebroid} here to maintain naming consistency with the rest of the thesis but we choose to adopt the well-established nomenclature of Jacobi algebroid found in the literature to avoid the cacophony of the name L-Lie.} when there is a L-vector bundle morphism covering the identity $\rho:A^L\to \Der L$, called the anchor, and a Lie bracket on sections $(\Sec{A},[,])$ for which $\rho$ induces a Lie algebra morphism and such that the following Leibniz identity holds
\begin{equation*}
    [a,f\cdot b]=f\cdot [a,b]+\delta\circ \rho(a)[f]\cdot b,
\end{equation*}
where $\delta: \Der L \to \Tan M$ is the anchor of the der bundle. Then, in direct analogy with tangent bundles being the trivial example of Lie algebroids, we find der bundles of arbitrary line bundles, with identity as anchor and with the commutator of derivations as Lie bracket, to be the trivial example of Jacobi algebroids. Then we can rewrite the der functor as mapping line bundles to Jacobi algebroids
\begin{equation*}
    \Der : \Line_\Man \to \textsf{Jacb}_\Man.
\end{equation*}
Note how the given definition of Jacobi algebroid required $\rho$ to induce a Lie algebra morphism from sections to derivations, which, again, is analogous to the anchor map of an ordinary Lie algebroid inducing a Lie algebra morphism from sections to vector fields. However, we remark that proposition \ref{SymbolLieAlgebra} ensured that anchors of Lie algebroids induce Lie algebra morphisms as a consequence of the natural definition of Lie algebroids as derivative Lie algebra structures on ordinary vector bundles. Indeed, if we were to define Jacobi algebroids as above but dropping the requirement of $\rho$ inducing a Lie algebra morphism  this will not result in an equivalent definition to the one given above, essentially because $\Ker{\delta}\neq 0$ generically. We accept this discrepancy here so that the definitions of this section recover the standard notions found in the existing literature where Jacobi algebroids are commonly defined as Lie algebroids together with a Lie algebroid representation on the derivations of some line bundle. This issue is related to the general problem of finding an analogue of the ring structure of functions in the unit-free generalization as sections of a line bundle. These topics will be analyzed in detail in chapter \ref{DimAlgebraGeometry}.\newline

Since the anchor of a Jacobi algebroid is a Lie algebroid morphism $\rho:A^L\to \Der L$, thus recovering what was called a Lie algebroid representation of $A$ on the vector (line) bundle $L$, it follows that there is a Lie algebroid cohomology complex with values in $L$, $(\Omega^\bullet_L(A),\wedge,d_A^L)$. Upon inspection of the definition given for $\Omega^\bullet_L(A)$ in Section \ref{LieAlgebroids}, we readily check that the Lie algebroid cohomology with values in $L$ of $A$ is nothing but the regular de Rham complex of the L-dual vector bundle $\Omega^\bullet_L(A)=\Sec{\wedge^\bullet A^{*L}}=:\Omega^\bullet(A^L)$. This further establishes the choice of definition of Jacobi algebroid as the line bundle equivalent of Lie algebroids. The complex $\Omega^\bullet(A^L)$ is called the \textbf{Jacobi algebroid cohomology (with trivial coefficients)} and the associated identities for the obvious definitions of inner product and Lie derivative are called the \textbf{Cartan calculus} of the Jacobi algebroid. Similarly to the Gerstenhaber algebra of a Lie algebroid, the Lie bracket of a Jacobi algebroid can be extended uniquely into a graded Jacobi bracket on the multisections $(\Sec{\wedge^\bullet A^L},[,])$. Having identified the de Rham complex of a Jacobi algebroid we naturally define \textbf{morphism of Jacobi algebroids} as a L-vector bundle morphism $F:A^L\to B^{L'}$ such that the induced pull-back on the de Rham complexes $F^*:\Omega^\bullet(B^{L'})\to \Omega^\bullet(A^L)$ is a morphism of differential graded algebras.\newline

An important result for Jacobi algebroids, which links them intimately with Jacobi structures, is the analogue of the one-to-one correspondence between linear Poisson structures and Lie algebroids proved in proposition \ref{LieAlgebroidLinearPoisson}. Fibre-wise linear Jacobi structures on L-vector bundles were defined in Section \ref{ContactGeometry} as Jacobi brackets compatible with the submodule of spanning sections on the total space of the L-vector bundle. Given a Jacobi algebroid $(\alpha:A^L\to M,\rho,[,])$, following an entirely analogous construction to that of proposition \ref{LieAlgebroidLinearPoisson}, we find a correspondence with a linear Jacobi structure on the L-dual bundle $A^{*L}$ realized by imposing the following identities on spanning sections
\begin{align*}
\{l_{a},l_{b}\} & = l_{[a,b]}\\
\{l_{a},\alpha^*s\} & = \alpha^*\rho(a)[s]\\
\{\alpha^*s,\alpha^*r\} & = 0
\end{align*}
for all $s,r\in\Sec{L}$ and $a,b\in\Sec{A^L}$. This construction is easily checked to be one-to-one from the fact that L-duals behave like ordinary duals, i.e. there is a canonical isomorphism of L-vector bundles $(A^{*L})^{*L}\cong A^L$ and, furthermore, functorial with respect to Jacobi maps on L-vector bundles restricting to L-vector bundle morphisms covering diffeomorphisms.\newline

In the same spirit of how Courant algebroids were motivated by Lie bialgebroids at the beginning of Section \ref{CourantAlgebroids}, let us introduce the notion of \textbf{Jacobi bialgebroid} as a pair of a L-vector bundle and its L-dual $(A^L,A^{*L})$ each with a Jacobi algebroid structure, denoted by $(A^L,\rho,[,])$ and $(A^{*L},\rho_*,[,]_*)$, satisfying the compatibility condition
\begin{equation*}
    d_{A^{*L}}^L[a,b]=[d_{A^{*L}}^La,b]+[a,d_{A^{*L}}^Lb]
\end{equation*}
for all $a,b\in\Sec{A^L}$, where $[,]$ denotes the Schouten-Jacobi bracket of the multisections $(\Sec{\wedge^\bullet A^L},[,])$ and $d_{A^{*L}}$ is the graded differential of $A^{*L}$ defined on the complex $\Sec{\wedge^\bullet (A^{*L})^{*L}}\cong \Sec{\wedge^\bullet A^L}$. Just like in the case of ordinary Lie bialgebroids, this compatibility condition allows to form Leibniz representations of $\Sec{A^L}$ and $\Sec{A^{*L}}$ on each other via the Cartan calculus identities thus inducing two semidirect product Courant algebras on $\Sec{A^L}\oplus\Sec{A^{*L}}$ and then forming the sum Courant algebra with bracket explicitly given by
\begin{equation*}
    [a\oplus\alpha ,b\oplus \beta]:= [a,b] + \LDer_\alpha b -i_\beta d_{A^{*L}}^La \oplus [\alpha,\beta]_* + \LDer_a\beta -i_bd_{A}^L\alpha.
\end{equation*}
This is then a Courant algebra structure on the sections of the L-vector bundle $A^L \oplus A^{*L}$ whose fibres are clearly L-Courant spaces with the obvious bilinear form induced by the dual pairing and the anchor map defined as $\rho\oplus \rho_*:A^L \oplus A^{*L}\to \Der L$.\newline

This construction motivates the definition of \textbf{L-Courant algebroid} as a L-vector bundle $E^L\to M$ together with a structure $(E^L,\langle,\rangle,\rho: E^L\to \Der L,[,])$ where $\langle,\rangle\in\Sec{\odot^2 E^{*L}}$ is a symmetric non-degenerate $L$-valued bilinear form and $\rho: (\Sec{E^L},[,])\to (\Der L,[,])$ is a Courant algebra structure on the module of sections over the Lie algebra of derivations such that
\begin{align*}
    [a,f\cdot b] &=f\cdot [a,b]+\delta\circ\rho(a)[f]\cdot b\\
    \rho(a)[\langle b , c \rangle] &=\langle [a,b] , c \rangle + \langle b , [a,c] \rangle\\
    [a,a] &=D\langle a , a \rangle
\end{align*}
for all $a,b\in\Sec{E^L}$, $f\in\Cin{M}$, where $D:=\sharp\circ \rho^* \circ j^1:\Sec{L}\to \Sec{E^L}$ is defined from the musical isomorphism $\sharp:E^{*L}\to E^L$ induced by the non-degenerate $L$-valued bilinear form and the jet extension map $j^1:\Sec{L}\to \Sec{\Jet^1 L}\cong\Sec{\Der L^{*L}}$. Notice how these defining axioms can be regarded as the line bundle analogues of the axioms of ordinary Courant algebroids given in Section \ref{CourantAlgebroids} simply by interchanging der bundles and tangent bundles and accounting for forms taking fibre-wise values on $L$ instead of $\Real$. It is in this sense that the notion of Courant algebra of Section \ref{CourantAlgebras} was general enough to capture the modules of sections of both Courant algebroids and L-Courant algebroids.\newline

A \textbf{L-Dirac structure} in a L-Courant algebroid $(E^L,\langle,\rangle,\rho: E^L\to \Der L,[,])$ is a maximally isotropic subbundle $L\subset E^L$ (possibly supported on a submanifold) with involutive sections $[\Sec{L},\Sec{L}]\subset \Sec{L}$. At this point, we could clearly carry out definitions and constructions for L-Courant algebroids and L-Dirac structures that mirror those given in our discussion of Section \ref{CourantAlgebroids} for ordinary Courant algebroids and Dirac structures. In particular, by interchanging the roles that $\Tan M$  and $\Cot M$ played for an ordinary Courant algebroid $E\to M$ with $\Der L$ and $\Jet^1 L$ for a L-Courant algebroid $E^L\to M$ and accounting for forms taking values in $L$ instead of $R$, we find the obvious definitions of \textbf{exact L-Courant algebroids}, \textbf{L-Courant tensor}, \textbf{product L-Courant algebroid} with the product $\boxplus$ of L-vector bundles, \textbf{morphisms of L-Courant algebroids} as L-Dirac structures supported on L-graphs of factors of line bundles and \textbf{L-Courant maps} as graphs of L-vector bundle morphisms that are L-Dirac structures in the product L-Courant algebroid. These then allow us to identify \textbf{the category of L-Courant algebroids} $\LCrnt_\Man$ and \textbf{the category of L-Dirac structures} $\LDir_\Man$.\newline

To close this section we prove that ordinary Courant algebroids and Dirac structures are recovered as particular cases of L-Courant algebroids and L-Dirac structures.

\begin{prop}[L-Courant Algebroids encompass Courant Algebroids]\label{LCourantCourant}
Let $\normalfont (E,\langle,\rangle,\rho: E\to \Tan M,[,])$ be a Courant algebroid, then the L-vector bundle $E^{\Real_M}$ is naturally a L-Courant algebroid. Furthermore, if $K\subset E$ is a Dirac structure, $K\subset E^{\Real_M}$ is a L-Dirac structure with respect to the L-Courant algebroid structure.
\end{prop}
\begin{proof}
This follows from the results about trivial line bundles presented at the end of Section \ref{UnitFreeManifolds}. In particular, recall that the der bundle of a trivial line bundle is $\Der \Real_M\cong \Tan M\oplus \Real_M$, the jet bundle is $\Jet^1\Real_M\cong \Cot M\oplus \Real_M$ and that the jet prolongation map $j^1:\Sec{\Real_M}\to \Jet^1 \Real_M$ is given by $j^1=d\oplus \Id_{\Real_M}$, where $d$ is the ordinary exterior derivative. The Courant bilinear form $\langle,\rangle$ is trivially promoted to a $\Real_M$-valued bilinear form and the anchor map defines a L-Courant anchor map simply by setting $\rho \oplus 0: E\to \Tan M\oplus \Real_M$. Since the symbol of the Lie bracket of derivations of a trivial line bundle becomes $\delta = \Proj_1 : \Tan M\oplus \Real_M \to \Tan M$ and we have $\Sec{\Real_M}\cong\Cin{M}$, it is clear that the three compatibility conditions of the Courant bracket on $E$ are equivalent to the same compatibility conditions of the L-Courant bracket on $E^{\Real_M}$. It is also easy to see that these make the defining conditions for a Dirac structure on $E$ equivalent to the defining conditions of a L-Dirac structure on $E^{\Real_M}$ with respect to the L-Courant structure defined above.
\end{proof}

\chapter{History and Foundations of Metrology and Phase Space} \label{FoundationsMetrologyPhaseSpace}

\section{A Brief History of Metrology} \label{BriefHistoryOfMetrology}

Consider the following school textbook problem: ``Question: If water flows out of a tap at $4.3\, \text{L}/\text{min}$, how long will it take for an empty $300\, \text{cm}^3$ cup to fill up? Answer: The flow through a pipe measures volume of fluid traversing the pipe per unit time, then setting $F=V/T$ with $F=4.3\, \text{L}/\text{min}$ and $V=300\, \text{cm}^3=0.3\, \text{L}$ we find $T=0.06\, \text{min}$ so the cup will fill up in approximately $4$ seconds''. In order to produce this result, a pupil will be required to confidently operate with the basic concepts of physical quantity and units of measurement. As illustrated by this example, the essential features of dimensional analysis (see \cite{barenblatt1996scaling}) can be summarized as follows: physical quantities are given by rational numerical values and are accompanied by units; different units may be used to express the same physical quantity, in which case conversion factors exist between them; arithmetical products of physical quantities induce multiplication of numerical values and the formation of composite units; arithmetical addition can only happen between physical quantities of the same kind, i.e. measuring the same characteristic. In this section we present an overview of the long trace of cultural and scientific ideas that has led to these modern conventions in metrology and dimensional analysis.\newline

\begin{figure}[h]
\centering
\includegraphics[scale=0.4]{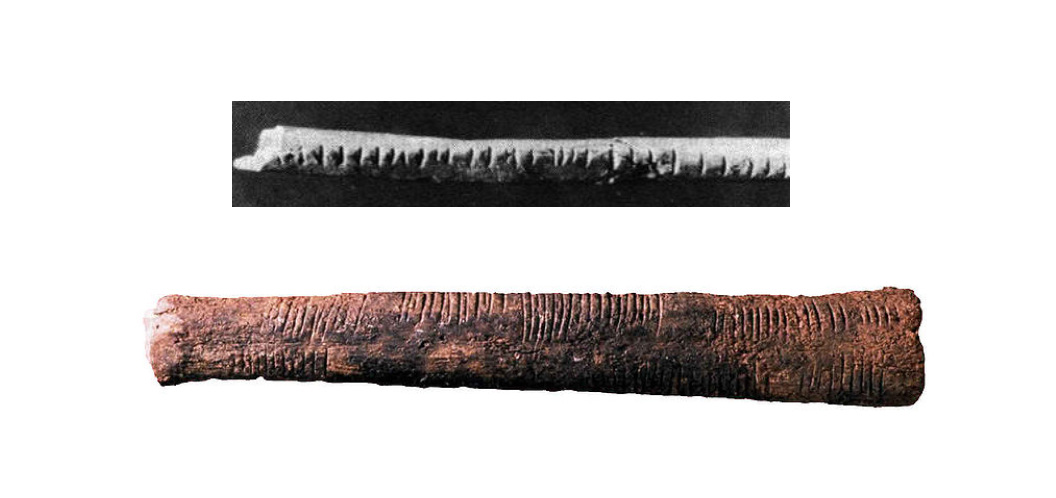}
\caption{The Lebombo bone $44,000\, \text{BC}$ (top) and the Ishango bone $20,000\, \text{BC}$ (bottom).}
\label{bones}
\end{figure}

The earliest application of measurement procedures in hominids likely predates the oldest available evidence of metrological artifacts made by \emph{Homo Sapiens}. It is hard to give a precise estimate of when, or by which species in the \emph{Homo} genus, primitive units of measurement were first used. To gain some insight on this question, we first look at the ample evidence of non-human animals with sufficiently developed nervous systems displaying relatively sophisticated analytical behaviour that can only be attributed to the use of rudimentary mathematical ideas. For instance, higher mammals such as elephants and whales, are known to use complex cause-and-effect connections between events in their environment to adapt creatively to new situations (see \cite{kuczaj2009intelligent} for a study on dolphin intelligence) and some avian species have been observed to have the ability to solve simple problems that involve counting small quantities and comparing quantities of different objects (see \cite{taylor2010complex} for evidence of inventive problem-solving skills in crows). Since no direct evidence of metrological activity in non-human animals has been found, we are led to conclude that the appearance of primitive metrology in early humans, who shared their basic cognitive mathematical prowess with other species, was largely due to their main distinguishing features: dexterous hands, which allowed for the ability to manipulate objects with precision and that, ultimately, resulted in sophisticated tool use, and speech, which gradually evolved to communicate intersubjective and environmental information with increasing levels of accuracy. The first direct evidence of metrological artifacts are the Lebombo and Ishango bones (depicted in figure \ref{bones}, see \cite[Chapter 2]{bangura2011african}), pieces of animal bone carved with what looks like tallies, found in central Africa and dated around $44,000\, \text{BC}$ and $20,000\, \text{BC}$, respectively. These are, however, most likely late examples of that kind of tool, as the technological sophistication and intellectual ability to produce them were present in hominids as early as $500,000\, \text{BC}$, the approximate dating of the first controlled use of fire (see \cite{balter1995did}).\newline

Metrological concepts and techniques will develop, like any other aspects of ancient human culture, in multiple complex streams of ideas throughout millennia. This only partially changed with the appearance of the first standardized systems for large civilizations, despite disconnected communities of humans continuing to use their particular units of measurement and metrological procedures independently to this day. This is evidenced by the diverse and multifaceted measurement and number systems of the few isolated tribes that have units of length based on body parts of living people and counting systems based on arboreal fruits (see \cite{walker2015protecting}).\newline

\begin{figure}[h]
\centering
\includegraphics[scale=0.7]{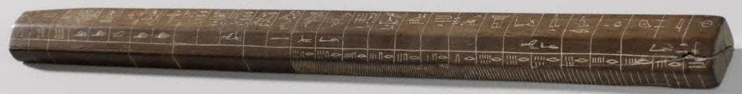}
\caption{A cubit measuring rod dating from $\text{c.}1330\, \text{BC}$.}
\label{cubit}
\end{figure}

There were two major driving forces behind the development of uniform systems of measurement within large human communities: firstly, planning and arrangement of spatio-temporal circumstances, such as keeping track of seasons, distributing land for agriculture or building stable structures, and secondly, and perhaps more importantly, commerce, which created the need to precisely quantify goods of varying nature for their efficient trade. Ancient civilizations of Mesopotamia and the Indus valley are known to have used units of length based on body parts of ruling figures and units of time based on astronomical cycles as early as $4000\, \text{BC}$ (see \cite{hodgkin2005history}). The first attested standard unit is the Egyptian royal cubit (a well preserved specimen is depicted in figure \ref{cubit}, see \cite{scott1942egyptian}), probably based on the length of the forearm of an unknown early Egyptian king as measured from the elbow to the tip of the middle finger. The first recorded use of the cubit appears in a hieroglyph describing the flood levels of the river Nile dated around $3000\, \text{BC}$ (see \cite{clagett1989ancient,blanco2017atlas}). Several measuring rods using the cubit, or the many variations that appeared later, have been linked to Mediterranean civilizations as late as the Romans with written references to the cubit appearing also in early biblical texts and the Talmud (see \cite{scott1959weights}). This shows that, despite the inevitable variance of exact length standard given by the limitations of the time, the abstract notion of a uniform and universal unit of measurement must have been well-established by approximately $2500\, \text{BC}$.
 
\begin{figure}[h]
\centering
\includegraphics[scale=0.6]{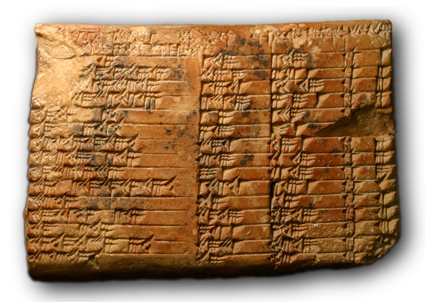}
\caption{The Plimpton 322 clay tablet, originating from Larsa (present-day southern Iraq) $\text{c.}1800\, \text{BC}$.}
\label{plimp}
\end{figure}

The early development of metrological theory and systems of units went, naturally, hand in hand with the development of general mathematical ideas about arithmetic and geometry. This is particularly evident in the case of Assyro-Babylonian mathematics. There are abundant samples of clay tablets produced by peoples of Mesopotamia between $2000\, \text{BC}$ and $600\, \text{BC}$ that show mathematical ideas of considerable sophistication and a notation that has many similarities with our modern conventions (see \cite{fowler1998square}). As a reference, we could consider the Plimpton 322 clay tablet (depicted in figure \ref{plimp}), which contains a list of Pythagorean triples expressed in cuneiform script. Upon careful interpretation, this tablet showcases a few remarkable features generally found in Assyro-Babylonian mathematics: arithmetic involving addition, multiplication, subtraction, division and primitive notions of powers and roots; a sexagesimal (base $60$) numeral system that used place-value notation that aided in the computation of fractions; formulas for the areas and volumes of several shapes involving approximations of irrational numbers such as $\sqrt{2}$ or $\pi$. This early mathematical tradition had far-reaching influences across Europe, northern Africa and India, particularly, early Greek mathematicians borrowed greatly from the Babylonians. This influence can be seen to this day in the $60$ subdivisions of hours and minutes or the $360$ degrees of a circumference.

\begin{figure}[h]
\centering
\includegraphics[scale=0.15]{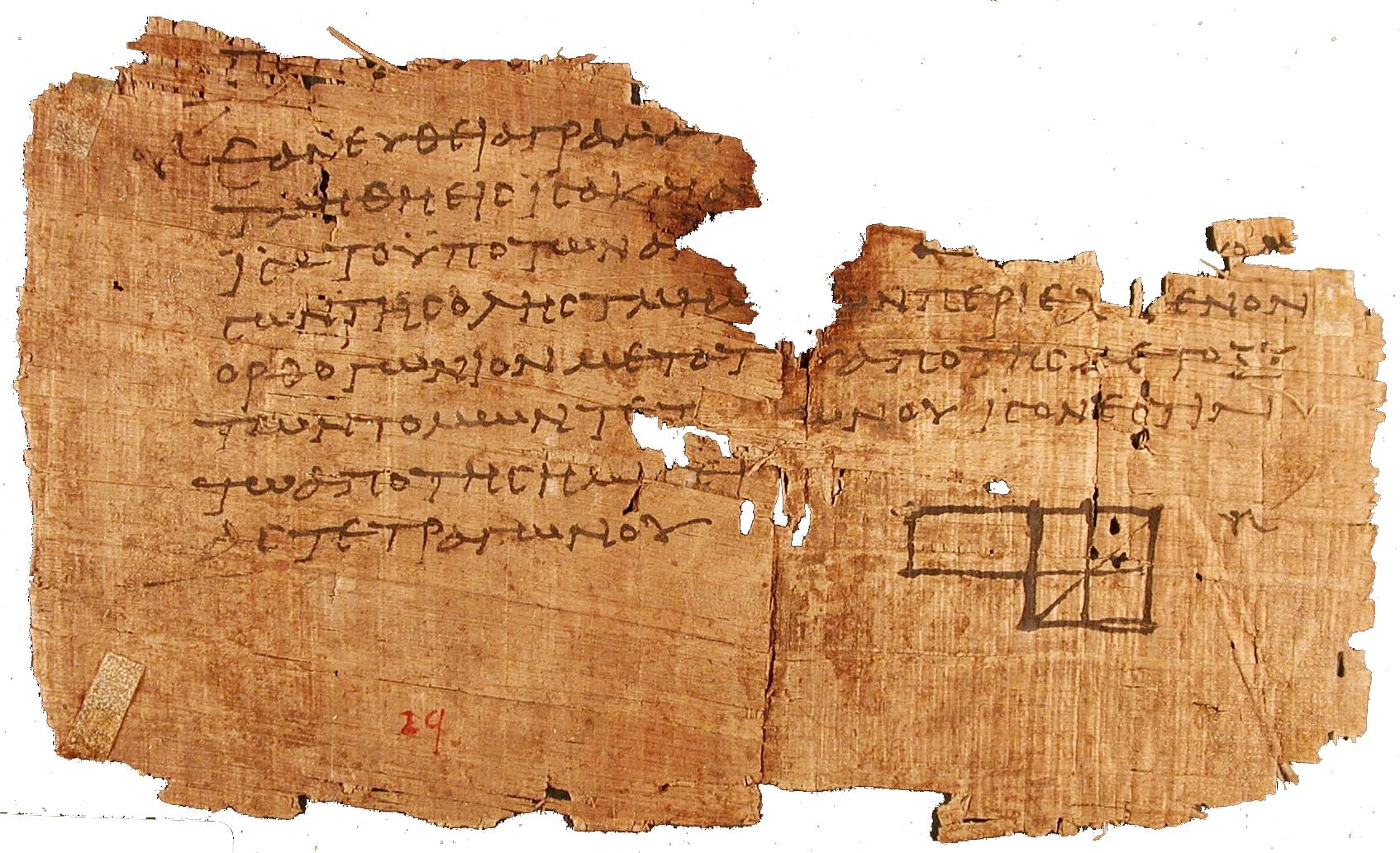}
\caption{The oldest surviving fragment of an edition of Euclid's \emph{Elements}, a papyrus from central Egypt dated $\text{c.}100$. Complete versions of the \emph{Elements} have reached us via transmission of the numerous translations and copies that were made by the civilizations that have risen and fallen since Euclid's work $\text{c.}300\, \text{BC}$.}
\label{elements}
\end{figure}

The works of the ancient Greek mathematicians, particularly those of the Hellenistic period, from $350\, \text{BC}$ to $50\, \text{BC}$, after the cross-pollination of ideas with the Egyptian and Babylonian civilizations, represent a crucial milestone for the development of modern mathematics and science (see \cite[Chapters 2,3]{hodgkin2005history}). Of singular importance is the treatise \emph{Elements} from $\text{c.}300\, \text{BC}$ (figure \ref{elements}, see \cite{heath1956thirteen}), where Euclid (Alexandria, $\text{c.}350-\text{c.}250\, \text{BC}$) compiles a breadth of results in planar geometry and number theory presented following a strict logical structure with theorems and proofs, thus laying the foundations of the contemporary format of mathematical texts (see \cite{mueller1981philosophy}). The \emph{Elements} are, furthermore, of particular importance in the history of metrological theory since they contain notes on how equations involving quantities of different kinds, such as length and area, should be subject to homogeneity conditions. This represents the first written evidence of a theory that abstractly treats physical quantities subject to algebraic rules in a mathematical framework, thus hinting to what will be later known as dimensional analysis. The deep influence of Euclid's \emph{Elements} in the history of science and mathematics is clear from the fact that they were used as the standard of mathematical tuition, particularly in the area of geometry, until the late $1800$s. It should be noted then, that the seminal works of the parents of modern science were directly influenced by this treaty, as evidenced by both the prominent role of geometry in the early mathematical formulation of physics and the preponderance of the theorem-proof style of presentation of arguments.\newline

The rise of empires such as the Roman $\text{c.}30\, \text{BC}$ or Chinese $\text{c.}200$, brought units of measurement that were uniformly used over ever greater expansions of territory. For example, the Roman pes (Latin for `foot') was a virtually universal unit of length for the peoples within the empire's domains (see \cite{duncan1980length}). So much so, that our modern unit of foot, which derives directly from it, only differs by $0.03\%$ in magnitude after $2000$ years. Another contributing factor was, of course, the refinement of tool-making technologies, which enabled for more accurate and consistent standards to be used for longer periods of time. For many centuries, there was no clear conceptual development of the basic ideas of metrology besides the gradual improvements in accuracy and the increasing levels of standardization.\newline

All the written evidence of the symbolic manipulation of units of measurement predating $1300$ shows that arithmetical operations were only performed with numbers quantifying the same kind of magnitude. In most instances, in fact, units were entirely omitted from the calculations as they were assumed from contextual information. This meant that the mathematical treatment of physical quantities was lacking one of the key aspects of our modern system: products and quotients of magnitudes of different kinds to describe derived magnitudes such as momentum, velocity or density. Although it is clear that Greek science involved concepts similar to the modern notions of velocity, density or torque, direct evidence of their mathematical implementation has not been found. The first known instance of a quantity involving the product of two magnitudes of different kinds appears in the astronomical treaties of Johannes de Muris (France, $1290 - 1355$) (see \cite{de1998arte}). Throughout the $1300$s and $1400$s, the schools of Oxford and Paris developed precise notions in kinematics and statics, many of them involving physical quantities that required products and quotients of measurements of time, space and weight (see \cite[Chapter 1]{nolte2018galileo}). These set the metrological foundations upon which Galileo would eventually base his empirical methods.

A good reference for the state of metrological methods and conventions around this time is the publication in $1585$ of a booklet by Simon Stevin (Belgium, $1548 - 1620$) entitled \emph{De Thiende}, that roughly translates as `the tenths', where a systematic decimal notation for fractions is introduced and shown in practice for several arithmetical operations (see \cite{stevin1965thiende}). From his writings we can infer that the notion of physical quantity as a rational numerical value together with a unit of measurement must have been well established among scientists and mathematicians by $1600$. In this booklet, a rudimentary version of the method of unit fractions for the conversion of units of measurement can be found as well as algebraic operations involving units of more than one kind.\newline 

The next qualitative advancement in the history of metrology came, perhaps unsurprisingly, with the Scientific Revolution of the $1500$s and $1600$s in Europe. As has already been mentioned, the profoundly influential works of Galileo Galilei (Italy, $1564 - 1642$) and Isaac Newton (England, $1642 - 1726$) in the field of mechanics made explicit use of the metrological standards of the time, thus further establishing them, but it was the contributions of Rene Descartes (France, Sweden, $1596 - 1650$) that will have deeper consequences for the development of modern metrology. On the one hand, he is credited for being the first to use the term `dimension' - borrowed from the context of geometry - in a sense that we would attribute today to the concept of physical dimension (see \cite{macagno1971historico}). His writings on the topic of mechanics contain several comments that essentially amount to dimensional arguments for the justification of algebraic expressions of physical laws. Descartes' work represents the first theoretical, although far from systematic, use of dimensional analysis. The second contribution by Descartes was the mathematical legacy of analytical geometry and the modern conventions and notations of algebra. These ideas will prove instrumental in the centuries to come for the development of theories that will involve operations typically performed in the context of geometry, such as changing scale via congruences in planar geometry, that are applied to non-geometric entities, such as the change of reference unit of time or mass.\newline

During the $1700$s, as part of the accelerated scientific, technological, artistic and political progress of the Age of Enlightenment, many advancements were made across several disciplines. In the field of physics, most of the efforts were concentrated either on the experimental side, exploiting the new technological possibilities enabled by the Industrial Revolution, or the theoretical side, discovering deeper mathematical patterns in the description of natural phenomena (see \cite{wolf2019history}). The basic notions of metrology and dimensional analysis were still rudimentary and had not changed much since their inception during the late Middle Ages. Although Leonhard Euler (Switzerland, Russia, $1707 - 1783$) includes a chapter on the subject of units and homogeneity in his $1765$ book \emph{Theoria motus corpum solidorum seu rigidorum}, his comments remained largely obscure for the scientific community at the time and it will not be until half a century later that a more comprehensive treatment of dimensional analysis will appear. Around this time, coinciding with the American and French revolutions, the metric system of natural units was developed in France as a first attempt at a truly universal standardization of measurement. This will later develop into the modern S.I. of units, in common use today by the international scientific community, and it can be regarded as the end of the millenniums-long history of complicated evolution of standards of measurement.\newline

It was Joseph Fourier (France, $1768 - 1830$), in the third edition of his treaty on thermodynamics \emph{Theorie analytique de la chaleur} \cite{baron1822theorie} in $1822$, who first identified the algebraic properties of physical quantities in an abstract sense, laying the foundations of modern dimensional analysis. All the basic ideas of the modern formalism used by scientists and engineers today can be found in Fourier's work. Another relevant feature of Fourier's writings is that they contain abundant comments on the importance of a systematic treatment of physical dimension and the usefulness of identifying dimensionless coefficients to characterize physical laws. In $1877$, Strutt Rayleigh (England, $1842 - 1919$) published his \emph{The Theory of Sound} \cite{strutt1945theory} containing a section entitled \emph{The Method of Dimensions} in which concrete applications of the ideas introduced by Fourier can be found alongside a partial proof of the fundamental theorem of dimensional analysis, a result that gives a systematic way to identify a minimal set of dimensionless coefficients characterizing a physical law involving several physical quantities of different dimensions. This marks the beginning of what we know today as dimensional analysis.\newline

A first complete proof of the fundamental theorem of dimensional analysis was given by Aime Vaschy (France, $1857 - 1899$) in his paper \emph{Sur les lois de similitude en physique} \cite{vaschy1892lois} of $1892$. This proof was, unfortunately, later forgotten by the international community. Some years later, starting in $1914$, Edgar Buckingham (United States, $1867 - 1940$) published a series of papers providing a new and more accessible proof which (see \cite{buckingham1915principle}), together with Rayleigh's introductory articles published in \emph{Nature} around that time, launched the subject of dimensional analysis into the mainstream of science and engineering. Thereafter, concepts and techniques of dimensional analysis have proved extremely powerful in many different fields, both theoretical and applied, specially before the raise of numerical methods aided by computers in the $1960$s. It is perhaps due to the more pragmatical nature of the vast majority of applications of dimensional analysis to physics that a revision of the foundations of the topic has not been conducted until recent years. This thesis aims to provide a humble contribution in this direction.\newline

The crystallization of the long history of metrology outlined in this section can be seen today in the standards and common practice kept by the International Bureau of Weights and Measures (BIPM, from the French acronym). It is safe to say that this organization, representing the metrological conventions and scientific language shared by all the global community, has played a major role in the erection of the technological edifice of modern civilization.

\section{The Measurand Space Formalism} \label{MeasurandFormalism}

In this section we aim to define a formal notion of physical quantity that naturally encompasses the standard operational properties that dimensional analysis seeks to exploit. Our mathematical characterization should not only capture the essential features\footnote{See the opening paragraph of Section \ref{BriefHistoryOfMetrology} for a summary of these.} of units of measurement as used in practical science and engineering but also provide with a clear metrological interpretation of any formal constructions found therein. One of the goals of an axiomatic revision of the basic notions of theoretical metrology is the development of a unit-free formulation of modern physical theories where physical quantities are more than just $\Real$-valued variables. Since our proposal will not be based in any particular exiting definition of physical quantity found in the literature, we will use the less common term ``measurand'' to refer to this concept. As it will become apparent from our discussion below, this choice of terminology is further motivated by the emphasis that we place on the physical property that is to be measured rather than the measurement outcome itself.  The common usage of this term in metrology, as seen from the entry of the BIPM-endorsed International Vocabulary of Metrology \cite{bimp2012metrology} quoted below, indeed justifies our choice:\newline

\emph{2.3 (2.6)\textbf{Measurand}: quantity intended to be measured. NOTE 1 The specification of a measurand requires knowledge of the kind of quantity, description of the state of the phenomenon, body, or substance carrying the quantity, including any relevant component, and the chemical entities involved.}\newline

In order to motivate the mathematical definitions of measurand and unit of measurement, let us first depict a scenario of very rudimentary science that, however simple and mundane as it may seem, illustrates quite clearly the main formal properties of physical quantities and the process of measurement. Suppose that one is outside in the countryside wandering near a riverbed where plenty of canes and flat cobblestones can be found. These objects, canes and flat cobblestones, are the matter of scientific inquiry as one sets out to develop an elementary theory of length, area and weight in this primitive scenario, not unlike that of the humans who used the Lebombo and Ishango bones in prehistoric times, as discussed in Section \ref{BriefHistoryOfMetrology}. Let us denote a cane as $C\in \mathcal{C}$ and a flat cobblestone as $F\in\mathcal{F}$, where the calligraphic letters denote the somewhat abstract set of all existing objects within reach. The collection of all these objects $\mathcal{D}:=\mathcal{C}\cup \mathcal{F}$ will be called the domain of observation. One is interested in quantitatively comparable characteristics of objects within the domain of observation; these shall be called measurands. What we mean by ``quantitatively comparable'' here is that one should be able to give a numerical correspondence to a practical construction involving a pair of objects in the domain of observation. The examples of measurands below will make this notion precise. For these we shall only assume basic manipulations that can be practically made with canes and cobblestones together with the inherent human ability to count and keep a consistent memory record.
\begin{itemize}
    \item \textbf{Length.} Take two canes $C_1,C_2\in\mathcal{C}$; lay them on the ground side-by-side with near ends matching; by marking where far ends fall on the ground, move the near end of the shorter cane to the marked position; repeat this process until the far end of the shorter cane exceeds the far end of the longer cane or the two ends meet\footnote{The ``meet'' condition essentially amounts to the accuracy of the measurement process being described here. Assuming that there is a minimum accuracy is not only more realistic but is, in fact, operationally necessary for the measurement algorithm described to terminate in a finite number of steps. We will not deal with the notion of measurement uncertainty in our discussion and simply note that a more realistic numerical output will be two (instead of one) rational numbers corresponding to the lower and upper bounds of the ``meet'' condition, essentially given in this example by what the naked human eye can discern and the shape of the canes' ends.}; if the ends didn't meet, move the longer cane to its far end; repeat until the far ends meet; count the number of markings made for far ends of each cane $n_1,n_2\in\mathbb{N}$ after the process has stopped. The length comparison of cane $C_1$ to $C_2$ is defined as the quotient of this two numbers:
    \begin{equation*}
        l(C_1,C_2):=\frac{n_1}{n_2}\in\Rat.
    \end{equation*}
    The property that makes this construction into the desired quantitative comparison is the fact that its transitivity is precisely modelled by the multiplicative structure of $\Rat$. To see this, consider a third cane $C_3\in\mathcal{C}$ and perform the length comparisons with $C_1$ and $C_2$ to verify the following experimental fact:
    \begin{equation} \label{2outof3}
        l(C_1,C_2)\cdot l(C_2,C_3)\cdot l(C_3,C_1)=1.
    \end{equation}
    We will refer to this transitivity condition as the \textbf{2-out-of-3 property}. In fact, it is this property what guarantees the consistency of  the formal assignment of a fraction (seen as an equivalence class of solutions to the basic multiplicative equation of integers) to the pair of naturals resulting from the length comparison algorithm described above. We can then symbolically summarize the length comparison procedure as an assignment of the form $l:\mathcal{C}\times \mathcal{C}\to \Rat$ satisfying the 2-out-of-3 property.
    \item \textbf{Area.} Via a procedure that now takes flat cobblestones and compares their surface area by superimposing them in a similar fashion to the algorithm of the comparison of lengths, we characterize the comparison of area via a map $a:\mathcal{F}\times \mathcal{F}\to \Rat$ satisfying the 2-out-of-3 property.
    \item \textbf{Weight.} With the simple materials at hand one can construct a very primitive model of scale that will allow for the balancing of pairs of objects on an arm of varying span. An entirely analogue construction of quantitative comparison thus follows but now for any object in our observation domain $\mathcal{D}=\mathcal{C}\cup \mathcal{F}$. We shall call this comparison of weight and denote it by $w:\mathcal{D}\times \mathcal{D}\to \Rat$.
\end{itemize}
For the reminder of this section we make a few remarks about these concrete cases of Length, Area and Weight; nevertheless, the formal constructions that follow shall illustrate general points about arbitrary measurands.\newline

The measurand ``length of canes'' is conceptually thought of as the unifying property of all canes (or any straight object, in general) that allows for the operational system described above, including subdomain of applicability, algorithmic description of the procedure and numerical assignment. In particular, we see that a measurand is inseparably tied to the concrete class of physical objects on which it is defined. If we fix an object to set it as a reference or standard, say we take a cane $C_0\in\mathcal{C}$, we arrive at the notion of \textbf{unit of measurement}: all canes are assigned a rational number given by their length comparison with the unit $C_0$. In other words, a choice of a unit $C_0\in\mathcal{C}$ gives a map defined by
\begin{align*}
l_{C_0}: \mathcal{C} & \to \Rat\\
C & \mapsto l(C_0,C).
\end{align*}
Simply having a map of this form is not particularly useful but the 2-out-of-3 property of $l$ (\ref{2outof3}) allows to identify each arbitrary comparison of length with a quotient of the rational values of the map. More explicitly, having fixed a unit of length so we have a map  $l_{C_0}: \mathcal{C} \to \Rat$, the length comparison of any arbitrary pair of canes $C_1,C_2\in\mathcal{C}$ is given by
\begin{equation*}
        l(C_1,C_2)=\frac{l_{C_0}(C_1)}{l_{C_0}(C_2)}.
\end{equation*}
The length comparison $l(C_1,C_2)$ was defined from the experimental procedure with no reference whatsoever to the notion of unit, we thus expect the above equation to be independent of the choice of unit. Indeed, again by the 2-out-of-3 property (\ref{2outof3}), choosing a different unit $\overline{C}_0$ will introduce the factor $l(\overline{C}_0,C_0)$ in both numerator and denominator of the fraction, hence giving the same value for the length comparison. It is then easy to see that this leads to the usual notion of ``intrinsic length of a cane'' as the abstract entity that only becomes concrete once a unit is chosen. The language often employed to refer to measurement outcomes, that we will enforce throughout our discussion, is now fully justified: a cane $C$ is said to have length $l_{C_0}(C)\in\Rat$ in units of $C_0$.\newline

From the elementary measurement process described above for Length it is clear that given two canes $C_1,C_2\in\mathcal{C}$ we can make the far end of one meet with the near end of the other to virtually (i.e. to all metrological effects) create a new cane $C_1\parallel C_2$, we call it the \textbf{combination} of the canes $C_1$ and $C_2$. It is then an experimental fact that for any other cane $C_3$ we have:
\begin{equation*}
    l(C_1\parallel C_2,C_3)=l(C_1,C_3)+l(C_2,C_3).
\end{equation*}
Choosing another pair of canes $\overline{C}_1$ and $\overline{C}_2$ whose comparisons with the previous pair, $l(C_1,\overline{C}_1)=q_1$ and $(C_2,\overline{C}_2)=q_2$, are known, it is also an experimental fact that
\begin{equation} \label{combination}
    l(\overline{C}_1\parallel \overline{C}_2,C_3)=q_1\cdot l(C_1,C_3)+q_2\cdot l(C_2,C_3),
\end{equation}
and thus we see that the combination construction satisfies a $\Rat$-linearity property. Under the choice of a unit $C_0$, this implies, in particular, an additivity property for the lengths of the canes in units of $C_0$:
\begin{equation*}
    l_{C_0}(C_1\parallel C_2)=l_{C_0}(C_1)+l_{C_0}(C_2).
\end{equation*}
It is not hard to see that any sensible definition of the combination construction for objects comparable by area or weight will satisfy a similar linearity property.\newline

Next we formalize the intuition that areas are products of lengths. Suppose that one finds a corner-shaped piece of wood or stone to act as a standard wedge $\vee$ (a form of rudimentary protractor that only measures one non-flat angle) so it is possible to consistently array two canes in two independent planar directions. Taking any two canes $C_1$ and $C_2$, with the use of the standard wedge one virtually (i.e. to all metrological effects) forms a planar shape directly comparable in area to flat cobblestones, we denote it by $C_1 \vee C_2$ and call it the \textbf{tile} of $C_1$ and $C_2$. We can thus compare the area of a generic flat cobblestone $F\in\mathcal{F}$ to the area of a tile via $a(F,C_1\vee C_2)\in\Rat$, and we call that value the area of $F$ in units of $C_1$by$C_2$. It is clear form the measurement procedures for length and area comparison that the following are experimental facts for any flat cobblestone $F\in\mathcal{F}$ and canes $C_1,C_2,C_3\in\mathcal{C}$:
\begin{align}\label{tile1}
    a(F,C_1\vee C_2)&=a(F,C_2\vee C_1), \\ 
    a(F,C_1\vee (C_2\parallel C_3))&=a(F,C_1\vee C_2)+a(F,C_1\vee C_3).
\end{align}
Consider now a different tile $\overline{C}_1\vee \overline{C}_2$ whose canes have length comparisons $l(C_1,\overline{C}_1)=q_1$ and $l(C_2,\overline{C}_2)=q_2$ with the canes of the tile $C_1\vee C_2$. It is then an experimental fact that, for any flat cobblestone $F\in\mathcal{F}$:
\begin{equation}\label{tile2}
    a(F,C_1\vee C_2)=q_1\cdot q_2\cdot a(F,\overline{C}_1\vee \overline{C}_2),
\end{equation}
which gives the usual identification of areas with products of lengths.\newline

Lastly, we briefly remark on a simple implementation of a theory of density enabled by the fact that the rudimentary weighting scale construction described above allows to give quantitative comparisons between any objects. In particular, if we consider canes $\mathcal{C}$, they can be compared by length or by weight. By simultaneously comparing both characteristics it is possible to construct a new practical procedure that defines a new measurand, what we may call Density, and which is experimentally shown to be given by:
\begin{align*}
d: \mathcal{C}\times \mathcal{C} & \to \Rat\\
(C_1,C_2) & \mapsto \frac{w(C_1,C_2)}{l(C_1,C_2)}
\end{align*}
Such a measurand will be recognized as an entity that measures pairs of canes forming what we call a \textbf{ratio} by weight and length, denoted by $C_1|^w_lC_2$. To make this notion more concrete, we fix a unit of weight $C_0^w$ and a unit of length $C_0^l$, the density of an arbitrary cane $C$ will be specified by the quotient of the comparison in weight $w(C_0^w,C)$ by the comparison in length $l(C_0^l,C)$, then we say that this is the value of the density of $C$ measured in units of $C_0^w$per$C_0^l$. The upshot of the construction of Density is that it provides a natural example of a ``quantitative comparison of quantitative comparisons'', thus defining what we call a \textbf{derived measurand}. These play a fundamental role in practical science, as physical laws are typically formulated in terms of derived measurands. In this example, Density is the common characteristic of canes that allows for a prediction of a comparison of weights with the knowledge of a comparison of length or \emph{vice versa}.\newline

So far we have explored the empirico-formal structure of some rudimentary examples of measurands in a fair amount of detail so as to motivate the purely mathematical definition that will be given below. Before we proceed, however, let us make two brief remarks about the use of rational numbers as the outcome of quantitative comparisons.\newline

Firstly, invoking what we could call the \textbf{Principle of Refinement}, we note that in real physics, particularly in the historical conventions of classical mechanics, one assumes that smaller units always exist so that, as far as the minimum accuracy of a measurement process allows, comparisons can be refined indefinitely. One possible way in which this principle could be implemented abstractly in our formal treatment of Length above is by assuming that the outcome of a length comparison is no longer a single rational value $l(C_1,C_2)\in\Rat$ but, in general, a sequence of rationals $\{l_i(C_1,C_2)\}_{i=1}^\infty$ with
\begin{equation*}
   \lim_{i,j\to\infty} |l_i(C_1,C_2)- l_j(C_1,C_2)|=0.
\end{equation*}
In other words, $\{l_i(C_1,C_2)\}_{i=1}^\infty$ is a Cauchy sequence in $\Rat$. Since from an operational point of view a length comparison needs to yield a single well-defined numerical value, the Cauchy sequence should be forced to be convergent to a single value that we would denote by $l(C_1,C_2)$. This is nothing but the usual construction of $\Real$ by completion of $\Rat$ as a metric space. Quantitative comparison outcomes will henceforth be considered real values and we shall update all the notions introduced so far in this section by replacing any instance of $\Rat$ by $\Real$. It is worth pointing out, however, that this change is only strictly necessary in order to connect with classical mechanics in later sections, where we work within the formalism of differential geometry, which the formal construction of the real numbers built into its core.  Although all the mathematical structure necessary to formalize the definitions of this section is already present in $\Rat$, for the sake of consistency with the rest of the thesis, we will only work over the field of reals $\Real$ from now on.\newline

Secondly, we have not given any physical (i.e. operationally constructive) meaning to the zero in the field of numbers, and therefore neither have we given any meaning to negative values as outcomes of quantitative comparisons. What this is really telling us is that only the affine structure of $\Rat$ (or $\Real$, as we will consider from now on) is relevant for a consistent definition of measurand. Although we could proceed by simply assuming the affine structure, we will keep a ``rooted'' approach, with the zero as a valid possible outcome of a measurement, and only explicitly use the linear structure of $\Rat$. Since all vector spaces are, in particular, affine spaces, this is a convenient choice that will facilitate our discussions involving vector bundles in other chapters.\newline

From the example of measuring canes and stones above we see that the abstract notion of measurand should be given abstractly by a set $\mathcal{M}$, whose elements are identified with the physical objects sharing some property that can be quantitatively compared, and a comparison map $\mathcal{P}:\mathcal{M}\times \mathcal{M}\to \Real$ satisfying the 2-out-of-3 property (\ref{2outof3}). Recall that proposition \ref{RatMaps} ensured that objects in the line category carry ratio maps satisfying a 2-out-of-3 identity, we are thus compelled to give the following definition: a \textbf{measurand} is identified abstractly with a 1-dimensional real vector space $L\in\Line$. Physical objects sharing the same quantitatively comparable quality are identified with non-zero elements of $L$. The \textbf{comparison map} is given canonically with the choice of $L$ by the ratio map $l:L\times L^{\times}\to \Real$, which was shown in (\ref{ratio1}) to satisfy the 2-out-of-3 identity
\begin{equation*}
    l_{ab}\cdot l_{bc}\cdot l_{ca}=1
\end{equation*}
as desired.\newline

Fixing a non-zero element of a measurand gives the notion of a \textbf{unit of measurement} $u\in L$. Since a choice of a non-zero element is indeed equivalent to a choice of a basis, a unit $u\in L$ is equivalently defined as a factor, i.e. an isomorphism of 1-dimensional vector spaces, of the form
\begin{equation*}
    u:L\to \Real.
\end{equation*}

The additive structure of a line $L\in\Line$ identifying a measurand is empirically motivated by the combination construction $\parallel$ introduced above. There we saw how the comparison map, and, consequently, the unit map satisfy a linearity property (\ref{combination}) with respect to the $\parallel$ operation. An obvious way to encode this property is to endow the set of physical objects themselves with a $\Real$-vector space structure, identifying $\parallel$ with usual vector addition $+$, and requiring maps to be linear. This then justifies the use of lines and factors in the category $\Line$ to abstractly represent measurands and the choice of units.\newline

With the proposed abstract definition of measurand in mind, let us revisit the examples of Length and Area. We identify the length and the area measurand with lines $L\in\Line$ and $A\in\Line$ respectively. The tile construction $\vee$ gave an assignment of pairs of lengths to areas satisfying equations (\ref{tile1}) and (\ref{tile2}), a closer look at them reveals that they are simply a bilinearity property of the $\vee$ construction. If one takes elements of the length measurand $b_1,b_2\in L$ representing two canes $C_1$ and $C_2$, respectively, it is a simple exercise to check that making the association $b_1\otimes b_2 \sim C_1 \vee C_2$ and using the definition of the ratio map for the tensor product line leads precisely to equations (\ref{tile1}) and (\ref{tile2}). We thus see that this has effectively created a new measurand out of $L$ by taking its tensor product $L\otimes L$. If the area measurand $A$ was defined independently, as in our example with canes and cobblestones, the fact that tiles of canes are comparable by area to flat cobblestones (ensured by the presence of a standard wedge $\vee$) will be represented by the existence of a canonical factor $w_\vee:L\otimes L\to A$.\newline

In a similar fashion, the notion of derived measurands gives an example of a measurand that is characterising the association of comparisons of two other measurands. Let us consider the Weight, Length and Density measurands as lines $W,L,D\in\Line$, respectively. In the discussion above, we observed that quantitative comparisons of Density correspond to associations of comparisons of Length with comparisons of Weight. It is therefore a simple check to see that setting $D:=\Hom{L,W}\cong L^*\otimes W$ will recover the desired properties of Density.\newline

In summary, we see that our example of a simple theory of Length, Area and Weight of canes and cobblestones can be mathematically formulated in terms of two basic measurands: ``length of a straight object'' $L\in\Line$ and ``weight of an object'' $W\in\Line$. Other measurands of empirical relevance, such as Area $A$ and Density $D$, are given as tensor products of the lines characterizing the basic measurands:
\begin{equation*}
    A=L\otimes L, \qquad D=L^*\otimes W.
\end{equation*}

These considerations lead us to the following general definition: having set a \textbf{domain of observation}, a notion that will not be formalized here but that conceptually corresponds to specifying the class of physical objects that will be subject to scientific enquiry, a choice of quantitatively comparable characteristics between discernible natural phenomena will correspond to a choice of a finite family of lines $\{L_1,\dots, L_k\}$, these will be called the \textbf{base measurands}. All other measurands that may be of scientific relevance will be objects in the \textbf{measurand space} $M$, defined as the potential\footnote{See Section \ref{PotentialFunctor} for the definition and discussion of the potential functor.} of the family of base measurands in the category of lines $M=(L_1\dots L_k)^\odot$.\newline

The measurand space of an arbitrary physical theory is thus shown to have the structure a dimensioned field of real numbers\footnote{See Section \ref{DimRingsModules} for a definition.}. In chapter \ref{DimAlgebraGeometry}, where dimensioned algebra will be systematically studied, the definition of dimensioned ring will be directly motivated by the properties of physical quantities encapsulated in the notion of measurand space described in this section.

\section{A Brief History of Phase Space and Derived Notions} \label{BriefHistoryOfPhaseSpace}

\begin{figure}[h]
\centering
\includegraphics[scale=0.35]{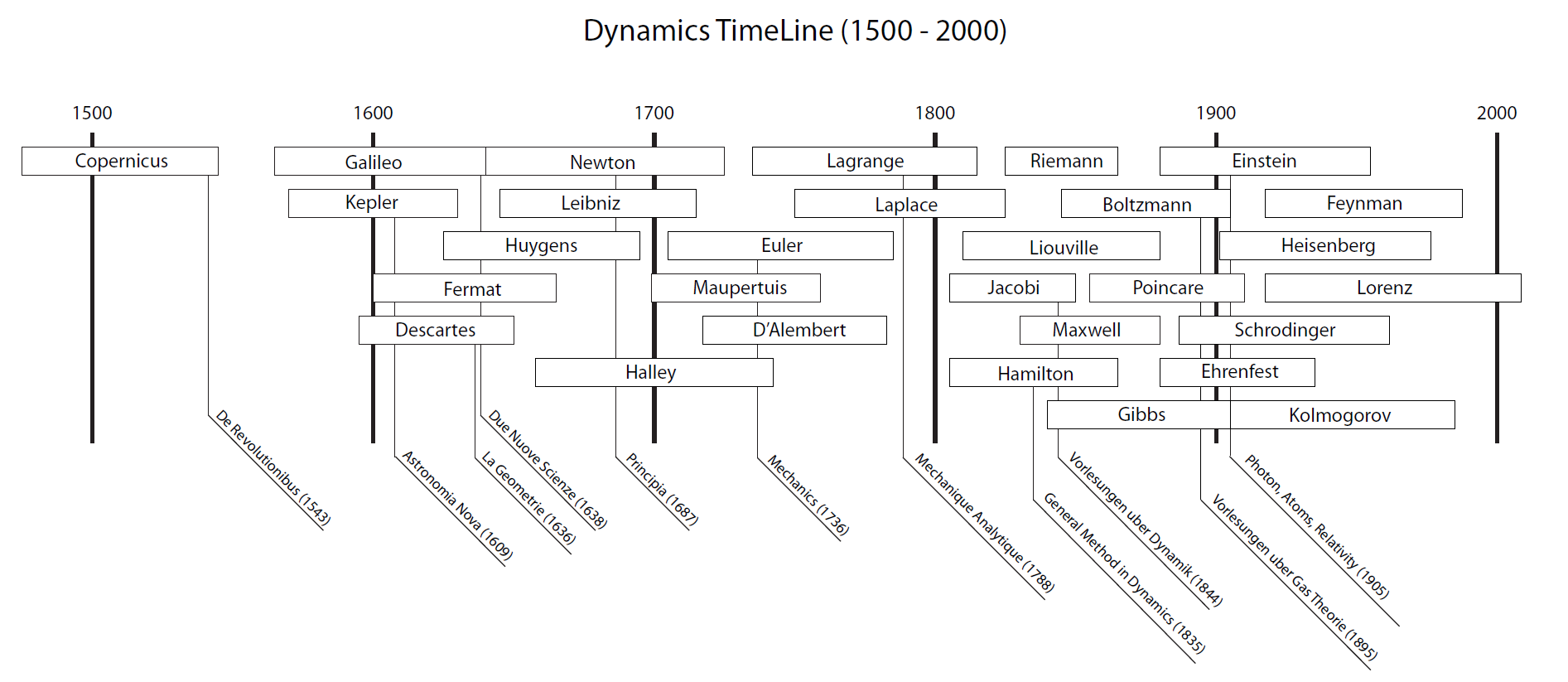}
\caption{Source: \emph{Galileo Unbound: A Path Across Life, the Universe and Everything} \cite{nolte2018galileo} by D. Nolte.}
\label{TimeLineDyn}
\end{figure}

The history of phase space is the history of the abstracted principles of mechanics cast into geometrical form and, as such, it is inextricably tied to the history of modern mathematics and physics. Any comprehensive historical account of the contemporary concept of phase space must necessarily involve a detailed discussion on the origins of the modern mathematical notions of space and dynamics and a genealogy of the multitude of theories that led to our current version of classical mechanics. An excellent reference in this direction is the recently published \emph{Galileo Unbound: A Path Across Life, the Universe and Everything} \cite{nolte2018galileo} by D. Nolte. Figure \ref{TimeLineDyn}, taken from this reference, offers an overview of the milestones and main contributors to the development of dynamics throughout the history of modern science, within which, the idea of phase space originated, developed and, ultimately, crystallized into its modern form. In this section we shall give an account of this course of events that, far from being comprehensive, is intended to draw a clear picture of the evolution of geometric mechanics and the sister mathematical fields of Poisson geometry and the theory of Lie groupoids and Lie algebroids.

The central idea of our modern understanding of phase space is that abstract states of a physical system are represented as single points of a mathematical set and change in the system is represented by a curve being traced on that set. Since the symbolic and mathematical study of static systems is directly tied to the classical discipline of geometry, the first evidence of the use of the idea of an imaginary space where temporal change is represented will come with the advent of modern mechanics.  Building upon a centuries-long tradition of early kinematic theories with roots in ancient Greek science and Babylonian mathematics, Isaac Newton (England, $1642 - 1726$) captured in his laws of motion, formulated in the $1687$ landmark publication \emph{Philosophiae Naturalis Principia Mathematica}, the first clear instance of the concept of kinematic state of a mechanical system. The mathematics necessary to fully articulate these notions precisely did not exist at the time and were indeed being actively developed by Newton himself and Gottfried Wilhelm Leibniz (Germany, $1646 - 1716$). Newton's fluxions and Leibniz's derivatives allowed for a precise mathematical characterization of the intuition of velocity of a moving object which, together with spatial positions, were recognized as the relevant initial data that determines future motion. The invention of Calculus at the end of the $1600$s thus marks the inception of the idea of an abstract set whose elements - positions and velocities - can be conceptually identified with states of a physical system (see \cite{arthur1995newton,grabiner1983changing}).\newline

\begin{figure}[h]
\centering
\includegraphics[scale=0.25]{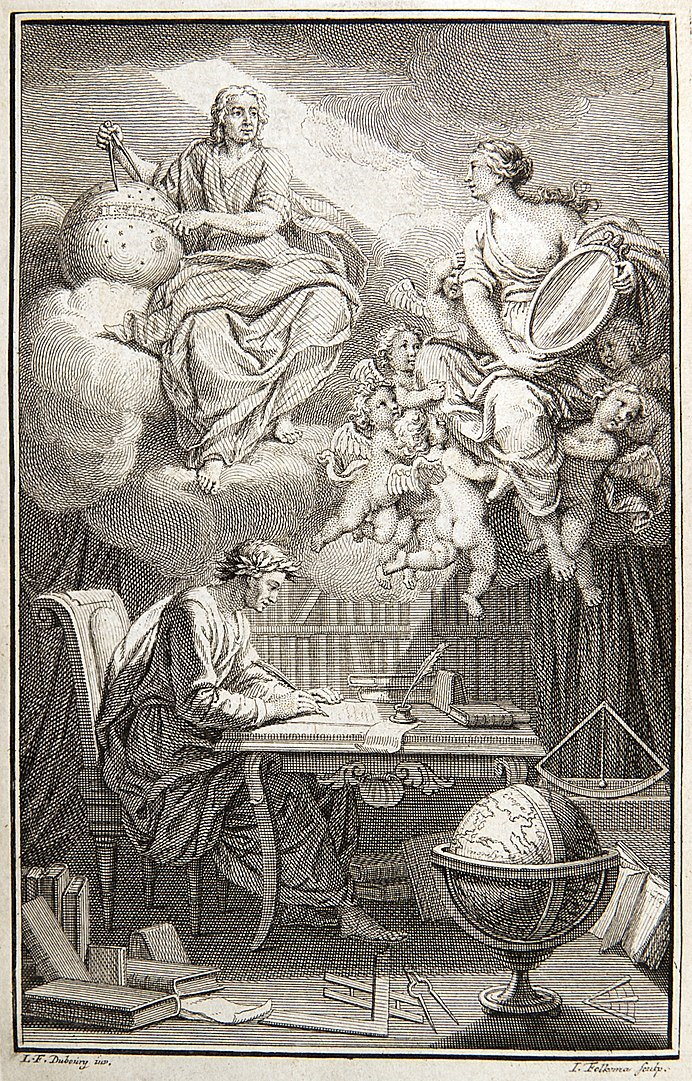}
\caption{Frontispiece of Voltaire's $1738$ \emph{Elemens de la philosophie de Newton}. Newton (top left) is depicted as casting celestial light that Chatelet (top right) reflects down to Voltaire (bottom left). With the hindsight of history, we can aptly replace Voltaire's figure with Lagrange's (as a representative of the French enlightenment school of mechanics), who was greatly influenced by Chatelet's French translation of the \emph{Principia} to make one of the most profound contributions to the history of phase space.}
\label{duChatelet}
\end{figure}

The development of calculus, which, as exemplified above, enabled the transition from kinematics to dynamics and thus kick-started the history of modern mechanics, was largely possible due to the seminal work of Rene Descartes (France, Sweden, $1596 - 1650$) and Pierre de Fermat (France, $1607 - 1665$) on analytic geometry (see \cite{nikulin2002matter}). The original ideas of analytic geometry, by which one could convert geometric problems into algebraic ones, will have profound consequences in the French school of mechanics during the $1700$s. An indispensable figure in this chapter of the history of mechanics is Emilie du Chatelet (France, $1706 - 1749$), whose commentated French translation of Newton's \emph{Principia} made the work a great deal more accessible, critically contributing to its eventual widespread adoption across Europe (see \cite[Chapter 4]{nolte2018galileo}). Chatelet also deserves a special mention in our account since she was one of the pioneers of the idea of kinetic energy and its conservation, thus paving the way for what will later become conservative mechanics. The influence of Chatelet's work in the next generation of French physicists will prove crucial, as we shall discuss below; in a curious historical coincidence, a book published by Voltaire (France, $1694 - 1778$), who was a close acquaintance of hers, had a frontispiece (figure \ref{duChatelet}) depicting what we can regard today as an allegory of Chatelet's role in the history of mechanics.

The French school of mechanics during the age of enlightenment (see \cite{goodman1996republic}) took the scientific insights of Newtonian physics and the mathematical devices of analytic geometry and calculus to essentially elevate mechanics to the subject we recognize today. During this time, not only in France but certainly best exemplified there, the theoretical treatises of the $1600$s and early $1700$s were revised and translated into vernacular languages, making them widely accessible. This meant that theoretical principles, which had been studied only by the intellectual elites, were increasingly understood by a general audience that had more practical goals in mind. Modern engineering, the discipline that uses the best theoretical tools available to solve the most challenging practical problems, was thus initiated, in turn enabling the Industrial Revolution of the $1800$s.\newline

\begin{figure}[h]
\centering
\includegraphics[scale=0.25]{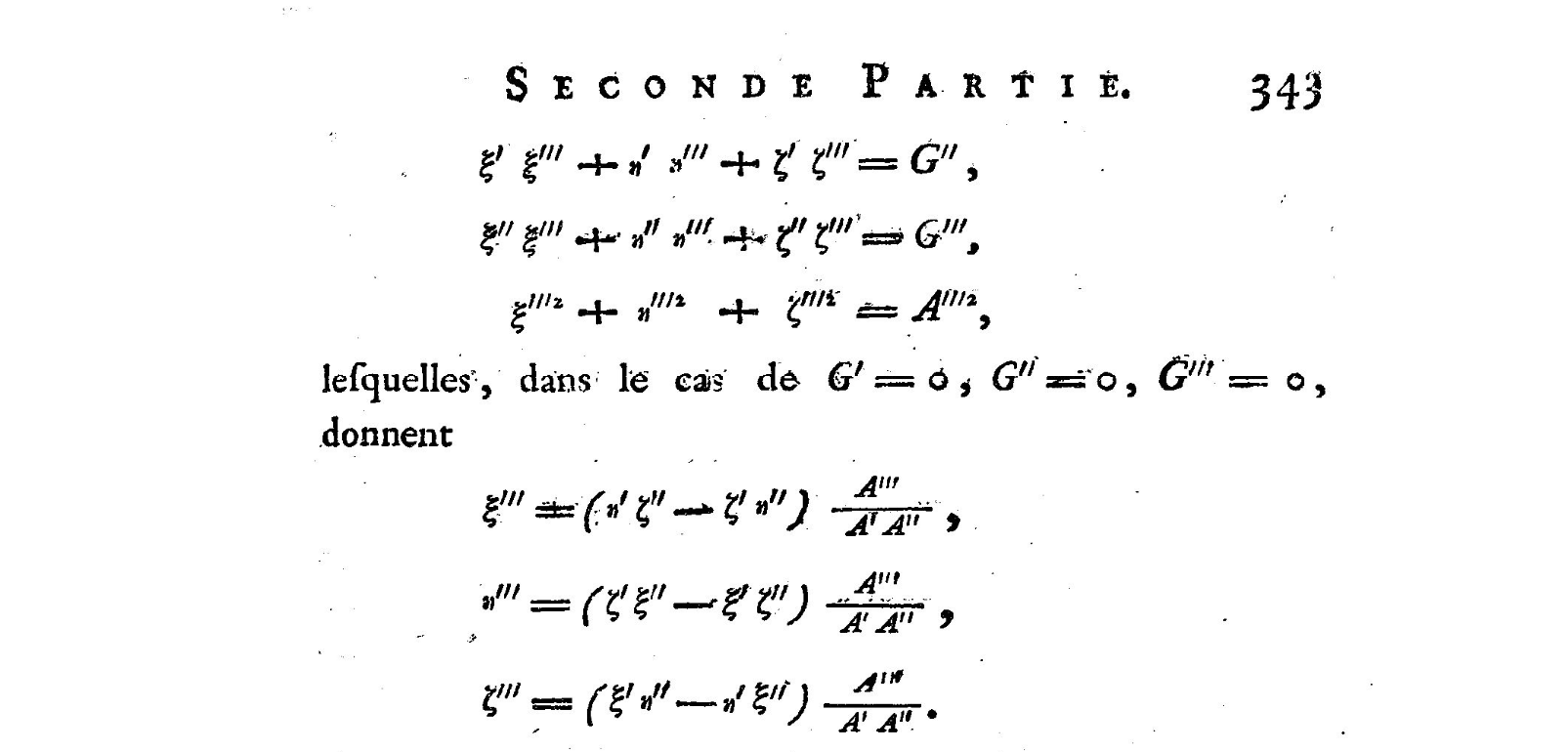}
\caption{An excerpt from \emph{Section VI: Sur la Rotation des Corps} in Lagrange's \emph{Mechanique Analytique} of $1789$. We clearly recognize the usual commutation relations of the Lie algebra of 3-dimensional rotations in the second set of equations.}
\label{LagrangeSO(3)}
\end{figure}

It is in the context of the French enlightenment that we find the first traces of the modern notion of phase space. Joseph-Louis Lagrange (Italy, Germany, France, $1736 - 1813$) was a self-proclaimed champion of the idea that mechanics should be reduced to algebra as he regarded his general methods for finding equations of motion in his \emph{Mechanique Analytique} of $1789$, what we know today as Lagrangian mechanics, as having removed geometry from physics - he is quoted as being proud of the fact that his landmark publication contained no diagrams, only equations (see \cite{pepe2014lagrange}). This is clearly exemplified by the introduction of what we know today as generalized coordinates on configuration space, which indeed bear little connection to any explicit geometric structure. We thus see how Lagrange's work deliberately impelled mechanics into ever increasing levels of abstraction. Aside from the monumental achievement that the \emph{Mechanique Analytique} represents in the history of mechanics as a pinnacle of a century of scientific exchanges and theoretical synthesis since Newton's \emph{Principia}, Lagrange is also credited with having made the first explicit use of two key mathematical structures of modern mechanics: symplectic structures, identified as brackets of certain constants of motion in his paper \emph{Memoire sur la theorie generale de la variation des constantes} of $1809$, and Lie algebras, found in his studies of rotating bodies (see figure \ref{LagrangeSO(3)}).

\begin{figure}[h]
\centering
\includegraphics[scale=0.29]{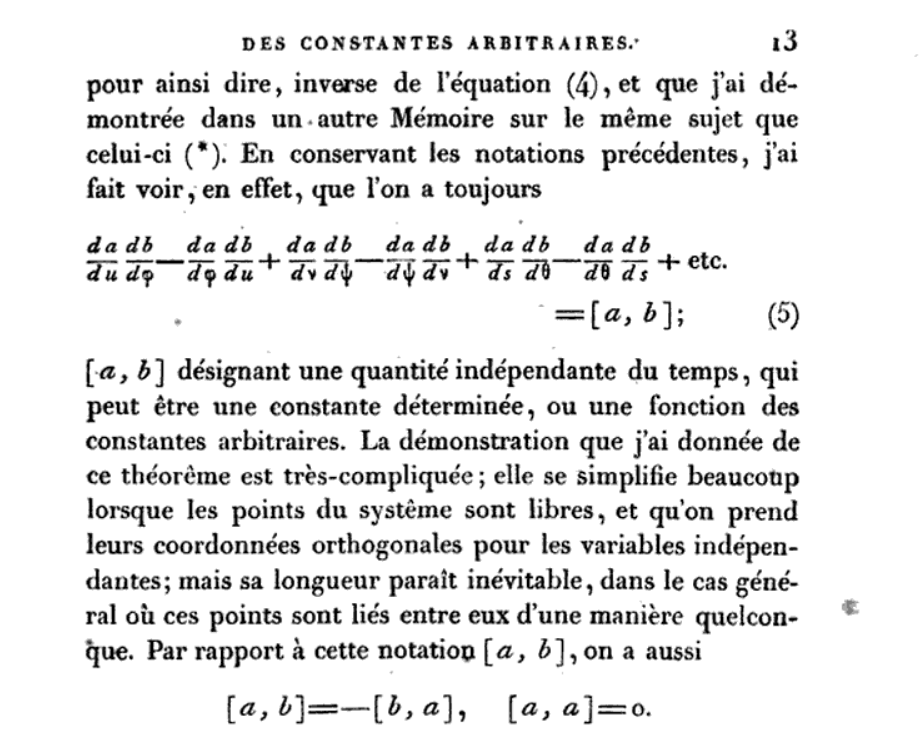}
\caption{An excerpt from Poisson's \emph{Sur la variation des constantes arbitraires dans les questions de mecanique} of $1809$. This entry matches our modern coordinate definition of Poisson bracket on phase space (top equations) precisely, furthermore, Poisson immediately notes (bottom equations) the key antisymmetry property of the bracket that becomes central in the modern treatment of conservation laws.}
\label{PoissonBracc}
\end{figure}

In the same year of $1809$, Simeon Denis Poisson (France, $1781 - 1840$), who had been Lagrange's student and with whom enjoyed a life-long friendship, published a paper aiming to improve on Lagrange's brackets entitled \emph{Sur la variation des constantes arbitraires dans les questions de mecanique} where the first instance of what we call today Poisson bracket explicitly appears (see figure \ref{PoissonBracc}). Both authors pointed out the skew-symmetry of their brackets and were aware of the many simplifications that occurred when using them to determine constants of motion, however, it was Carl Gustav Jacob Jacobi (Germany, $1804 - 1851$) who, a few years later, noticed the property that bears his name today.

\begin{figure}[h]
\centering
\includegraphics[scale=0.3]{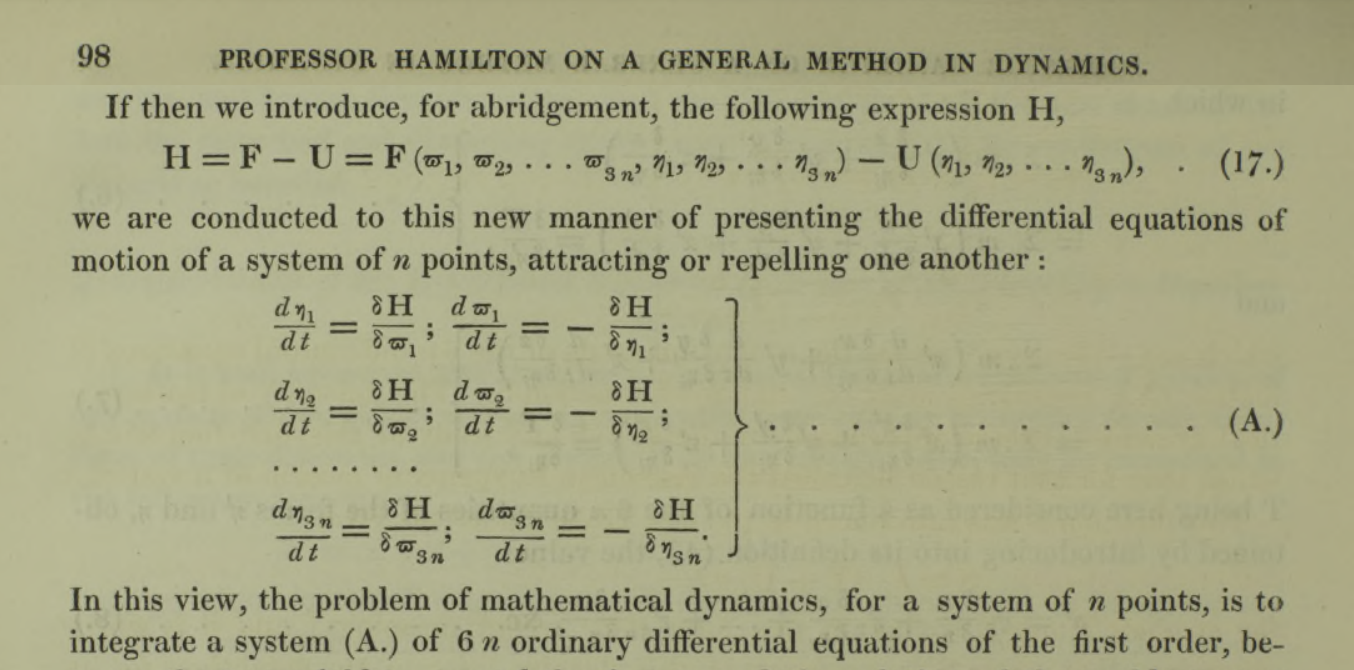}
\caption{An excerpt from Hamiltons's \emph{Second essay on a general method in dynamics} \cite{hamilton1835vii} of $1835$. Here we find Hamilton's equations of motion in its modern form for the conjugate coordinates $(\eta_i,\varpi_i)$.}
\label{HamiltonEq}
\end{figure}

Extending and refining the analytical mechanics of Lagrange and Poisson, the work of William Rowan Hamilton (Ireland, $1805 - 1865$) pushed the abstract understanding of mechanics closer to its modern embodiment in two different ways: on the one hand, his reformulation of Euler-Lagrange equations via conjugate momenta (see figure \ref{HamiltonEq}) constitutes the first example of a set of ordinary differential equations giving the coordinate expression of what is called today a Hamiltonian vector field; and on the other, the method of characteristic functions and the development of the formalism given the name of Hamilton-Jacobi today, makes extensive use of functions dependent on position and momenta variables thus conferring an even more tangible status to the abstract set of position and momenta as the space representing the kinematic states of a mechanical system.

By $1850$ it is clear that the parameter space of position and momenta variables is gaining increasing conceptual weight, however, the mathematical treatment of mechanics was still firmly rooted in Lagrange's view as it was mostly concerned with the algebraic manipulation of equations of motion in search for analytic solutions and constants of motion. The missing piece of mathematical technology required for a precise understanding of the set of position and momenta variables as a space was intrinsic geometry and the theory of differentiable manifolds. These came, motivated by a long history of problems in curved surfaces and geodesy, by the hand of Carl Friedrich Gauss (Germany, $1777 - 1855$) and his doctoral student Bernhard Riemann (Germany, Italy, $1826 - 1866$). In his lecture of $1854$ \emph{Uber die Hypothesen, welche der Geometrie zu Grunde liegen} delivered in $1854$ at the University of Gottingen, Riemann defined the notion of manifold as an abstract space that could be locally described by an indefinite number of variables. This lecture essentially marks the beginning of modern differential geometry and it will be from this point that the idea of space of arbitrary dimensions, in particular, greater than 3, will gradually find its way into the mainstream scientific vernacular.\newline

During the second half of the $1800$s, while the German school of early differential geometers laid the foundations of most of the notions that are familiar to us today, e.g. groups of transformations, non-euclidean geometries, homogeneous spaces, geodesics etc., significant advances in the kinetic interpretation of thermodynamics were being made by Rudolf Clausius (Germany, $1822 - 1888$) and James Clerk Maxwell (Scotland, England, $1831 - 1879$). This interpretations involved thinking of phenomena such as heat dissipation or gas expansion as the average effect of the complicated motion of vast numbers of small constituents of matter, indeed exploiting the recently developed ideas of molecular theory (see \cite{rosenfeld1971men}). Ludwig Boltzmann (Austria, Italy, $1844 - 1906$) devoted his scientific career to the task of deriving the laws of thermodynamics, which had been known in some capacity for more than a century before him, entirely from the first principles of kinetic theory, establishing, in the process, the field that we know today as statistical mechanics. The great technical challenge of this endeavour lies in the fantastically large amount of coupled differential equations one would need to solve to exactly calculate the motion of the gas molecules (specially in a time when computers were far from being a reality), hence approximate methods or simplifying assumptions were required. The conceptual breakthrough that allowed Boltzmann to make his most significant contributions was the idea that one should consider regions of the space of positions and momenta to represent statistical aggregates of kinematic states of the molecules of a gas, as was clearly formulated in his \emph{Lectures on Gas Theory} \cite{boltzmann2012lectures} of $1896$. In order to illustrate his ergodic hypothesis, Boltzmann used an analogy between the kinematic state of a molecule and a Lissajous figure, resulting from the plotting of two phases, in the sense of angle of rotation, on a plane. Borrowing from this analogy, the term `phase point', used to designate the state of a gas at a given instant in time, is found throughout Boltzmann's writings; it is here that the term `phase', referring to an abstract state of a physical system, first enters the language of mechanics.\newline

Although Boltzmann himself did not embrace the geometric nature of his new idea, it was after the notion of intrinsic geometric spaces of arbitrary dimension had been firmly mathematically established and physicists had begun adopting it, due, in part, to the great success of the papers of Albert Einstein (Germany, United States $1879 - 1955$) on relativity theory, that, when Paul Ehrenfest (Austria, Netherlands, $1880 - 1933$) and Paul Hertz (Germany, United States, $1881 - 1940$) collaborated in writing an encyclopedia entry about Boltzmann's `phase point' in $1910$, the term `phase space' was explicitly used to refer to the set of all possible `phase points' and the interpretation of the time evolution of the system as a single trajectory in this `phase space' is directly alluded to. This constitutes the first explicit appearance of the modern notion of phase space, both in conception as well as in name, thus reuniting the science of motion with its geometric origins and setting the stage for the development of contemporary geometric mechanics during the $1900$s.\newline

Upon seeing how well suited Riemann's intrinsic geometry was for the description of mechanics, mathematicians were quick to tackle many long-standing problems that were then made approachable from different angles by the new abstractions that followed from the introduction of manifold theory into analytical mechanics. This is best exemplified by the work of Henri Poincare (France, $1854 - 1912$) on the three-body problem around $1882$, in which the idea of the evolution of the system tracing a single trajectory in phase space proved crucial for a qualitative treatment of what he will later show to be one of the first examples of a chaotic system. Perhaps the single most influential contribution to the eventual adoption of the formalism of differentiable manifolds as a core component of modern theoretical physics was Einstein's general relativity, a theoretical framework first encapsulated in his seminal paper \emph{Die Feldgleichungen der Gravitation} \cite{einstein1915feldgleichungen} of $1915$, which paved the way for some of the most profound scientific breakthroughs in the history of science to this day.\newline

By the time of the widespread adoption of general relativity and of the increase in popularity of manifold theory among physicists in the early $1900$s, mathematicians had been elevating the initial ideas of differential geometry, which were closely linked to analytical mechanics in their origins, to ever more abstract heights. The work of Sophus Lie (Norway, $1842 - 1899$) on groups of continuous transformations and the integrability of distributions on manifolds, which established all the notions that bear his name familiar to us today, is of particular relevance to the present mathematical embodiment of phase spaces. Not only did Lie groups and Lie algebras play a major role in the study and implementation of symmetry in modern physics, but his writings from the $1890$s also contain the first explicit formal uses of the mathematical objects that we recognize today as symplectic manifolds, Poisson structures, moment maps or Lie algebroids (see \cite{fritzsche1999sophus,weinstein1983sophus}).\newline

Although clearly originating in early analytical mechanics and the Erlangen program, the history of the modern notions mentioned in the previous paragraph throughout the $1900$s is a complicated one since they were independently discovered in different contexts with different flavours by mathematicians and physicists alike. The term `symplectic', which stems from the Greek word for `complex', was first used by Hermann Weyl (Germany, Switzerland, $1885 - 1955$) as a reference to the similarity between symplectic structures and complex structures. It wasn't until the $1970$s that Andre Lichnerowicz (France, $1915 - 1998$) isolated the axioms for our modern definition of Poisson manifold, also coining the term in so doing (see \cite{flato1976deformations}). The first general definition of groupoid was given by Heinrich Brandt (Germany, $1886 - 1954$) in $1926$ with Charles Ehresmann (Germany, France, $1905 - 1979$) adding smooth structures on them soon after and thus marking the beginning of the use of what was later known as Lie groupoids as a tool in differential topology (see \cite{weinstein1996groupoids}). The work of Michael Atiyah (England, United States, Scotland, $1929 - 2019$) on principal bundles (see \cite{atiyah1988michael}) in the $1950$s contains, perhaps ironically due to his personal dislike of the term `algebroid', one of the first uses of Lie algebroids in a modern sense. These were defined and identified as the infinitesimal counterparts of Lie groupoids by Jean Pradines (France, $-$), a doctoral student of Ehresmann, in the $1960$s. Around this time, Bertram Kostant (United States, $1928 - 2017$) and Jean-Marie Souriau (France, $1922 - 2012$) independently discovered the concept of moment map in symplectic manifolds, with the latter author coining its current name (in fact, the word `moment' has been kept from his original terminology in French for `momentum mapping') and studying their properties extensively (see \cite{souriau1970structure}).\newline

By the late $1970$s, the aggregated works of several mathematicians and physicists, perhaps best exemplified by the extensive contributions of Jerrold Eldon Marsden (United States, $1942 - 2010$) and Vladimir Igorevich Arnold (Russia, $1937 - 2010$), led to the establishment of what we know today as geometric mechanics (see \cite{abraham1978foundations,arnold2013mathematical}). The birth of this discipline, which has since been developed by mathematical physicists, geometers and engineers alike, represents a return to the physical origins of all the geometric abstractions that had been developed since the times of Hamilton. Geometric mechanics has proven instrumental for both the theoretical understanding of fundamental physics, with the Hamiltonian formalism being built into the early formulations of quantum mechanics by Paul Adrien Maurice Dirac (England, United States, $1902 - 1984$), and for powerful applications in control theory and mechanical engineering.\newline

Despite Dirac's scientific career focusing mostly on quantum mechanics and particle physics, it was his novel approach to constrained Hamiltonian systems, clearly presented in his book \emph{Lectures on Quantum Mechanics} \cite{dirac1966lectures} of $1966$, what would shape the fields of Poisson and Symplectic geometry in the following decades. In the later decades of the $1900$s, the incredible experimental accuracy of predictions from quantum field theory together with the fact that a solid mathematical foundation for this formalism was missing, launched a golden era of mathematical innovation within the theoretical physics community that resulted in the birth of entire industries of mathematical physics such as string theory, M-theory or holography, that were pushing the boundaries of what was considered as a physical phase space into ever more exotic geometric realms.\newline

The modern generalizations of the notion of phase space that are suggested by the developments in mathematical physics of the second half of the $1900$s mentioned above find a natural mathematical home in the fields of Dirac geometry, Lie Groupoid and Lie algebroid theory and non-commutative geometry. Although still a vividly active area of research today with several hundreds of mathematicians and physicists around the world working on it, there are three key moments in the recent history of Dirac geometry that help us understand its current state. Firstly, motivated by the deformations of Poisson algebras and their connection with non-commutative geometry, there was the realization by several mathematicians that Lichnerowicz's Poisson manifolds and Pradines' Lie algebroids, which had been developed independently, were intimately related. Alan David Weinstein (United States, $-$), in addition to his momentous contributions to the theory of reduction in geometric mechanics, played a key role in the early days of Dirac geometry when, together with collaborators, he showed in $1987$ that the base manifold of a symplectic groupoid naturally inherits a Poisson structure from a Lie algebroid bracket on its cotangent bundle. The second important event was the discovery, independently by Theodore James Courant (United States, $-$) and Irina Jakovlevna Dorfman (Russia, $1948-1994$), of a bracket on the sum of the tangent and cotangent bundle. Courant's thesis \emph{Dirac Manifolds} of $1990$ provides a precise mathematical framework for the theory of constraints proposed by Dirac some $40$ years earlier, from which the name originates, and lays the foundations for the research on this area to this day (see \cite{kosmann2013courant}). Lastly, we find the effort to extend Lie's theorems of integrability to Lie algebroids culminating in the work of Rui Loja Fernandes (Portugal, United States, $-$) and Marius Crainic (Netherlands, $-$) as collected in their paper \emph{Integrability of Lie Brackets} \cite{crainic2003integrability} of $2003$.\newline

As a closing remark linking back to our account of the history of metrology of Section \ref{BriefHistoryOfMetrology}, we note that by the time the core principles of dimensional analysis were identified by Fourier in $1822$, the mathematical theory of analytical mechanics was well underway. Seeing how the topic of dimensional analysis was only seriously considered by the physics community in the late $1800$s, it is no surprise that contemporary geometric mechanics shows no explicit influence of the formal treatment of physical dimension.

\begin{figure}
    \centering
    \begin{timeline}{1687}{1915}{2cm}{2.5cm}{13cm}{20cm}
        \entry{1687}{Newton's \emph{Philosophiae Naturalis Principia Mathematica}}
        \entry{1695}{Mathematical notion of kinematic states via derivatives in early calculus}
        \entry{1736}{Euler's \emph{Mechanica}}
        \entry{1749}{Chatelet's commentated translation of the \emph{Principia}}
        \entry{1770}{Early uses of conservation of mechanical energy}
        \entry{1780}{First explicit occurrence of a Lie algebra in Lagrange's work on rotating bodies}
        \entry{1789}{Lagrange's \emph{Mechanique Analytique}}
        \entry{1809}{First occurrences of symplectic structures in the works of Lagrange and Poisson}
        \entry{1822}{Fourier's \emph{Theorie Analytique de la Chaleur}, beginning of dimensional analysis}
        \entry{1825}{Jacobi notes the identity of Poisson brackets that bears his name}
        \entry{1835}{Hamilton's equations appear in his \emph{General Methods in Dynamics}}
        \entry{1850}{Clausius postulates the first two laws of thermodynamics}
        \entry{1854}{Riemann defines manifold, beginning of modern differential geometry}
        \entry{1865}{Boltzmann uses the term `phase point' to refer to the state of a gas}
        \entry{1877}{Rayleigh's \emph{Theory of Sound}, further results in dimensional analysis}
        \entry{1882}{Poincare's tackles the three-body problem}
        \entry{1888}{Lie's work contains modern versions of Poisson manifolds and moment maps}
        \entry{1892}{Vaschy gives a first proof of the fundamental theorem of dimensional analysis}
        \entry{1896}{Botzmann's \emph{Lectures on Gas Theory} establish statistical mechanics}
        \entry{1910}{First explicit use of the term phase space by Ehrenfest}
        \entry{1915}{Widespread adoption of dimensional analysis}
        \entry{1915}{Einstein's \emph{Die Feldgleichungen der Gravitation}}
    \end{timeline}
\label{PhaseSpaceTimeline}
\end{figure}

\section{The Hamiltonian Phase Space Formalism} \label{FoundationsHamiltonianPhaseSpace}

As we saw in Section \ref{BriefHistoryOfPhaseSpace}, there is a long and rich history of mathematical ideas and scientific breakthroughs that has led to the modern formulations of classical mechanics. This has the downside of making it quite hard to find a single source where the fundamental theoretical principles are clearly laid out and are expressed using contemporary mathematical language. This section is aimed at providing such a source.\newline

\begin{figure}[h]
\centering
\includegraphics[scale=0.4]{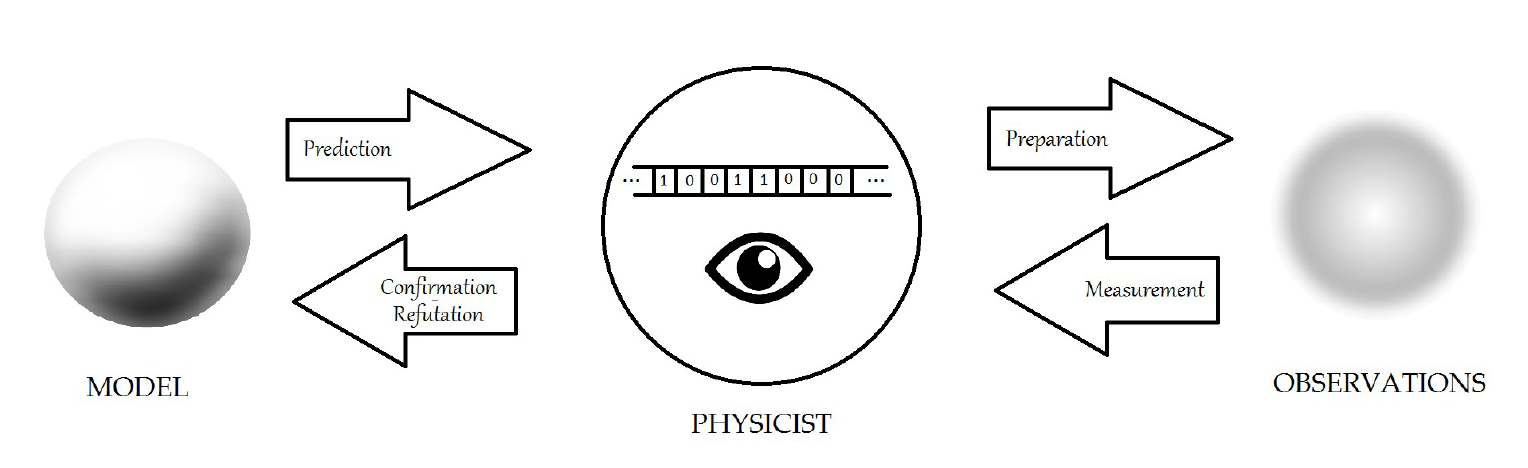}
\caption{A general sketch of the scientific activity.}
\label{ModPhysObs}
\end{figure}

In order to fix some ideas before we begin our discussion of the physical principles underlying Hamiltonian mechanics, let us consider a general picture of a working physicist who is studying a temporal series of measurements, represented by the tape of binary values in figure \ref{ModPhysObs}. In that diagram, metrological operations, e.g the practical measurement with a flowmeter of the volume of water traversing a pipe, are summarized as the right arrows and the mathematical formulation of predictions and the notion of compatibility within a theoretical framework, e.g. the mathematical definition of flow and Bernoulli's principle, are summarized as the left arrows. Although we will be dealing mostly with the ``Model'' side in our discussion below, it will be helpful to have this general picture in mind throughout the reminder of this section.\newline

In a practical experimental situation, the temporal series of measurements will be a discrete sequence of rational numbers $\{m_i\in\Rat\}$, each value being the outcome of a quantitative comparison as described in Section \ref{MeasurandFormalism}, e.g. positions of a particle at different times. Following the Principle of Refinement of Section \ref{MeasurandFormalism}, the discrete temporal series of measurements will be replaced by a smooth dependence of real values $m(t)\in\Real$ with a real \textbf{time} parameter $t\in I\subset \Real$. An orientation of the observer's time interval $I$ is chosen so that the \textbf{past} $\{t<t_0\}$ and \textbf{future} $\{t>t_0\}$ of a given \textbf{present} $t_0\in I$ are defined. The task of a theory of mechanics is to find a mathematical model that takes as inputs the mathematical parameters corresponding to a given experimental preparation and produces as output a temporal series $m(t)$. What follows is a list of the main physical intuitions behind the formal content of a dynamical theory in Hamiltonian phase spaces. Note that the Principle of Refinement will be implicitly invoked every time we use smooth manifolds or smooth maps to represent physical notions.
\begin{itemize}
    \item \textbf{Principle of Realism.} A physical system exists independent of the observer studying it. Discernible configurations of the system, sets of equivalent experimental preparations and outcomes, are called \textbf{states} $s$ and are identified with points in a smooth manifold $P$, called the \textbf{phase space}.
    \item \textbf{Principle of Characterization.} Properties or characteristics of a system are smooth assignments of measurement outcome values to each point of the phase space. In the case of conventional classical mechanics this gives the usual definition of \textbf{observable} as a real-valued function $f\in \Obs{P}:=\Cin{P}$. Note that observables $\Obs{P}$ form a ring with the usual operations of point-wise addition and multiplication.
    \item \textbf{Principle of Kinematics.} The observer studying the system exists simultaneously with the system. In the same way that the observer's memory state is mapped uniquely into the time interval used to array the experimental temporal series, the physical system is thought to be in a single state corresponding to each of the observer's time parameter values. A \textbf{motion} is defined as a smooth curve in phase space $c:I\to P$ parameterized by the observer's time $t\in I\subset \Real$. Phase spaces are assumed to be path-connected so that any state is  connected to any other state by a motion, at least virtually, not necessarily physically. A family of motions $\{c(t)\}$ is called an \textbf{evolution} $\mathcal{E}$ on the phase space $P$ if all the states are included in the path of some curve. More precisely, we call $\mathcal{E}=\{c(t)\}$ an evolution if
    \begin{equation*}
        \forall s_0\in P \quad \exists c(t)\in \mathcal{E} \quad | \quad c(t_0)=s_0 \text{ for some } t_0\in I.
    \end{equation*}
    \item \textbf{Principle of Observation.} The observed time series of measurement outcomes are the result of observables taking values\footnote{We will not deal with non-trivial implementations of measurement outcomes, such as the collapse mechanism in the usual Hilbert space formulation of quantum mechanics, but we note that the principles of the Hamiltonian phase space formalism described in this section, excluding the Principle of Observation, generally apply to quantum mechanics, as will be discussed in Section \ref{QuantumPhaseSpace}} along a particular motion. More concretely, given a motion $c(t)$ and an observable $f\in\Obs{P}$, the predicted temporal series is simply given by
    \begin{equation*}
        m(t)=(f\circ c)(t).
    \end{equation*}
    \item \textbf{Principle of Reproducibility.} Similar experimental preparations of a physical system should give similar observational outcomes\footnote{By ``observational outcome'' here we mean aggregates of experimental results on which statistical analysis is necessary. The usual classical and quantum measurement paradigms fit into this description.}. An evolution $\mathcal{E}$ implementing this principle will satisfy the following property for all pairs of motions $c(t),c'(t)\in\mathcal{E}$:
    \begin{equation*}
        t_0\in I,\quad c(t_0)=c'(t_0) \quad \Rightarrow \quad c(t)=c'(t) \quad \forall t\in I,
    \end{equation*}
    in other words, $\mathcal{E}$ is a family of non-intersecting curves parameterized by $t\in I$.
    \item \textbf{Principle of Dynamics.} Future states of a physical system are completely determined\footnote{We are careful not to call this determinism since it is the abstract state of a system, not the measurement outcomes, what are assumed to evolve deterministically. Even with measurement paradigms such as a collapse mechanism that forces the definition of a pre-measurement state and post-measurement state, if the theory relies on ordinary differential equations for the modelling of time evolution, the principle of dynamics will be implicitly used, at least to some degree.} by any given present state, at least locally. Enforcing this condition on motions for arbitrary small time intervals leads to the condition that a motion must be uniquely specified by the values of its tangent everywhere. In other words, a motion satisfying the Principle of Dynamics must be an integral curve of some smooth vector field, thus justifying the chosen name for this principle. Indeed, existence and uniqueness of ODEs implies that a smooth vector field $X\in\Sec{\Tan P}$ gives (at least locally) an evolution defined by the family of its integral curves $\mathcal{E}^X$ automatically satisfying the Principle of Reproducibility and the Principle of Dynamics. Arguing from a different angle, conjunction of the Principle of Observation with the Principle of Dynamics leads, via the local notion of directional derivative, to predicted temporal series being the integration of some derivation on observables. We thus conclude that an evolution on a phase space $P$ satisfying all the principles stated above is given by a choice of \textbf{dynamics} $X\in\Dyn{P}$, which is equivalently understood as a vector field on the phase space or a derivation of observables
    \begin{equation*}
        \Sec{\Tan P}\cong : \Dyn{P} :\cong \Dr{\Cin{P}}.
    \end{equation*}
    \item \textbf{Principle of Conservation.} Conserved quantities are a fundamental building block of experimental mechanics: one can only study time-dependent phenomena effectively when enough variables can be assumed to be constants to the effects of the experiment at hand. We could promote this to the more general and abstract requirement that any observable $f\in\Obs{P}$ has an associated evolution $\mathcal{E}^f$, given by some vector field $X_f$, along which the predicted time series are constant. We call such an assignment a \textbf{Hamiltonian map}
    \begin{equation*}
        \eta : \Obs{P}\to \Dyn{P},
    \end{equation*}
    which is required to satisfy
    \begin{equation*}
        \eta(f)[f]=0 \quad \forall f\in\Obs{P}.
    \end{equation*}
    In conventional Hamiltonian mechanics this is assumed to be given by the slightly stronger structure of a Poisson bracket on observables $(\Obs{P},\{,\})$, thus making the phase space into a Poisson manifold $(P,\pi)$. The Hamiltonian map is then $\eta:=\pi^\sharp \circ d$ and, following from the Jacobi identity of the Poisson bracket, it is a Lie algebra morphism.
    \item \textbf{Principle of Reductionism.} The theoretical description of a physical system specified as a subsystem of a larger system must be completely determined by the theoretical description of the larger system and the information of how the smaller system sits inside. This principle is implemented by demanding that knowledge about submanifolds in a phase space allows us to construct new phase spaces. In conventional Hamiltonian mechanics this corresponds the problem of reduction in Poisson manifolds.
    \item \textbf{Principle of Combination.} The theoretical description of a system formed as a combination of other two systems must be completely determined by the theoretical descriptions of each of the parts and the information of how they interact. This is implemented by demanding that there is a \textbf{combination product} construction for phase spaces
    \begin{equation*}
        \between :  (P_1,P_2)\mapsto P_1 \between P_2.
    \end{equation*}
    In the case of conventional Hamiltonian mechanics, this is simply given by the usual product of Poisson manifolds.
    \item \textbf{Principle of Symmetry.} A theoretical description of a system containing states that are physically indistinguishable should contain all the information to form a faithful theoretical description of the system. Physically indistinguishable states are commonly regarded to be orbits of some Lie group action on the phase space, thus an implementation of this principle will require that from the information of a Lie group action preserving some existing structure on a phase space, a new phase space is constructed whose states are classes of physically-indistinguishable states. In Hamiltonian mechanics this is implemented via the theory of Poisson group actions and equivariant moment maps.
\end{itemize}

Blending all these principles together and casting them into categorical form, we introduce the general notion of a \textbf{theory of phase spaces} consisting of the following categorical data: a \textbf{category of phase spaces} $\textsf{Phase}$, which can be identified with some subcategory of smooth manifolds $\textsf{Phase}\subset\Man$, carrying natural notions of subobjects and quotients, thus also inducing a notion of reduction, and a categorical product; a \textbf{category of observables} $\textsf{Obs}$ whose objects carry local\footnote{Here local is used in the sense introduced in the beginning of Section \ref{DerivativeAlgebras}.} algebraic structures reflecting the measurement paradigm in which the theory will fit; a \textbf{category of dynamics} $\textsf{Dyn}$ whose objects carry local algebraic structures reflecting the time evolution of the quantities to be measured; and three functors fitting in the following commutative diagram:
\begin{equation*}
\begin{tikzcd}[row sep=small]
& \textsf{Obs} \arrow[dd,"\text{Evl}"] \\
\textsf{Phase} \arrow[ur,"\text{Obs}"] \arrow[dr,"\text{Dyn}"'] & \\
& \textsf{Dyn} 
\end{tikzcd}
\end{equation*}
where the \textbf{observable functor} $\text{Obs}$ represents the assignment of measurable properties to a given system, the \textbf{dynamics functor} $\text{Dyn}$ represents the correspondence between motions and smooth curves on phase spaces and the \textbf{evolution functor} $\text{Evl}$ represents how motions should induce the time evolution of observables infinitesimally. A theory of phase spaces is called \textbf{Hamiltonian} if for all phase spaces $P\in\textsf{Phase}$ there exists a canonical choice of Hamiltonian map
\begin{equation*}
    \eta_P: \Obs{P}\to \Dyn{P}
\end{equation*}
which is compatible with the algebraic structures and that captures the notion of conservative evolution. We could, at this point, discuss in great detail the physical significance of the interaction of the categorical structure on $\textsf{Phase}$ and the observable and dynamics functors, however we postpone such matters to the specific cases of phase space theories that are studied in detail in the chapters to follow.

\chapter{Symplectic Phase Spaces} \label{SymplecticPhaseSpaces}

In this chapter we begin by presenting standard classical mechanics on cotangent bundles with an emphasis on its categorical features so as to identify the first example of a theory of phase spaces in the sense of Section \ref{FoundationsHamiltonianPhaseSpace}. The main traits of this formulation will be summarized in the notion of a Hamiltonian functor, which will be latter used as a template for several generalizations of classical mechanics. In light of the technical results of sections \ref{SymplecticGeometry} and \ref{DiracGeometry}, the generalization from symplectic manifolds to Poisson, presymplectic or Dirac structures will be discussed in relation to the Hamiltonian functor. Lastly, a series of results are proved that seem to give a distinct Poisson flavour to all the natural candidates for generalizations of symplectic phase spaces.

\section{Canonical Symplectic Phase Spaces and the Cotangent Functor} \label{CanonicalSymplectic}

Classical mechanics is erected upon the fundamental notion of \textbf{configuration space}: the set of static states of a physical system, such as the possible spatial positions of moving particles or the possible shapes of a vibrating membrane. Invoking the Principle of Refinement formulated in Section \ref{MeasurandFormalism} by which mathematical objects representing physical entities are assumed to be continuous and smooth, our definition of a configuration space will be simply that of a smooth manifold whose points $q\in Q$ are identified with the different static states of a given physical system. We then propose the definition of \textbf{the category of configuration spaces} simply to be the category of smooth real manifolds $\Man$. Under the understanding that smooth manifolds represent sets of static states of a physical system, we now give physical interpretation to the natural categorical structure present in $\Man$:
\begin{itemize}
    \item The measurable static properties of a physical system with configuration space $Q\in\Man$, what we call \textbf{static observables}, are simply the smooth real-valued functions $\Obs{Q}:=\Cin{Q}$. More generally, we see that the assignment of observables to configuration spaces is given by a contravariant functor $\text{C}^\infty:\Man\to\Ring$, which we regard as the categorical version of the Principle or Characterization for configuration spaces.
    \item A subsystem is characterised by restricting possible positions of a larger system, that is, by an inclusion of an embedded submanifold $i:S\to Q$. We then see that monomorphisms in the category $\Man$ implement the Principle of Reductionism for configuration spaces.
    \item Physically-indistinguishable static states in a configuration space $Q\in\Man$ are related by equivalence relations $\sim$ that have quotients faithfully characterizing the physical system, that is, there is a surjective submersion $p:Q\to Q/\sim$. In particular, a free and proper action of a Lie group $G\Acts Q$ gives an example implementing the Principle of Symmetry via $p:Q\to Q/G$. These are epimorphisms in the category of configuration spaces $\Man$.
    \item Given two configuration spaces $Q_1,Q_2\in\Man$ representing the possible positions of two physical systems, the combined system resulting from considering them together will have static states given by all the possible pairs of static states in each of the two systems. The categorical manifestation of the Principle of Combination for configuration spaces is then simply the presence of the Cartesian product of smooth manifolds $Q_1\times Q_2\in\Man$ as a categorical product.
    \item A temporal series of static states will be called a \textbf{path} of the physical system. Paths in a configuration space $Q\in\Man$ will be given by smooth curves $r:I\to Q$ parameterized by the observer's time parameter $t\in I\subset \Real$. Following the discussion on the Principle of Dynamics of Section \ref{FoundationsHamiltonianPhaseSpace}, dynamics on the static states of a configuration space $Q\in\Man$ are identified with smooth vector fields, i.e. $\Dyn{Q}=\Sec{\Tan Q}$.
\end{itemize}

We thus see that the category of configuration spaces is the first natural example of a theory of phase spaces where the observable functor is simply the assignment of the ring of smooth functions to a manifold $\text{Obs}=\text{C}^\infty:\Man \to \Ring$, the dynamics functor is the tangent functor $\text{Dyn}=\Tan : \Man \to \Lie_\Man$ and the evolution functor is given by taking the vector bundle of smooth ring derivations $\text{Evl}=\text{Der}: \Ring \to \Lie_\Man$. Indeed, these three functors fit in the phase space theory commutative diagram
\begin{equation*}
\begin{tikzcd}[row sep=small]
& \Ring \arrow[dd,"\text{Der}"] \\
\Man \arrow[ur,"\text{C}^\infty"] \arrow[dr,"\Tan"'] & \\
& \Lie_\Man 
\end{tikzcd}
\end{equation*}
as a direct consequence of the standard isomorphism between smooth vector fields and derivations of the ring of smooth functions:
\begin{equation*}
    \Sec{\Tan Q}\cong \Dr{\Cin{Q}}.
\end{equation*}
Note, however, that, since the categorical information of each object of $\Man$ is strictly a smooth manifold with no canonical choice of extra structure on it, the category of configuration spaces does not provide an example of a Hamiltonian theory of phase spaces.\newline

The mathematical implementation of the conventional \textbf{canonical formalism of classical mechanics} can be understood quite simply, when approached from this categorical angle, as the search for a category of phase spaces associated to the category of configuration spaces that forms a Hamiltonian theory of phase spaces in a natural or \emph{canonical} way, hence the name. As shown below, this will be achieved by the identification of the canonical symplectic structures on the cotangent bundles of smooth manifolds.\newline

Another, more physically-informed, approach to the canonical formalism of classical mechanics is the search for a category of phase space whose objects are the sets that encompass both the static and dynamics states of a physical system as motivated by the basic postulates of Newtonian mechanics, where positions and velocities are the initial data for the deterministic evolution of the system. More precisely, a phase space $P$ associated to a configuration space $Q$ should carry a space of observables naturally containing the static observables of $Q$ and its dynamics $\Obs{Q},\Dyn{Q}\subset\text{Obs}(P)$. We see that the cotangent bundle $P=\Cot Q$ appears, again, as the natural or \emph{canonical} choice for such a phase space since it is clear that the static observables $\Obs{Q}=\Cin{Q}$ and the dynamics $\Dyn{Q}=\Sec{\Tan Q}$ are recovered as the fibre-wise constant and fibre-wise linear functions of the cotangent bundle $\Cin{\Cot Q}$.\newline

In order to formally identify the categorical properties of what will become the category of canonical phase spaces we recall some well-known facts about the symplectic geometry of cotangent bundles.\newline

The cotangent bundle of any smooth manifold $Q$ carries a canonical symplectic structure
\begin{equation*}
    \begin{tikzcd}
    \Cot Q \arrow[d, "\pi_Q"]\\
    Q
    \end{tikzcd}
    \qquad \qquad \omega_Q:=-d\theta_Q, \quad \text{ with }\quad  \theta_Q|_{\alpha_q}(v):=\alpha_q(\Tan_{\alpha_q}\pi_Q (v)) \quad \forall v\in \Tan_{\alpha_q}(\Cot Q).
\end{equation*}
The non-degenerate Poisson structure $(\Cin{\Cot Q},\cdot,\{,\}_Q)$, with the natural inclusions of vector fields and functions on $Q$ as linear and basic functions on $\Cot Q$, respectively:
\begin{equation*}
    l_Q:\Sec{\Tan Q}\to \Cin{\Cot Q}, \qquad \pi_Q^*:\Cin{Q}\to \Cin{\Cot Q},
\end{equation*}
is clearly shown to be fibre-wise linear:
\begin{align*}
    \{l_Q(X),l_Q(Y)\}_Q & =  l_Q([X,Y])\\
    \{l_Q(X),\pi_Q^*f\}_Q & =  \pi_Q^* X [f]\\
    \{\pi_Q^*f,\pi_Q^*g\}_Q & = 0.
\end{align*}
The cotangent bundle of a Cartesian product is canonically symplectomorphic to the vector bundle product of cotangent bundles with the induced product symplectic forms:
\begin{equation*}
    (\Cot (Q_1\times Q_2),\omega_{Q_1\times Q_2})\cong (\Cot Q_1\boxplus \Cot Q_2, \Proj_1^*\omega_{Q_1} \oplus \Proj_2^*\omega_{Q_2}).
\end{equation*}
A smooth map $\varphi:Q_1\to Q_2$ induces a Lagrangian relation $\Cot \varphi\subset \Cot Q_1 \times \overline{\Cot Q_2}$ called the cotangent lift of $\varphi$  and defined by
\begin{equation*}
   \Cot \varphi:=\{(\alpha_q,\beta_p)| \quad \varphi(q)=p, \quad \alpha_q=(\Tan_q\varphi)^*\beta_p\}.
\end{equation*}
Here $\overline{\Cot Q_2}$ denotes $(\Cot Q_2,-\omega_{Q_2})$. When $\varphi:Q_1\to Q_2$ is a diffeomorphism, its cotangent lift becomes the graph of a symplectomorphism $\Cot \varphi :\Cot Q_2\to \Cot Q_1$. Composition of cotangent lifts as Lagrangian relations in the Weinstein symplectic category $\Symp_\Man$ is always strongly transversal, i.e giving a well-defined smooth Lagrangian submanifold, and thus the cotangent bundle construction can be seen as a contravariant functor
\begin{equation*}
    \Cot : \Man \to \Symp_\Man.
\end{equation*}
An embedded submanifold $i:S\to Q$, the intrinsic definition of holonomic constraints in standard classical mechanics, gives the coisotropic submanifold $\Cot Q|_S\subset (\Cot Q,\omega_Q)$ inducing the coisotropic reduction diagram
\begin{equation*}
    \begin{tikzcd}[sep=small]
        \Cot Q|_S \arrow[r, hookrightarrow] \arrow[d, twoheadrightarrow] & (\Cot Q,\omega_Q) \\
        (\Cot S,\omega_S).
    \end{tikzcd}
\end{equation*}
Note that the surjective submersion results from quotienting by the foliation given by the conormal bundle of $S$, in other words
\begin{equation*}
        \Cot S\cong \Cot Q|_S/(\Tan S)^0
\end{equation*}
as vector bundles over $S$. Lastly, for a free and proper group action $\phi: G\times Q\to Q$, with infinitesimal action $\psi:\mathfrak{g}\to \Sec{\Tan Q}$ and smooth orbit space $\Tilde{Q}:=Q/G$, the cotangent lift gives an action by symplectomorphisms $\Cot \phi : G\times \Cot Q\to \Cot Q$ with equivariant (co)moment map given by
\begin{equation*}
        \overline{\mu}:=l\circ \psi : \mathfrak{g}\to \Cin{\Cot Q},
\end{equation*}
this induces the symplectic reduction diagram
\begin{equation*}
    \begin{tikzcd}[sep=small]
    \mu^{-1}(0) \arrow[r, hookrightarrow] \arrow[d, twoheadrightarrow] & (\Cot Q,\omega_Q) \\
    (\Cot \Tilde{Q},\omega_{\Tilde{Q}}).
    \end{tikzcd}
\end{equation*}

In light of these results, we are compelled to define the \textbf{category of canonical symplectic phase spaces} simply as the image of the category of configuration spaces under the cotangent functor $\Cot (\Man)$. This is clearly a theory of phase spaces with notions of observable, dynamics and evolution functors as in the case of configuration spaces. Furthermore, the presence of a canonical symplectic structure on each phase space allows to define Hamiltonian maps simply by assignment of the Hamiltonian vector field to a function
\begin{equation*}
    \eta_Q:=\omega_Q^\sharp \circ d = \text{ad}_{\{,\}_Q}: \Obs{\Cot Q}\to \Dyn{\Cot Q}.
\end{equation*}
This makes the category $\Cot (\Man)$ into a Hamiltonian theory of phase spaces, thus achieving the motivating goal of finding a Hamiltonian theory of phase spaces canonically associated with the category of configuration spaces.\newline 

We are now in the position to argue that the structural content of standard canonical Hamiltonian mechanics can be encapsulated in the cotangent functor, which indeed preserves the categorical structures that are naturally identified with the defining principles of Hamiltonian phase spaces outlined in Section \ref{FoundationsHamiltonianPhaseSpace}. The cotangent functor is the \textbf{Hamiltonian functor} for the theory of phase spaces of configuration spaces and its action on the different formal notions with direct physical interpretation is summarized in the following list of correspondences:
\begin{center}
\begin{tabular}{ c c c }
Category of Configuration Spaces & Hamiltonian Functor  & Category of Phase Spaces \\
\hline
 $\Man$ & $\begin{tikzcd}\phantom{A} \arrow[r,"\mathbb{H}_\Man"] & \phantom{B} \end{tikzcd}$ & $\Symp_\Man$ \\ 
 $Q$ & $\begin{tikzcd} \phantom{Q} \arrow[r, "\Cot ",mapsto] & \phantom{Q} \end{tikzcd}$ & $(\Cot Q,\omega_Q)$ \\ 
 $\Obs{Q}$ & $\begin{tikzcd} \phantom{Q} \arrow[r,hookrightarrow, "\pi_Q^*"] & \phantom{Q}  \end{tikzcd}$ & $\Obs{\Cot Q}$ \\
 $\Dyn{Q}$ & $\begin{tikzcd} \phantom{Q} \arrow[r,hookrightarrow, "l_Q"] & \phantom{Q}  \end{tikzcd}$ & $\Obs{\Cot Q}$ \\
 $Q_1\times Q_2$ & $\begin{tikzcd} \phantom{Q} \arrow[r,mapsto, "\Cot "] & \phantom{Q}  \end{tikzcd}$ & $\Cot Q_1\boxplus \Cot Q_2$ \\
 $\varphi:Q_1\to Q_2$ & $\begin{tikzcd} \phantom{Q} \arrow[r,mapsto, "\Cot "] & \phantom{Q}  \end{tikzcd}$ & $\Cot \varphi\subset \Cot Q_1\boxplus \overline{\Cot Q_2}$ Lagrangian\\
 $i:S\hookrightarrow Q$ & $\begin{tikzcd} \phantom{Q} \arrow[r,mapsto, "\Cot"] & \phantom{Q}  \end{tikzcd}$ & $\Cot Q|_S\subset \Cot Q$ coisotropic \\
 $G\Acts Q$ & $\begin{tikzcd} \phantom{Q} \arrow[r,mapsto, "\Cot"] & \phantom{Q}  \end{tikzcd}$ & $G\Acts \Cot Q$ Hamiltonian action \\
\end{tabular}
\end{center}
The Hamiltonian functor connects the category of configuration spaces and the category of canonical symplectic phase spaces regarded as theories of phase spaces as summarized by the following diagram:
\begin{equation*}
\begin{tikzcd}[row sep=small]
 & \Ring \arrow[dd,"\text{Der}"']& & & \textsf{PoissAlg} \arrow[dd,"\text{Der}"] \\
 & & \Man \arrow[ul, "\text{C}^\infty"']\arrow[dl,"\Tan"]\arrow[r,"\Cot"] & \Symp_\Man \arrow[ur,"\text{C}^\infty"] \arrow[dr,"\Tan"'] & \\
 & \Lie & & & \Lie
\end{tikzcd}
\end{equation*}

Once a physical system with configuration space $Q$ is canonically assigned a phase space $\Cot Q$ as above, the only remaining task left for the physicist is to determine which choice of dynamics will produce evolutions of the system that match experimental temporal series of measurements. The Poisson structure present in $\Cot Q$ reduces this problem to the choice of an observable $h\in\Obs{\Cot Q}$ which, under the Hamiltonian map, gives a choice of conservative dynamics $\eta_Q(h)\in\Dyn{\Cot Q}$. This distinguished observable is often called the \textbf{energy} of the system. It generates the time evolution of the system and is, itself, a fundamental conserved quanity. Mathematically, this is a trivial fact by construction since
\begin{equation*}
    \eta_Q(h)[h]=\{h,h\}_Q=0
\end{equation*}
by antisymmetry of the Lie bracket. We see this as the formal implementation of the Principle of Conservation of Section \ref{FoundationsHamiltonianPhaseSpace} for the category of symplectic phase spaces. Given two systems with a choice of energy $(Q_1,h_1)$ and $(Q_2,h_2)$, where $h_i\in\Cin{\Cot Q_i}$, their combined product has a canonical choice of energy given by the sum of pull-backs via the canonical projections $h_1+h_2:=\Proj_1^*h_1+\Proj_2^*h_2\in\Cin{\Cot Q_1\times \Cot Q_2}$. This gives an extra line of assignments to the Hamiltonian functor:
\begin{center}
\begin{tabular}{ c c c }
Category of Configuration Spaces & Hamiltonian Functor  & Category of Phase Spaces \\
\hline
  $(Q_1,h_1)\times (Q_2,h_2)$ & $\begin{tikzcd} \phantom{Q} \arrow[r,mapsto, "\Cot "] & \phantom{Q}  \end{tikzcd}$ & $h_1+h_2 \in\Cin{\Cot Q_1\times \Cot Q_2}$ \\
\end{tabular}
\end{center}

The phase space formalism described so far in this section is general enough to account for a vast class of mechanical systems, however, this generality comes at a price: the Hamiltonian functor above fails to select a preferred choice of energy observable for a given configuration space. Turning to one of the earliest examples of mechanics we find inspiration to redefine the category of configuration spaces in order to account for some extra physical intuitions. In Newtonian mechanics, configuration spaces are submanifolds of Euclidean space and Cartesian products thereof, thus always carrying a Riemannian metric, encoding the physical notion of distance and angle; often also with a choice of potential, which is a function on the configuration space. This motivates us to refine our definition of configuration space and define the \textbf{category of Newtonian configuration spaces} $\Man_{\textsf{Newton}}$ whose objects are triples $(Q,g,V)$, where $Q\in\Man$, $g\in \Sec{\odot ^2 \Cot Q}$ Riemannian metric and $V\in\Cin{Q}$, and whose morphisms $\varphi:(Q_1,g_1,V_1)\to (Q_2,g_2,V_2)$ are smooth maps $\varphi:Q_1\to Q_2$ such that $g_1-\varphi^*g_2$ is positive semi-definite and $V_1=\varphi^*V_2$. When $\varphi$ is a diffeomorphism, a morphism in this category is, in particular, an isometry between $Q_1$ and $Q_2$. Note that a metric defines a quadratic form on tangent vectors $K_g:\Tan Q\to \Real$ given simply by $K_g(v):=\tfrac{1}{2}g(v,v)$, then, given the equivalence between quadratic forms and positive-definite bilinear forms, we see that the datum of a Newtonian configuration space $(Q,g,V)$ is equivalent to a choice of function on static states $V\in\Cin{Q}$ and a choice of quadratic function on velocities $K\in\Cin{\Tan Q}$, indeed the two fundamental notions of classical mechanics. The Cartesian product of configuration spaces gets updated to a categorical product in $\Man_{\textsf{Newton}}$ by setting
\begin{equation*}
    (Q_1,g_1,V_1) \times (Q_2,g_2,V_2) := (Q_1\times Q_2,g_1+g_2,V_1+V_2)
\end{equation*}
where
\begin{align*}
    g_1+g_2 &:= \Proj^*_1g_1 \oplus \Proj_2^*g_2\in \Sec{\odot ^2 (\Tan Q_1\boxplus \Tan Q_2))}\cong\Sec{\odot ^2 \Tan (Q_1\times Q_2)}, \\
    V_1+V_2 &:= \Proj^*_1V_1+\Proj_2^*V_2\in \Cin{Q_1\times Q_2}.
\end{align*}
A metric $g$ on $Q$ gives the usual musical isomorphism
\begin{equation*}
        \begin{tikzcd}
        \Tan Q \arrow[r, "\flat",yshift=0.7ex] & \Cot Q \arrow[l,"\sharp",yshift=-0.7ex]
        \end{tikzcd},
\end{equation*}
which can be used to regard the quadratic function $K$ identified with the metric $g$ as a quadratic function on the cotangent bundle by pull-back: $\sharp^*K_g\in\Cin{\Cot Q}$. Once a configuration space $(Q,g,V)$ is fixed, we now see how the phase space formalism for Newtonian configuration spaces does provide a canonical choice of energy by
\begin{equation*}
    E_{g,V}:=\sharp^*K_g+\pi_Q^*V,
\end{equation*}
what we call the \textbf{Newtonian energy}. This name is further justified by the fact that a direct computation shows that solving the Hamiltonian dynamics of this observable is equivalent to solving Newton's equations on a Riemannian manifold background $(Q,g)$ and with force field $F=-\sharp(dV)$. For two Newtonian configuration spaces $(Q_1,g_1,V_1)$ and $(Q_2,g_2,V_2)$, the categorical product construction above gives the following additivity property of Newtonian energy:
\begin{equation*}
    E_{g_1+g_2,V_1+V_2}=E_{g_1,V_1}+E_{g_2,V_2}.
\end{equation*}

Newtonian configuration spaces $\Man_{\textsf{Newton}}$ form a subcategory of $\Man$, thus the Hamiltonian functor must operate in a similar way but with the updated categorical product and the added assignment of the Newtonian energy as a canonical choice of observable to generate the dynamical evolution. We shall then update the Hamiltonian functor given above for general configuration spaces by adding the extra correspondence of Newtonian energy:
\begin{center}
\begin{tabular}{ c c c }
Category of Configuration Spaces & Hamiltonian Functor  & Category of Phase Spaces \\
\hline
 $\Man_{\textsf{Newton}}$ & $\begin{tikzcd}\phantom{A} \arrow[r,"\mathbb{H}_{\Man_{\textsf{Newton}}}"] & \phantom{B} \end{tikzcd}$ & $\Symp_\Man$ \\ 
 $(Q,g,V)$ & $\begin{tikzcd} \phantom{Q} \arrow[r, "\Cot ",mapsto] & \phantom{Q} \end{tikzcd}$ & $E_{g,V}\in\Obs{\Cot Q}$ \\
 $(Q_1,g_1,V_1)\times (Q_2,g_2,V_2)$ & $\begin{tikzcd} \phantom{Q} \arrow[r, "\Cot ",mapsto] & \phantom{Q} \end{tikzcd}$ & $(\Cot (Q_1\times Q_2),E_{g_1,V_1}+E_{g_2,V_2})$
\end{tabular}
\end{center}
Note that this result in the category of Newtonian configuration spaces motivates the additivity of energy of the combination product of two general phase spaces $(\Cot Q_1\boxplus \Cot Q_2,h_1+h_2)$ discussed above.

\section{Poisson, Presymplectic and Dirac Generalizations} \label{GeneralizationSymplectic}

The canonical symplectic structures found on cotangent bundles, seen in Section \ref{BriefHistoryOfPhaseSpace} as the historical origin of the notion of phase space, have been firmly established in Section \ref{CanonicalSymplectic} above as the natural mathematical objects encapsulating the standard formalism of classical Hamiltonian mechanics. In this section we present some mathematically-motivated generalizations of cotangent symplectic manifolds and discuss to what degree they retain the desirable features of phase spaces as summarized by the list of principles in Section \ref{FoundationsHamiltonianPhaseSpace}.\newline

An immediate generalization of cotangent bundles are, of course, general symplectic manifolds which can be regarded as a kind of global counterpart to cotangent bundles since Darboux theorem ensures that all symplectic manifolds are locally symplectomorphic to cotangent bundles of open subsets. The obvious notions that are compromised by taking general symplectic manifolds as phase spaces are configuration spaces and the Hamiltonian functor. Indeed, in symplectic cotangent manifolds, the zero section and the fibres are transversal Lagrangian submanifolds that have clear physical interpretation as the static states and the path dynamics of the system, i.e. the classic \emph{position and momenta} $(q,p)$, however, in general symplectic manifolds it is not possible to canonically identify concrete Lagrangian submanifolds representing these notions. Such Lagrangian submanifolds certainly exist locally in great generality, but the data of a general symplectic manifold does not specify a choice that allows for a clear physical interpretation. Nonetheless, it follows from the results in Section \ref{SymplecticGeometry} that the category\footnote{For the reminder of this chapter categories and functors are considered ``up to smoothness issues'' with the intent to illustrate the overarching patterns between the different geometric structures involved in the generalizations of phase space but keeping in mind that one may need to restrict to ``regular'' subclasses of objects and morphisms for the composition and functors to be well-defined in the standard categorical sense.} of general symplectic manifolds forms a Hamiltonian theory of phase spaces:
\begin{equation*}
\begin{tikzcd}[row sep=small]
& \textsf{PoissAlg} \arrow[dd,"\text{Der}"] \\
\Symp_\Man \arrow[ur,"\text{C}^\infty"] \arrow[dr,"\Tan"'] & \\
& \Lie_\Man 
\end{tikzcd}
\end{equation*}
where the Hamiltonian maps $\eta_M:\Cin{M}\to \Sec{\Tan M}$ are given simply by the Hamiltonian vector field map of the non-degenerate Poisson structure. This is called \textbf{the theory of symplectic phase spaces}. Clearly, all the principles listed in Section \ref{FoundationsHamiltonianPhaseSpace} are accounted for in this theory of phase spaces. We remark that, even if static states and path dynamics don't find a natural mathematical embodiment in this theory, this does not constitute an obstacle for practical physics, since what will be contrasted with experiment are the temporal evolution of the values of observables. Then, more generally, as long as values of observables have a clear physical interpretation, any theory of phase spaces can, in principle, be put to work to account for experimental data.\newline

The first obvious generalization that departs from the category of symplectic manifolds the relaxation of the non-degeneracy condition of symplectic structures regarded as 2-forms. Presymplectic manifolds appear naturally as the pull-back of symplectic structures to arbitrary submanifolds, particularly, energy surfaces, i.e. level sets of a choice of energy in symplectic phase spaces, are always strictly presymplectic manifolds since they are odd-dimensional. It follows from the discussion about admissible functions and Hamiltonian vector fields of Section \ref{SymplecticGeometry} that presymplectic manifolds form a theory of phase spaces
\begin{equation*}
\begin{tikzcd}[row sep=small]
& \textsf{PoissAlg} \arrow[dd,"\text{Der}"] \\
\textsf{PreSymp}_\Man \arrow[ur,"\text{C}^\infty_{\text{adm}}"] \arrow[dr,"\Tan|_{\text{Ham}}"'] & \\
& \Lie_\Man 
\end{tikzcd}
\end{equation*}
where $\Tan|_{\text{Ham}}$ stands for the involutive tangent distribution of Hamiltonian vector fields. We call this \textbf{the theory of presymplectic phase spaces}. The fact that Hamiltonian vector fields of admissible functions are defined up to sections of the characteristic distribution of a presymplectic manifold means that there is no canonical way identify Hamiltonian maps $\eta_M:\text{C}^\infty_{\text{adm}}(M)\to \Sec{\Tan M}$ and thus the theory of presymplectic phase spaces fails to be Hamiltonian. Therefore, given that the only cause for the non-uniqueness of Hamiltonian vector fields is the degeneracy of the 2-form, we see that by moving from symplectic to presymplectic phase spaces the Principle of Conservation is the only physical interpretation that is compromised from the list of Section \ref{FoundationsHamiltonianPhaseSpace}.\newline

Similarly, we may consider a departure from symplectic phase spaces by relaxing the non-degeneracy condition of the symplectic structures regarded as a bivectors. Poisson structures appear naturally as coisotropic reductions of symplectic manifolds, more concretely, non-holonomic constraints and symmetries of symplectic phase spaces give several examples such as, notably, the identification of the kinematic states of a rigid body with the linear Poisson manifold $\mathfrak{so}(3)^*$. The results in Section \ref{SymplecticGeometry} ensure that Poisson manifolds form a theory of Hamiltonian phase spaces:
\begin{equation*}
\begin{tikzcd}[row sep=small]
& \textsf{PoissAlg} \arrow[dd,"\text{Der}"] \\
\Poiss_\Man \arrow[ur,"\text{C}^\infty"] \arrow[dr,"\Tan"'] & \\
& \Lie_\Man 
\end{tikzcd}
\end{equation*}
with Hamiltonian maps given by the Hamiltonian vector field map, as in the case of symplectic phase spaces. We call this \textbf{the theory of Poisson phase spaces}. Note that the physical interpretation of this theory is essentially identical to that of general symplectic manifolds. In this respect, the only added generality of Poisson manifolds, aside from the mathematical fact that the lack of a non-degeneracy condition allows for a far broader class manifolds including any smooth manifold with the zero bracket on functions, comes from the absence of a Darboux theorem that allows to locally identify submanifolds representing static states and path dynamics of a physical system.\newline

The integrability condition of a presymplectic manifold $(M,\pi)$, $d\omega=0$, or a Poisson manifold $(M,\pi)$, $[\pi,\pi]=0$, can be relaxed by the introduction of a twisting 3-form $H\in\Omega_\text{cl}^3(M)$ and setting
\begin{equation*}
    d\omega=H,\qquad [\pi,\pi]=\tfrac{1}{2}\pi^\sharp(H).
\end{equation*}
As mentioned in Section \ref{BriefHistoryOfPhaseSpace}, this was first proposed in the context of string theory and it has led, among other things, to the mathematical formulation of quasi-Hamiltonian reduction with group-valued moment maps. The presence of a twisting constitutes a departure from the category of Poisson algebras as the natural setting for the algebraic structure of observables, since the induced brackets on functions have Jacobiators controlled by the twisting 3-form. Moreover, although Hamiltonian vector fields exist, they fail to preserve the twisted brackets and thus the usual physical interpretation of the integrability theory of conserved quantities is compromised. In the two propositions below we show explicitly that twisted presymplectic and twisted Poisson structures are, in fact, particular cases of Dirac structures as defined in Section \ref{DiracGeometry}, then we delay our comments on the viability of these as candidates for phase spaces for the same question formulated in the more general context of Courant algebroids.

\begin{prop}[Twisted Presymplectic Manifolds are Dirac Structures] \label{TwistedPreSymplecticDirac}
 A 2-form $\omega\in \Omega^2(M)$ is a presymplectic form twisted by $\normalfont H\in\Omega_\text{cl}^3(M)$ iff its graph $\normalfont \Grph{\omega^\flat}\subset \Tan M\oplus \Cot M$ is a Dirac structure in the $H$-twisted standard Courant algebroid $\mathbb{T}M_H$.
\end{prop}
\begin{proof}
Elements of the graph of the musical map associated to a 2-form $\omega\in \Omega^2(M)$ are clearly spanned by the following subspace of sections of the standard Courant algebroid
\begin{equation*}
    \{X\oplus \omega^\flat(X), X\in \Sec{\Tan M} \}\subset \Sec{\mathbb{T}M_H}.
\end{equation*}
These span the subbundle of the fibre-wise graphs of $\omega^\flat$ so, in particular, $\dimm \Grph{\omega^\flat_x}=\dimm \Tan_xM$ for all $x\in M$, then $\Grph{\omega^\flat}$ is clearly a subbundle of half rank in $\mathbb{T}M_H$. The definition of the standard symmetric bilinear form via the dual pairing shows that isotropy of $\Grph{\omega^\flat}$ is equivalent to the antisymmetry of the 2-form $\omega$, then it only remains to show that the integrability condition is equivalent to the involutivity of $\Grph{\omega^\flat}$ with respect to $H$-twisted Courant bracket. We check this directly, for any $X,Y\in\Sec{\Tan M}$, we have:
\begin{align*}
[X\oplus i_X\omega,Y\oplus i_Y\omega]_H&=[X,Y]\oplus \mathcal{L}_Xi_Y\omega +\mathcal{L}_Yi_X\omega +\tfrac{1}{2}d(i_Xi_Y\omega -i_Yi_X\omega) +i_Yi_XH\\
&=[X,Y]\oplus \mathcal{L}_Xi_Y\omega +\mathcal{L}_Yi_X\omega +di_Xi_Y\omega+i_Yi_XH\\
&=[X,Y]\oplus \mathcal{L}_Xi_Y\omega +\mathcal{L}_Yi_X\omega +\mathcal{L}_Xi_Y \omega -i_X\mathcal{L}_Y-i_Xi_Yd\omega+i_Yi_XH\\
&=[X,Y]\oplus i_X\mathcal{L}_Y\omega +\mathcal{L}_Yi_X\omega +i_Xi_Yd\omega+i_Yi_XH\\
&=[X,Y]\oplus i_{\mathcal{L}_XY}\omega +i_Xi_Yd\omega+i_Yi_XH\\
&=[X,Y]\oplus i_{[X,Y]}\omega +i_Xi_Y(d\omega-H)
\end{align*}
where, exploiting the fact that $\Grph{\omega^\flat}$ is an isotropic subbundle, we have used the antisymmetrized Dorfman bracket. Then it clearly follows that
\begin{equation*}
    [X+\omega^\flat(X),Y+\omega^\flat(Y)]_H=[X,Y]+\omega^\flat([X,Y]) \quad \Leftrightarrow \quad d\omega = H
\end{equation*}
as desired.
\end{proof}
\begin{prop}[Twisted Poisson Manifolds are Dirac Structures] \label{TwistedPoissonDirac}
 A bivector $\pi\in \Sec{\wedge^2\Tan M}$ is a Poisson bivector twisted by $\normalfont H\in\Omega_\text{cl}^3(M)$ iff its graph $\normalfont \Grph{\pi^\sharp}\subset \Tan M\oplus \Cot M$ is a Dirac structure in the $H$-twisted standard Courant algebroid $\mathbb{T}M_H$.
\end{prop}
\begin{proof}
We proceed in an entirely analogous way to the proof of proposition \ref{TwistedPreSymplecticDirac} above by identifying the graph of a bivector as the subbundle spanned by the following subset of sections:
\begin{equation*}
    \{\pi^\sharp(\alpha)\oplus \alpha, \alpha\in \Sec{\Cot M} \}\subset \Sec{\mathbb{T}M_H},
\end{equation*}
and by directly computing for any two $\alpha,\beta \in \Sec{\Cot M}$:
\begin{align*}
[\pi^\sharp(\alpha)\oplus \alpha, \pi^\sharp(\beta)\oplus \beta]_H&=[\pi^\sharp(\alpha),\pi^\sharp(\beta)]\oplus \mathcal{L}_{\pi^\sharp(\alpha)}\beta -\mathcal{L}_{\pi^\sharp(\beta)} \alpha +\tfrac{1}{2}d(i_{\pi^\sharp(\alpha)}\beta -i_{\pi^\sharp(\beta)}\alpha) +i_{\pi^\sharp(\beta)}i_{\pi^\sharp(\alpha)}H\\
&=[\pi^\sharp(\alpha),\pi^\sharp(\beta)]\oplus  \mathcal{L}_{\pi^\sharp(\alpha)}\beta -\mathcal{L}_{\pi^\sharp(\beta)} \alpha +d\pi(\alpha,\beta) +i_{\pi^\sharp(\beta)}i_{\pi^\sharp(\alpha)}H\\
&:=[\pi^\sharp(\alpha),\pi^\sharp(\beta)]\oplus [\alpha,\beta]_{\pi,H}\\
&= \pi^\sharp([\alpha,\beta]_{\pi,H}) - \mathcal{L}_{\pi^\sharp(\alpha)}\beta +\mathcal{L}_{\pi^\sharp(\beta)} \alpha-\pi^\sharp(i_{\pi^\sharp(\beta)}i_{\pi^\sharp(\alpha)}H) \oplus [\alpha,\beta]_{\pi,H}\\
&= \pi^\sharp([\alpha,\beta]_{\pi,H}) +2[\pi,\pi](\alpha,\beta)-\pi^\sharp(i_{\pi^\sharp(\beta)}i_{\pi^\sharp(\alpha)}H) \oplus [\alpha,\beta]_{\pi,H}
\end{align*}
where the intrinsic definition of the Schouten bracket has been used. It follows that
\begin{equation*}
    [\pi^\sharp(\alpha)\oplus \alpha, \pi^\sharp(\beta)\oplus \beta]_H=\pi^\sharp([\alpha,\beta]_{\pi,H})\oplus [\alpha,\beta]_{\pi,H}  \quad \Leftrightarrow \quad [\pi,\pi]=\tfrac{1}{2}\pi^\sharp(H)
\end{equation*}
as desired.
\end{proof}

It follows as a corollary of these two propositions that Poisson and presymplectic manifolds are recovered as particular cases of Dirac structures on (untwisted) standard Courant algebroids. Now, since proposition \ref{SeveraClass} identified the twisting 3-forms of standard Courant algebroids with the Severa class of an exact Courant algebroid, we are compelled to consider Dirac manifolds, defined as manifolds carrying exact Courant algebroids with a maximally isotropic involutive subbundle as in Section \ref{CourantAlgebroids}, as the natural generalized notion of phase space that will encompass all the previously discussed generalizations.\newline

From our discussion in Section \ref{DiracGeometry} it follows that Dirac manifolds form a category\footnote{This category entails similar subtleties to those of the Weinstein category of symplectic manifolds.} $\Dir_\Man$ with morphisms covering smooth maps between the underlying manifolds. Dirac morphisms covering embeddings can be regarded as subobjects, the natural categorical product for exact Courant algebroids of Section \ref{CourantAlgebroids} extends to a categorical product of Dirac structures in an obvious way and the results about extended actions on Courant algebroids and the Reduction of Dirac structures of Section \ref{CourantReduction} ensure that there are natural notions of surjections and subquotients in the category $\Dir_\Man$. If we recall the motivating physical principles that are desirable in a theory of phase spaces listed in Section \ref{FoundationsHamiltonianPhaseSpace}, we see that the category of Dirac manifolds carries the necessary categorical structure to encode the Principles of Reductionism, Combination and Symmetry, in addition addition to the Principles of Realism, Characterization, Observation and Dynamics common to all phase spaces modelled on smooth manifolds.\newline

The main advantage of identifying \textbf{the category of Dirac phase spaces} $\Dir_\Man$ is that many constructions that involved different geometric structures that may be regarded, in their own right, as generalized phase spaces, are now reinterpreted as natural categorical constructions of Dirac manifolds. For instance, general submanifolds of Poisson or presymplectic phase spaces can be uniformly understood as Dirac embeddings. Another important and illustrative example is coisotropic reduction: consider a symplectic manifold $(M,\omega)$ with a coisotropic $i:S \hookrightarrow M$ submanifold that induces a reduction scheme
\begin{equation*}
    \begin{tikzcd}[sep=small]
        S \arrow[r, hookrightarrow] \arrow[d, twoheadrightarrow] & (M,\omega) \\
        (M',\pi').
    \end{tikzcd}
\end{equation*}
such that $(M',\pi')$ is a Poisson manifold. In this situation, three different notions of phase space appear in a diagram of smooth maps: the symplectic and Poisson manifolds but also the coisotropic submanifold, which is an intermediate object that we identify as a presymplectic manifold with the restricted symplectic form $(S,i^*\omega)$. If we regard these three manifolds as Dirac structures of the corresponding (untwisted) standard Courant algebroids, the reduction condition simply becomes the fact that the above diagram is a diagram of well-defined Dirac morphisms.\newline

General Dirac manifolds, however, do not allow for an obvious physical interpretation of the Hamiltonian vector fields and the algebraic structures on the admissible functions as observables, since the Dirac bracket no longer satisfies the Jacobi identity and it is unclear what notion of derivation should be used to account for the dynamics and evolution functors.

\section{Poissonization Functors} \label{PoissonizationFunctors}

In this section we present some results that seem to indicate that all the generalizations presented in Section \ref{GeneralizationSymplectic} above share a fundamental Poisson flavour. This will be argued by showing that there are functors mapping different generalizations of phase spaces into the category of Poisson phase spaces and that these functors, furthermore, respect subobjects, products and reduction.\newline

The main strategy in constructing Poissonization functors is to exploit the one-to-one correspondence between Lie algebroids and linear Poisson structures. We begin by considering general Poisson phase spaces. The identification of Poisson manifolds as Dirac structures that follows from proposition \ref{TwistedPoissonDirac} and the fact that Dirac structures are, in particular, Lie algebroids shows that the datum of a Poisson manifold $(M,\omega)$ can be equivalently encoded in a Lie algebroid structure on the cotangent bundle $(\Cot M,\pi^\sharp,[,]_\pi)$ with the bracket being defined by
\begin{equation*}
    [\alpha,\beta]_\pi:=\mathcal{L}_{\pi^\sharp(\alpha)}\beta -\mathcal{L}_{\pi^\sharp(\beta)} \alpha +d\pi(\alpha,\beta)
\end{equation*}
for all $\alpha,\beta\in\Sec{\Cot M}$. Note that this Lie bracket is characterized from the bracket on differentials of pairs of functions $f,g\in\Cin{M}$
\begin{equation*}
    [df,dg]_\pi=d\{f,g\}
\end{equation*}
by Leibniz extension with the anchor $\pi^\sharp$. We call $(\Cot M,\pi^\sharp,[,]_\pi)$ the \textbf{Poisson-Lie algebroid}. Now, in virtue of proposition \ref{LieAlgebroidLinearPoisson}, the Poisson-Lie algebroid induces a linear Poisson structure on the tangent bundle that we denote by $(\Tan M,p^\pi)$. Since differentials of functions span all sections of the cotangent bundle, the linear Poisson constructed on the dual bundle of a Lie algebroid of proposition \ref{LieAlgebroidLinearPoisson} gives the following characterization of the linear Poisson bracket on $\Tan M$:
\begin{align*}
\{l_{df},l_{dg}\}_{\Tan M} & = l_{d\{f,g\}_M}\\
\{l_{df},\tau^*g\}_{\Tan M} & = \tau^*\{f,g\}_M\\
\{\tau^*f,\tau^*g\}_{\Tan M} & = 0
\end{align*}
for all $f,g\in \Cin{M}$ and were $\tau: \Cot M\to M$ is the tangent bundle projection. We thus find the tangent functor as the first example of a Poissonization functor.

\begin{prop}[The Tangent Functor as a Poissonization Functor] \label{TangentPoissonizationFunctor}
The assignment of the linear Poisson structure $\normalfont (\Tan M,p^\pi)$ to a general Poisson manifold $(M,\pi)$ is functorial. Furthermore, products of Poisson manifolds are mapped to products of linear Poisson structures on vector bundles and coisotropic submanifolds of general Poisson structures are mapped to coisotropic subbundles of the linear Poisson structures.
\end{prop}
\begin{proof}
Functoriality follows from the following basic facts that result from the structure of the spanning functions on tangent bundles:
\begin{equation*}
    (\Tan \varphi)^*l_{df}=l_{d\varphi^*f}, \qquad (\Tan \varphi)^*\tau_M^*f=\tau_N^*\varphi^*f
\end{equation*}
where $f\in\Cin{M}$ and $\varphi:N\to M$ is an arbitrary smooth function. It is then clear that if $\varphi : (N,\pi_N) \to (M,\pi_M)$ is a Poisson map then so is $\Tan \varphi : (\Tan N,p^{\pi_N}) \to (\Tan M,p^{\pi_M})$. It follows from the discussion on the tangent functor for general smooth manifolds of Section \ref{TangentFunctor} that $\Tan$ sends Cartesian products of manifolds into vector bundle products of tangent bundles. It follows from the defining condition of the product Poisson structure $(M_1\times M_2,\pi_{12})$ of two Poisson manifolds $(M_i,\pi_i)$, $i=1,2$:
\begin{equation*}
    \Proj_1^*\{f_1,g_1\}_1=\{\Proj_1^*f_1,\Proj_1^*g_1\}_{12} \qquad \Proj_1^*\{f_2,g_2\}_2=\{\Proj_2^*f_2,\Proj_2^*g_2\}_{12}
\end{equation*}
for all  $f_i,g_i\in\Cin{M_i}$, by taking tangents of the projection maps that the linear Poisson structures are related by
\begin{equation*}
    p^{\pi_{12}}=(\Tan \Proj_1)^*p^{\pi_1}+ (\Tan \Proj_2)^*p^{\pi_2}.
\end{equation*}
Lastly, we use proposition \ref{CoisotropicSubPoisson} to characterize a coisotropic submanifold $C\subset M$ by its vanishing ideal as a coisotrope in the Poisson algebra of functions $I_C\subset (\Cin{M},\cdot,\{,\})$ and note that the tangent subbundle $\Tan C\subset \Tan M$ has vanishing ideal generated by functions $l_df$ and $\tau^*f$ for $f\in I_C$. The defining conditions of the linear Poisson bracket then clearly show that these functions form a Lie subalgebra that is, furthermore, closed under multiplication, thus giving a coisotrope in the linear Poisson manifold $(\Tan M,p^\pi)$.
\end{proof}

We shall prove now that the tangent functor preserves reduction schemes. We begin by considering a general Poisson reduction scheme as defined in Section \ref{SymplecticGeometry}.

\begin{thm}[Poisson Reduction induces Linear Poisson Reduction]\label{LinearPoissonReduction}
Let $(M,\pi)$ be a general Poisson structure, $i:C\to M$ an immersed \emph{constraint} submanifold, $D\subset \Tan M$ a tangent distribution that becomes a subbundle of $\normalfont \Tan M|_C$ when restricted to $C$ and let $((M,\pi),C,D)$ form a reduction data triple with reduced Poisson structure $(\tilde{M},\tilde{\pi})$, then the tangent counterparts $\normalfont ((\Tan M,p^\pi),\Tan C,\Tan D)$  are a reduction data triple which reduces to the linear Poisson structure $\normalfont (\Tan\tilde{M},p^{\tilde{\pi}})$. 
\end{thm}
\begin{proof}
The first two integrability conditions for the tangent triple are direct implications of the integrability of the base triple; we simply observe that $\Tan \Tan C\cap \Tan D = \Tan K$ with $K=\Tan C\cap D$ the integrable subbundle on $i:C \to M$ with foliation $\mathcal{F}$, which defines the leaf space $q:C\to \tilde{M}$, and that the integrable subbundle $\Tan K\subset \Tan C$ has regular foliation the family of tangent bundles of the leaves $\Tan \mathcal{F}$. The space of leaves is a smooth manifold since when the foliation of the tangent bundles of the leaves is examined closely, we notice that:
\begin{equation*}
\Tan C/\Tan \mathcal{F}=\frac{\Tan C}{K}/\mathcal{F}\cong \Tan \tilde{M},
\end{equation*}
so it is naturally identified with the tangent bundle of the reduced manifold. The third compatibility condition of the triple and reducibility involve the Poisson structure and, since a bivector field is specified by the values on all covectors at every cotangent space, we shall consider families of spanning functions on $\Tan  M$ and $\Tan \tilde{M}$ whose differentials span the cotangent bundle. As discussed in the proof of proposition \ref{TangentPoissonizationFunctor} this is done by considering subspaces of functions of the form
\begin{equation*}
f=\tau^*\phi + l_{d\phi'},
\end{equation*}
for $\phi,\phi'\in \Cin{M}, \Cin{\tilde{M}}$. We see that, since $D\subset \Tan M$  is a subbundle when restricted to $C$, and, simply using the vector bundle structure, we have:
\begin{align*}
&\Tan \tau_M |_{\Tan D}:\Tan D\to D\\
&dl_\alpha|_{\Tan M|_C}\in(\Tan D)^0 \Leftrightarrow \alpha|_C \in D^0
\end{align*}
for any $\alpha\in \Gamma(\Cot M)$. We can now explicitly show that the third compatibility condition for the reducible triple $((M,\pi),C,D)$ implies the analogous condition for $((\Tan M,p^\pi),\Tan C,\Tan D)$ by writing:
\begin{equation*}
\{f,g\}_{\Tan M}=\{\tau_M^*\phi + l_{d\phi'},\tau_M^*\gamma + l_{d\gamma'}\}_{\Tan M}=\tau_M^*(\{\phi',\gamma\}_M+\{\phi,\gamma'\}_M)+l_{d\{\phi',\gamma'\}_M}
\end{equation*}
where we have used the definition of the linear Poisson structure on $\Tan M$ from the Poisson structure on $M$. Using the two properties listed above it is immediate to see now that $d\phi|_C,d\gamma|_C\in D^0\Rightarrow d\{\phi,\gamma\}_M|_C\in D^0$ is equivalent to $df|_{\Tan C},dg|_{\Tan C}\in (\Tan D)^0\Rightarrow d\{f,g\}_{\Tan M}|_{\Tan C}\in (\Tan D)^0$. The reducibility condition for the linear Poisson structure $p^\pi$ into $p^{\tilde{\pi}}$ reads:
\begin{equation*}
(\Tan q)^*\{f,g\}_{\Tan \tilde{M}}=(\Tan i)^*\{F,G\}_{\Tan M}
\end{equation*}
for all $f,g\in \Cin{\Tan \tilde{M}}$ and $F,G \in \Cin{\Tan M}$ extensions of $(\Tan q)^*f,(\Tan q)^*g$ satisfying $dF|_{\Tan C},dG|_{\Tan C}\in (\Tan D)^0$. It will suffice check the above condition for spanning functions and, recalling the arguments above, it should be noted that spanning extensions with differentials in $(\Tan D)^0$ for functions of the form $(\Tan q)^*(\tau_{\tilde{M}}^*\phi + l_{d\phi'})$ can be given by $\tau_M^*\Phi + l_{d\Phi'}$, where $\Phi,\Phi'\in \Cin{M}$ are extensions of $q^*\phi,q^*\phi'$ satisfying $d\Phi,d\Phi'\in D^0$, without loss of generality. The reducibility condition is then rewritten as:
\begin{align*}
&(\Tan q)^*\tau_{\tilde{M}}^*(\{\phi',\gamma\}_{\tilde{M}}+\{\phi,\gamma'\}_{\tilde{M}}+l_{d\{\phi',\gamma'\}_{\tilde{M}}})=\\
=&(\Tan i)^*\tau_M^*(\{\Phi',\Gamma\}_M+\{\Phi,\Gamma'\}_M+l_{d\{\Phi',\Gamma'\}_M})
\end{align*}
for all $\phi,\phi',\gamma,\gamma'\in \Cin{\tilde{M}}$ and $\Phi,\Phi',\Gamma,\Gamma'\in \Cin{M}$ extensions as above, which, using the commutativity of the tangent projection and tangent map diagrams, readily becomes:
\begin{align*}
&\tau_C^*(q^*\{\phi',\gamma\}_{\tilde{M}}+q^*\{\phi,\gamma'\}_{\tilde{M}})+l_{dq^*\{\phi',\gamma'\}_{\tilde{M}}}=\\
=&\tau_C^*(i^*\{\Phi',\Gamma\}_M+i^*\{\Phi,\Gamma'\}_M)+l_{di^*\{\Phi',\Gamma'\}_M}.
\end{align*}
The above equality obviously holds from the reducibility condition for $((M,\pi),C,D)$, thus proving reducibility for the tangent triple.
\end{proof}

In the next two propositions we prove that the tangent Poissonization functor is furtheremore compatible with two special reduction schemes of particular relevance in mechanics: the case of coisotropic constraints, where a coisotropic submanifold $C\subset M$ is defined as the zero locus of a map $\Phi:M\to \mathbb{R}^l$ with $0\in \mathbb{R}^l$ a regular value and whose components satisfy $\{\phi_i,\phi_j\}=f^k_{ij}\phi_k$ for some $f_{ij^k}\in \Cin{M}$, and the case of a Hamiltonian group action with equivariant momentum map $\mu:M\to \mathfrak{g}^*$.

\begin{prop}(Coisotropic Reduction induces Linear Coisotropic Reduction)\label{LinearCoistropicReduction}
Let a coisotropic reduction data triple $((M,\pi),\Phi^{-1}(0),\pi^\sharp(d\Phi))$ defined by a set of coisotropic constraints $\Phi:M\to \mathbb{R}^l$ which reduces to the Poisson manifold $(\tilde{M},\tilde{\pi})$ then the associated linear reduction data triple given by $\normalfont ((\Tan M,p^\pi),\Tan \Phi^{-1}(0),\Tan \pi^\sharp(d\Phi))$ is coisotropic and reduces to the linear Poisson structure $\normalfont (\Tan \tilde{M},p^{\tilde{\pi}})$.
\end{prop}
\begin{proof}
The reducibility part of the statement above is a direct consequence of theorem \ref{LinearPoissonReduction} so we only need to show that the linear counterpart of the coisotropic submanifold $\Tan \Phi^{-1}(0)$ is a coisotropic submanifold of $(\Tan M,p^\pi)$. This can be shown by the explicit construction of the following coisotropic constraint function:
\begin{equation*}
\Lambda_{\Phi}:=l_{d\Phi}\times \tau^*\Phi : \Tan M \to \mathbb{R}^l\times \mathbb{R}^l
\end{equation*}
whose components are of the form $\{l_{d\phi_i},\tau^*\phi_j\}$. Both $l$ and $\tau^*$ are injective maps thus $0\in \mathbb{R}^{2l}$ will be a regular value of $\Lambda_\Phi$. More explicitly, we can directly show
\begin{equation*}
\Lambda_\Phi^{-1}(0)=\Tan \Phi^{-1}(0).
\end{equation*}
To complete the proof we show that these constraints are coisotropic by directly computing the pertinent Poisson brackets:
\begin{align*}
\{l_{d\phi_i},l_{d\phi_j}\}_{\Tan M} & = \tau^*(f^k_{ij})l_{d\phi_i}+l_{df_{ij}^k}\tau^*\phi_k\\
\{l_{d\phi_i},\tau^*\phi_j\}_{\Tan M} & = \tau^*(f^k_{ij})\tau^*\phi_k\\
\{\tau^*\phi_i,\tau^*\phi_j\}_{\Tan M} & = 0
\end{align*}
and noting that all three brackets give $\Cin{M}$-linear combinations of the constraint functions.
\end{proof}

\begin{prop}[Hamiltonian Reduction induces Linear Hamiltonian Reduction]\label{LinearHamiltonianReduction}
Let a Hamiltonian reduction data triple $((M,\pi),\mu^{-1}(0),\psi(\mathfrak{g}))$ defined by a Poisson $G$-action on $(M,\pi)$ and a moment map $\mu:M\to \mathfrak{g}^*$ which reduces to the Poisson manifold $(\tilde{M},\tilde{\pi})$, then the linear reduction data triple $\normalfont ((\Tan M,p^\pi),\Tan \mu^{-1}(0),\Tan \psi(\mathfrak{g}))$ is Hamiltonian for a Poisson $\normalfont \Tan G$-action on $\normalfont (\Tan M,p^\pi)$ and reduces to the linear Poisson structure $\normalfont (\Tan \tilde{M},p^{\tilde{\pi}})$.
\end{prop}
\begin{proof}
Again, the reducibility part of the statement is a consequence of theorem \ref{LinearPoissonReduction} but that a Hamiltonian $G$-action on $(M,\pi)$ gives a Hamiltonian $\Tan G$-action on $(\Tan M,p^\pi)$ remains to be shown. We first note that the tangent bundle of a Lie group $\Tan G$ has a natural Lie group structure inherited from that of $G$ as shown in detail in Section \ref{TangentFunctor}. We can use the left-invariant canonical 1-form $\theta^L\in\Omega^1(G,\mathfrak{g})$ to make the usual diffeomorphism $\Tan G\cong G\times \mathfrak{g}$ into the following homomorphism of Lie groups:
\begin{equation*}
(\Tan G,\cdot) \cong G\times_{\text{Ad}}\mathfrak{g}
\end{equation*}
where $\mathfrak{g}$ is seen as an abelian group with the addition of vectors. The Lie algebra of this Lie group is then easily shown to be:
\begin{equation*}
\mathfrak{tg} \cong \mathfrak{g}\oplus_{\text{ad}}\mathfrak{g}
\end{equation*}
where the second component is taken as the abelian Lie algebra. Furthermore, given a smooth action on a manifold $\phi: G\times M \to M$ we can use the fact that $\Tan (G\times M)\cong \Tan G \times \Tan M$ to define an action as follows:
\begin{align*}
\Phi(\cong \Tan \phi): \Tan G\times \Tan M & \to \Tan M\\
v_g+w_x & \mapsto \Tan _g\phi(\cdot,x)v_g+\Tan _x\phi(g,\cdot)w_x.
\end{align*}
This is easily checked to be an action with respect to the product defined above $(\Tan G,\cdot)$. The infinitesimal action will be denoted by $\Psi: \mathfrak{tg}\to \Gamma(\Tan \Tan M)$ and it suffices to specify how the infinitesimal generators of arbitrary elements of the Lie algebra act as derivations on spanning functions, which follows from a direct computation using the definitions above:
\begin{align*}
\Psi(\xi + \eta)[\tau^*f] & = \tau^*\psi(\xi)[f]\\
\Psi(\xi + \eta)[l_{df}] & = \tau^*\psi(\eta)[f]+l_{d\psi(\xi)[f]}.
\end{align*}
Now in the specific case when the manifold is Poisson with a (co)momentum map $\mu:\mathfrak{g}\to \Cin{M}$ it is immediate to check that the above infinitesimal action becomes Hamiltonian with (co)moment map given by:
\begin{align*}
\Lambda_\mu: \mathfrak{tg} \cong \mathfrak{g}\oplus_{\text{ad}}\mathfrak{g} & \to \Cin{\Tan M}\\
\xi + \eta & \mapsto l_{d\mu_\xi}+\tau^*\mu_\eta.
\end{align*}
Finally it is then clear that this (co)moment map corresponds to the dual of the constraint map in proposition \ref{LinearCoistropicReduction} with $\Phi=\mu$, therefore the result follows.
\end{proof}

With these results at hand one may be tempted to think that linear Poisson structures on tangent bundles faithfully characterize general Poisson structures on the base manifolds. This is, however, not the case. To understand this, note that, indeed, all Poisson structures can be seen as linear Poisson structures via the one-to-one correspondence of linear Poisson structures with Lie algebroids, but if we recall that this was motivated by the fact that Poisson structures are particular cases of Dirac structures, it is easy to see that general reduction schemes of linear Poisson structures correspond to general reduction schemes of Poisson-Dirac structures which may produce non-Poisson Dirac structures in general.\newline

Let us now turn our attention to presymplectic phase spaces. The identification of presymplectic manifolds with the graphs of the 2-forms as Dirac structures in the standard Courant algebroid $\Tan  M \oplus \Cot M$ allows us to regard the datum of a presymplectic manifold as the standard Lie algebroid $\Tan M$ together with an isotropic embedding given by a musical map $\omega^\flat$. Note that the one-to-one correspondence between Lie algebroids and linear Poisson structures allows us to recognize the linear Poisson structure associated to the canonical symplectic structure on the cotangent bundle $(\Cot M,\Omega)$ as the dual of the standard Lie algebroid structure on $\Tan M$. This suggests that we make use of the cotangent functor to to assign a non-degenerate spanning Poisson structure to an arbitrary presymplectic manifold.

\begin{prop}[The Cotangent Functor as a Poissonization Functor] \label{CotangentPoissonizationFunctor}
Let $(M,\omega)$ be a presymplectic manifold, then there exists a unique symplectic structure the cotangent bundle $\normalfont (\Cot M,\Omega^\omega)$ such that its associated non-degenerate Poisson bracket preserves spanning functions and satisfies
\begin{equation*}
\{l_X,l_Y\}=l_{[X,Y]}-\tau^*\omega(X,Y).
\end{equation*}
for all $\normalfont X,Y\in\Sec{\Tan  M}. $Furthermore, the assignment $\normalfont (M,\omega) \mapsto (\Cot M,\Omega^\omega)$ is functorial, preserves products and sends isotropic submanifolds on the base into coisotropic submanifolds on the cotangent bundle.
\end{prop}
\begin{proof}
We see existence and uniqueness by explicitly constructing the symplectic structure as
\begin{equation*}
\Omega^\omega:=\Omega+\tau^*\omega
\end{equation*}
where $\Omega\in\Omega^2(\Cot M)$ denotes the canonical 2-form and $\tau:\Cot M\to M$ is the standard vector bundle projection. This 2-form is clearly closed and expressing it in canonical coordinates gives the matrix:
\begin{equation*}
[\Omega^\omega]=\begin{pmatrix} [\omega] & I\\ -I & 0 \end{pmatrix}
\end{equation*}
where $[\omega]$, the coordinate expression of the presymplectic form, is a general antisymmetric matrix. This is clearly a non-degenerate matrix with inverse
\begin{equation*}
[\Omega^\omega]^{-1}=\begin{pmatrix} 0 & -I\\ I & [\omega] \end{pmatrix}.
\end{equation*}
These matrices allow us to write the local expression for the corresponding non-degenerate bivector and compute explicitly all the spanning brackets:
\begin{align*}
\{l_X,l_Y\} & = l_{[X,Y]}-\tau^*\omega(X,Y)\\
\{l_X,\tau^*f\} & = \tau^*X[f]\\
\{\tau^*f,\tau^*g\} & = 0
\end{align*}
thus showing the required conditions for spanning Poisson structure. To see functoriality we simply note that the cotangent lift of a presymplectic map $\varphi:(M_1,\omega_1)\to (M_2,\omega_2)$ is a Lagrangian submanifold of $\Cot  M_1 \times \overline{\Cot M_2}$ with respect to the canonical symplectic form $\Omega_1 \oplus -\Omega_2$ from the general argument made for the cotangent functor $\Cot : \Man \to \Symp_\Man$ and its becoming Lagrangian with respect to the symplectic form $\Omega_1^{\omega_1} \oplus -\Omega_2^{\omega_2}$ is easily seen to be equivalent to $\varphi^*\omega_2=\omega_1$, indeed the condition of $\varphi$ to a presymplectic map. The usual isomorphism of symplectic manifolds $(\Cot (M_1\times M_2),\Omega_{12})\cong (\Cot M_1,\Omega_1)\times (\Cot M_2,\Omega_2)$ gives the following identity
\begin{equation*}
    \Omega^{\omega_{12}}_12=\Proj_1^* \Omega_1^{\omega_1}+ \Proj_2^* \Omega_2^{\omega_2}
\end{equation*}
from the fact that $\tau_{M_1\times M_2}=\tau_{M_1} \times \tau_{M_2}$, thus products of presymplectic manifolds induce products of their associated symplectic structures on the cotangent bundle. An isotropic submanifold $i:C\hookrightarrow M$ satisfies $i^*\omega=0$ then it is clear that the conormal bundle $\Tan  C^0\subset \Cot M$, which is always a coisotropic submanifold with respect to the canonical symplectic structure $\Omega$, is also coisotropic with respect to the symplectic form $\Omega^\omega$.
\end{proof}

We now prove that this Poissonization functor is compatible with presymplectic reduction.

\begin{thm}[Presymplectic Reduction induces Symplectic Reduction] \label{PresymplecticSymplecticReduction}
Let a presymplectic structure $(M,\omega)$ that reduces to $(\tilde{M},\tilde{\omega})$ via a Hamiltonian $G$-action with equivariant moment map $\mu$ as in proposition \ref{PresymplecticReduction}, then $\normalfont (\Cot M,\Omega^\omega)$ reduces to the symplectic manifold $\normalfont (\Cot \tilde{M},\tilde{\Omega}^{\tilde{\omega}})$ via a Hamiltonian $\normalfont \Tan G$-action and equivariant moment map
\begin{equation*}
 \lambda:=l\circ \psi + \tau^*\mu.    
\end{equation*}
\end{thm}
\begin{proof}
Let us denote the presymplectic group action on $(M,\omega)$ as $\phi:G\times M\to M$ and the infinitesimal action as $\psi:\mathfrak{g}\to \Gamma(\Tan M)$. If we denote the cotangent lift action by $\Phi:G\times \Cot M\to \Cot M$ with $\Phi:=(\Tan \varphi^{-1})^*$, we can define the following $\Tan G$ action:
\begin{align*}
\Phi^\omega : \Tan G\times \Cot M & \to \Cot M\\
((g,\xi),\alpha) & \mapsto \Phi_g(\alpha+i_{\psi(\xi)}\omega), 
\end{align*}
which is well-defined since $\phi$ is a presymplectic action and where we have used the standard isomorphism $\Tan G \cong G\times_{\text{Ad}}\mathfrak{g}$. It is then a simple computation to show that this is a symplectic action, $(\Phi^\omega_{(g,\xi)})^*\Omega^\omega=\Omega^\omega$. The corresponding infinitesimal action $\Psi^\omega:\mathfrak{g}\oplus_{\text{ad}} \mathfrak{g} \to \Gamma(\Tan \Cot M)$ can be computed explicitly on spanning functions:
\begin{align*}
\Psi(\xi \oplus \eta)[\tau^*f] & = \tau^*\psi(\xi)[f]\\
\Psi(\xi \oplus \eta)[l_X] & = l_{[\psi(\xi),X]}-\tau^*\omega(\psi(\xi+\eta),X).
\end{align*}
A routine computation shows that this is a Hamiltonian action for the comoment map:
\begin{align*}
\lambda : \mathfrak{g}\oplus_{\text{ad}} \mathfrak{g} & \to \Cin{\Cot M}\\
(\xi \oplus \eta) & \mapsto l_{\psi(\xi)}+\tau^*\mu_\eta. 
\end{align*}
We have thus shown that the $\Tan G$ action above is symplectic. It remains to verify that the reduced manifold is indeed $(\Cot \tilde{M},\tilde{\Omega}^{\tilde{\omega}})$. To see this, note that the constraint submanifold is $I:\lambda^{-1}(0)=\psi(\mathfrak{g})^0|_C\to \Cot M$. The fact that the $\Tan G$ action is a fibre-bundle action (although not linear) on $\Cot M$ implies that
\begin{equation*}
\psi(\mathfrak{g})^0|_C/G\times \mathfrak{g} \cong \frac{\psi(\mathfrak{g})^0|_C}{\omega^\flat(\psi(\mathfrak{g}))|_C}/G.
\end{equation*}
The RHS above is simply the orbit space of the cotangent lift action of $G$ that descends to the quotient by equivariance of $\mu$ and the fact that the original group action preserves the presymplectic form. The $\Tan G$ action is then seen to be free and proper on $\psi(\mathfrak{g})^0|_C$ if the original action is so on $C$, the corresponding surjective submersion will be denoted by $Q:\psi(\mathfrak{g})^0|_C\to \psi(\mathfrak{g})^0|_C/\Tan G$. By construction, $C$ is the constraint manifold given by the momentum map, so $\omega^\flat(\psi(\mathfrak{g}))|_C=(\Tan C)^0|_C$, hence we find the desired reduced manifold
\begin{equation*}
\frac{\psi(\mathfrak{g})^0}{(\Tan C)^0}\bigg|_C/G \cong \Cot (C/G).
\end{equation*}
It is then straightforward to show that the reducibility condition $Q^*\tilde{\Omega}^{\tilde{\omega}}=I^*\Omega^\omega$ is a direct implication of the reducibility condition for the presymplectic forms and the fact that the canonical forms reduce under cotangent lifts of group actions.
\end{proof}

The two Poissonization functors introduced so far have been motivated by the interpretation of Poisson and presymplectic manifolds as particular cases of Dirac structures, it is then natural to consider extending such construction to the category of Dirac manifolds. In Section \ref{DiracGeometry} it was shown that when an exact Courant algebroid $E$ admits two transversal Lagrangian subbundles $L$ and $L'$ then we have $E\cong L\oplus L'$ and the non-degenerate bilinear for on $E$ allows to identify the dual of each Lagrangian subbundle with its complement $L^*\cong L'$. Since Dirac structures are Lagrangian subbundles that become Lie algebroids with the restricted Courant bracket, following the examples set by Poisson and presymplectic manifolds, we suspect that a Poissonization functor exists for Dirac structures.

\begin{conj}[Poissonization Functor for Dirac Structures] \label{PoissonizationDirac}
Let $E$ be a Courant algebroid and $L\subset E$ a Dirac structure admitting a Lagrangian complement $L'\subset E$, then there is a spanning Poisson structure on the Lagrangian complement seen as a vector bundle $(L',p^L)$. Furthermore, the datum of the Dirac structure $L$ seen as a Lie algebroid with an isotropic embedding into $E$ is equivalent to the spanning Poisson structure $(L',p^L)$.
\end{conj}

Such a construction will constitute a generalization of the one-to-one correspondence between Lie algebroids and linear Poisson structures of proposition \ref{LieAlgebroidLinearPoisson}. Indeed, given any Lie algebroid $(A,\rho,[,])$ we can form a bialgebroid $(A,A^*)$ trivially by equipping $A^*$ with the zero anchor and zero bracket, thus making $A\oplus A^*$ into a Courant algebroid where $A \oplus 0$ and $0\oplus A^*$ are transversal Dirac structures. The spanning Poisson structure $(0\oplus A^*, p^L)$ given by conjecture \ref{PoissonizationDirac} will simply be the linear Poisson structure induced on $A^*$ by proposition \ref{LieAlgebroidLinearPoisson} reflecting the fact that the isotropic embedding $A\hookrightarrow A\oplus 0$ contains no more information than $A$ itself. Assuming the conjecture holds, at least in some generality within the category of Dirac manifolds, it is easy to see that the assignment will be functorial, sending Dirac morphisms to coisotropic relations between the spanning Poisson manifolds, and that products are preserved.\newline

Theorems \ref{LinearPoissonReduction} and \ref{PresymplecticSymplecticReduction} set the example for the following general conjecture that, in case of being proven true in some generality, will establish the assignment of a spanning Poisson manifold to a Dirac structure in conjecture \ref{PoissonizationDirac} as a Poissonization functor.

\begin{conj}[Dirac Reduction induces Spanning Poisson reduction] \label{DiracPoissonReduction}
Let $E$ be a Courant algebroid with a pair of transversal Dirac structures $L,L'\subset E$, assume that $E$ reduces to $E_r$ as in proposition \ref{CourantAlgebroids} and that $L,L'\subset E$ reduce to the Dirac structures $L_r,L_r'\subset E_r$, then the spanning Poisson structure $(L',p^L)$ reduces to the spanning Poisson structure $(L_r',p^{L_r})$.
\end{conj}

Figure \ref{GeneralizationDiagramSymplectic} captures pictorially the generalizations of canonical symplectic phase spaces proposed in this section.

\chapter{Contact Phase Spaces} \label{ContactPhaseSpaces}

In this chapter we present a generalization of the canonical phase space formalism incorporating the notion of physical dimension. This will be done by replacing ordinary configuration spaces, seen in Section \ref{CanonicalSymplectic} to be smooth manifolds, by unit-free manifolds i.e. line bundles over smooth manifolds. A series of theorems are proved in order to show that the category of line bundles admits a Hamiltonian functor into a category of contact manifolds, thus showing that, from a categorical perspective, the proposed generalization of configuration spaces produces a similar Hamiltonian phase space theory. Generalizations of contact manifolds to Jacobi, precontact and L-Dirac as well as the possibility of Jacobization functors are briefly discussed in relation to the Hamiltonian functor.

\section{Introducing Physical Dimension into Geometric Mechanics} \label{DimensionGeometricMechanics}

The introduction and discussion around the concept of measurand spaces in Section \ref{MeasurandFormalism} in the context of metrology established lines, i.e. 1-dimensional real vector spaces, as the natural objects that abstractly represent physical quantities in a truly unit-free manner. On the other hand, line bundles were argued in Section \ref{UnitFreeManifolds} to encapsulate the notion of unit-free free manifolds. Then, if one wishes to see physical dimension being incorporated into the framework of classical mechanics in an organic way, it follows from these considerations that an obvious possibility is to take phase spaces to be unit-free manifolds. Although this generalization will be first motivated from physical principles in this section, where we introduce line bundles as the enhanced version of configuration spaces, the degree of its success will be measured by whether we can recover a canonical theory of phase spaces via a Hamiltonian functor and whether dynamical equations recover the equations of motion of conventional classical mechanics. In Section \ref{CanonicalContact} below we prove that a Hamiltonian functor indeed exists and in chapter \ref{SpeculationsDynIntQuant} we discuss how the standard Hamilton's equations are recovered in the context of line bundle phase spaces.\newline

Following the principles of Realism, Characterization, Kinematics and Observation of Section \ref{FoundationsHamiltonianPhaseSpace} and guided by the intuition that the characteristics to be measured of a physical system should be implemented by elements of a measurand space, we are compelled to define the \textbf{category of unit-free configuration spaces} as the category of line bundles over smooth manifolds $\Line_\Man$. As usual in this thesis, objects in this category will be indistinctly denoted by $\lambda:L\to Q$, $L_Q$ or $L$. Having all the technical results of Section \ref{CategoryOfLineBundles} at hand, we give a physical interpretation for the categorical structure naturally present in $\Line_\Man$ in direct analogy with the list of physical intuitions behind the category of conventional configuration spaces discussed at the beginning of Section \ref{CanonicalSymplectic}:
\begin{itemize}
    \item The physical interpretation of the points of the base manifold of a unit-free configuration space $L_Q$ is exactly the same as the points of a configuration space, they represent the \textbf{static states} of the physical system. In this sense, the space of static states of a physical system is independent of the particular dimensions of the physical quantities that will be measured from it, as is indeed the case in conventional geometric mechanics.
    \item The measurable static properties of the same kind of a physical system are identified with the smooth sections of some fixed line bundle $\Sec{L_Q}$. The line bundle $L_Q$ will be appropriately called the \textbf{configuration space $Q$ of dimension $L$}. We thus identify the collection of all possible measurable properties of a fixed physical dimension with a choice of unit-free configuration space $L_Q$. We call these the \textbf{static observables of dimension $L$} of the configuration space $Q$ and denote them by $\Obs{L_Q}:=\Sec{L_Q}$. Properties of factor pull-backs in the category of line bundles then ensure that we have an observable contravariant functor $\text{Obs}:\Line_\Man\to\textsf{RMod}_1$.
    \item We recover the notion of ``unit-less'' observable of ordinary configuration spaces via the notion of local unit, which we recall is a local non-vanishing section $u\in\Sec{L_Q^\times}$ on some open subset $U\subset Q$. Restricting to the open subset $U$ allows us to see any other static observable $s\in\Obs{L_Q}$, i.e. the restriction of some arbitrary section to the open subset, as a local real-valued function $\Tilde{s}_u$ determined uniquely by the equation $s=\Tilde{s}_u\cdot u$. Thus, a local unit $u$ allows (locally) for a functorial assignment of the form
    \begin{equation*}
        u:\Obs{L_Q}\to \Obs{Q}.
    \end{equation*}
    \item A subsystem is characterised by restricting possible positions of a larger system, that is, by an inclusion of an embedded submanifold $i:S\to Q$. Our discussion about submanifolds of line bundles and embedding factors in Section \ref{CategoryOfLineBundles} ensures that this situation is indeed equivalent to the monomorphisms of the category of unit-free configuration spaces.
    \item Physically-indistinguishable static states must produce measurement outcomes that are indistinguishable as elements of a line bundle over the configuration space and thus should give a basic quotient of line bundles. Our discussion about submersion factors in Section \ref{CategoryOfLineBundles} shows that this notion precisely corresponds to the epimorphisms of the category of unit-free configuration spaces.
    \item Given two unit-free configuration spaces $L_{Q_1}$ and $L_{Q_2}$, the line product construction of Section \ref{CategoryOfLineBundles} gives the direct analogue of the Cartesian product of conventional configuration spaces. We thus regard $L_{Q_1}\utimes L_{Q_2}$ as the categorical implementation of the Principle of Combination for unit-free configuration spaces.
    \item Paths of a physical system, i.e. temporal series of static states, are simply recovered as smooth curves on the base space of a unit-free configuration space. In this manner, conventional dynamics $\Dyn{Q}$ are simply recovered as the vector fields on the base manifold. However, the extra structure introduced by the presence of the line bundle induces a new dynamical aspect of configuration spaces. Given a unit-free configuration space $L_Q$, all the non-zero fibre elements over a point $L_q$ represent different choices of unit for the same type of physical quantity. Then, any measurement performed on a system moving along a path $c(t)$ passing through $q$ at a time $t_0$ will have to be unit-compatible with any measurement performed at a later time $t>t_0$. This means that the choice of unit should be preserved along the motion of a path. Considering a unit $u$ as a local section, this is ensured locally by construction, however, taking all the possible arbitrary choices of local unit around the point $q$, forces the existence of a $1$-parameter family of fibre-wise isomorphisms covering the smooth curve $c(t)$. These are nothing but smooth families of line bundle automorphisms, which are given infinitesimally by line bundle derivations, and thus we identify the \textbf{unit-free dynamics} of a unit-free configuration space $L_Q$ as the der bundle of the line bundle $\Dyn{L_Q}:=\Der L_Q$. The anchor of the der bundle $\delta:\Der L_Q\to \Tan Q$ allows to connect dimensioned dynamics with ordinary dynamics via the surjective map
    \begin{equation*}
        \delta:\Dyn{L_Q}\to \Dyn{Q}.
    \end{equation*}
\end{itemize}

Similarly to the case of ordinary configuration spaces, we thus see that unit-free configuration spaces provide an example of a theory of phase spaces where the observable functor is the assignment of sections of line bundles $\text{Obs}=\Gamma :\Line_\Man \to \textsf{RMod}_1$, the dynamics functor is the der functor $\text{Dyn}=\Der : \Line_\Man\to \textsf{Jacb}_\Man$ and the evolution functor is given by taking the vector bundle of module derivations of the spaces of sections $\text{Evl}=\text{Der}:\textsf{RMod}_1\to \textsf{Jacb}_\Man$. The identification of module derivations as the sections of the der bundle
\begin{equation*}
    \Sec{\Der L}\cong \Dr{L},
\end{equation*}
implies that these functors fit in the following commutative diagram
\begin{equation*}
\begin{tikzcd}[row sep=small]
& \textsf{RMod}_1 \arrow[dd,"\text{Der}"] \\
\Line_\Man \arrow[ur,"\Gamma"] \arrow[dr,"\Der"'] & \\
& \textsf{Jacb}_\Man 
\end{tikzcd}
\end{equation*}
thus making the category of unit-free configuration spaces into a theory of phase spaces. Since no additional structure on the line bundle is specified by objects in the category of unit-free configuration spaces, again in analogy with the case of conventional configuration spaces, this theory of phase spaces is not Hamiltonian.

\section{Canonical Contact Phase Spaces and the Jet Functor} \label{CanonicalContact}

Our aim in this section is to construct a Hamiltonian functor for the category of unit-free configuration spaces hoping to find an analogous categorical structure to that described in Section \ref{CanonicalSymplectic} for the category of canonical symplectic phase spaces. The analogue of the cotangent construction in the case of line bundles will be given by the bundle of first jets, then, before we identify the Hamiltonian functor in full, we prove a series of results about the contact geometry of jet bundles.\newline

Firstly, we identify the \textbf{canonical contact structure} found in the jet bundles of general line bundles.

\begin{prop}[The Canonical Contact Manifold Associated to a Line Bundle]\label{JetCanonicalContactStructure}
Let $\lambda:L\to Q$ be a line bundle and $\normalfont \pi:\Jet^1L \to Q$ its jet bundle, then there is a canonical contact structure $\normalfont (\Jet^1 L,\text{H}_L)$ such that the contact line bundle is isomorphic to the pull-back of the line bundle on the base manifold:
\begin{equation*}
\normalfont
    \Tan (\Jet^1 L)/\text{H}_L\cong \pi^*L.
\end{equation*}
We denote this line bundle by $\normalfont L_{\Jet^1 L}$. Furthermore, the non-degenerate Jacobi structure induced in the line bundle $\normalfont (L_{\Jet^1 L},\{,\}_L)$ is fibre-wise linear and it is completely determined by the algebraic structure of derivations acting on sections
\begin{align*}
    \{l_a,l_b\}_L & =  l_{[a,b]}\\
    \{l_a,\pi^*s\}_L & =  \pi^* a [s]\\
    \{\pi^*s,\pi^*r\}_L & = 0
\end{align*}
for all $s,r\in\Sec{L}$, $\normalfont a,b\in\Sec{\Der L}$ and where $\normalfont l:\Sec{\Der L}\to \Sec{L_{\Jet^1 L}}$ is the inclusion of derivations as fibre-wise linear sections of $\normalfont L_{\Jet^1 L}$.
\end{prop}
\begin{proof}
The hyperplane distribution $\text{H}_L\subset \Tan (\Jet^1 L)$ is the usual Cartan distribution defined in general jet bundles of vector bundles, however, in our case it can be regarded as the kernel of the \textbf{canonical contact form} $\theta\in\Sec{\Tan^{*(\pi^*L)}(\Jet^1 L)}$ defined as the line bundle analogue of the Liouville $1$-form on the cotangent bundle. More precisely, denoting by $\varpi:\Jet^1 L\to L$ the surjective bundle map of the jet sequence, the canonical contact form is explicitly given at any point $j^1_xu\in\Jet^1 L$ of the jet bundle by
\begin{equation*}
    \theta_{j^1_xu}:=\Tan_{j^1_xu}\varpi - \Tan_x u\circ \Tan_{j^1_xu}\pi
\end{equation*}
where we note that the explicit use of $\Tan_xu$ is well-defined from the fact that $j^1u$ is defined as the equivalence class of all sections agreeing in value and tangent map at $x\in Q$. The map above is mapping tangent spaces of the vector bundles $\theta_{j^1_xu}:\Tan_{j^1_xu}\Jet^1 L\to \Tan_{u(x)}L$, in order to make it into a $L$-valued $1$-form we will use the fact that vertical subspaces of the total space of a vector bundle are canonically isomorphic to the fibres, $\Tan_{u(x)}^{\text{Vert}}L=\Ker{\Tan_{u(x)} \lambda}\cong L_x$ and that the image of $\theta_{j^1_xu}$ is always vertical
\begin{equation*}
    \Tan_{u(x)} \lambda\circ \theta_{j^1_xu}=\Tan_{j^1_xu}(\lambda\circ \varpi)-\Tan_{x}(\lambda\circ u)\circ \Tan_{j^1_xu}\pi=\Tan_{j^1_xu}\pi - \Tan_{j^1_xu}\pi =0,
\end{equation*}
where we have used the fact that sections and jet prolongations fit in the commutative diagram
\begin{equation*}
\begin{tikzcd}
\Jet^1 L \arrow[rr, "\varpi"] \arrow[dr, "\pi"] & & L \arrow[dl, "\lambda"'] \\
 & Q \arrow[ul,bend left,"j^1u"]\arrow[ur,bend right,"u"'] & 
\end{tikzcd}
\end{equation*}
We then define $\text{H}_L:=\Ker{\theta}$, which is shown to be a hyperplane distribution from simple point-wise dimension counting. Applying the first isomorphism theorem for vector spaces fibre-wise, we find $\Tan (\Jet^1 L)/\text{H}_L\cong \pi^*L$, as desired. The non-degenerate Jacobi structure on $\pi^*L$ appears as the $L$-dual to the Jacobi algebroid structure present in $\Der L$ analogously to the canonical symplectic structure being the linear Poisson structure on $\Cot Q$ dual to the Lie algebroid $\Tan Q$. This can be seen  explicitly via the natural inclusions of fibre-wise constant sections $\pi^*:\Sec{L}\hookrightarrow \Sec{\pi^*L}$ and fibre-wise linear sections $l:\Sec{\Der L}\hookrightarrow \Sec{\pi^*L}$ satisfying
\begin{equation*}
    \pi^*(f\cdot u)=\pi^*f\cdot \pi^*u \qquad l_{f\cdot a}=\pi^*f\cdot l_a
\end{equation*}
for all $f\in\Cin{Q}$, $u\in\Sec{L}$ and $a\in\Sec{\Der L}$, which uniquely determine a fibre-wise linear $\Real$-linear bracket $\{,\}$ determined by the defining identities in the proposition and whose Jacobi identity follows from the Jacobi identity of the commutator bracket of derivations. We can show that this bracket extends to a Jacobi structure by identifying a symbol and squiggle pair $(X,\Lambda)$ for an extension by symbol argument. A spanning subspace of functions for $\Cin{\Jet^1 L}$ is given by the natural inclusions of fibre-wise constant functions $\pi^*:\Cin{Q}\hookrightarrow \Cin{\Jet^1 L}$ and fibre-wise linear functions $l:\Sec{\Der L\otimes L^*}\hookrightarrow \Cin{\Jet^1 L}$. Then, we can explicitly define the symbol and squiggle acting on spanning sections and functions:
\begin{align*}
    X_{l_a}[l_{b\otimes \sigma}]&=l_{[a,b]\otimes \sigma}   & X_{l_a}[\pi^*f]&=\pi^*\delta(a)[f]\\
    X_{\pi^*s}[l_{b\otimes \sigma}]&=-\pi^*\sigma(b[s]) & X_{\pi^*s}[\pi^*f]&=0\\
    \Lambda^\sharp(dl_{a\otimes \rho}\otimes l_b)[l_{c\otimes \sigma}]&= l_{[a,c]\otimes \rho}l_{b\otimes \sigma} & \Lambda^\sharp(dl_{a\otimes \rho}\otimes l_b)[\pi^*f]&= \pi^*\delta(a)[f]l_{b\otimes \rho}\\
    \Lambda^\sharp(d\pi^*f\otimes l_a)[l_{b\otimes \rho}]&= -\pi^*\delta(b)[f]l_{a\otimes \sigma} & \Lambda^\sharp(d\pi^*f\otimes l_a)[\pi^*g]&= 0
\end{align*}
for all $f,g\in\Cin{Q}$, $s\in\Sec{L}$, $a,b,c\in\Sec{\Der L}$ and $\rho, \sigma\in\Sec{L^*}$. Direct computation using the compatibility formulas between basic and linear sections above, as well as the usual isomorphism $\text{End}(L)=\Real_Q$, shows that $(X,\Lambda)$ defined by these identities satisfy the extension by symbol axioms of proposition \ref{ExtensionBySymbol} and thus define a Jacobi structure on the line bundle $L_{\Jet^1 L}$.
\end{proof}

The following result establishes the interaction between the jet bundle construction and taking line products of line bundles.

\begin{thm}[Jet Bundle of a Line Product]\label{JetLineProduct}
Let $\lambda_1:L_1\to Q_1$, $\lambda_2:L_2\to Q_2$ two line bundles and denote their line product by $L_1\utimes L_2$, then there is a canonical factor $W$ giving the isomorphism of line bundles:
\begin{equation*}
\normalfont
\begin{tikzcd}
L_{\Jet^1 (L_1\utimes L_2)} \arrow[r, "W"] \arrow[d] & L_{\Jet^1 L_1} \utimes L_{\Jet^1 L_2} \arrow[d] \\
\Jet^1 (L_1\utimes L_2) \arrow[r, "w"'] & \Jet^1 L_1\dtimes \Jet^1 L_2
\end{tikzcd}
\end{equation*}
Furthermore, the factor $W$ is an isomorphism of Jacobi manifolds.
\end{thm}
\begin{proof}
Let us first construct the isomorphism factor $W$ explicitly. Recall from proposition \ref{DerLineProduct} that $\Der(L_1\utimes L_2)\cong \Der L_1 \boxplus \Der L_2$, the fact that the jet bundle is the L-dual of the der bundle in the category of L-vector bundles allows us to write:
\begin{equation*}
    \Jet^1 (L_1\utimes L_2)\cong(p_1^*\Der L_1\oplus p_2^*\Der L_2)^*\otimes (L_1\utimes L_2),
\end{equation*}
then, using the swapping isomorphism factor for the second term, we find the following isomorphism of vector bundles
\begin{equation*}
    t:\Jet^1 (L_1\utimes L_2)\to p_1^*\Jet^1 L_1\oplus p_2^*\Jet^1 L_2.
\end{equation*}
covering the identity map on $Q_1\dtimes Q_2$. Now, by the definition of line product, it is clear that we can define the following map
\begin{align*}
z: p_1^*\Jet^1 L_1\oplus p_2^*\Jet^1 L_2 & \to \Jet^1 L_1\dtimes \Jet^1L_2\\
(\alpha_{q_1}\oplus \beta_{q_2})_{B_{q_1q_2}} & \mapsto C_{\alpha_{q_1}\beta_{q_2}}
\end{align*}
where $C_{\alpha_{q_1}\beta_{q_2}}=B_{q_1q_2}$, which is well-defined since the line bundle over a jet bundle is defined as the pull-back line bundle from the base. Conversely, given a fibre-wise invertible map between pull-back line bundles we can project to a fibre-wise map between the base line bundles, denote this projection by $\overline{\pi}:\Jet^1 L_1 \dtimes \Jet^1 L_2 \to Q_1 \dtimes Q_2$, so there is an obvious inverse for the map above
\begin{equation*}
    z^{-1}(C_{\alpha_{q_1}\beta_{q_2}})=(\alpha_{q_1}\oplus \beta_{q_2})_{\overline{\pi}(C_{\alpha_{q_1}\beta_{q_2}})}.
\end{equation*}
We thus find the desired diffeomorphism
\begin{equation*}
    w:=z\circ t:\Jet^1 (L_1\utimes L_2) \to \Jet^1 L_1\dtimes \Jet^1L_2.
\end{equation*}
Let us write the line product commutative diagram for the line bundles over the jets as
\begin{equation*}
\begin{tikzcd}
L_{\Jet^1L_1} \arrow[d, "\mu_1"'] & L_{\Jet^1L_1}\utimes L_{\Jet^1L_2} \arrow[l,"R_1"']\arrow[d, "\mu_{12}"]\arrow[r,"R_2"] & L_{\Jet^1L_2} \arrow[d,"\mu_2"] \\
\Jet^1L_1 & \Jet^1L_1 \dtimes \Jet^1L_2 \arrow[l,"r_1"]\arrow[r,"r_2"'] & \Jet^1L_2
\end{tikzcd}
\end{equation*}
and denote the compositions $o_i:=\pi_i\circ r_i:\Jet^1L_1\dtimes \Jet^1L_2\to Q_i$. It then follows by construction that
\begin{equation*}
    L_{\Jet^1L_1}\utimes L_{\Jet^1L_2}=o_1^*L_1
\end{equation*}
and
\begin{equation*}
    L_{\Jet^1(L_1\utimes L_2)}\cong w^*o_1^*L_1
\end{equation*}
thus showing
\begin{equation*}
    L_{\Jet^1(L_1\utimes L_2)}\cong w^*(L_{\Jet^1L_1}\utimes L_{\Jet^1L_2}).
\end{equation*}
This is, of course, tantamount to there being a factor $W:L_{\Jet^1(L_1\utimes L_2)}\to L_{\Jet^1L_1}\utimes L_{\Jet^1L_2}$ covering the diffeomorphism $w$ that essentially acts as the fibre-wise identity on $L_1$. In order to show that this factor is, in fact, a Jacobi map we can write the pull-backs of (non-zero) brackets of spanning sections directly
\begin{align*}
    W^*\{R_i^*l_{a_i},R_i^*l_{b_i}\}_{12}&=W^*R_i^*\{l_{a_i},l_{b_i}\}_i=(R_i\circ W)^*l_{[a_i,b_i]}\\
    W^*\{R_i^*l_{a_i},R_i^*\pi_i^*u_i\}_{12}&=W^*R_i^*\{l_{a_i},\pi_i^*u_i\}_i=(R_i\circ W)^*\pi_i^*a_i[u_i]
\end{align*}
for $u_i\in\Sec{L_i}$, $a_i,b_i\in\Sec{\Der L_i}$, $i=1,2$, and where the definition of product Jacobi structure has been used to write the RHS expressions. Note that the construction of $W$ above is such that
\begin{equation*}
    (R_i\circ W)^*l_{a_i}=\overline{l}_{k_i(a_i)}\qquad (R_i\circ W)^*\pi_i^*u_i=\overline{\pi}^*P_i^*u_i
\end{equation*}
where $\overline{\pi}:\Jet^1(L_1\utimes L_2)\to Q_1\dtimes Q_2$ is the jet bundle projection for the line product, $k_i:\Sec{\Der L_i}\to \Sec{\Der (L_1\utimes L_2)}$ are the natural Lie algebra injections of derivations from proposition \ref{DerivationsLineProduct} and $\overline{l}:\Sec{\Der (L_1\utimes L_2)}\to L_{\Jet^1(L_1\utimes L_2)}$ is the inclusion of fibre-wise linear sections. A direct computation then gives the brackets of pull-backs of spanning sections
\begin{align*}
    \{W^*R_i^*l_{a_i},W^*R_i^*l_{b_i}\}&=\overline{l}_{[k_i(a_i),k_i(b_i)]}=(R_i\circ W)^*l_{[a_i,b_i]}\\
    \{W^*R_i^*l_{a_i},W^*R_i^*\pi_i^*u_i\}&=\overline{\pi}^*k_i(a_i)[P_i^*u_i]=\overline{\pi}^*P_i^*a_i[u_i]=(R_i\circ W)^*\pi_i^*a_i[u_i]
\end{align*}
which agree with the pull-backs of brackets above, thus completing the proof.
\end{proof}

Consider a factor between some pair of line bundles $B:L_1\to L_2$ covering a smooth map $b:Q_1\to Q_2$. Recall that the der map always gives a well-defined morphism of L-vector bundles $\Der B:\Der L_1\to \Der L_2$, however only factors covering diffeomorphisms can be dualized to give a well-defined map between the jet bundles. We can, nevertheless, define the \textbf{jet lift} of the factor $B$ as the following subamanifold of the base product of jet bundles
\begin{equation*}
    \Jet^1 B:=\{C_{\alpha_{q_1}\beta_{q_2}}\in \Jet^1 L_1\dtimes \Jet^1 L_2|\quad q_2=b(q_1),\quad (\Der_{q_1}B)^{*L}\beta_{q_2}=\alpha_{q_1}, \quad C_{\alpha_{q_1}\beta_{q_2}}=B_{q_1}\}
\end{equation*}
where the last condition is understood using the fact that the line bundles over the jet bundles are simply the pull-back bundles by the canonical projection to the base. Note that when the factor $B$ covers a diffeomorphism, the L-dual of the der map $\Der B$ is a well-defined isomorphism of L-vector bundles and so the jet lift becomes simply the graph $\Jet^1 B=\LGrph{(\Der B)^{L_2*L_1}}$. We now show that the jet lift of a factor is, in fact, a Legendrian submanifold of the product Jacobi manifold, thus defining a Legendrian relation.

\begin{prop}[Jet Lift of a Factor]\label{JetLiftFactor}
Let $B:L_1\to L_2$ be a factor between two line bundles, then its jet lift
\begin{equation*}
\normalfont
    \Jet^1 B\subset \Jet^1 L_1 \dtimes \Jet^1 L_2
\end{equation*}
is a Legendrian submanifold of the base of the Jacobi product $\normalfont L_{\Jet^1 L_1} \utimes \overline{L}_{\Jet^1 L_2}$.
\end{prop}
\begin{proof}
We first show that $\Jet^1B$ is a coisotropic submanifold by identifying a set of generating sections of its vanishing submodule $\Gamma_{\Jet^1B}\subseteq \Sec{L_{\Jet^1(L_1\utimes L_2)}}$ and showing that they indeed form a Lie subalgebra. Firstly, given $u_i\in\Sec{L_i}$, $a_i\in\Sec{\Der L_i}$, $i=1,2$, let us define the following sections of the line product of the canonical contact line bundles on the jet bundles $\Jet^1L_i$ following the notation introduced in proposition \ref{JetLineProduct}:
\begin{align*}
    l_{a_1a_2}&:=R_1^*l_{a_1}-R_2^*l_{a_2}\\
    \pi^*_{u_1u_2}&:= R_1^*\pi_1^*u_1-R_2^*\pi_2^*u_2.
\end{align*}
Consider a point on the jet lift $C_{\alpha_{q_1}\beta_{q_2}}\in \Jet^1 B$ so that $q_2=b(q_1)$, $(\Der_{q_1}B)^{*L}\beta_{q_2}=\alpha_{q_1}$ and $C_{\alpha_{q_1}\beta_{q_2}}=B_{q_1}$. Let us evaluate the defined sections on it:
\begin{align*}
l_{a_1a_2}(C_{\alpha_{q_1}\beta_{q_2}}) & = R_1|^{-1}_{C_{\alpha_{q_1}\beta_{q_2}}}l_{a_1}(\alpha_{q_1})-R_2|^{-1}_{C_{\alpha_{q_1}\beta_{q_2}}}l_{a_2}(\beta_{q_2})\\
& = l_{a_1}(\alpha_{q_1})-B^{-1}_{q_1}(l_{a_2}(\beta_{q_2}))\\
& = B^{-1}_{q_1}\beta_{q_2}(\Der_{q_1}B(a_1|_{q_1})-a_2|_{q_2})\\
& = B^{-1}_{q_1}\beta_{q_2}(\Der B \circ a_1 - a_2\circ b)(q_1)
\end{align*}
and
\begin{align*}
\pi^*_{u_1u_2} & = R_1|^{-1}_{C_{\alpha_{q_1}\beta_{q_2}}}u_1(q_1)-R_2|^{-1}_{C_{\alpha_{q_1}\beta_{q_2}}}u_2(q_2)\\
& = u_1(q_1)-B_{q_1}^{-1}u_2(q_2)\\
& = u_1(q_1)-B_{q_1}^{-1}u_2(b(q_1))\\
& = (u_1-B^*u_2)(q_2).
\end{align*}
It is then obvious that these spanning sections vanish on the jet lift $\Jet^1 B$ iff the derivations are $B$-related and the sections are mapped by the pull-back $B^*$, i.e.
\begin{align*}
    \Gamma_{\Jet^1 B}\ni l_{a_1a_2} &\Leftrightarrow   a_1\sim_B a_2\\
    \Gamma_{\Jet^1 B}\ni \pi^*_{u_1u_2}  &\Leftrightarrow  u_1=B^*u_2.
\end{align*}
These are the generating vanishing sections so it will suffice to check that evaluations of brackets among these on an arbitrary point of the jet lift $C_{\alpha_{q_1}\beta_{q_2}}\in \Jet^1 B$ vanish. From the defining conditions of a the linear Jacobi we see that $\{\pi^*_{u_1u_2},\pi^*_{u_1'u_2'}\}=0$, so we are left with the two other possible brackets. For the bracket of fibre-wise linear sections we compute explicitly using the defining properties of the product Jacobi bracket:
\begin{equation*}
    \{l_{a_1a_2} ,l_{a_1'a_2'} \}= R^*_1\{l_{a_1},l_{a_1'}\}_1-R^*_2\{l_{a_2},l_{a_2'}\}_2=R^*_1l_{[a_1,a_1']}-R^*_2l_{[a_2,a_2']}
\end{equation*}
for $a_i,a_i'\in\Sec{\Der L_i}$, $i=1,2$. Then if $l_{a_1a_2}$ and $l_{a_1'a_2'}$ are vanishing sections, the derivations are $B$-related and, by virtue of proposition \ref{DerFunctor}, where $\Der B$ is shown to be a Lie algebroid morphism, the brackets are also $B$-related $[a_1,a_1']\sim_B [a_2,a_2']$ making the expression above into a vanishing section. Only the cross bracket $\{l_{a_1a_2} ,\pi^*_{u_1u_2} \}= R^*_1\pi_1^*a_1[u_1]_1-R^*_2\pi_2^*a_2[u_2]$ remains, for which we evaluate on a point of the jet lift and show it vanishes by directly computing using $u_1=B^*u_2$ and $a_1\sim_B a_2$:
\begin{align*}
\{l_{a_1a_2} ,\pi^*_{u_1u_2} \}(C_{\alpha_{q_1}\beta_{q_2}}) & = a_1|_{q_1}(u_1)-B_{q_1}^{-1}a_2|_{q_2}(u_2)\\
& = a_1|_{q_1}(u_1)-B_{q_1}^{-1}\Der_{q_1}B(a_1|_{q_1})(u_2)\\
& = a_1|_{q_1}(B^*u_2)-B_{q_1}^{-1}B_{q_1}a_1|_{q_1}(B^*u_2)\\
& = a_1|_{q_1}(B^*u_2)-a_1|_{q_1}(B^*u_2)\\
& = 0.
\end{align*}
To show that $\Jet^1 B$ is maximally isotropic we simply do a dimension count. By definition, it is clear that $\dimm \Jet^1 B=\dimm Q_1\dtimes Q_2$, but it follows from proposition \ref{JetLineProduct} that $\dimm \Jet^1 L_1\dtimes \Jet^1 L_2 = 2(\dimm Q_1\dtimes Q_2)+1$, then clearly $\Jet^1 B\subset \Jet^1 L_1\dtimes \Jet^1 L_2$ is maximally isotropic and thus Legendrian.
\end{proof}

With this last proposition at hand, we can now see the jet bundle construction as a functor from the category of line bundles.

\begin{thm}[The Jet Functor] \label{JetFunctor}
The assignment of jet bundles to line bundles is a contravariant functor
\begin{equation*}
\normalfont
    \Jet^1 : \Line_\Man \to \Cont_\Man.
\end{equation*}
\end{thm}
\begin{proof}
A factor $B:L_1\to L_2$ gives a Legendrian relation $\Jet^1B:\Jet^1L_2\dashrightarrow \Jet^1L_1$ in virtue of proposition \ref{JetLiftFactor}, then it only remains to check functoriality. It is obvious by definition that the identity factor $\Id_L:L\to L$ gives the L-diagonal relation $\text{Ldiag}(L_{\Jet^1L})\subset \Jet^1 L\dtimes \Jet^1L$, where the L-diagonal is the natural subset of fibre-wise identity maps in a line product $L\utimes L$. Consider now two factors $B:L_1\to L_2$ and $B':L_2\to L_3$ covering the smooth maps $b:Q_1\to Q_2$ and $b':Q_2\to Q_3$. By definition of composition of relations, we find
\begin{equation*}
    \Jet^1 B\circ \Jet^1B':=\{C_{\alpha_{q_1}\gamma_{q_3}}\in \Jet^1 L_1\dtimes \Jet^1 L_3| q_3=b'(b(q_1)), (\Der_{q_1}B)^{L_2*L_1}(\Der_{b(q_1)}B')^{L_3*L_2}\gamma_{q_3}=\alpha_{q_1}, C_{\alpha_{q_1}\gamma_{q_3}}=B_{q_1}\circ B_{b(q_1)}'\}
\end{equation*}
Then, it follows from the fact that L-duality is a contravariant autofunctor of L-vector spaces, as shown in Section \ref{CategoryOfLVectorSpaces}, that the jet lift of factors is contravariant with respect to composition of relations
\begin{equation*}
    \Jet^1 B\circ \Jet^1B'=\Jet^1(B'\circ B).
\end{equation*}
\end{proof}

There are two alternative ways to regard the jet bundle of the line bundle induced on a submanifold of the base. On the one had, we could take the line bundle as a an embedded subbundle and construct its jet bundle from the ambient jet bundle. On the other, we could simply regard the restricted line bundle as an intrinsic line bundle and canonically construct its jet bundle. We now prove that these two constructions are equivalent.

\begin{thm}[Canonical Coisotropic Reduction in Jet Bundles]\label{JetCoisotropicReduction}
Let $i:S\hookrightarrow Q$ be a submanifold of a line bundle $L$, then the restriction of the ambient jet bundle to the submanifold $\normalfont (\Jet^1 L)|_S$ is a coisotropic submanifold with respect to the canonical contact structure on $\normalfont \Jet^1 L$. Furthermore, there is a submersion factor covering the surjective submersion $\normalfont z:(\Jet^1 L)|_S\twoheadrightarrow \Jet^1 L_S$ given by the fibre-wise quotient:
\begin{equation*}
\normalfont
    \Jet^1_qL/(\Der L_S)^{0L}\cong \Jet^1_q L_S\qquad q\in S
\end{equation*}
so that the canonical contact structure on $\normalfont \Jet^1 L$ Jacobi reduces to the canonical contact structure on $\normalfont \Jet^1 L_S$.
\end{thm}
\begin{proof}
Let us first prove that the submersion $ z:(\Jet^1 L)|_S\twoheadrightarrow \Jet^1 L_S$ fits in a Jacobi reduction scheme. It is a direct implication of proposition \ref{DerBundleLineSubmanifold} and the basic properties of L-vector spaces applied fibre-wise, that we have the following diagram of line bundle morphisms
\begin{equation*}
\begin{tikzcd}[sep=tiny]
L_{\Jet^1L|_{S}} \arrow[rr,"\iota"] \arrow[dd,"\zeta"'] \arrow[dr]& & L_{\Jet^1L} \arrow[dr] & \\
& \Jet^1L|_S \arrow[rr,"i"', hook] \arrow[dd, "z",twoheadrightarrow] & & \Jet^1L \\
L_{\Jet^1L_S} \arrow[dr] & & & \\
 & \Jet^1L_S &  & 
\end{tikzcd}
\end{equation*}
where $\iota$ denotes, abusing notation, the embedding factor induced by the submanifold $i:S\hookrightarrow Q$ and $\zeta$ is defined as the fibre-wise identity of the pull-back line bundles covering the point-wise linear submersion
\begin{equation*}
    z_q:\Jet^1_q L\twoheadrightarrow \Jet^1_qL/(\Der L_S)^{0L}\cong \Jet^1_q L_S, \qquad q\in S.
\end{equation*}
Following proposition \ref{DerivationsLineSubmanifold}, we find natural isomorphisms $\Sec{L_S}\cong\Sec{L}/\Gamma_S$ and $\Dr{L_S}\cong\text{Der}_S(L)/\text{Der}_{0S}(L)$, and denoting the natural inclusions of spanning sections on $L_{\Jet^1 L_S}$ by $\overline{l}$ and $\overline{\pi}$, we can write the fibre-wise linear Jacobi structure on $\Jet^1L_S$ equivalently as
\begin{align*}
    \{\overline{l}_{\overline{a}},\overline{l}_{\overline{b}}\}_S & =  \overline{l}_{\overline{[a,b]}}\\
    \{\overline{l}_{\overline{a}},\overline{\pi}^*\overline{u}\}_S & =  \overline{\pi}^*\overline{a[u]}\\
    \{\overline{\pi}^*\overline{u},\overline{\pi}^*\overline{v}\}_S & = 0
\end{align*}
with $\overline{u},\overline{v}\in\Sec{L}/\Gamma_S$ and $\overline{a},\overline{b}\in\text{Der}_S(L)/\text{Der}_{0S}(L)$, which is well-defined precisely from the description of derivations as a subquotient of Lie algebras. The submersion factor $\zeta:L_{\Jet^1L|_{S}}\to L_{\Jet^1L_S}$ covering the quotient map $z:\Jet^1L|_{S}\to \Jet^1L_S$ has been defined such that it is the point-wise counterpart to the isomorphisms used above to rewrite the linear Jacobi bracket. Pull-backs via these factors satisfys the following identities by construction
\begin{equation*}
    \zeta^*\overline{\pi}^*\overline{u} =\iota^*\pi^*u \qquad \zeta^*\overline{l}_{\overline{a}} = \iota^* l_{a}
\end{equation*}
for all $u\in\Sec{L}$ and $a\in\text{Der}_S(L)$. This now clearly implies the reduction condition for all spanning sections
\begin{align*}
    \zeta^*\{\overline{l}_{\overline{a}},\overline{l}_{\overline{b}}\}_S & = \iota^* l_{[a,b]}=\iota^*\{l_a,l_b\}_L\\
    \zeta^*\{\overline{l}_{\overline{a}},\overline{\pi}^*\overline{u}\}_S & = \iota^* \pi^*a[u]=\iota^* \{l_a,\pi^*u\}_L\\
    \zeta^*\{\overline{\pi}^*\overline{u},\overline{\pi}^*\overline{v}\}_S & = 0 = \iota^* \{\pi^*u,\pi^*v\}_L
\end{align*}
thus showing that the linear Jacobi $(L_{\Jet^1L},\{,\}_L)$ reduces to $(L_{\Jet^1L_S},\{,\}_S)$. Lastly, it is easy to see that the vanishing sections of $\Jet^1L|_S$ seen as a submanifold of the jet bundle are precisely those of the form $l_{\text{Der}_{0S}(L)}$ and $\pi^*\Gamma_S$. It follows again from proposition \ref{DerivationsLineSubmanifold} that these form a Lie subalgebra of the linear Jacobi structure $(L_{\Jet^1L},\{,\}_L)$, thus explicitly showing that $\Jet^1L|_S\subseteq \Jet^1L$ is a coisotropic submanifold.
\end{proof}

Consider now a line bundle action $G\Acts L$. Recall that, in the case of a free and proper action, the orbit space is canonically a line bundle, denoted by $L/G$, and that there is a natural submersion factor $\sigma:L\to L/G$. The following proposition relates the canonical contact structures associated to these two line bundles.

\begin{thm}[Canonical Hamiltonian Reduction in Jet Bundles]\label{JetHamiltonianReduction}
Let $\Phi:G\times L\to L$ be a free and proper line bundle action of a connected Lie group $G$ on $L$, then the canonical contact structure on $\normalfont \Jet^1 L$ Jacobi reduces to the canonical contact structure on $\normalfont \Jet^1 (L/G)$. This reduction is, in fact, Hamiltonian: the jet lift of the line bundle action $\normalfont G\Acts \Jet^1 L$ preserves the canonical contact structure and has a natural comoment map given by
\begin{align*}
\overline{\mu}: \mathfrak{g} &\normalfont \to \Sec{L_{\Jet^1 L}}\\
\xi & \mapsto l_{\Psi(\xi)},
\end{align*}
where $\normalfont \Psi:\mathfrak{g}\to \Dr{L}$ is the infinitesimal line bundle action and $\normalfont l:\Dr{L}\to \Sec{L_{\Jet^1 L}}$ is the natural inclusion of derivations as fibre-wise linear sections on the jet bundle.
\end{thm}
\begin{proof}
Before we address the specific case of a group action, let us discuss the jet lift of a general diffeomorphic factor, i.e. a fibre-wise invertible line bundle morphism $B:L\to L$ covering a diffeomorphism $b:Q\to Q$. The \textbf{jet lift} of $B$ is a factor of the canonical contact line bundle on the jet bundle $L_{\Jet^1 L}$ defined by
\begin{equation*}
    \Jet^1 B:(\alpha_q,l_q)\mapsto (\Jet^1 b(\alpha_q),B^{-1}_{q}l_q)
\end{equation*}
where
\begin{equation*}
    \Jet^1 b(\alpha_q)=(\Der_qB)^{L_q*L_{b^{-1}(q)}}(\alpha_q)\in\Jet_{b^{-1}(q)}^1L.
\end{equation*}
The cotangent lift then induces the following commutative diagram of vector bundle morphisms
\begin{equation*}
    \begin{tikzcd}[sep=small]
     L_{\Jet^1L}\arrow[r,"\Jet^1B"]\arrow[d] & L_{\Jet^1L}\arrow[d] \\
     \Jet^1L\arrow[r,"\Jet^1 b"] \arrow[d] & \Jet^1L \arrow[d] \\
     Q\arrow[r,"b^{-1}"] & Q
    \end{tikzcd}
\end{equation*}
For any two diffeomorphic factors $B,F:L\to L$ and the identity factor $\Id_L:L\to L$ it is a simple check to show that
\begin{equation*}
    \Jet^1(B\circ F)=\Jet^1F\circ \Jet^1B, \qquad \Jet^1(\Id_L)=\Id_{L_{\Jet^1L}}.
\end{equation*}
The spanning sections $l_a,\pi^*u\in\Sec{L_{\Jet^1L}}$ transform under pull-back by a jet lift of a diffeomorphic factor according to the following expressions
\begin{equation*}
    \Jet^1B^*l_a=l_{B_*a}, \qquad \Jet^1B^*\pi^*u=\pi^*(B^{-1})^*u,
\end{equation*}
then it follows that the jet lift of a diffeomorphic factor $\Jet^1B$ is indeed a Jacobi map of the canonical contact structure on the jet bundle:
\begin{align*}
    \Jet^1B^*\{l_a,l_b\}_L & = \Jet^1B^*l_{[a,b]}=l_{[B_*a,B_*b]}=\{l_{B_*a},l_{B_*b}\}_L=\{\Jet^1B^*l_a,\Jet^1B^*l_b\}_L \\
    \Jet^1B^*\{l_a,\pi^*u\}_L & = \Jet^1B^*\pi^*a[u]=\pi^*(B_*a)[(B^{-1})^*u]=\{l_{B_*a},\pi^*(B^{-1})^*u\}_L=\{ \Jet^1B^*l_a, \Jet^1B^* \pi^*u\}_L \\
    \Jet^1B^*\{\pi^*u,\pi^*v\}_L & = 0 = \{\pi^*(B^{-1})^*u,\pi^*(B^{-1})^*v\}_L=\{\Jet^1B^*\pi^*u,\Jet^1B^*\pi^*v\}_L.
\end{align*}
The jet lift of the group action $G\Acts \Jet^1 L$ is defined by the jet lifts of the diffeomorphic factors corresponding to each group element
\begin{equation*}
    (\Jet^1 \Phi)_g:=\Jet^1\Phi_g,
\end{equation*}
which, in light of the above results for general jet lifts of diffeomorphic factors, is readily checked to be a group action that acts via Jacobi maps. This is a Hamiltonian action with comoment map simply given by $\overline{\mu}:=l\circ \Psi :\mathfrak{g}\to \Sec{L_{\Jet^1L}}$. Observe that the zero locus of the moment map is naturally identified with the L-annihilator of the subspace of derivations spanned by the infinitesimal generators regarded as a subbundle of the jet bundle
\begin{equation*}
    \mu^{-1}(0)=\Psi(\mathfrak{g})^{0L}\subseteq \Jet^1L.
\end{equation*}
Note that $G$-equivariance of the infinitesimal action $\Psi$ implies that the jet lifted action of $G$ restricts to a $G$-action on $\Psi(\mathfrak{g})^{0L}$, indeed we check for any $j_q^1u\in\Psi(\mathfrak{g})^{0L}$ and $\xi\in\mathfrak{g}$
\begin{equation*}
    \Jet^1\phi_g(j_q^1 u)(\Psi(\xi))=\Psi(\xi)[\Phi^*_gu]=(\Phi_g)_{b^{-1}(q)}\Der_q\Phi_g\Psi(\xi)=(\Phi_g)_{b^{-1}(q)}j_q^1u(\Psi(\text{Ad}_g(\xi))=0.
\end{equation*}
Using Proposition \ref{DerBundleGroupAction} and simple linear algebra of L-vector bundles we find the following point-wise isomorphism
\begin{equation*}
    \Jet^1_{[q]}(L/G):=(\Der_{[q]}(L/G))^{*L/G}\cong(\Der_q/\Psi(\mathfrak{g})_q)^{*L_q}\cong \Psi(\mathfrak{g})_q^{0L_q}.
\end{equation*}
This allows us to write the following factor reduction diagram
\begin{equation*}
\begin{tikzcd}[sep=tiny]
L_{\Psi(\mathfrak{g})^{0L}} \arrow[rr,"\iota"] \arrow[dd,"\zeta"'] \arrow[dr]& & L_{\Jet^1L} \arrow[dr] & \\
& \Psi(\mathfrak{g})^{0L} \arrow[rr,"i"', hook] \arrow[dd, "z",twoheadrightarrow] & & \Jet^1L \\
L_{\Jet^1(L/G)} \arrow[dr] & & & \\
 & \Jet^1(L/G) &  & 
\end{tikzcd}
\end{equation*}
where $\iota$ is the embedding factor for the L-annihilator $\Psi(\mathfrak{g})^{0L}$ seen as a subbundle (submanifold) of the jet bundle and $\zeta$ is the submersion factor induced by the jet lifted action restricted $\Psi(\mathfrak{g})^{0L}$, which is clearly free and proper. Since $G$ is connected, recall that the sections and derivations of the quotient line bundle can be equivalently regarded as
\begin{equation*}
    \Sec{L/G}\cong \Sec{L}^\mathfrak{g}, \qquad \Dr{L/G}\cong \Dr{L}^{\mathfrak{g}}/\Psi(\mathfrak{g}),
\end{equation*}
thus the linear Jacobi structure determined by the spanning sections $\overline{l}_{\overline{a}},\overline{\pi}^*\overline{u}\in\Sec{L_{\Jet^1(L/G)}}$ can be fully characterised under these isomorphisms. The factors constructed above are such that we have the following explicit treatment of extensions of spanning sections
\begin{align*}
    \zeta^*\overline{\pi}^*\overline{u}=\iota^*\pi^*u \quad &\Leftrightarrow \quad u\in\Sec{L}^\mathfrak{g}\\
    \zeta^*\overline{l}_{\overline{a}}=\iota^*l_a \quad &\Leftrightarrow \quad \overline{a}=a+\Psi(\mathfrak{g}), a\in\Dr{L}^\mathfrak{g}
\end{align*}
The reduction condition is now easily checked
\begin{align*}
    \zeta^*\{\overline{l}_{\overline{a}},\overline{l}_{\overline{b}}\}_{L/G} & = \iota^* l_{[a,b]}=\iota^*\{l_a,l_b\}_L\\
    \zeta^*\{\overline{l}_{\overline{a}},\overline{\pi}^*\overline{u}\}_{L/G} & = \iota^* \pi^*a[u]=\iota^* \{l_a,\pi^*u\}_L\\
    \zeta^*\{\overline{\pi}^*\overline{u},\overline{\pi}^*\overline{v}\}_{L/G} & = 0 = \iota^* \{\pi^*u,\pi^*v\}_L
\end{align*}
for all spanning sections $\overline{l}_{\overline{a}},\overline{l}_{\overline{b}},\overline{\pi}^*\overline{u},\overline{\pi}^*\overline{v}\in\Sec{L_{\Jet^1(L/G)}}$ and extensions $l_a,l_b,\pi^*u,\pi^*v\in\Sec{L_{\Jet^1L}}$, thus concluding the proof.
\end{proof}

We are now in the position to define \textbf{the category of canonical contact phase spaces} as the image of the category of unit-free configuration spaces under the jet functor, $\Jet^1(\Line_\Man)$, which, together with the same notions of observable functor $\text{Obs}=\Gamma: \Cont_\Man\to \textsf{RMod}_1$, dynamics functor $\text{Dyn}=\Der:\Cont_\Man\to \textsf{Jacb}_\Man$ and evolution functor $\text{Evl}=\text{Der}:\textsf{RMod}_1\to \textsf{Jacb}_\Man$, forms a theory of phase spaces. Clearly, the presence of the canonical contact structure on jet bundles ensures that $\Cont_\Man$ is, furthermore, a Hamiltonian theory of phase spaces, with Hamiltonian maps given by the Hamiltonian derivation of the non-degenerate Jacobi structures:
\begin{equation*}
    \eta_L:=\Lambda_L^\sharp\circ j^1=\text{ad}_{\{,\}_L}:\Obs{\Jet^1 L}\to \Dyn{\Jet^1 L}.
\end{equation*}

The results proved in this section can be encapsulated in the \textbf{Hamiltonian functor} for unit-free configuration spaces, which is given by the following list of correspondences:

\begin{center}
\begin{tabular}{ c c c }
Unit-Free Configuration Spaces & Hamiltonian Functor  & Unit-Free Phase Spaces \\
\hline
 $\Line_\Man$ & $\begin{tikzcd}\phantom{A} \arrow[r,"\mathbb{H}_{\Line_\Man}"] & \phantom{B} \end{tikzcd}$ & $\Cont_\Man$ \\ 
 $L_Q$ & $\begin{tikzcd} \phantom{Q} \arrow[r, "\Jet^1 ",mapsto] & \phantom{Q} \end{tikzcd}$ & $(\Jet^1 L_Q,\text{H}_{L_Q})$ \\ 
 $\Obs{L_Q}$ & $\begin{tikzcd} \phantom{Q} \arrow[r,hookrightarrow, "\pi^*"] & \phantom{Q}  \end{tikzcd}$ & $\Obs{L_{\Jet^1 L_Q}}$ \\
 $\Dyn{L_Q}$ & $\begin{tikzcd} \phantom{Q} \arrow[r,hookrightarrow, "l"] & \phantom{Q}  \end{tikzcd}$ & $\Obs{L_{\Jet^1 L_Q}}$ \\
 $L_1\utimes L_2$ & $\begin{tikzcd} \phantom{Q} \arrow[r,mapsto, "\Jet^1 "] & \phantom{Q}  \end{tikzcd}$ & $\Jet^1 L_1\dtimes \Jet^1 L_2$ \\
 $B:L_1\to L_2$ & $\begin{tikzcd} \phantom{Q} \arrow[r,mapsto, "\Jet^1 "] & \phantom{Q}  \end{tikzcd}$ & $\Jet^1 B\subset \Jet^1 L_1\dtimes \Jet^1 L_2$ Legendrian\\
 $\iota:L_S\hookrightarrow L_Q$ & $\begin{tikzcd} \phantom{Q} \arrow[r,mapsto, "\Jet^1"] & \phantom{Q}  \end{tikzcd}$ & $(\Jet^1 L_Q)|_S\subset \Jet^1 L_Q$ coisotropic \\
 $G\Acts L$ & $\begin{tikzcd} \phantom{Q} \arrow[r,mapsto, "\Jet^1"] & \phantom{Q}  \end{tikzcd}$ & $G\Acts \Jet^1 L$ Hamiltonian action \\
\end{tabular}
\end{center}
This correspondence is clearly categorically analogous to the Hamiltonian functor for ordinary configuration spaces of Section \ref{CanonicalSymplectic}, hence, canonical contact manifolds appear as valid unit-free generalizations of conventional, ``unit-less'', phase spaces. Furthermore, a similar diagram connecting the category of unit-free configuration spaces and the category of canonical contact phase spaces as theories of phase spaces is given by the jet functor:
\begin{equation*}
\begin{tikzcd}[row sep=small]
 & \textsf{RMod}_1 \arrow[dd,"\text{Der}"']& & & \textsf{LocLieAlg} \arrow[dd,"\text{Der}"] \\
 & & \Line_\Man \arrow[ul, "\Gamma"']\arrow[dl,"\Der"]\arrow[r,"\Jet^1"] & \Cont_\Man \arrow[ur,"\Gamma"] \arrow[dr,"\Der"'] & \\
 & \textsf{Jacb}_\Man & & & \textsf{Jacb}_\Man 
\end{tikzcd}
\end{equation*}

Similarly to ordinary symplectic phase spaces, the notion of \textbf{L-energy} as a choice of observable $h\in\Sec{L_{\Jet^1 L}}$ to determine the dynamics of a physical system appears naturally in the category of canonical contact phase spaces. The fact that a Jacobi structure is a Lie algebra allows for the interpretation of $h$ as a fundamental conserved quantity of the evolution of the system, $\eta_L(h)[h]=\{h,h\}_L=0$, and the construction of the categorical product of line bundles gives a natural additivity property for the L-energies of two physical systems under the Hamiltonian functor
\begin{center}
\begin{tabular}{ c c c }
Unit-Free Configuration Spaces & Hamiltonian Functor  & Unit-Free Phase Spaces \\
\hline
 $(L_1,h_1), (L_2,h_2)$ & $\begin{tikzcd} \phantom{Q} \arrow[r,mapsto, "\Jet^1 "] & \phantom{Q}  \end{tikzcd}$ & $(L_{\Jet^1 L_1}\utimes L_{\Jet^1 L_2},P_1^*h_1 + P_2^*h_2)$ 
\end{tabular}
\end{center}

This additivity property of energy can be motivated from the natural unit-free generalization of the notion of Riemannian metric. A \textbf{L-metric} on a unit-free configuration space $\lambda:L\to Q$ is a non-degenerate symmetric bilinear form $\gamma:\Der L\odot \Der L\to L$. The presence of a L-metric induces a musical isomorphism
\begin{equation*}
        \begin{tikzcd}
        \Der L \arrow[r, "\gamma^\flat",yshift=0.7ex] & \Jet^1 L \arrow[l,"\gamma^\sharp",yshift=-0.7ex]
        \end{tikzcd},
\end{equation*}
which, in turn, allows for the definition of \textbf{L-kinetic energy} as a quadratic section $K_\gamma\in \Sec{L_{\Jet^1L}}$ via
\begin{equation*}
    K_\gamma(\alpha):=\gamma(\gamma^\sharp(\alpha),\gamma^\sharp(\alpha)).
\end{equation*}
With a choice of a section of the line bundle $v\in \Sec{L}$, that we identify as the \textbf{L-potential}, we define the \textbf{Newtonian L-energy} as
\begin{equation*}
    \Sec{L_{\Jet^1L}}\ni E_{\gamma, v}:= K_\gamma +\pi^*v.
\end{equation*}
It follows from proposition \ref{DerLineProduct} that, given two line bundles with choices of Newtonian L-energy $(L_1,E_{\gamma_1, v_1})$ and $(L_2,E_{\gamma_2, v_2})$, their line product carries a natural choice of Newtonian L-energy
\begin{equation*}
    (L_1\utimes L_2,E_{\gamma_1, v_1}+E_{\gamma_2, v_2})
\end{equation*}
where
\begin{equation*}
    E_{\gamma_1, v_1}+E_{\gamma_2, v_2}:=P_1^*E_{\gamma_1, v_1}+P_2^*E_{\gamma_2, v_2}=E_{p_1^*\gamma_1 \oplus p_2^*\gamma_2, p_1^*v_1+p_2^*v_2}.
\end{equation*}
In Section \ref{ProgramDimDyn} we discuss via a toy model to what extent unit-free configuration spaces endowed with a L-metric and a choice of L-potential can be interpreted as the valid unit-free generalization of the category of ordinary Newtonian configuration spaces.

\section{Generalizations and Jacobization Functors} \label{GeneralizationContact}

In Section \ref{CanonicalContact} we have established contact manifolds as the natural setting for a unit-free treatment of canonical Hamiltonian mechanics. Contact manifolds were also found in Section \ref{ContactGeometry} as the line bundle analogue of symplectic manifolds that unified the non-degenerate cases of Jacobi and precontact manifolds, which were, in turn, presented as the direct line bundle analogues of Poisson and presymplectic manifolds, respectively. This means that a discussion about the natural generalizations of canonical contact phase spaces would mirror that of Section \ref{GeneralizationSymplectic}, where symplectic phase spaces were generalized to Poisson, presymplectic, twisted and Dirac phase spaces.\newline

In light of this complete analogy, generalizations of unit-free phase spaces will then appear in the form of \textbf{the theory of Jacobi phase spaces} and \textbf{the theory of precontact phase spaces}. Similar limitations to those found in the theories of Poisson and precontact phase spaces apply to these together with the added caveat that observables are now unit-free and thus represented by sections of a line bundle. Twistings can be introduced in an entirely analogous way by considering closed 3-forms in the de Rham complex of the line bundle. Lastly, following from our discussion on L-Courant algebroids and L-Dirac structures of Section \ref{LDiracGeometry}, the line bundle analogue of propositions \ref{TwistedPreSymplecticDirac} and \ref{TwistedPoissonDirac} will establish \textbf{the theory of L-Dirac phase spaces} as the natural general unit-free framework encompassing Jacobi phase spaces and precontact phase spaces.\newline

In the same vein of the results of Section \ref{PoissonizationFunctors}, the identification of Jacobi manifolds as unit-free Poisson manifolds clearly shows that the der functor and the jet functor provide direct line bundle analogues of propositions \ref{TangentPoissonizationFunctor} and \ref{CotangentPoissonizationFunctor}, respectively, thus allowing us to regard Jacobi and precontact manifolds as spanning Jacobi structures on vector bundles. The assignments of such spanning Jacobi structures are similarly shown to be functorial and to respect subobjects, products and reduction, thus we call them \textbf{Jacobization functors}. Lastly, in analogy with conjectures \ref{PoissonizationDirac} and \ref{DiracPoissonReduction}, we are compelled to conjecture that a functorial assignment of spanning Jacobi structures to L-Dirac structures generalizing the jet and der Jacobization functors may exist.\newline

Figure \ref{GeneralizationDiagramContact} captures pictorially the generalizations of canonical contact phase spaces proposed in this section.

\chapter{Landscape of Hamiltonian Phase Spaces} \label{LandscapePhaseSpaces}

In this brief chapter we present a visual guide to the breadth of geometric structures introduced in chapters \ref{SymplecticPhaseSpaces} and \ref{ContactPhaseSpaces} during the course of the discussion on the generalization of the concept of phase space. The diagrams displayed below are aimed at providing a clear picture of the available options of generalized phase spaces and how they overlap with each other.\newline

\begin{figure}[h!]
\centering
\begin{tikzpicture}

\draw[thick] (0,0) ellipse (8.5 and 4.4);
\node at (0,-4) {$\textsf{Dir}$};

\draw (-2.5,0) ellipse (4.5 and 3.5);
\node at (-4,-2.6) {$\textsf{twPoiss}$};
\draw (2.5,0) ellipse (4.5 and 3.5);
\node at (4,-2.6) {$\textsf{twPreSymp}$};

\draw[deepblue,thick] (-2.2,0) ellipse (3.7 and 2);
\node at (-3.2,-1.4) {$\textsf{Poiss}$};
\draw[deepred,thick] (2.2,0) ellipse (3.7 and 2);
\node at (3.2,-1.4) {$\textsf{PreSymp}$};

\draw (-2,0.3) ellipse (3 and 1);
\node at (-2.6,-0.3) {$\textsf{sPoiss}$};

\draw[deepgreen,thick] (0,0.3) ellipse (0.75 and 0.5);
\node at (0,0.3) {$\Cot \Man$};
\node at (0,-1) {$\Symp$};
\node at (0,-2.3) {$\textsf{twSymp}$};

\end{tikzpicture}

\caption{Generalizations of \textcolor{deepgreen}{canonical symplectic phase spaces} $\Cot \Man$. Non-degenerate structures are placed in the centre of the diagram, with structures generalizing the symplectic 2-form to the right and structures generalizing the symplectic bivector to the left. Spanning Poisson structures $\textsf{sPoiss}$ are highlighted for their role in the Poissonization functors of Section \ref{PoissonizationFunctors}. In virtue of propositions \ref{TwistedPreSymplecticDirac} and \ref{TwistedPoissonDirac}, twisted Poisson structures $\textsf{twPoiss}$ and twisted presymplectic structures $\textsf{twPreSymp}$ form proper subclasses of Dirac structures $\Dir$ since generic involutive tangent distributions or generic Dirac reductions provide examples of Dirac structures that do not correspond to graphs of bivectors or 2-forms. \textcolor{deepblue}{Poisson manifolds} $\Poiss$ and \textcolor{deepred}{presymplectic manifolds} $\textsf{PreSymp}$ are simply the untwisted subclasses of these.}

\label{GeneralizationDiagramSymplectic}
\end{figure}
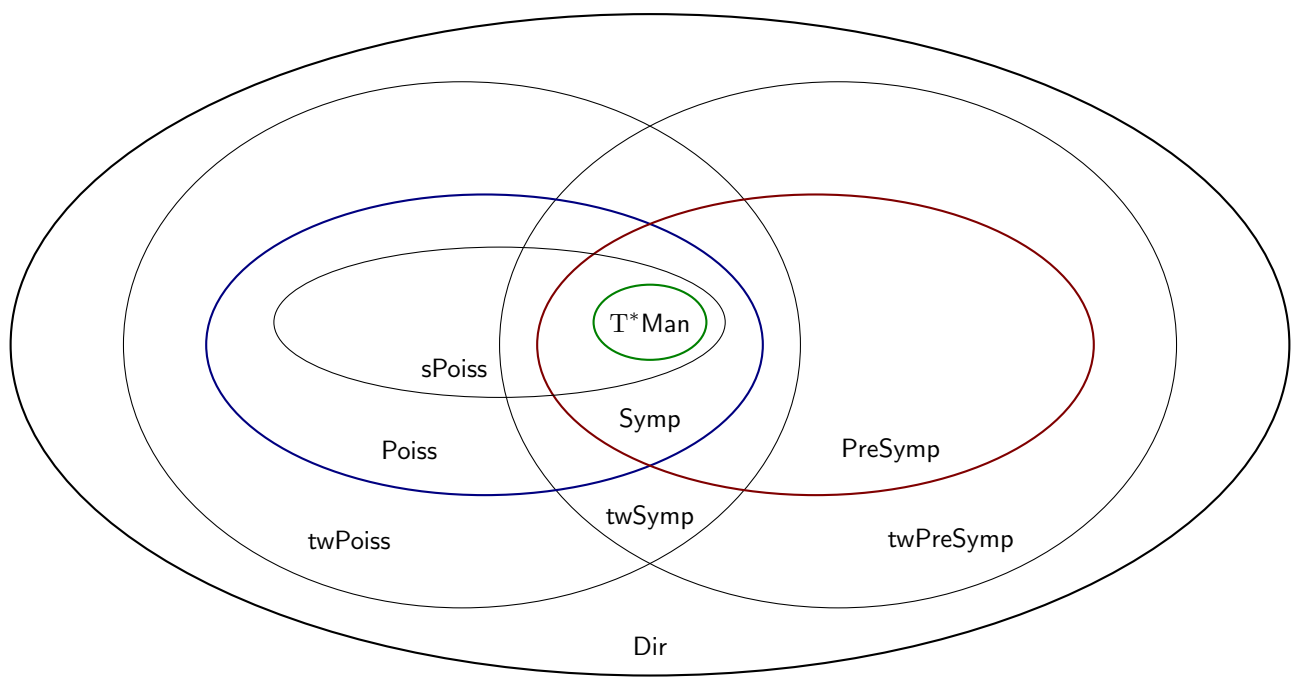

\begin{figure}[h!]
\centering
\begin{tikzpicture}
\draw[thick] (0,0) ellipse (8.5 and 4.4);
\node at (0,-4) {$\textsf{LDir}$};

\draw (-2.5,0) ellipse (4.5 and 3.5);
\node at (-4,-2.6) {$\textsf{twJac}$};
\draw (2.5,0) ellipse (4.5 and 3.5);
\node at (4,-2.6) {$\textsf{twPreCont}$};

\draw[deepblue,thick] (-2.2,0) ellipse (3.7 and 2);
\node at (-3.2,-1.4) {$\textsf{Jac}$};
\draw[deepred,thick] (2.2,0) ellipse (3.7 and 2);
\node at (3.2,-1.4) {$\textsf{PreCont}$};

\draw (-2,0.3) ellipse (3 and 1);
\node at (-2.6,-0.3) {$\textsf{sJac}$};

\draw[deepgreen,thick] (0,0.3) ellipse (0.75 and 0.5);
\node at (0,0.3) {$\Jet^1 \Line$};
\node at (0,-1) {$\Cont$};
\node at (0,-2.3) {$\textsf{twCont}$};
\end{tikzpicture}
\caption{Generalizations of \textcolor{deepgreen}{canonical contact phase spaces} $\Jet^1 \Line$. This is the unit-free analogue of figure \ref{GeneralizationDiagramSymplectic}.}
\label{GeneralizationDiagramContact}
\end{figure}
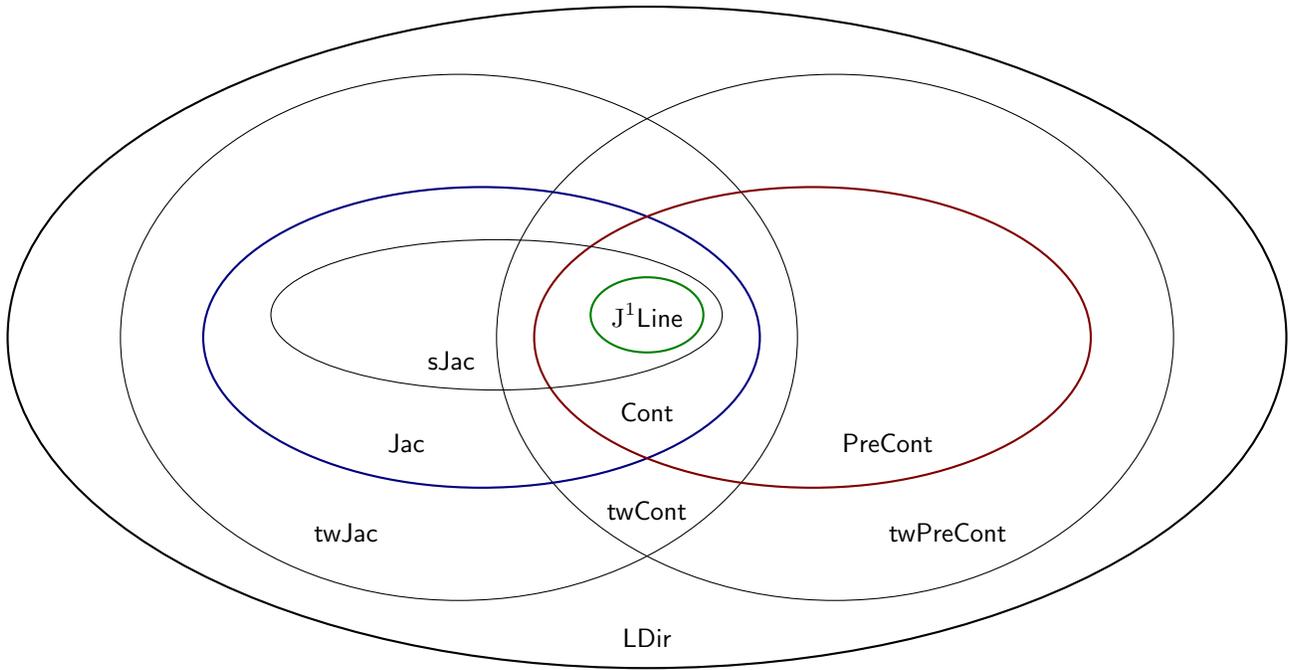

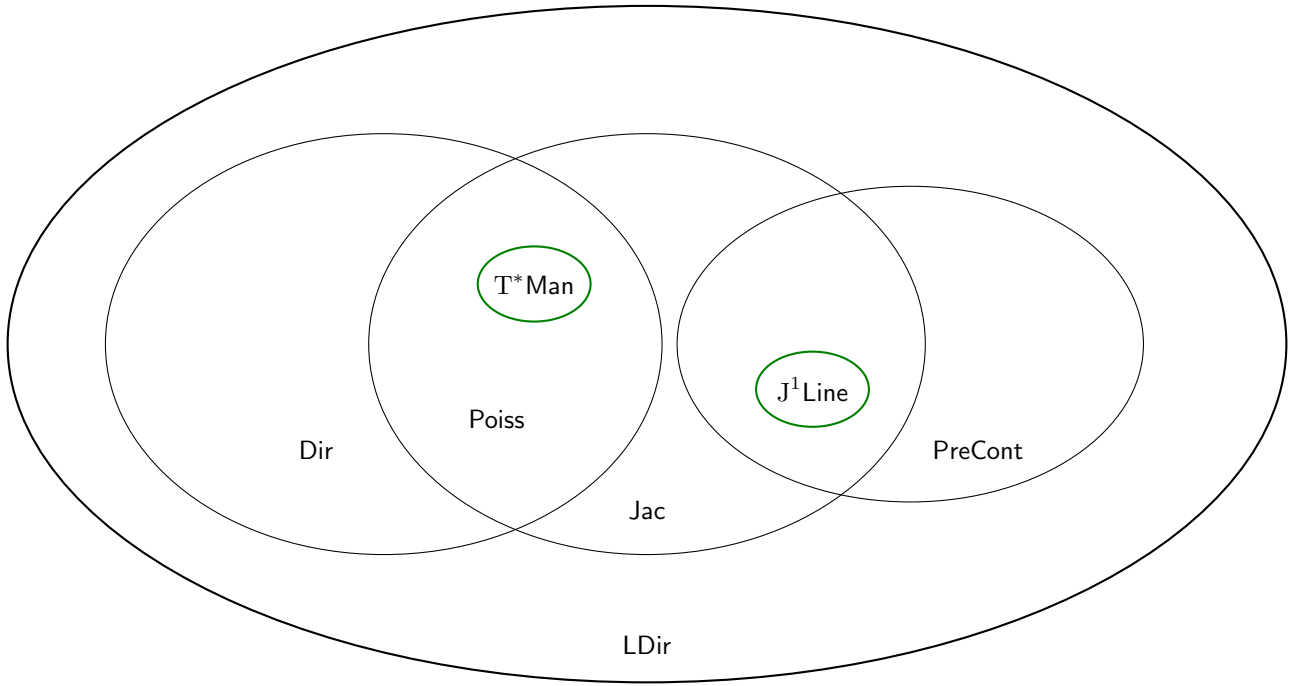
\begin{figure}[h!]
\centering
\begin{tikzpicture}
\draw[thick] (0,0) ellipse (8.5 and 4.5);
\node at (0,-4) {$\textsf{LDir}$};

\draw (3.5,0) ellipse (3.1 and 2.1);
\node at (4.4,-1.4) {$\textsf{PreCont}$};

\draw (0,0) ellipse (3.7 and 2.8);
\node at (0,-2.2) {$\textsf{Jac}$};

\draw (-3.5,0) ellipse (3.7 and 2.8);
\node at (-4.4,-1.4) {$\textsf{Dir}$};
\node at (-2,-1) {$\textsf{Poiss}$};

\draw[deepgreen,thick] (2.2,-0.6) ellipse (0.75 and 0.5);
\node at (2.2,-0.6) {$\Jet^1 \Line$};

\draw[deepgreen,thick] (-1.5,0.8) ellipse (0.75 and 0.5);
\node at (-1.5,0.8) {$\Cot \Man$};

\end{tikzpicture}
\caption{The landscape of Hamiltonian phase spaces. This diagram is drawn in this way since proposition \ref{LocalUnitsJacobi} implied that Poisson manifolds $\Poiss$ can be regarded as special cases of Jacobi structures $\textsf{Jac}$ on trivial line bundles, and also because proposition \ref{LCourantCourant} showed that L-Dirac structures $\textsf{LDir}$ encompass Dirac structures $\Dir$.}
\label{LandscapeHamiltonianPhaseSpaces}
\end{figure}

\chapter{Dimensioned Algebra and Geometry} \label{DimAlgebraGeometry}

As a motivation for the mathematical definitions that will be introduced in this chapter, let us begin by briefly discussing the use of physical quantities in classical thermodynamics. It can be determined empirically that there are four parameters that uniquely characterise the state of a gas, its pressure $P$, volume $V$, amount of matter $N$ and temperature $T$, and that they can all be measured independently by comparison to some reference unit, for example an atmosphere atm, a litre L, a mol and a Kelvin K, respectively. The gas under study will then have a state specified by $(P=p\, \text{atm},V=v\,\text{L},N=n\,\text{mol},T=t\,\text{K})$ with $p,v,n,t$ $\Real$-valued variables. By assuming, for instance, the ideal gas law to model the equilibrium behaviour of the gas, one can derive other physical quantities such as the energy $U=PVN^{-1}$, carrying units atm$\cdot$L$/$mol, or the entropy $E=PVN^{-1}T^{-1}$, carrying units atm$\cdot$L$/$mol$/$K. In general, any measurable physical quantity defined for the state of the gas will be of the form $Q=q\,P^aV^bN^cT^d$ for some $q\in\mathbb{R}$ and $a,b,c,d\in\mathbb{Z}$. The product of two physical quantities $Q$ and $Q'$ is then given by
\begin{equation*}
    Q\cdot Q'= qq'\,P^{a+a'}V^{b+b'}N^{c+c'}T^{d+d'},
\end{equation*}
and two such quantities can be added together to form $Q+Q'$ only when $a=a'$, $b=b'$, $c=c'$ and $d=d'$.\newline

The goal of this chapter is to develop a precise mathematical formalism that encapsulates the formal properties of physical quantities as commonly used in practical science and engineering (see \cite[Chapter 1]{barenblatt1996scaling} for a standard modern reference in dimensional analysis and \cite{hart2012multidimensional} for an attempt at extending the notions of dimensional analysis to linear algebra and calculus). It was argued in our historical note of Section \ref{BriefHistoryOfPhaseSpace} that the algebraic building blocks of the standard formulation of modern physical theories, i.e. fields, vector spaces, groups, etc., were originally conceived without any consideration of the natural structure of physical quantities and with an emphasis on endowing abstract sets with total binary operations, i.e. algebraic operations defined for all pairs of elements. As exemplified in the previous paragraph, the two distinguishing formal features of physical quantities are:
\begin{itemize}
    \item The presence of a set of dimensions or, more concretely, units of measurement, indexing the set of all physical quantities which, furthermore, carries a partial operation given by addition of physical quantities of homogeneous dimensions.
    \item The set of physical quantities also carries a total operation, given by multiplication, and the set of dimensions carries an abelian group structure, in the example above simply $\Int^4$, with a natural compatibility condition between the two, seen explicitly in the formula for $Q\cdot Q'$ in the example above.
\end{itemize}
It is no surprise, then, that the structures of conventional algebra, based on single abstract sets endowed with total operations, prove insufficient to capture the essential algebraic properties of physical quantities. The definitions in the sections to follow are directly motivated by these two characteristic features of the formal manipulation of physical quantities.

\section{Dimensioned Groups, Rings and Modules} \label{DimRingsModules}

Let us begin by introducing the general notion of \textbf{dimensioned set} as an abstract set $A$ together with a surjective map $\alpha:A\to D$ onto another abstract set $D$, called the \textbf{dimension set} of $A$. We will commonly use the notation $A_D$ for a dimensioned set $A$ with dimensions $D$. The surjective map induces a partition of $A$ by preimages $A_d:=\alpha^{-1}(d)\subset A$ that are called the subsets of \textbf{homogeneous dimension}. An element $a\in A_d$ is said to have \textbf{dimension} $d\in D$. A map between dimensioned sets $\Phi:A_D\to B_E$ is called a \textbf{morphism of dimensioned sets} if subsets of homogeneous dimension of $A$ are mapped into subsets of homogeneous dimension of $B$, in other words, there exists a map between the dimension sets $\phi:D\to E$ such that the following diagram commutes
\begin{equation*}
\begin{tikzcd}
A \arrow[r, "\Phi"] \arrow[d, "\alpha"'] & B \arrow[d, "\beta"] \\
D \arrow[r, "\phi"'] & E
\end{tikzcd}
\end{equation*}
the map $\phi$ will be called the \textbf{dimension map}. Clearly, dimensioned sets with the notion of morphism introduced above form a category, we call it  \textbf{the category of dimensioned sets} and denote it by $\textsf{DimSet}$.\newline

Groups are the first class of algebraic structures that we generalize to the dimensioned setting. A dimensioned set $A_D$ is called a \textbf{dimensioned group} if there is a partial binary operation $\cdot$ on $A$ that, upon restriction, induces a group structure on each of the homogeneous subsets $(A_d,\cdot|_{A_d})$. We shall denote dimensioned groups by $(A_D,\cdot_D)$. A dimensioned morphism between dimensioned groups $\Phi:A_D\to B_E$ is called a \textbf{dimensioned group morphism} when the restrictions to the homogeneous subsets $\Phi|_{A_d}:A_d\to B_{\phi(d)}$ are group morphisms for all $d\in D$. The notion of dimensioned group is a direct generalization of the ordinary notion of group as one recovers the defining axioms of group when $D$ is a set with a single element. It clearly follows from our definition that dimensioned groups with dimensioned group morphisms form a subcategory of $\textsf{DimSet}$ that we call \textbf{the category of dimensioned groups} and denote it by $\textsf{DimGrp}$.\newline

Let $(A_D,\cdot_D)$ be a dimensioned group, then the subset $0_D:=\{0_d\in (A_d,\cdot_d),\quad d\in D\}$ is called the \textbf{zero} of $A_D$. A subset $S\subset A_D$ is called a \textbf{dimensioned subgroup} when $S\cap A_d\subset (A_d,\cdot_d)$ are subgroups for all $d\in D$. A dimensioned subgroup $S\subset A_D$ is clearly a dimensioned group with dimension set given by $\alpha(S)$. We can define the \textbf{kernel} of a dimensioned group morphism $\Phi:A_D\to B_E$ in the obvious way
\begin{equation*}
    \Ker{\Phi}:=\{a_d\in A_D|\quad \Phi(a_d)=0_{\phi(d)}\},
\end{equation*}
then clearly $\Ker{\Phi}_D\subset A_D$ is a dimensioned subgroup. A dimensioned subgroup $S\subset A_D$ whose homogeneous intersections $S\cap A_d\subset (A_d,\cdot_d)$ are normal subgroups also induces a natural notion of \textbf{quotient}:
\begin{equation*}
    A_D/S:=\bigcup_{d\in \alpha(S)} A_d/(S\cap A_d)
\end{equation*}
which has an obvious dimensioned group structure with dimension set $\alpha(S)$. There is also a natural notion of \textbf{product} of two dimensioned groups $A_D$, $B_E$ given by
\begin{equation*}
    \alpha\times \beta : A\times B\to D\times E
\end{equation*}
with the partial multiplication defined in the obvious way
\begin{equation*}
    (a_d,b_e)\cdot_{(d,e)} (a_d',b_e'):=(a_d\cdot_da_d',b_e\cdot_e b_e').
\end{equation*}

A dimensioned group is called \textbf{abelian} when all its homogeneous subsets are abelian groups. As an example of a familiar class of objects that displays this kind of structure, we note that vector bundles can be seen as dimensioned groups where the fibre-wise vector addition gives the partial abelian group multiplication and the base manifold is the set of dimensions. An abelian dimensioned group $A$ with set of dimensions $D$ will be denoted by $(A_D,+_D)$. These clearly form a subcategory of $\textsf{DimGrp}$ that we call \textbf{the category of dimensioned abelian groups} and denote it by $\textsf{DimAb}$. Abelian dimensioned groups display structures analogous to those of ordinary abelian groups when we fix a dimension set $D$ and we consider \textbf{dimension-preserving morphisms}, i.e. dimensioned group morphisms $\Phi:A_D\to B_D$ for which the induced map on the dimension sets is the identity $\Id_D:D\to D$. Abelian dimensioned groups over a fixed dimension set $D$ together with dimension-preserving morphisms form a subcategory $\textsf{DimAb}_D\subset\textsf{DimAb}$ that, in addition to the notions of subgroup, kernel and quotient of the general category $\textsf{DimGrp}$, admits a \textbf{direct sum} defined as $A_D\oplus B_D:=(A\ftimes{\alpha}{\beta}B)_D$ with partial multiplication given in the obvious way
\begin{equation*}
    (a_d,b_d)+_d (a_d',b_d'):=(a_d+_d a_d',b_d +_d b_d').
\end{equation*}
It is easy to prove that this direct sum operation on $\textsf{DimAb}_D$ acts as a product and coproduct for which the notions of kernel and quotient identified in the general category $\textsf{DimGrp}$ satisfy the axioms of abelian category. We call $\textsf{DimAb}_D$ \textbf{the category of $D$-dimensioned abelian groups}.\newline

Another similarity between abelian dimensioned groups and ordinary abelian groups is that the sets of morphisms carry a natural dimensioned group structure. Indeed, if we consider a pair of dimensioned group morphisms between abelian dimensioned groups $\Phi, \Psi:(A_D,+_D)\to (B_E,+_E)$ with dimension maps $\phi,\psi: D\to E$, we could attempt to define
\begin{equation*}
    (\Phi + \Psi) (a_d) :=\Phi(a_d)+_e \Phi(a_d)
\end{equation*}
but we will have to choose $e\in E$ in a way that is consistent with the dimension maps. This is clearly achieved by setting $\phi(d)=e=\psi(d)$, thus addition of dimensioned morphisms is only defined for pairs of dimensioned morphism covering the same dimension map. If we denote by $\text{Map}(D,E)$ the set of maps between the dimension sets, we clearly see that the set of dimensioned group morphisms between $A_D$ and $B_E$ has a structure of a dimensioned group with dimension set given by the set of maps between the dimension sets $D$ and $E$. This will be called the \textbf{dimensioned group of morphisms} and we will denote it by $(\Dim{A_D,B_E}_{\text{Map}(D,E)},+_{\text{Map}(D,E)})$ or for a single dimensioned group $\Dim{A_D}_{\text{Map}(D)}:=\Dim{A_D,A_D}_{\text{Map}(D,D)}$.\newline

Taking inspiration from the definition of a conventional ring as an abelian group with a compatible multiplication operation, we define \textbf{dimensioned ring} as a dimensioned abelian group $R_G$ whose dimension set $G$ carries a monoid structure (with multiplication denoted by juxtaposition) and with a total binary operation, called the dimensioned multiplication $\cdot : R_G\times R_G\to R_G$, satisfying the following axioms
\begin{itemize}
    \item[1)] $a\cdot (b\cdot c)=(a\cdot b)\cdot c$,
    \item[2)] $\exists!\quad 1\in R_G|\quad 1\cdot a = a = a\cdot 1$,
    \item[3)] $\rho(a\cdot b)=\rho(a)\rho(b)$,
    \item[4)] $(a+a')\cdot b=a\cdot b+a'\cdot b$,
\end{itemize}
for all $a,a',b,c\in R_G$ with $\rho(a)=\rho(a')$ and where $\rho: R\to G$ denotes the surjective dimension map. Note that in order to demand this list of axioms of a multiplication operation $\cdot : R_G\times R_G\to R_G$ in consistency with the dimensioned structure, the presence of a monoid structure on $G$ is necessary. A dimensioned ring will be denoted by $(R_G,+_G,\cdot)$ and we will employ the subindex notation $a_g\in R_G$ to keep track of the dimensions of elements. Using the explicit index notation, axiom 4) reads
\begin{equation*}
    (a_g+_ga_g')\cdot b_h=a_g\cdot b_h+_{gh}a_g'\cdot b_h.
\end{equation*}
It follows immediately from these axioms that the dimensioned zero $0_G$ acts as an absorbent subset in the following sense
\begin{equation*}
    0_g\cdot a_h=0_{gh},
\end{equation*}
that the multiplicative identities are mapped to each other
\begin{equation*}
    \rho(1)=e\in G,
\end{equation*}
and that the homogeneous subset over the monoid identity forms an ordinary ring with the restricted operations $(R_e,+_e,\cdot|_{R_e})$, we call it the \textbf{dimensionless ring} of $R_G$. We see that a dimensioned ring can be simply understood as a dimensioned abelian group with a monoid structure projecting to the dimension set and that is distributive with respect to the dimensioned addition defined on homogeneous subsets. A dimensioned ring is called \textbf{commutative} when the monoid structures are commutative. For the reminder of this thesis dimensioned rings will be assumed to be commutative unless otherwise stated.\newline

Let $(R_G,+_G,\cdot)$ and $(P_H,+_H,\cdot)$ be two dimensioned rings, a dimensioned group morphism $\Phi:R_G\to P_H$ is called a \textbf{morphism of dimensioned rings} when
\begin{equation*}
    \Phi(a\cdot b)=\Phi(a)\cdot \Phi(b), \qquad \Phi(1_{R_G})=1_{P_H}
\end{equation*}
for all $a,b\in R_G$. The map between the dimension monoids $\phi:G\to H$ is, then, necessarily a monoid morphism. Note that, from the general definition of dimensioned morphism, homogeneous subsets $R_G$ are mapped into homogeneous subsets of $P_H$:
\begin{equation*}
    \Phi(R_g)\subset P_{\phi(g)}.
\end{equation*}
In particular, a dimensioned ring morphism induces a morphism of ordinary rings over the identity
\begin{equation*}
    \Phi|_{R_e}:(R_e,+_e,\cdot|_{R_e})\to (P_{\phi(e)},+_{\phi(e)},\cdot|_{P_{\phi(e)}}).
\end{equation*}
We thus identify dimensioned rings with these morphisms as \textbf{the category of dimensioned rings} and denote it by \textsf{DimRing}. We note again that, as in the general case of dimensioned groups above, setting the dimension monoids to the trivial monoid in all the above instances the usual basic definitions and results for ordinary rings are recovered. It is clear then that we can regard the category of ordinary rings as a subcategory of dimensioned rings $\textsf{Ring}\subset \textsf{DimRing}$.\newline

An abelian dimensioned subgroup $S\subset R_G$ is called a \textbf{dimensioned subring} when $S\cdot S\subset S$ and $1\in S$. A dimensioned subring $I\subset R_G$ is called a \textbf{dimensioned ideal} if for all elements $a_g\in R_G$ and $i_h\in I$ we have
\begin{equation*}
    a_g\cdot i_h\in I\cap R_{gh}.
\end{equation*}
This condition ensures that the general construction of quotient by an abelian dimensioned subgroup applied to the case of an ideal gives the \textbf{dimensioned quotient ring} $R_G/I$ in a natural way. If we denote $I_g:=I\cap R_g$, the dimensioned quotient ring multiplication is explicitly checked to be well-defined:
\begin{equation*}
    (a_g+_gI_g)\cdot (b_h+_hI_h)=a_g\cdot b_h +_{gh} a_g\cdot I_h +_{gh} b_h\cdot I_g +_{gh} I_g\cdot I_h=a_g\cdot b_h +_{gh} I_{gh}.
\end{equation*}

A \textbf{choice of units} $u$ in a dimensioned ring $R_G$ is a splitting of the dimension projection
\begin{equation*}
\begin{tikzcd}[row sep=small]
R  \arrow[d, "\rho"'] \\
G \arrow[u, "u"', bend right=50] 
\end{tikzcd}\quad \rho \circ u =\Id_G,
\qquad \text{ such that } \quad
    u_{gh}=u_g\cdot u_h, \qquad u_g\neq 0_g
\end{equation*}
for all $g,h\in G$. In other words, a choice of units is a splitting of monoids $u:G\to R$ with non-zero image. Choices of units can be regarded as the dimensioned generalization of the notion of non-zero element of a ring with the caveat that they may not exist due to the non-vanishing condition being required for all of $G$. It was noted above that vector bundles give examples of dimensioned rings, then, considering the Moebius band as a dimensioned ring with dimension set the abelian group $\text{U}(1)$ and the zero multiplication operation, we find an explicit example of a dimensioned ring that does not admit choices of unit.\newline

As first examples of dimensioned rings we have already mentioned ordinary rings and vector bundles with the zero multiplication. Another important example is given by pairs of ordinary rings and monoids: let $R$ be a ring and $G$ a monoid, then the Cartesian product $R\times G$ carries a natural dimensioned ring structure, called the \textbf{trivial dimensioned ring} $R$ with dimensions in $G$, defined in the obvious way
\begin{equation*}
    \Proj_2:R\times G\to G \qquad (a,g)+_g(b,g):=(a+b,g), \qquad (a,g)\cdot (b,h):=(a\cdot b,gh).
\end{equation*}

A dimensioned ring $R_G$ is called a \textbf{dimensioned field} when
\begin{equation*}
    \forall a\notin 0_G\quad \exists a^{-1}|\quad a\cdot a^{-1}=1=a^{-1}\cdot a,
\end{equation*}
note that for this requirement to be consistent with the dimension projection the monoid structure of $G$ must be a group structure, for this reason $G$ will be called the dimension group of the dimensioned field. It is easy to see that multiplicative inverses of a dimensioned field are mapped into inverses of the dimension group $\rho(a^{-1})=\rho(a)^{-1}$.\newline

A direct consequence of the defining condition of dimensioned field is that non-zero elements induce group isomorphisms between homogeneous subsets. Indeed, for a non-zero element $0_g\neq a_g\in R_G$ we have the following maps
\begin{align*}
a_g\cdot : R_h & \to R_{gh}\\
b_h & \mapsto a_g\cdot b_h,
\end{align*}
which are group morphisms from axiom 4) of dimensioned rings and are invertible with inverse given by $a_g^{-1}\cdot$. These maps allow to prove a general result that confers a role to choices of unit on dimensioned fields similar to that of a trivialization of a fibre bundle.

\begin{prop}[Choices of Units in Dimensioned Fields] \label{UnitsDimFields}
Let $(R_G,+_G,\cdot)$ be a dimensioned field, then a choice of units $u:G\to R$ induces an isomorphism with the trivial dimensioned field $R_e$ with dimensions in $G$:
\begin{equation*}
    R_G\cong R_e\times G.
\end{equation*}
\end{prop}
\begin{proof}
We can explicitly construct the following map
\begin{align*}
\Phi^u: R_e\times G & \to R_G\\
(r,g) & \mapsto u_g\cdot r,
\end{align*}
which is a clearly bijective morphism of dimensioned abelian groups from the fact that it is constructed with the group isomorphisms $u_g\cdot$ for all the values of the choice of units $u$. It only remains to check that it is indeed a dimensioned ring morphism, this follows directly by construction and the fact that $u$ is a morphism of monoids:
\begin{align*}
    \Phi^u((r_1,g)\cdot (r_2,h))&= \Phi^u((r_1\cdot r_2,gh))=u_{gh}\cdot r_1\cdot r_2=u_g\cdot u_h\cdot r_1\cdot r_2 =\\
    &=(u_g\cdot r_1)\cdot (u_h\cdot r_2)=\Phi^u(u_g\cdot r_1)\cdot \Phi^u(u_h\cdot r_2).
\end{align*}
\end{proof}
This last proposition shows that the dimensioned fields for which choices of units exist are non-canonically isomorphic to the trivial dimensioned fields $\mathbb{F}\times G$ with $\mathbb{F}$ an ordinary field and $G$ an abelian group.\newline

Let us consider again a general a dimensioned ring $(R_G,+_G,\cdot)$. Recall from our discussion above that the dimensioned morphisms from $R_G$ into itself form an abelian dimensioned group $(\Dim{R_G}_{\text{Map}(G)},+_{\text{Map}(G)})$ where $\text{Map}(G)$ denotes the set of maps from $G$ into itself. The presence of the dimensioned ring multiplication allows for the definition of the following module-like structure
\begin{equation*}
    \cdot :R_G\times \Dim{R_G} \to \Dim{R_G}
\end{equation*}
defined via
\begin{equation*}
    (a_g\cdot \Phi)(b_h):=a_g\cdot \Phi(b_h).
\end{equation*}
We note that $a_g\cdot \Phi$ is a well-defined dimensioned morphism from the fact $G$ acts naturally on $\text{Map}(G)$ by composition with the left action of $G$ on itself, indeed the if $\phi:G\to G$ is the dimension map of $\Phi$, then $a_g\cdot \Phi$ has dimension map $L_g\circ \phi:G\to G$. It follows directly from the axioms of the dimensioned ring multiplication that that this operation satisfies the usual linearity properties of the conventional notion of $R_G$-module with the caveat that addition is only partially defined.\newline

Recall now that products and direct sums of general abelian dimensioned groups can be taken so, considering the abelian dimensioned group part of a dimensioned ring $(R_G,+_G,\cdot)$, we can form the product $R_G\times R_G$, which is an abelian dimensioned group with dimension set $G\times G$, or the direct sum $R_G \oplus_G R_G$, which is a dimensioned abelian group with dimension set $G$. In both cases we can form module-like maps by setting
\begin{equation*}
    a_g\cdot (b_h,c_k):=(a_g\cdot b_h,a_g\cdot c_k), \qquad a_g\cdot (b_h\oplus c_h):= a_g\cdot b_h \oplus a_g\cdot c_h.
\end{equation*}
These module-like actions preserve the dimensioned structure from the fact that, in the first case, $G$ acts diagonally on $G\times G$ and, in the second case, $G$ acts on itself by left multiplication. Furthermore, from the defining axioms of dimensioned ring, these maps satisfy the usual linearity properties of the conventional notion of $R_G$-module with the caveat that addition is only partially defined.\newline

These examples motivate the definition of dimensioned modules: let $(R_G,+_G,\cdot)$ be a dimensioned ring and $(A_D,+_D)$ a dimensioned abelian group, then $A_D$ is called a \textbf{dimensioned $R_G$-module} if there is a map 
\begin{equation*}
    \cdot :R_G\times A_D \to A_D
\end{equation*}
that is compatible with the dimensioned structure via a monoid action $G\times D\to D$ (denoted by juxtaposition) in the following sense
\begin{equation*}
    r_g\cdot a_d=(r \cdot a)_{gd}
\end{equation*}
and that satisfies the following axioms
\begin{itemize}
    \item[1)] $r_g\cdot (a_d+b_d)=r_g\cdot a_d + r_g\cdot b_d$,
    \item[2)] $(r_g+p_g)\cdot a_d=r_g\cdot a_d + p_g \cdot a_d$,
    \item[3)] $(r_g\cdot p_h)\cdot a_d=r_g\cdot (p_h\cdot a_d)$,
    \item[4)] $1\cdot a_d=a_d$
\end{itemize}
for all $r_g,p_h\in R_G$ and $a_d,b_d\in A_D$. Note that these four axioms for a map $\cdot :R_G\times A_D \to A_D$ can only be demanded in consistency with the dimensioned structure in the presence of a monoid action $G\times D\to D$. With this definition at hand, we recover the motivating examples for a dimension ring $R_G$ introduced above: the \textbf{dimensioned module of dimensioned morphisms} $\Dim{R_G}$ is clearly a dimensioned $R_G$-module with dimension set $\text{Map}(G)$; the product $R_G\times R_G$ is a dimensioned $R_G$-module with dimension set $G\times G$; and the direct sum $R_G \oplus_G R_G$ is a dimensioned $R_G$-module with dimension set $G$.\newline

Let $(A_D,+_D)$ and $(B_E,+_E)$ be two dimensioned $R_G$-modules, a morphism of abelian dimensioned groups $\Phi:A_D\to B_E$ is called \textbf{$R_G$-linear} if
\begin{equation*}
    \Phi(r_g\cdot a_d)=r_g\cdot \Phi(a_d)
\end{equation*}
for all $r_g\in R_G$ and $a_d\in A_D$. Note that this condition forces the dimension map $\phi:D\to E$ to satisfy
\begin{equation*}
    \phi(gd)=g\phi(d)
\end{equation*}
for all $g\in G$ and $d\in D$, in other words, the dimension map $\phi$ must be $G$-equivariant with respect to the monoid actions of the dimension sets $D$ and $E$. Let us denote the set of $G$-equivariant dimension maps as
\begin{equation*}
    \text{Map}^G(D,E):=\{\phi:D\to E|\quad \phi\circ g = g \circ \phi, \quad \forall g\in G\},
\end{equation*}
then it follows that the dimensioned group of morphisms $\Dim{A_D,B_E}_{\text{Map}(D,E)}$ contains a dimensioned subgroup of morphisms covering $G$-equivariant dimension maps for which the following dimensioned module map can be defined
\begin{equation*}
    (r_g\cdot \Phi)(a_d):=r_g\cdot \Phi(a_d)=\Phi(r_g\cdot a_d).
\end{equation*}
This clearly makes $\Dim{A_D,B_E}_{\text{Map}^G(D,E)}$ into a dimensioned $R_G$-module, we shall call this the \textbf{dimensioned module of $R_G$-linear morphisms}.\newline

Let $(A_D,+_D)$ be a dimensioned $R_G$-module, an abelian dimensioned subgroup $S\subset A_D$ is called a \textbf{dimensioned submodule} if
\begin{equation*}
    r_g\cdot s_d\in S\cap A_{gd}
\end{equation*}
for all $r_g\in R_G$ and $s_d\in S$. All the notions introduced at the beginning of this section for general abelian dimensioned groups, e.g. direct sums, products, quotients, etc., apply to dimensioned $R_G$-modules in particular. Furthermore, given a dimensioned submodule $S\subset A_D$ that is a dimensioned $I$-module for $I\subset R_G$ a dimensioned ideal, there is a natural notion of \textbf{dimensioned quotient module} $A_D/S$ with dimensioned ring given by the dimensioned quotient ring $R_G/I$.\newline

As in the case or ordinary modules, Dimensioned modules admit a tensor product construction: Let $(A_D,+_D)$ and $(B_E,+_E)$ be two dimensioned $R_G$-modules, then we define their \textbf{dimensioned tensor product} as
\begin{equation*}
    A_D\otimes_{R_G} B_E := R_G\bullet (A_D \times B_E)/\sim
\end{equation*}
where $R_G\bullet (A_D \times B_E)$ denotes the free abelian dimensioned group of pairs $(a_d,b_e)$ with coefficients in $R_G$ and $\sim$ denotes taking a quotient with respect to the following relations within the free abelian dimensioned group
\begin{equation*}
    (a_d+a_d',b_e)=(a_d,b_e)+(a_d',b_e), \qquad (a_d,b_e+b_e')=(a_d,b_e)+(a_d,b_e'), \qquad (r_g\cdot a_d,b_e)=(a_d,r_g\cdot b_e).
\end{equation*}
Note that these are the same relations used to define the tensor product of ordinary modules with the added caveat that addition is only partially defined. This construction clearly makes $A_D\otimes_{R_G} B_E$ into a dimensioned $R_G$-module with dimension set $D\times E$ and monoid $G$-action given by the diagonal action. The dimensioned tensor product so defined makes the abelian category of dimensioned $R_G$-modules into a monoidal category with the tensor unit given by $R_G$. Particularly, this definition makes the dimensioned tensor product distributive with respect to the dimensioned direct sum since it easy to check directly from the definition that, for three dimensioned $R_G$-modules $A_D$, $B_D$ and $C_E$, there is a canonical $R_G$-linear isomorphism
\begin{equation*}
    (A_D \oplus_D B_D)\otimes_{R_G} C_E \cong A_D\otimes_{R_G} C_E \oplus_{D\times E} B_D\otimes_{R_G} C_E.
\end{equation*}

\section{Dimensioned Algebras} \label{DimAlgebras}

Let us motivate our discussion on the dimensioned generalization of the notion of algebra for ordinary rings and modules by considering the dimensioned morphisms of a dimensioned ring $R_G$. In Section \ref{DimRingsModules} above it was shown that $\Dim{R_G}_{\text{Map}(G)}$ is a dimensioned $R_G$-module, we note that the dimension set given by the maps from $G$ into itself $\text{Map}(G)$ carries a natural monoid structure given by composition of maps. Denoting three dimensioned morphisms by $\Phi_\phi,\Theta_\phi, \Psi_\psi:R_G\to R_G$ where $\phi,\psi:G\to G$ are the dimension maps, it follows directly from the defining properties of dimensioned rings that the composition of the dimensioned morphisms is consistent with the monoid structure of the dimension set $\text{Map}(G)$
\begin{equation*}
    \Phi_\phi \circ \Psi_\psi = (\Phi \circ \Psi)_{\phi\circ \psi},
\end{equation*}
and that it interacts with the $R_G$-module structure of $\Dim{R_G}$ as a bilinear operation
\begin{equation*}
    (\Phi_\phi +_\phi \Theta_\phi )\circ \Psi_\psi= \Phi_\phi\circ \Psi_\psi +_{\phi\circ \psi} \Theta_\phi \circ \Psi_\psi \qquad r_g\cdot (\Phi_\phi \circ \Psi_\psi)= (r_g\cdot \Phi_\phi) \circ \Psi_\psi=\Phi_\phi \circ (r_g\cdot \Psi_\psi)
\end{equation*}
for all $r_g\in R_G$. This shows that $(\Dim{R_G},\circ)$ gives a prime example of a bilinear associative operation on a dimensioned module and prompts us to give the following general definition.\newline

Let $(A_D,+_D)$ be a dimensioned $R_G$-module, a map $M:A_D\times A_D\to A_D$ is called a \textbf{dimensioned bilinear multiplication} if it satisfies
\begin{align*}
    M(a_d \, +_d\, b_d,c_e)&=M(a_d,c_e) \,+_{\mu(d,e)} \, M(b_d,c_e)\\
    M(a_d,b_e \, +_e \, c_e)&=M(a_d,b_e) \, +_{\mu(d,e)} \, M(a_d,c_e)\\
    M(r_g\cdot a_d,s_h\cdot b_e)&=r_g\cdot s_h\cdot M(a_d, b_e)
\end{align*}
for all $a_d,b_d,b_e,c_e\in A_D$, $r_g,s_h\in R_G$ and for a \textbf{dimension map} $\mu:D\times D\to D$ which is $G$-equivariant in both entries, i.e.
\begin{equation*}
    \mu(gd,he)=gh\mu(d,e)
\end{equation*}
for all $g,h\in G$ and $d,e\in D$. When such a map $M$ is present in a dimensioned $R_G$-module $A_D$, the pair $(A_D,M)$ is called a \textbf{dimensioned $R_G$-algebra}. The notion of dimensioned tensor product given a the end of Section \ref{DimRingsModules} above allows to reformulate the definition of a dimensioned bilinear multiplication $M:A_D\times A_D\to A_D$ as a dimensioned $R_G$-linear morphism
\begin{equation*}
    M:A_D\otimes_{R_G} A_D \to A_D.
\end{equation*}
Note that the dimension set of the tensor product $A_D\otimes_{R_G} A_D$ is $D\times D$ with the diagonal $G$-action induced from the $R_G$-module structure, then we see that the double $G$-equivariant condition of $\mu$ is reinterpreted now as ordinary $G$-equivariance with respect to the natural monoid actions.\newline

The natural notions of morphisms and subalgebras of ordinary algebras extend naturally to the dimensioned case. Let $(A_D,M)$ and $(B_E,N)$ be two dimensioned $R_G$-algebras, a $R_G$-linear morphism $\Phi:A_D\to B_E$ is called a \textbf{morphism of dimensioned algebras} if
\begin{equation*}
    \Phi(M(a,a'))=N(\Phi(a),\Phi(a')),
\end{equation*}
for all $a,a'\in A_D$. A submodule $S\subset A_D$ such that $M(S,S)\subset S$ is called a \textbf{dimensioned subalgebra}.\newline

The dimension map $\mu$ of a dimensioned bilinear multiplication in a dimensioned $R_G$-algebra $(A_D,M_\mu)$ is naturally regarded as an binary operation on the set of dimensions $D$. In a general sense, dimension sets of dimensioned algebras carry the most basic algebraic structures, commonly known as magmas. However, if one wishes to demand specific algebraic properties, such as commutativity or associativity, the algebraic structure present in the dimension magma becomes richer. Let $(A_D,M_\mu)$ be a dimensioned $R_G$-algebra, we say that it is \textbf{symmetric} or \textbf{antisymmetric} if
\begin{equation*}
    M(a_d,b_e)=M(b_e,a_d), \qquad M(a_d,b_e)=-M(b_e,a_d)
\end{equation*}
for all $a_d,b_e\in A_D$, respectively. The dimension magmas of symmetric or antisymmetric dimensioned algebras are necessarily commutative, i.e. $\mu(d,e)=\mu(e,d)$ for all $d,e\in D$. The usual 3-element-product properties of ordinary algebras can be demanded for dimensioned algebras in an analogous way, in particular $(A_D,M_\mu)$ is called \textbf{associative} or \textbf{Jacobi} if
\begin{equation*}
    \text{Ass}_M(a_d,b_e,c_f)=0,\qquad \text{Jac}_M(a_d,b_e,c_f)=0
\end{equation*}
for all $a_d,b_e,c_f\in A_D$, respectively. The dimension magmas of associative or Jacobi dimensioned algebras are necessarily associative, i.e. $\mu(\mu(d,e),f)=\mu(d,\mu(e,f))$ for all $d,e,f\in D$, making them into semigroups. Returning to the motivating example presented at the beginning of this section, we now see that the dimensioned morphisms of a dimensioned ring $R_G$ give the prime example of dimensioned associative algebra $(\Dim{R_G},\circ)$.\newline

In parallel with the definitions of ordinary algebras, we define \textbf{dimensioned commutative algebra} as a symmetric and associative dimensioned algebra and a \textbf{dimensioned Lie algebra} as an antisymmetric and Jacobi dimensioned algebra. Note that dimensioned commutative and dimensioned Lie algebras necessarily carry dimension sets that are commutative semigroups.\newline

In keeping with the general philosophy to continue to scrutinize the natural algebraic structure present in the dimensioned module of dimensioned morphisms of a dimensioned ring $R_G$, let us attempt to find the appropriate dimensioned generalization of the notion of derivations of a ring. Working by analogy, a dimensioned derivation will be a dimensioned morphism $\Delta\in \Dim{R_G}$ covering a dimension map $\delta:G\to G$ satisfying a Leibniz identity with respect to the dimensioned ring multiplication
\begin{equation*}
    \Delta(r_g\cdot s_h)= \Delta(r_g)\cdot s_h + r_g\cdot\Delta(s_h),
\end{equation*}
for all $r_g,s_h\in R_G$, however, for the right-hand-side to be well-defined, both terms must be of homogeneous dimension, which means that the dimension map must satisfy
\begin{equation*}
    \delta(gh)=\delta(g)h=g\delta(h)
\end{equation*}
for all $g,h\in G$. Since $G$ is a monoid, this condition is equivalent to the dimension map being given by multiplication with a monoid element, i.e. $\delta=L_d$ for some element $d\in G$. Following from this observation, we see that there is a natural dimensioned submodule of the dimensioned module of dimensioned morphisms $\Dim{R_G}_G\subset \Dim{R_G}_{\text{Map(G)}}$ given by the dimensioned morphisms whose dimension maps are specified by multiplication with a monoid element. Recall that dimensioned rings are assumed to be commutative and, thus, the dimension monoid has commutative binary operation. This allows for the identification of the first natural example of dimensioned Lie algebra: consider the commutator of the associative dimensioned composition
\begin{equation*}
    [\Delta,\Delta']:=\Delta \circ \Delta'-\Delta' \circ \Delta,
\end{equation*}
it is easy to check that this bracket is indeed antisymmetric and Jacobi, thus making $(\Dim{R_G}_G,[,])$ into the \textbf{dimensioned Lie algebra of dimensioned morphisms} of a dimensioned ring $R_G$. Notice that this bracket can only be defined on the dimensioned submodule $\Dim{R_G}_G\subset \Dim{R_G}_{\text{Map(G)}}$ since the two terms of the right-hand-side for general dimensioned morphisms will have dimensions given by the composition of maps from $G$ into itself which is a non-commutative binary operation in general. It is then clear that the Leibniz condition proposed above can be demanded in consistency with the dimensioned structure of dimensioned morphisms within the Lie algebra of dimensioned morphisms, so we see the dimensioned Lie algebra of \textbf{derivations of a dimensioned ring} $R_G$ as the natural dimensioned Lie subalgebra of the dimensioned morphisms
\begin{equation*}
    \Dr{R_G}\subset (\Dim{R_G}_G,[,]).
\end{equation*}
Derivations covering the monoid identity, i.e. those with dimension map $\Id_G:G\to G$, are called \textbf{dimensionless derivations} and it is clear by definition that they form an ordinary Lie algebra with the commutator bracket $(\Dr{R_G}_0, [,])$. Restricting their action to elements of the dimensionless ring $R_e\subset R_G$ we recover the ordinary Lie algebra of ring derivations, in other words, there is a surjective map of Lie algebras
\begin{equation*}
    (\Dr{R_G}_0,[,])\to (\Dr{R_e},[,]).
\end{equation*}

The example of the dimensioned Lie algebra of derivations of a dimensioned ring illustrates the case of a dimensioned algebra whose space of dimensions is a (commutative) monoid and whose dimension map is simply given by the monoid multiplication. For the reminder of this chapter, the dimension sets of dimensioned modules will be assumed to carry a commutative monoid structure (with multiplication denoted by juxtaposition of elements) unless stated otherwise. Let $A_G$ be a dimensioned module, a dimensioned algebra multiplication $M:A_G\times A_G\to A_G$ is said to be \textbf{homogeneous of dimension $m$} if the dimension map $\mu:G\times G\to G$ is given by monoid multiplication with the element $m\in G$, i.e. $\mu(g,h)=mgh$ for all $g,h\in G$. Assuming a monoid structure on the dimension set of a dimensioned module and considering dimensioned algebra multiplications of homogeneous dimension is particularly useful in order to study several algebra multiplications coexisting on the same set. Indeed, given two homogeneous dimensioned algebra multiplications $(A_G,M_1)$ and $(A_G,M_2)$ with dimensions $m_1\in G$ and $m_2\in G$, respectively, the fact that the monoid operation is assumed to be associative and commutative, allows for consistently demanding properties of the interaction of the two dimensioned multiplications involving expressions of the form $M_1(M_2(a,b),c)$ without any further requirements.\newline

Let $A_G$ be a dimensioned $R_H$-module and let two dimensioned algebra multiplications $*:A_G\times A_G\to A_G$ and $\{,\}:A_G\times A_G\to A_G$ with homogeneous dimensions $p\in G$ and $b\in G$, respectively, the triple $(A_G,*_p,\{,\}_b)$ is called a \textbf{dimensioned Poisson algebra} if
\begin{itemize}
    \item[1)] $(A_G,*_p)$ is a dimensioned commutative algebra,
    \item[2)] $(A_G,\{,\}_b)$ is a dimensioned Lie algebra,
    \item[3)] the two multiplications interact via the Leibniz identity
    \begin{equation*}
        \{a,b*c\}=\{a,b\}*c+b*\{a,c\},
    \end{equation*}
    for all $a,b,c\in A_G$.
\end{itemize}
Note that the Leibniz condition can be consistently demanded of the two dimensioned algebra multiplications since the dimension projections of each of the terms of the Leibniz identity for $\{a_g,b_h*c_k\}$ are: 
\begin{equation*}
    bgphk,\quad pbghk,\quad phbgk,
\end{equation*}
but they are indeed all equal from the fact that the monoid binary operation is associative and commutative.\newline

A morphism of dimensioned modules between dimensioned Poisson algebras $\Phi:(A_G,*_p,\{,\}_b)\to (B_H,*_r,\{,\}_c)$ is called a \textbf{morphism of dimensioned Poisson algebras} if $\Phi:(A_G,*_p)\to (B_H,*_r)$ is a morphism of dimensioned commutative algebras and also $\Phi:(A_G,\{,\}_b)\to (B_H,\{,\}_c)$ is a morphism of dimensioned Lie algebras. A submodule $I\subset A_G$ that is a dimensioned ideal in $(A_G,*_p)$ and that is a dimensioned Lie subalgebra in $(A_G,\{,\}_b)$ is called a \textbf{dimensioned coisotrope}.

\begin{prop}[Dimensioned Poisson Reduction] \label{DimensionedPoissonReduction}
Let $(A_G,*_p,\{,\}_b)$ be a dimensioned Poisson algebra and $I\subset A_G$ be a coisotrope, then there is a dimensioned Poisson algebra structure induced in the subquotient
\begin{equation*}
    (A'_G:=N(I)/I,*'_p,\{,\}'_b)
\end{equation*}
where $N(I)$ denotes the dimensioned Lie idealizer of $I$ regarded as a submodule of the dimensioned Lie algebra.
\end{prop}
\begin{proof}
We assume without loss of generality that the dimension projection of $I$ is the whole of $G$, the intersections with the homogeneous subsets are denoted by $I_g:=I\cap A_g$. The dimensioned Lie idealizer is defined in the obvious way
\begin{equation*}
    N(I):=\{n_g\in A_G|\quad \{n_g,i_h\}\in I_{bgh} \quad \forall i_h\in I\}.
\end{equation*}
We clearly see that $N(I)$ is the smallest dimensioned Lie subalgebra that contains $I$ as a dimensioned Lie ideal. The Leibniz identity implies that $N(I)$, furthermore, is a dimensioned commutative subalgebra with respect to $*_p$ in which $I$ sits as a dimensioned commutative ideal, since it is a commutative ideal in the whole $A_G$. It follows that we can form the dimensioned quotient commutative algebra $(N(I)/I,*')$ in an entirely analogous way to the construction of the dimensioned quotient ring presented in Section \ref{DimRingsModules}. The only difference with that case is that commutative multiplication covers a dimension map that is given by the monoid multiplication with a non-identity element $p\in G$, but this has no effect on the quotient construction. To obtain the desired quotient dimensioned Lie bracket we set:
\begin{equation*}
    \{n_g+I_g,m_h+I_h\}':=\{n_g,m_h\}+I_{bgh}
\end{equation*}
which is easily checked to be well-defined and that inherits the antisymmetry and Jacobi properties directly from dimensioned Lie bracket $\{,\}$ and the fact that $I\subset N(I)$ is a dimensioned Lie ideal.
\end{proof}

\section{The Potential Functor} \label{PotentialFunctor}

In this section we will find a close link between dimensioned algebras and the category of lines $\Line$ introduced in Section \ref{CategoryOfLines}. The results presented in this section will be instrumental for section  \ref{JacobiManifoldsRevisited}, where we argue that the general notions of dimensioned algebra introduced in sections \ref{DimRingsModules} and \ref{DimAlgebras} provide the natural language to express constructions of Jacobi manifolds algebraically. Furthermore, these results provide the direct connection of the notion of dimensioned ring with the motivating example of physical quantities with units and the abstract notion of measurand space introduced in Section \ref{MeasurandFormalism}.\newline

Recall that in Section \ref{CategoryOfLines} it was shown that the categorical structure of 1-dimensional vector spaces allowed for the construction of abelian groups of tensor powers. More concretely, the \textbf{potential} of a line $L\in\Line$ was defined as the set of all tensor powers
\begin{equation*}
    L^\odot:=\bigcup_{n\in \Int} L^n
\end{equation*}
where $L^n$ denotes the tensor powers of $L$ for $n>0$, the tensor powers of the dual line $L^*$ for $n<0$ and the patron line $\Real$ for $n=0$. This set has than an obvious dimensioned set structure with dimension set $\Int$:
\begin{equation*}
    \pi:L^\odot\to \Int.
\end{equation*}
Since homogeneous subsets are precisely the tensor powers $L^n$, they carry a natural $\Real$-vector space structure thus clearly making the potential of $L$ into an abelian dimensioned group $(L^\odot_\Int,+_\Int)$. The next proposition shows that the ordinary $\Real$-tensor product of vector spaces endows $L^\odot$ with a dimensioned field structure.

\begin{prop}[Dimensioned Ring Structure of the Potential of a Line] \label{DimRingPotential}
Let $L\in\Line$ be a line and $(L^\odot_\Int,+_\Int)$ its potential, then the $\Real$-tensor product of elements induces a dimensioned multiplication
\begin{equation*}
    \odot: L^\odot \times L^\odot\to L^\odot
\end{equation*}
such that $(L^\odot_\Int,+_\Int,\odot)$ becomes a dimensioned field.
\end{prop}
\begin{proof}
The construction of the dimensioned ring multiplication $\odot$ is done simply via the ordinary tensor product of ordinary vectors and taking advantage of the particular properties of 1-dimensional vector spaces. The two main facts that follow from the 1-dimensional nature of lines are: firstly, that linear endomorphisms are simply real numbers
\begin{equation*}
    \text{End}(L)\cong L^*\otimes L\cong \Real
\end{equation*}
which, at the level of elements, means that
\begin{equation*}
    \text{End}(L)\ni\alpha \otimes a =\alpha(a)\cdot \Id_{L}
\end{equation*}
as it can be easily shown by choosing a basis; and secondly, that the tensor product becomes canonically commutative, since, using the isomorphism above, we can directly check
\begin{equation*}
    a\otimes b(\alpha,\beta)=\alpha(a)\beta(b)=\alpha(b)\beta(a)=b\otimes a(\alpha,\beta),
\end{equation*}
thus showing
\begin{equation*}
    a\otimes b=b\otimes a \in L\otimes L=L^2.
\end{equation*}
The binary operation $\odot$ is then explicitly defined for elements $a,b\in L=L^1$, $\alpha,\beta\in L^*=L^{-1}$ and $r,s\in \Real=L^0$ by
\begin{align*}
    a\odot b &:= a\otimes b\\
    \alpha\odot \beta &:= \beta \otimes \alpha\\
    r\odot s &:= r\otimes s=rs\\
    r\odot a &:= ra\\
    r\odot \alpha &:= r\alpha\\
    \alpha \odot a &:=\alpha(a) = a(\alpha) =: a \odot \alpha
\end{align*}
Products of two positive power tensors $a_1\otimes \cdots \otimes a_q$, $b_1\otimes \cdots \otimes b_p$ and negative powers $\alpha_1\otimes \cdots \otimes \alpha_q$, $\beta_1\otimes \cdots \otimes \beta_p$ are defined by
\begin{align*}
    (a_1\otimes \cdots \otimes a_q)\odot(b_1\otimes \cdots \otimes b_p) &:= a_1\otimes \cdots \otimes a_q\otimes b_1\otimes \cdots \otimes b_p\\
    (\alpha_1\otimes \cdots \otimes \alpha_q)\odot (\beta_1\otimes \cdots \otimes \beta_p) &:= \alpha_1\otimes \cdots \otimes \alpha_n \otimes \beta_1\otimes \cdots \otimes \beta_m
\end{align*}
and extending by $\Real$-linearity. This clearly makes the dimensioned ring product satisfy, for $q,p>0$,
\begin{equation*}
    \odot: L^q \times L^p \to L^{q+p}, \qquad \odot: L^{-q} \times L^{-p} \to L^{-q-p}, \qquad \odot: L^0 \times L^0 \to L^0.
\end{equation*}
For products combining positive power tensors $a_1\otimes \cdots \otimes a_q$ and negative power tensors $\alpha_1\otimes \cdots \otimes \alpha_p$ we critically make use of the isomorphism $L^*\otimes L\cong \Real$ to define, without loss of generality for $p>q>0$,
\begin{equation*}
    (a_1\otimes \cdots \otimes a_q) \odot (\alpha_1\otimes \cdots \otimes \alpha_p) := \alpha_1( a_1) \cdots \alpha_q( a_q) \alpha_{p-q}\otimes \cdots \otimes \alpha_p.
\end{equation*}
It is then clear that the multiplication $\odot$ satisfies, for all $m,n\in \Int$,
\begin{equation*}
    \odot: L^m \times L^n \to L^{m+n}
\end{equation*}
and so it is compatible with the dimensioned structure of $L^\odot_\Int$. The multiplication $\odot$ is clearly associative and bilinear with respect to the addition on each homogeneous subset from the fact that the ordinary tensor product is associative and $\Real$-bilinear. Then it follows that $(L^\odot_\Int,+_\Int,\odot)$ is a commutative dimensioned ring. It only remains to show that non-zero elements of $L^\odot$ have multiplicative inverses. Note that a non-zero element corresponds to some non-vanishing tensor $0\neq h\in L^n$, but, since $L^n$ is a 1-dimensional vector space for all $n\in \Int$, we can find a unique $\eta\in (L^n)^*=L^{-n}$ such that $\eta(h)=1$. It follows from the above formula for products of positive and negative tensor powers that, in terms of the dimensioned ring multiplication, this becomes
\begin{equation*}
    h\odot \eta =1,
\end{equation*}
thus showing that all non-zero elements have multiplicative inverses, making the dimensioned ring $(L^\odot_\Int,+_\Int,\odot)$ into a dimensioned field.
\end{proof}

We now prove that the construction of the potential dimensioned field of a line is, in fact, functorial.

\begin{thm}[The Potential Functor for Lines] \label{PotentialIsAFunctorLine}
The assignment of the potential construction to a line is a functor
\begin{equation*}
\normalfont
    \odot: \Line \to \textsf{DimRing}.
\end{equation*}
Furthermore, a choice of unit in a line $L\in \Line$ induces a choice of units in the dimensioned field $(L^\odot_\Int,+_\Int,\odot)$ which, since $L^0=\Real$, then gives an isomorphism with the trivial dimensioned field
\begin{equation*}
    L^\odot \cong \Real \times \Int.
\end{equation*}
\end{thm}
\begin{proof}
To show functoriality we need to define the potential of a factor of lines $B:L_1\to L_2$
\begin{equation*}
    B^\odot:L_1^\odot \to L_2^\odot.
\end{equation*}
This can be done explicitly in the obvious way, for $q>0$
\begin{align*}
    B^\odot|_{L^q} &:= B\otimes \stackrel{q}{\cdots} \otimes B :L_1^q\to L_2^q\\
    B^\odot|_{L^0} &:= \Id_{\Real}:L_1^0\to L_2^0\\
    B^\odot|_{L^{-q}} &:= (B^{-1})^*\otimes \stackrel{q}{\cdots} \otimes (B^{-1})^* :L_1^{-q}\to L_2^{-q}
\end{align*}
where we have crucially used the invertibility of the factor $B$. By construction, $B^\odot$ is compatible with the $\Int$-dimensioned structure and since $B$ is a linear map with linear inverse, all the tensor powers act as $\Real$-linear maps on the homogeneous sets, thus making $B^\odot:L_1^\odot \to L_2^\odot$ into a morphism of abelian dimensioned groups. Showing that $B^\odot$ is a dimensioned ring morphism follows easily by the explicit construction of the dimensioned ring multiplication $\odot$ given in proposition \ref{DimRingPotential} above. This is checked directly for products that do not mix positive and negative tensor powers and for mixed products it suffices to note that
\begin{equation*}
    B^\odot (\alpha)\odot B^\odot (a)=(B^{-1})^*(\alpha) \odot B(a)= \alpha (B^{-1}(B(a)))= \alpha(a) = \Id_{\Real}(\alpha(a)) =B^\odot (\alpha\odot a).
\end{equation*}
It follows from the usual properties of tensor products in vector spaces that for another factor $C:L_2\to L_3$ we have
\begin{equation*}
    (C\circ B)^\odot=C^\odot \circ B^\odot, \qquad (\Id_L)^\odot=\Id_{L^\odot},
\end{equation*}
thus making the potential assignment into a functor. Recall that a choice of unit in a line $L\in \Line$ is simply a choice of non-vanishing element $u\in L^\times$. In proposition \ref{DimRingPotential} we saw that $L^\odot$ is a dimensioned field, so multiplicative inverses exist, let us denote them by $u^{-1}\in (L^*)^\times$. Using the notation for $q>0$
\begin{align*}
    u^q &:= u\odot \stackrel{q}{\cdots} \odot u\\
    u^0 &:= 1\\
    u^{-q} &:= u^{-1}\odot \stackrel{q}{\cdots} \odot u^{-1},
\end{align*}
it is clear that the map
\begin{align*}
u: \Int & \to L^\odot\\
n & \mapsto u^n
\end{align*}
satisfies
\begin{equation*}
    u^{n+m}=u^n\odot u^m.
\end{equation*}
By construction, all $u^n\in L^n$ are non-zero, so $u:\Int \to L^\odot$ is a choice of units in the dimensioned field $(L^\odot_\Int,+_\Int,\odot)$. The isomorphism of dimensioned fields $L^\odot \cong \Real \times \Int$ follows from proposition \ref{UnitsDimFields} and the observation that, by definition, $(L^\odot)_0=L^0=\Real$.
\end{proof}

Given a collection of lines $L_1,\dots,L_k\in\Line$, the above constructions generalize to the following notion of \textbf{potential}:
\begin{equation*}
    (L_1,\dots,L_k)^\odot:=\bigcup_{n_1,\dots n_k\in \Int} L_1^{n_1}\otimes \cdots \otimes L_k^{n_k},
\end{equation*}
which has a natural abelian dimensioned group structure given by $\Real$-linear addition and has dimension group $\Int^k$. The dimensioned filed structure generalizes in the obvious way:
\begin{equation*}
    (a_1\otimes \dots \otimes a_k) \odot (b_1\otimes \dots \otimes b_k):= a_1\odot b_1 \otimes \dots \otimes a_k\odot b_k
\end{equation*}
thus making $((L_1,\dots,L_k)^\odot_{\Int^k},+_{\Int^k},\odot)$ into a dimensioned field. Note that the potentials of each individual line $L_i$ can be found as dimensioned subfields $L_i^\odot\subset (L_1,\dots,L_k)^\odot$ since they are simply the dimensional preimages of the natural subgroups $\Int\subset \Int^k$. Furthermore, a choice of unit in each of the individual lines $u_i\in L_i^\times$ naturally induces a choice of units for the potential in a natural way
\begin{align*}
U: \Int^k & \to (L_1,\dots,L_k)^\odot\\
(n_1,\dots,n_k) & \mapsto u_1^{n_1}\odot \cdots \odot u_k^{n_k}.
\end{align*}

L-vector spaces were introduced in Section \ref{CategoryOfLines} as the natural generalization of vector spaces when lines are interpreted as ``unit-free'' fields of numbers. That interpretation, however, was only partial, since a proper generalization of notion of vector space would have to include a module structure with respect to the generalization of the field of numbers. In Section \ref{DimRingsModules} it was argued that the notion of dimensioned field indeed captures this generalization and dimensioned modules were introduced in a natural way. The technology developed in this chapter so far allows us to give a proper generalization of vector space in this sense via the notion of \textbf{dimensioned vector space over a dimensioned field} defined simply as a dimensioned module over a dimensioned field. These, together with dimensioned morphisms, form a category, $\textsf{DimVect}$, and the general notions of dimensioned modules introduced in Section \ref{DimRingsModules} apply. In particular, if a dimensioned field $F_G$ is fixed, the subcategory of dimensioned vector spaces over it, $\textsf{DimVect}_{F_G}$, becomes an abelian monoidal category, in complete analogy with the category of ordinary vector spaces over a fixed field.\newline

Our original claim that L-vector bundles represented a valid line generalization of ordinary vector spaces is fully justified by the fact that the datum of a L-vector space $V^L\in\LVect$ gives a dimensioned vector space. This is accomplished explicitly by the construction of the \textbf{potential} of $V^L\in\LVect$ which is defined as an abelian dimensioned group with dimensions in $\Int$ in a natural way:
\begin{equation*}
    V^{\odot L}:=\bigcup_{n\in \Int} L^n\otimes V.
\end{equation*}
This abelian dimensioned group carries an obvious dimensioned $L^\odot$-module structure that can be defined explicitly by
\begin{equation*}
    a\cdot (b\otimes v) := (a\odot b)\otimes v 
\end{equation*}
for all $a\in L^n$, $b\in L^m$ and $v\in V$, and then extending by linearity. The next proposition shows that. much like in the case of the potential of lines, the potential construction of L-vector spaces is functorial.

\begin{prop}[The Potential Functor for L-Vector Spaces] \label{PotentialIsAFunctorLVect}
The assignment of the potential construction to a L-vector space is a monoidal functor
\begin{equation*}
\normalfont
    \odot:\LVect\to \textsf{DimVect}
\end{equation*}
compatible with duality. Furthermore, fixing a line $L\in\Line$, the potential assignment
\begin{equation*}
\normalfont
    \odot:\LVect_L\to \textsf{DimVect}_{L^\odot}
\end{equation*}
becomes an abelian functor.
\end{prop}
\begin{proof}
To prove functoriality we give the explicit construction of the potential of a L-vector space morphism $\psi^B:V_1^{L_1}\to V_2^{L_2}$ as follows:
\begin{equation*}
    \psi^{\odot B} (a\otimes v):=B^\odot(a)\otimes \psi(v)
\end{equation*}
for all $a\in L^n$ and $v\in V$ and extending by linearity. This clearly clearly makes $\psi^{\odot B}$ into an abelian dimensioned group morphism with identity dimension map
\begin{equation*}
\begin{tikzcd}
V_1^{\odot L_1} \arrow[r, "\psi^{\odot B}"] \arrow[d] & V_2^{\odot L_2} \arrow[d] \\
\Int \arrow[r, "\Id_\Int"'] & \Int
\end{tikzcd}
\end{equation*}
which, by construction, interacts with the dimensioned module morphism as follows
\begin{equation*}
    \psi^{\odot B} (b\cdot (a\otimes v))=B^\odot(b)\cdot \psi^{\odot B} (a\otimes v).
\end{equation*}
Note that this last expression is the natural generalization of the $\Real$-linearity of ordinary vector spaces where the dimensioned ring isomorphism $B^\odot$ will be replaced by the particular case of the identity map. Functoriality of $\psi^{\odot B}$ then simply follows by the functoriality of $B^\odot$ proved in theorem \ref{PotentialIsAFunctorLine} and the usual composition of linear maps between ordinary vector spaces. Recall that the L-tensor product is defined as $V_1^{\odot L_1} \otimes V_1^{\odot L_1} :=(V_1\otimes V_2)^{L_1\otimes L_2}$. We can take the potential of the two lines $(L_1L_2)^\odot$ and define
\begin{equation*}
    (V_1\otimes V_2)^{\odot(L_1\otimes L_2)}:=\bigcup_{n_1,n_2\in \Int} L_1^{n_1} \otimes L_2^{n_2}\otimes V
\end{equation*}
which is a dimensioned vector space over the dimensioned ring $(L_1L_2)^\odot$. Tensor products of L-vector bundle morphisms are clearly sent to tensor products of dimensioned morphism via
\begin{equation*}
    \odot: \psi^B \otimes \varphi^C \mapsto (\psi \otimes \varphi)^{\odot(B\otimes C)},
\end{equation*}
and thus we see that the potential functor is indeed monoidal. Recall that the L-dual of a L-vector space is defined as $V^{*L}:=(V^*\otimes L)^L$, then we observe that, after applying the usual canonical isomorphisms to reorder tensor products, the potential of the L-dual will have homogeneous dimension sets shifted by +1. This precisely corresponds to elements of $(V^{\odot L})^*$ being $L^\odot$-linear maps of the form $V^{\odot L}\to L^\odot$. Exploiting once more the invertibility of line factors, we easily see that the potential of the dual of a L-vector bundle morphism $\psi^{*B}:V_2^{*L_2}\to V_1^{*L_1}$ is an isomorphism of dimensioned vector spaces $\psi^{\odot(*B)}:(V_2^*)^{\odot L_2}\to (V_1^*)^{\odot L_1}$, then we see that the potential construction is compatible with L-duality. Lastly, when we fix a line $L\in\Line$ L-direct sums, subobjects and kernels become well-defined, it is then clear from the $\Real$-linearity of all the maps involved, these are preserved under the potential construction, thus showing that it is an abelian functor.
\end{proof}

\section{Measurand Spaces Revisited} \label{MeasurandSpacesRevisited}

Let us now return to the original question of the formal description of physical quantities and units of measurement but now with the machinery of dimensioned algebra developed in sections \ref{DimRingsModules}, \ref{DimAlgebras} and \ref{PotentialFunctor} at our disposal. We begin by pointing out that, it should have become obvious by now, the notion of dimensioned field indeed captures the general formal structure of physical quantities. Any concrete example of a scientific model involving quantitative measurements, such us the example of classical thermodynamics presented in the opening of this chapter (\ref{DimAlgebraGeometry}), will involve some finite number $k$ of basic units so the set of all theoretically possible physical quantities form a dimensioned ring $R_{\Int^k}$ with $R_0=\Real$.\newline

In light of the results of Section \ref{PotentialFunctor}, we have not only successfully recovered the algebraic structure of physical quantities, but also given a complete mathematical legitimization to the empirically-motivated definition of measurand spaces proposed in Section \ref{MeasurandFormalism}. Indeed, a physical theory will consist of a collection of basic measurable properties, what we called \textbf{base measurands}, that are mathematically identified with a collection of lines $L_1,\dots,L_k\in \Line$. The potential functor $\odot:\Line\to \textsf{DimRing}$ now gives a mathematically precise meaning to what was defined as the \textbf{measurand space} $M$ of the physical theory:
\begin{equation*}
    M=((L_1,\dots,L_k)^\odot_{\Int^k},+_{\Int^k},\odot)
\end{equation*}
which, following from theorem \ref{PotentialIsAFunctorLine}, has the structure of a dimensioned field. We are now in the position to give a precise mathematical definition of \textbf{physical quantity} $Q$ simply as an element in the potentil:
\begin{equation*}
    Q\in M=(L_1,\dots,L_k)^\odot_{\Int^k}.
\end{equation*}

The term choice of units for a monoid splitting of the dimension projection of a dimensioned ring $u:G\to R_G$ was introduced in Section \ref{DimRingsModules} in anticipation of the metrological interpretation that we now give. Units of measurement in applied science an engineering serve as the reference scale for all the measurements of a physical quantity of the same kind. Mathematically, this is conventionally represented by assigning the numerical value $1$ to the measurement of the physical quantity applied to the unit of measurement itself, hence the name \emph{unit} of measurement. Since a choice of units assigns an non-zero element to each set of homogeneous dimension, which in the case of a measurand space are simply lines, they form a basis for that set and their component expression is conventionally also given by the numerical value $1$. We have thus connected the abstract notion of \textbf{choice of units} in a potential of some collection of lines with a \textbf{set of units of measurement} that will be used in a physical theory that takes those lines as base measurands.\newline

In dimensional analysis, it is common to recombine basic measurands of a physical theory and express them in terms of products of other measurands, such as in the example of Section \ref{MeasurandFormalism} where area was expressed as the product of lengths. In our formalism, where a physical theory is identified with a measurand space $M=(L_1,\dots,L_k)^\odot_{\Int^k}$, this is simply encapsulated by the notion of dimensioned isomorphism of the measurand space $\Psi: M\to M$ which, in general, will have non-trivial dimension map $\psi:\Int^k\to \Int^k$, corresponding to the recombinations of basic measurands. More concretely, changes in choices of units for the same measurand space are, in virtue of theorem \ref{PotentialIsAFunctorLine}, completely characterized by dimensioned ring isomorphisms of the trivial dimensioned ring
\begin{equation*}
    \Real\times \Int^k \to \Real\times \Int^k.
\end{equation*}

\section{Unit-Free Manifolds Revisited} \label{UnitFreeRevisited}

In Section \ref{UnitFreeManifolds} we argued that line bundles could be understood as a ``unit-free'' analogue of ordinary manifolds if one conceptually replaces the role played by the ring of smooth functions with the module of sections of a line bundle. The many results that were proved in that section showing the parallels with the ordinary theory of manifolds attest to the adequacy of this interpretation. However, using the module of sections of a line bundle as the algebraic analogue of the ring of smooth functions meant that the ring multiplication did not have a direct analogue. In this section we will show that this analogue appears as a dimensioned ring structure when we generalize the potential construction of Section \ref{PotentialFunctor} to line bundles.\newline

Let $\lambda:L\to M$ be a line bundle, in Section \ref{CategoryOfLineBundles} we saw that there is a monoidal structure in the restricted category of line bundles over the same base $M$. We can use the tensor product in this category in a completely analogous way to the tensor product in the category of lines $\Line$ and thus form positive and negative powers of the line bundle $L^n$, $n\in \Int$. The \textbf{potential} of the line bundle $L$ is then defined in a natural way
\begin{equation*}
    \Sec{L}^\odot:=\bigcup_{n\in \Int} \Sec{L^n}.
\end{equation*}
This set carries an obvious dimensioned structure with dimension set $\Int$ and the usual module structure on sections for each power $(\Sec{L^n},+_n)$ clearly makes $\Sec{L}^\odot$ into an abelian dimensioned group. Furthermore, the construction of the dimensioned ring product $\odot$ detailed in the proof of proposition \ref{DimRingPotential} can be reproduced in this case verbatim, thus making the potential of a line bundle into a dimensioned ring $(\Sec{L}^\odot_\Int,+_\Int,\odot)$. We note that this dimensioned ring encapsulates the usual algebraic structures found in sections of line bundles: indeed, the dimensionless ring of $\Sec{L}^\odot$ is the ordinary ring of functions of the base manifold $\Sec{L^0}=\Sec{\Real_M}\cong\Cin{M}$ and, for $f\in\Sec{L^0}=\Cin{M}$, $s\in\Sec{L^1}=\Sec{L}$ and $\sigma\in\Sec{L^{-1}}=\Sec{L^*}$, the dimensioned products
\begin{equation*}
    f \odot s = f\cdot s, \qquad \sigma \odot s = \sigma(s)
\end{equation*}
amount to the $\Cin{M}$-module map and the duality pairing, respectively. We now show that, as was the case for lines, the potential construction of line bundles is functorial.

\begin{prop}[The Potential Functor for Line Bundles] \label{PotentialFunctorLineBundles}
The assignment of the potential construction to a line bundle is a contravariant functor
\begin{equation*}
\normalfont
    \odot : \Line_\Man \to \textsf{DimRing}.
\end{equation*}
\end{prop}
\begin{proof}
Let us first define the potential of a factor between line bundles $B:L_1\to L_2$ covering a smooth map $b:M_1\to M_2$. We aim to define a dimensioned ring morphism of the form
\begin{equation*}
    B^\odot: (\Sec{L_2}^\odot_\Int,+_\Int,\odot)\to (\Sec{L_1}^\odot_\Int,+_\Int,\odot),
\end{equation*}
our definition will be, furthermore, of a dimensionless morphism, in the sense that it will cover the identity on the dimension group $\Id_\Int:\Int\to \Int$, so that it will suffice to provide a collection of maps between the sections of all the tensor powers $B^\odot_n:\Sec{L_2^n}\to \Sec{L_1^n}$. The datum provided by the line bundle factor $B$ allows to define three maps
\begin{align*}
    b^* &:\Cin{M_2}\to \Cin{M_1}\\
    B^* &:\Sec{L_2}\to \Sec{L_1}\\
    B^* &:\Sec{L_2^*} \to \Sec{L_1^*}
\end{align*}
where the first is simply the pull-back of the smooth map between base manifolds, the second is the pull-back of sections induced by a factor of line bundles defined point wise by
\begin{equation*}
    B^*(s_2)(x):=B_x^{-1}(s_2(b(x)))
\end{equation*}
for all $s_2\in\Sec{L_2}$, and the third is the usual pull-back of dual forms on general vector bundles, defined point-wise for a general bundle map by
\begin{equation*}
    B^*\sigma_2(s_1)(x):=\sigma_2(b(x))(B_x(s_1(x)))
\end{equation*}
for all $\sigma_2\in \Sec{L_2^*}$, $s_1\in\Sec{L_1}$. The maps $B^\odot_n$ are then defined simply as the tensor powers of these pull-backs. Contravariance then follows directly from contravariance of the pull-backs. It is then clear by construction that $B^\odot$ so defined acts as a dimensioned ring morphism for products of positive or negative tensor powers, then it only remains to show that it also acts as such for mixed products of tensor powers. This is readily checked by considering the following observation for sections $s_2\in\Sec{L_2}$ and $\sigma_2\in\Sec{L_2^*}$:
\begin{align*}
    B^\odot\sigma_2\odot B^\odot s_2(x) &=B^*\sigma_2\odot B^*s_2(x) =\\ &=\sigma_2(b(x))(B_xB_x^{-1}(s_2(b(x))))=\sigma_2(b(x))(s_2(b(x)))=b^*(\sigma_2(s_2))(x)=B^\odot(\sigma_2\odot s_2)(x).
\end{align*}
\end{proof}

This last proposition provides the key result for the legitimization of the interpretation of line bundles as unit-free manifolds since we notice the similarity of the potential functor above with the ordinary contravariant functor given by the assignment of the ring of smooth functions to a manifold
\begin{equation*}
    \text{C}^\infty:\Man \to \Ring.
\end{equation*}
The potential functor $\odot$ is, in fact, a direct generalization of $\text{C}^\infty$ as ordinary rings can be regarded as dimensioned rings with dimension set the trivial monoid.\newline

Due to possible topological constraints, the notion of unit of a line, i.e. a non-vanishing element, can be recovered only locally for line bundles. Let $\lambda:L\to M$ be a line bundle and $U\subset M$ an open subset, the potential construction is clearly natural with respect to restrictions since the same prescription used for global sections can be used to define $\Sec{L|_U}^\odot$. Defining the positive and negative powers of $u$ as it was done for the line case, it is clear that a local unit induces a choice of units for the local potential
\begin{equation*}
    u:\Int \to \Sec{L|_U}^\odot.
\end{equation*}
It then follows from the second part of theorem \ref{PotentialIsAFunctorLine} that a local unit $u$ induces an isomorphism of the local potential with the trivial dimensioned ring of local functions with dimension set $\Int$:
\begin{equation*}
    \Sec{L|_U}^\odot\cong \Cin{U}\times \Int.
\end{equation*}

To further establish the idea that that potentials of line bundles generalize rings of functions of manifolds, we will discuss the generalization of two aspects of ordinary manifolds for which the ring of functions proves a convenient algebraic tool: derivations and vanishing ideals of submanifolds.\newline

Recall from Section \ref{DimRingsModules} that the derivations of a dimensioned ring $\Dr{R_G}$ form a dimensioned module with dimension set $G$ that contains the derivations of the dimensionless ring $\Dr{R_e}$ as a Lie subalgebra of dimensionless derivations. In the case of the potential of a line bundle $\Sec{L}^\odot$, this implies that derivations of smooth functions, or, equivalently, vector fields, are recovered as a Lie subalgebra of the dimensionless derivations
\begin{equation*}
    \Dr{\Cin{M}}\cong\Sec{\Tan M}\subset \Dr{\Sec{L}^\odot}.
\end{equation*}
In Section \ref{UnitFreeManifolds} we argued that the line bundle generalization of the tangent bundle was the der bundle whose sections are the line bundle derivations. The next proposition shows that the derivations of the potentil of a line bundle naturally include the line bundle derivations.

\begin{prop}[Dimensionless Potential Derivations] \label{DimensionlessPotentialDerivations}
Let $\lambda:L\to M$ be a line bundle and $\Sec{L}^\odot$ its potential, then there is an isomorphism of Lie algebras
\begin{equation*}
\normalfont
    \Dr{L}\cong \Dr{\Sec{L}^\odot}_0.
\end{equation*}
\end{prop}
\begin{proof}
We can give the isomorphism by explicitly specifying two maps. The first sends a dimensionless derivation to its restriction on the homogeneous subsets of dimension $0$ and $1$
\begin{equation*}
    \Dr{\Sec{L}^\odot}_0\ni P\mapsto (P|_{L^0},P|_{L^1})=:(X,D)
\end{equation*}
Clearly, from the fact that $P$ is a $\odot$-derivation these satisfy
\begin{align*}
    X(fg) &:=X(f\odot g)=X(f)\odot g+f\odot X(g)=X(f)g+fX(g)\\
    D(f\cdot s) &:= D(f\odot s)=X(f)\odot s+f\odot D(s)=X(f)\cdot s + f\cdot D(s),
\end{align*}
thus showing that $D$ is a line bundle derivation with symbol $X$. Conversely, given a line bundle derivation $D\in \Dr{L}$ with symbol $X$, we need to define a dimensionless derivation $\Dr{\Sec{L}^\odot}_0$. This is accomplished by extending $D$ as a $\odot$-derivation for non-negative tensor powers of $L$ the following basic identities
\begin{align*}
    P(f\odot g) &:= X[fg]=X[f]g+fX[g]\\
    P(f\odot s) &:= D(f\cdot s)=X[f]\cdot s + f\cdot D(s)\\
    P(s\odot r) &:= P(s)\odot r + s\odot P(s).
\end{align*}
To account for negative tensor powers we use proposition \ref{DerBundleDual}, which gives an isomorphism $\Dr{L}\cong \Dr{L^*}$ thus defining a derivation $\Delta\in\Dr{L^*}$ from $D$. This derivation $\Delta$ is extended as a $\odot$-derivation for non-positive tensor powers in a similar wat to $D$. To complete the extension of $D$ to $P$ as a $\odot$-derivation, it only remains to consider products mixing positive and negative tensor powers. This case is accounted for by the following consistency formula that follows from the definition of the dimensioned product in the potentil and the explicit isomorphism $\Dr{L}\cong \Dr{L^*}$ in the proof of proposition \ref{DerBundleDual}:
\begin{equation*}
    P(\sigma \odot s) = P(\sigma(s))=X(\sigma(s))= D^*(\sigma)(s)+\sigma(D(s))=D^*(\sigma)\odot s + \sigma \odot D(s)=P(\sigma)\odot s+ \sigma \odot P(s).
\end{equation*}
These two maps, which are clearly Lie algebra morphisms since the bracket is simply the commutator, are readily checked to be inverses of each other.
\end{proof}
We remark that dimensionless derivations do not determine all the derivations of a dimensioned ring. General derivations of the potentil of a line bundle are given by collections of differential operators between all the different tensor powers fitting consistently with the $\Int$-dimensioned structure.\newline

Consider now a submanifold $i:S\hookrightarrow M$ of a line bundle $\lambda:L\to M$. We saw in Section \ref{CategoryOfLineBundles} that a line bundle is induced on $S$ by pull-back with an inclusion factor covering the embedding $\iota:L_S\to L$. There, the set of vanishing sections on $S$, defined formally as the kernel of $\iota$, was shown to be a submodule of the sections of the ambient line bundle $\Gamma_S\subset \Sec{L}$ that can be seen as (locally) generated by the ideal of vanishing functions $I_S\subset\Cin{M}$. The following proposition shows that these two algebraic manifestation of a submanifold in a line bundle fit nicely into the potentil picture.

\begin{prop}[Vanishing Dimensioned Ideal of a Submanifold] \label{DimVanishingIdeal}
Let $\lambda:L\to M$ be a line bundle and $i:S\hookrightarrow M$ a submanifold carrying the restricted line bundle $L_S$. Let us denote the line bundle potentials by $\Sec{L}^\odot$ and $\Sec{L_S}^\odot$, then the submanifold defines a dimensioned ideal $I_S\subset \Sec{L}^\odot$ that allows to characterize (depending on the embedding $i$, perhaps only locally) the restricted potential as a quotient of dimensioned rings
\begin{equation*}
    \Sec{L_S}^\odot\cong \Sec{L}^\odot/I_S.
\end{equation*}
\end{prop}
\begin{proof}
The vanishing  dimensioned ideal is simply defined as the set of sections of all the tensor powers that vanish when restricted to $S$, that is
\begin{equation*}
    I_S:=\{a\in \Sec{L^n}|\quad a(x)=0\in L^n_x \quad \forall x\in S\}.
\end{equation*}
Note that this definition is indeed equivalent to the kernel of the potential of the inclusion factor $I_S=\Ker{\iota^\odot}$, the dimensioned ideal condition $\Sec{L}^\odot\odot I_S\subset I_S$ then follows:
\begin{equation*}
    \iota^\odot(r_n\odot a_m)= (\iota^*)^n r_n \odot (\iota^*)^m a_m=(\iota^*)^n r_n \odot 0=0
\end{equation*}
for all $r_n\in \Sec{L^n}$, $a_m\in I_S$, where the functoriality of the potential construction proved in proposition \ref{PotentialFunctorLineBundles} has been used. It is clear that, by construction, the ordinary ideal of vanishing functions is the dimensionless component of $I_S$ and that
\begin{equation*}
    \Gamma_S=I_S\cap \Sec{L^1}.
\end{equation*}
Using a local argument we can see that, similarly to the submodule of vanishing sections, the subsets of homogeneous dimension of $I_S$ can all be generated by elements in the dimensionless component. Just as in the case of the vanishing submodule of sections, it is clear that the quotient $\Sec{L}^\odot/I_S$ represents the identification of sections of the tensor powers of $L_S$ with extensions in $L$ that differ by a vanishing section of the corresponding tensor power, which then gives the desired the result.
\end{proof}

\section{Jacobi Manifolds Revisited} \label{JacobiManifoldsRevisited}

In Section \ref{ContactGeometry} Jacobi manifolds were presented as the natural generalization of Poisson structures on unit-free manifolds. In light of the results of Section \ref{UnitFreeRevisited} above connecting unit-free manifolds to dimensioned algebra, it is natural to ask whether Jacobi structures on line bundles are somehow reflected on the dimensioned rings associated to them under the potential functor. In this section we will see that this is indeed the case as we are now in the position to prove the main theorems of this chapter connecting Jacobi manifolds to dimensioned Poisson algebras in a natural way.\newline

We begin by identifying a dimensioned Poisson algebra structure on the potential of a Jacobi manifold.

\begin{thm}[Dimensioned Poisson Algebra associated with a Jacobi Manifold] \label{DimPoissonAlgebraJacobi}
Let $\lambda:L\to M$ be a line bundle and $(\Sec{L},\{,\})$ a Jacobi structure, then there exists a unique dimensioned Poisson algebra of dimension $-1$ on the potential dimensioned ring $(\Sec{L}^\odot_\Int,+_\Int,\odot_0,\{,\}_{-1})$ such that the brackets combining elements of dimensions $+1$, $0$, and $-1$ are determined by the Jacobi bracket $\{,\}$, its symbol $X$ and its squiggle $\Lambda$.
\end{thm}
\begin{proof}
We give an explicit construction of the dimensioned Poisson algebra on the potential $(\Sec{L}^\odot_\Int,+_\Int,\odot_0)$, that we regard here as a dimensioned commutative algebra over the real numbers with dimension set $\Int$ and dimensionless commutative multiplication $\odot_0$. Since the Jacobi bracket maps pairs of sections into sections $\{,\}:\Sec{L}\times \Sec{L}\to \Sec{L}$, we aim to extend it to all the tensor powers of the potential as a dimensioned algebra bracket of dimension $-1\in\Int$:
\begin{equation*}
    \{,\}_{-1}:\Sec{L^n}\times \Sec{L^m}\to \Sec{L^{n+m-1}}.
\end{equation*}
It is clear that we obtain a partial Lie bracket for all positive tensor powers simply by extending the Jacobi bracket as $\odot$-derivations in each argument, i.e. setting $\{a,b\}_{-1}:=\{a,b\}$ and generating all the brackets between higher powers from the basic identity:
\begin{equation*}
    \{a,b\odot c\}_{-1}:=\{a,b\}\odot c +b\odot \{a,c\}
\end{equation*}
for all $a,b,c\in \Sec{L^1}=\Sec{L}$. Note that this is analogous to using the isomorphism $\Dr{L}\cong\Dr{\Sec{L}^\odot}_0$ of proposition \ref{DimensionlessPotentialDerivations} to regard the Hamiltonian derivation of the Jacobi bracket $D_a$ as a dimensionless derivation of the potential. The symbol-squiggle identity of the Jacobi bracket written in terms of the potential dimensioned multiplication reads
\begin{equation*}
    \{f\odot a ,g\odot b\}_{-1}=f\odot g \odot \{a,b\} + f\odot X_a[g]\odot b - g\odot X_b[f]\odot a + \Lambda(df\otimes a, dg \otimes b)
\end{equation*}
for $f,g\in \Sec{L^0}\cong\Cin{M}$, $a,b\in \Sec{L^1}=\Sec{L}$, then we can extract the definition of the dimensioned Poisson bracket for non-negative tensor powers by reading off the above formula interpreted as a Leibniz rule of the $\odot$ multiplication:
\begin{align*}
    \{a,f\}_{-1} &:= X_a[f] = -\{f,a\}_{-1}\\
    \{f,g\}_{-1}(a) &:= \Lambda^\sharp( df\otimes a)[g] = -\{g,f\}_{-1}(a)
\end{align*}
Note that the second bracket has been defined on a generic argument since
\begin{equation*}
    \{,\}_{-1}:\Sec{L^0}\times \Sec{L^0}\to \Sec{L^{-1}}=\Sec{L^*}.
\end{equation*}
To define the bracket on negative tensor powers we first use the isomorphism $R:\Dr{L}\to\Dr{L^*}$ proved in proposition \ref{DerBundleDual} to define the Hamiltonian derivation on dual sections $\Delta_a:=R(D_a)\in\Dr{L^*}$ and set
\begin{equation*}
    \{a,\alpha\}_{-1}:= \Delta_a(\alpha)=-\{\alpha,a\}_{-1}.
\end{equation*}
Note that this definition is consistent with the previous definitions of brackets of non-negative tensor powers as we readily check that it acts as a $\odot$-derivation in both arguments:
\begin{equation*}
    \{a,\alpha \odot b\}_{-1}=\{a,\alpha(b)\}_{-1}=X_a[\alpha(b)]=\Delta_a(\alpha)(b)+\alpha(D_a(b))=\{a,\alpha\}_{-1}\odot b+\alpha \odot \{a,b\}_{-1}.
\end{equation*}
With the brackets defined so far for non-negative tensor powers and the mixed bracket above, we can expand the expression $\{f\odot a,\alpha \odot b\}$ by $\odot$-derivations (full details of the computation shown in appendix \ref{CalculationsDimensionedPoissonBracket}) to find the only non-yet defined bracket:
\begin{equation*}
    \{f,\alpha,\}_{-1}(a,b):=\Lambda^\sharp(df\otimes a)[\alpha(b)]+X_b[f]\alpha(a).
\end{equation*}
Similarly, expanding the bracket $\{\alpha \odot a,\beta \odot b\}$ (again, full details in appendix \ref{CalculationsDimensionedPoissonBracket}) we find:
\begin{equation*}
    \{\alpha,\beta\}_{-1}(a,b,c):=\Lambda^\sharp(d\alpha(a)\otimes b)[\beta(c)]+X_c[\alpha(a)]\beta(b)-\alpha(b)X_a[\beta(c)]+\alpha(b)\beta(\{a,c\}).
\end{equation*}
With these partial brackets we can now define the brackets of combinations of positive and negative tensor powers via extension as $\odot$-derivations. Clearly, following from the observation made at the end of the proof of proposition \ref{PotentialFunctorLineBundles}, the Jacobi identity of bracket for the negative tensor powers so defined will be directly dependent on the Jacobi identity for the bracket at dimensions $+1$ and $0$. In appendix \ref{CalculationsDimensionedPoissonBracket} the Jacobi identities for brackets of all the combinations of the tensor powers $+1$ and $0$ are shown to follow directly from the basic identities satisfied by the symbol $X$ and squiggle $\Lambda$ of the original Jacobi structure.
\end{proof}

The next theorem shows that, much like how the functor $\text{C}^\infty:\Man \to \Ring$ characterizes Poisson manifolds as a subcategory of Poisson algebras, the potential functor allows to regard the category of Jacobi manifolds as a subcategory of dimensioned Poisson algebras.

\begin{thm}[The Potential Functor for Jacobi Manifolds] \label{PotentialFunctorJacobi}
The assignment of the potential of a line bundle restricted to the category of Jacobi manifolds with Jacobi maps gives a contravariant functor
\begin{equation*}
\normalfont
    \odot: \textsf{Jac}_\Man\to \textsf{DimPoissAlg}.
\end{equation*}
\end{thm}
\begin{proof}
In proposition \ref{PotentialFunctorJacobi} it was shown that a line bundle factor $B:L\to L'$ covering a smooth map $b:M\to M'$. is mapped to a dimensioned ring morphism $B^\odot :\Sec{L'}^\odot\to \Sec{L}^\odot$ under the potential contravariant functor, then it will suffice to show that when $B$ is a Jacobi map, i.e.
\begin{equation*}
    B^*\{a,b\}'=\{B^*a',B^*b'\}
\end{equation*}
for all $a,b\in \Sec{L'}$, then the potential map is a dimensioned Lie algebra morphism, i.e. 
\begin{equation*}
    B^\odot\{s,r\}_{-1}'=\{B^\odot s,B^\odot r\}_{-1}
\end{equation*}
for all $s,r\in \Sec{L'}^\odot$. Note that $B^\odot$ was defined in proposition \ref{PotentialFunctorLineBundles} as the tensor powers of the pull-backs of sections of $L'$ , its dual and the smooth functions, then it is clear from the fact that the dimensioned Lie brackets $\{,\}_{-1}$, $\{,\}_{-1}'$ are defined by extension as $\odot$-derivations that the dimensioned Lie algebra morphism condition for brackets of positive and negative tensor powers is dependent on the same condition for all the bracket combinations of elements in dimensions $+1$, $0$ and $-1$. These conditions are checked directly using the definitions. For the bracket of a pair of elements of dimension $+1$, the condition of dimensioned Lie algebra morphism for $B^\odot$ is precisely the condition that $B$ is a Jacobi map. By considering a bracket of the form $\{a,f\cdot b\}'$ we can see that the fact that $B$ is a Jacobi map and the basic properties of the pull-backs of line bundle factors imply the morphism condition for the bracket of elements of dimension $+1$ and $0$:
\begin{equation*}
    B^\odot\{a,f\}_{-1}'=b^*X_a'[f]= X_{B^*a}[b^*f]=\{B^*a,b^*f\}_{-1}=\{B^\odot a,B^\odot f\}_{-1}.
\end{equation*}
From similar considerations for a bracket of the form $\{f\cdot a,g\cdot b\}$, it follows that
\begin{equation*}
    B^\odot\{f,g\}_{-1}'(c,d)=B^*\Lambda'(df \otimes c,dg\otimes d)=\Lambda(db^*f \otimes c,db^*g\otimes d)=\{b^*f,b^*g\}_{-1}(c,d)=\{B^\odot f,B^\odot g\}_{-1}(c,d)
\end{equation*}
for all $c,d\in \Sec{L}$. To account for brackets containing elements of dimension $-1$ we first consider the defining formula of the isomorphism $\Dr{L}\cong \Dr{L^*}$ under pull-back
\begin{align*}
B^*(\Delta_a'(\alpha)(c))&= b^*X_a'[\alpha(c)]- b^*\alpha(\{a,c\}')\\
&= X_{B^*a}[b^*\alpha(c)]-B^*\alpha(B^*\{a,c\})\\
&= X_{B^*a}[B^*\alpha(B^*c)]-B^*\alpha(\{B^*a,B^*c\})\\
&= \Delta_{B^*a}(B^*\alpha)(B^*c).
\end{align*}
Which clearly implies, in particular, the dimensioned Lie morphism condition for the bracket of mixed tensor powers
\begin{equation*}
    B^\odot\{a,\alpha\}_{-1}'=B^*\Delta_a'(\alpha)=\Delta_{B^*a}(B^*\alpha)=\{B^*a,B^*\alpha\}_{-1}=\{B^\odot a,B^\odot\alpha\}_{-1}.
\end{equation*}
In the proof of theorem \ref{DimPoissonAlgebraJacobi} it was shown that the brackets $\{f,\alpha\}_{-1}$ and $\{\alpha,\beta\}_{-1}$ were determined by extending the previously defined brackets between elements of dimensions $+1$, $0$ and $-1$ as $\odot$-derivations, thus the dimensioned Lie algebra morphism condition for these follows from the fact that $B^\odot$ is defined as the tensor powers of $B^*$.
\end{proof}

In our discussion of the problem of reduction for Jacobi manifolds in Section \ref{ContactGeometry}, vanishing submodules of coisotropic submanifolds were recognized to be analogous to the coisotropes of algebraic Poisson reduction, however, the lack of a commutative multiplication structure on the module of sections of a line bundle made it impossible to make the analogy more precise. Our identification of the dimensioned Poisson algebra on the potential of a Jacobi manifold does this precisely as the next proposition shows that coisotropic submanifolds induce dimensioned coisotropes.

\begin{prop} [Coisotropic Submanifolds induce Dimensioned Coisotropes] \label{CoisotropicSubmanifoldsDimensionedCoisotropes}
Let $\lambda:L\to M$ be a line bundle and $(\Sec{L},\{,\})$ a Jacobi structure, then the vanishing dimensioned ideal of a coistropic submanifold $i:S \hookrightarrow M$ is a dimensioned coisotrope of the dimensioned Poisson algebra on the potential
\begin{equation*}
    I_S\subset (\Sec{L}^\odot_\Int,+_\Int,\odot_0,\{,\}_{-1}).
\end{equation*}
\end{prop}
\begin{proof}
Proposition \ref{DimVanishingIdeal} shows that $I_S\subset \Sec{L}^\odot$ is a dimensioned $\odot$-ideal for any submanifold $S$, then it will suffice to show that $I_S$ is a dimensioned Lie subalgebra. The vanishing ideal is generated by the $\odot$-products of elements of dimension $+1$, $0$ and $-1$, then, by the Leibniz identity of the dimensioned Poisson bracket, it suffices to check the dimensioned Lie subalgebra conditions for elements of those dimensions. For this we will use characterizations 2, 3 and 4 of coisotropic submanifolds of a Jacobi manifold given in proposition \ref{CoisotropicSubmanifoldsJacobi}. Clearly the condition on the brackets of positive powers $\{a,b\}_{-1}$ for $a,b\in\Sec{L}$ is the fact that the vanishing submodule of sections of $S$ forms a Lie subalgebra of the Jacobi structure, characterization 3. For the bracket $\{a,f\}_{-1}=X_a[f]$ it is the fact that Hamiltonian vector fields of vanishing sections are tangent to the submanifold, i.e. $X_{\Gamma_S}[I_S]\subset I_S$, characterization 4. For the bracket $\{f,g\}_{-1}$ it is the condition $\Lambda^\sharp(dI_S\otimes L)\subset \Tan S$, characterization 2. From the observation that
\begin{equation*}
    a\in I_S\cap \Sec{L^1} \Rightarrow \alpha \odot a=\alpha(a)\in I_S
\end{equation*}
for any $\alpha\in\Sec{L^*}$, we check the dimensioned Lie subalgebra condition for the brackets $\{f,\alpha\}_{-1}$ and $\{\alpha,\beta\}_{-1}$ by writing the explicit defining formulas presented in the proof of theorem \ref{DimPoissonAlgebraJacobi} and using characterizations 2 and 4 from proposition \ref{CoisotropicSubmanifoldsJacobi} combined. 
\end{proof}

Lastly, when a coisotropic submanifold furthermore fits in a reduction scheme of Jacobi manifolds, the associated dimensioned Poisson algebras fit in a dimensioned algebra reduction scheme, again in direct analogy with the ordinary Poisson case.

\begin{thm}[Coisotropic Reduction induces Dimensioned Poisson Reduction] \label{CoisotropicReductionDimPoissonReduction}
Let $\lambda:L\to M$ be a line bundle with a Jacobi structure $(\Sec{L},\{,\})$ and let $i:S \hookrightarrow M$ be a coisotropic submanifold satisfying the assumptions of proposition \ref{CoisotropicReductionJacobi} so that there is a a reduced Jacobi structure $(\Sec{L'},\{,\}')$ fitting in the reduction diagram:
\begin{equation*}
    \begin{tikzcd}[sep=tiny]
    L_S \arrow[rr,"\iota"] \arrow[dd,"\pi"'] \arrow[dr]& & L \arrow[dr] & \\
    & S \arrow[rr,"i"', hook] \arrow[dd, "p",twoheadrightarrow] & & M \\
    L' \arrow[dr] & & & \\
    & M' &  & 
    \end{tikzcd}
\end{equation*}
then there is an isomorphism of dimensioned Poisson algebras between the potential of the reduced Jacobi structure and the algebraic dimensioned Poisson reduction by the vanishing dimensioned coisotrope:
\begin{equation*}
    \Sec{L'}^\odot\cong N(I_S)/I_S.
\end{equation*}
\end{thm}
\begin{proof}
Recall from Section \ref{ContactGeometry} that the Jacobi reducibility condition was given explicitly in terms of the brackets as
\begin{equation*}
    \pi^*\{a_1,a_2\}'=\iota^*\{A_1,A_2\}
\end{equation*}
for all $a_i\in \Sec{L'}$ and $A_i\in\Sec{L}$ extensions satisfying $\pi^*a_i=\iota^*A_i$. The definition of the potential of a line bundle factor of proposition \ref{PotentialFunctorLineBundles} and the explicit definition of the dimensioned bracket $\{,\}_{-1}$ clearly show that the reducibility condition translates into the potential setting verbatim as one finds that the dimensioned Poisson brackets on the potentials of $L$ and $L'$ are related by the following condition
\begin{equation*}
    \pi^\odot\{a_1,a_2\}'_{-1}=\iota^\odot\{A_1,A_2\}_{-1}
\end{equation*}
for all $a_i\in \Sec{L'}^\odot$ and $A_i\in\Sec{L}^\odot$ extensions satisfying $\pi^\odot a_i=\iota^\odot A_i$. We aim to relate the dimensioned Lie idealizer of the vanishing dimensioned ideal $N(I_S)$ to the submersion factor $\pi:L_S\to L'$ in a natural way. This will follow by the compatibility condition assumed in proposition \ref{CoisotropicReductionJacobi} for the coistropic submanifold:
\begin{equation*}
    \delta(\Ker{\Der \pi})=\Lambda^\sharp((\Tan S)^{\text{0}L})
\end{equation*}
which, exploiting the jet sequence of the Jacobi structure, can be rewritten as
\begin{equation*}
    \Ker{\Tan p}=(\tilde{\Lambda}^\sharp\circ i)(\Tan^0S \otimes L_S).
\end{equation*}
This equation gives the point-wise condition that the $p$-fibration on $S$ is a foliation integrating the tangent distribution of Hamiltonian vector fields of the vanishing sections. It follows from the observation made at the end of the proof of proposition \ref{CoisotropicSubmanifoldsDimensionedCoisotropes} that the Lie idealizer $N(I_S)$ is generated by the brackets of elements of dimension $+1$ and $0$, then it suffices to identify the elements satisfying the dimensioned Lie idealizer defining condition of these dimensions. These will be $f\in\Sec{L^0}\cong \Cin{M}$ and $s\in\Sec{L^1}=\Sec{L}$ such that
\begin{equation*}
    \{f,g\}_{-1}\in I_S, \qquad \{s,g\}_{-1}\in I_S, \qquad \{s,a\}_{-1}\in I_S
\end{equation*}
for all $g\in I_S\cap \Sec{L^0}$ and $a\in I_S\cap \Sec{L^1}$. From the explicit formulas given for the dimensioned bracket $\{,\}_{-1}$ given in the proof of theorem \ref{DimPoissonAlgebraJacobi} we clearly see that the compatibility condition of the coisotropic submanifold with the submersion factor gives a point-wise identification of elements in the idealizer and the infinitesimal description of the $p$-fibration with the restricted line bundles on $S$. Since, by assumption, the submersion factor fits in a reduction scheme of (smooth) Jacobi manifolds this infinitesimal identification carries over globally to allow the identification of $N(I_S)$ with the line bundle $\pi$-fibration. Quotienting by $I_S$, as seen in proposition \ref{DimVanishingIdeal}, amounts to restricting to the submanifold $S$, thus we find the equivalent description of the reduced bracket $\{,\}_{-1}'$ as the canonical dimensioned Poisson bracket on the quotient $N(I_S)/I_S$.
\end{proof}

\section{Dimensioned Manifolds} \label{DimManifolds}

In sections \ref{UnitFreeRevisited} and \ref{JacobiManifoldsRevisited} it was shown that the language of dimensioned algebra is particularly suited for the description of the natural algebraic structures associated with line bundles and Jacobi structures, particularly when these are understood as unit-free generalizations of ordinary manifolds and Poisson structures. Although, the results presented in these sections are indeed clear manifestations of what we could call dimensioned geometry in general, the notion of dimensioned vector space introduced in Section \ref{PotentialFunctor} constitutes another avenue for the discussion of dimensioned geometry. In this section we briefly touch upon the subject of finding a geometric implementation of the general ideas of dimensionality as motivated at the start of this chapter \ref{DimAlgebraGeometry}.\newline

With the machinery of dimensioned vector spaces developed in sections \ref{DimRingsModules} and \ref{PotentialFunctor} at our disposal, it would be possible, in principle, to develop a theory of calculus between dimensioned vector spaces in analogy with ordinary calculus on vector spaces. Natural notions of \textbf{dimensioned curves} and \textbf{dimensioned functions} can be defined as the ordinary analogues but taking into account that maps need to be consistent with the dimensioned structures. Indeed, ordinary curves and functions are recovered as \textbf{dimensionless curves} and \textbf{dimensionless functions} with dimensionless domain and dimensionless range, respectively. The notions of directional derivative, differential and smoothness of a dimensioned map between dimensioned vector spaces then follow naturally.\newline

An obvious approach to a definition of a \textbf{dimensioned manifold} is to consider a topological space $\mathcal{M}$ for which a cover of charts $(U_i,\varphi_i)$ exists such that the local coordinate maps $\varphi_i:U_i\to V_G$ are homeomorphisms into a model dimensioned vector space $V_G$ that give smooth invertible dimensioned maps of dimensioned vector spaces $\varphi_j\circ \varphi_i^{-1}: V_G\to V_G$ on overlaps $U_i\cap U_j\subset \mathcal{M}$. The usual notions of calculus such as curves, functions, smooth maps, submanifolds, etc. carry over from the calculus on dimensioned vector spaces just like in ordinary manifolds but with the added nuance of keeping compatibility with the dimensioned structure. We point out that in this general definition, the topological space of a dimensioned manifold need not be a dimensioned set itself, as the transition functions may have non-identity dimension maps.\newline

A detailed study of the theory of dimensioned manifolds is well beyond the scope of the present thesis, however, we do emphasize that unit-free manifolds can be recovered as a particular case of the proposed notion of dimensioned manifold. Consider a dimensioned manifold $\mathcal{M}$ with model dimensioned vector space given by the potential of a L-vector space $V^{\odot L}_\Int$ and whose transition functions $\varphi_j\circ \varphi_i^{-1}: V^{\odot L}_\Int\to V^{\odot L}_\Int$ are required to have identity dimension map $\Id_\Int:\Int\to \Int$. Hence, the local charts have transition functions that, when restricted to the dimensionless elements, give, in particular, a collection of ordinary transition functions $\varphi_j\circ \varphi_i^{-1}|_{V^{L^0}}: V\to V$ for the topological space $\mathcal{M}$, thus making $\mathcal{M}$ into an ordinary manifold. The dimensioned differential maps of the transition functions
\begin{equation*}
    \Tan_x (\varphi_j\circ \varphi_i^{-1}|_{V^{L^0}}): V^{\odot L}_\Int\to V^{\odot L}_\Int
\end{equation*}
are clearly dimensionless L-vector bundle isomorphisms for all $x\in U_i\cap U_j$. Now, since proposition \ref{PotentialIsAFunctorLVect} showed that the dimensionless morphisms of a potential of a L-vector space are given precisely by the potentials of of L-vector space morphisms, the differentials of the transition functions are specified by some L-vector space isomorphism $(\Phi_{ij}^{B^{ij}})_x:V^L\to V^L$, which in particular carries a line factor $B_x^{ij}:L\to L$. Factors are invertible maps on lines, then it follows from our discussion in Section \ref{CategoryOfLines} that they must satisfy the 2-out-of-3 identity point-wise over the intersection of three charts:
\begin{equation*}
    B_x^{ij}\circ B_x^{jk} \circ (B_x^{ik})^{-1}=B_x^{ij}\circ B_x^{ij} \circ B_x^{ki}=\Id_L
\end{equation*}
for all $x\in U_i\cap U_j\cap U_k$. Recognizing the cocycle condition in the identity above, we see that the $\GL{L}$-valued functions $B^{ij}$ clearly form the transition functions of a line bundle $\lambda: \mathcal{L}\to \mathcal{M}$ with typical fibre $L$. Since in Section \ref{UnitFreeManifolds} line bundles were argued to correctly encapsulate the notion of unit-free manifolds, we have thus shown that our proposed notion of dimensioned manifold encompasses unit-free manifolds.

\section{Dimensioned Mechanics} \label{DimMechanics}

In chapter \ref{ContactPhaseSpaces} line bundles were shown to adequately implement the notion of unit-free Hamiltonian phase space with an argument that relied on the categorical similarities between the cotangent functor and the jet functor, both regarded as examples of Hamiltonian functors from a category of configuration spaces into a category of Hamiltonian phase spaces. There is, nonetheless, a clear difference between the two when it comes to the algebraic manifestation of observables: in ordinary phase spaces observables are identified with smooth functions, thus carrying a ring structure naturally, whereas in unit-free phase spaces observables are identified with sections of a line bundle which form modules over rings of functions that carry no natural ring-like multiplication operation. We should emphasize at this point that this is not a minor mathematical detail, indeed, the ring structure of ordinary observables is crucial for the development of any realistic model in classical mechanics as the multiplication of physical quantities is at the foundation of how most mechanical theories are built, e.g. kinetic energy is a quadratic function on momenta, potentials depend on powers of position coordinates, magnetic terms involve products of velocities and position coordinates, etc. Furthermore, the Leibniz identity of the canonical Poisson bracket with respect to the ring multiplication is a key algebraic asset for the analysis of the integrability and dynamics of Hamiltonian systems. It may seem, then, that by attempting to introduce physical dimension into geometric mechanics via line bundles we are sacrificing one of its core features.\newline

If we are to claim to have successfully generalized conventional Hamiltonian mechanics into a unit-free framework, it is essential that clear unit-free counterparts of the ring multiplication of observables and the Leibniz rule of the canonical Poisson bracket are found. In light of the results of this chapter and, in particular, the identification of the dimensioned Poisson algebra on the potential of a Jacobi structure in theorem \ref{DimPoissonAlgebraJacobi}, this is achieved by taking the potential functor for line bundles as the observable functor. By replacing the section functor with the potential functor, the theory of canonical contact phase spaces defined in Section \ref{CanonicalContact} is now updated to give the following general categorical picture of classical mechanics incorporating physical dimension:
\begin{equation*}
\begin{tikzcd}[row sep=small]
 & \textsf{DimRing} \arrow[dd,"\text{Der}"']& & & \textsf{DimPoissAlg} \arrow[dd,"\text{Der}"] \\
 & & \Line_\Man \arrow[ul, "\odot"']\arrow[dl,"\Der"]\arrow[r,"\Jet^1"] & \Cont_\Man \arrow[ur,"\odot"] \arrow[dr,"\Der"'] & \\
 & \textsf{DimLie} & & & \textsf{DimLie} 
\end{tikzcd}
\end{equation*}
where the dynamics functor $\Der$ is regarded as the natural inclusion of the derivations of a line bundle in the dimensioned Lie algebra of the potential derivations, as seen in proposition \ref{DimensionlessPotentialDerivations}.\newline

The notion of unit-free configuration space as a line bundle $\lambda:L\to Q$ was introduced as the mathematical manifestation of the possible assignments of values of a physical quantity, i.e. a measurand in the sense of Section \ref{MeasurandFormalism} represented by the typical fibre of the line bundle $L$, to the set of static states of a mechanical system $Q$. In a realistic scenario, several measurands will be relevant for a complete description of a physical system and thus a collection of line bundles over the same configuration phase space $L_1,\dots,L_k$ must be considered. The construction of the potential of a collection of lines of Section \ref{PotentialFunctor} allows to naturally generalize the formalism of canonical contact mechanics to this more realistic setting. Although we do not investigate the matter in full detail in this thesis, we point out that it is possible to show that the potential functor for a collection of lines commutes with the jet functor that assigns to each tensor power of lines the associated canonical contact manifold. It is in this context that the Leibniz identity for the dimensioned Poisson bracket of observables in the unit-free framework is fully recovered.\newline

Another approach to the subject of introducing physical dimension into geometric mechanics will be to consider a category of configuration spaces whose objects are dimensioned manifolds in the sense of Section \ref{DimManifolds}. In that section we discussed how unit-free manifolds can be regarded as special cases of dimensioned manifolds, therefore, this direct approach to dimensioned mechanics is also motivated by the results proved for the theory of canonical contact phase spaces and their associated dimensioned algebraic structures. Although a detailed study of this possibility exceeds the scope of this thesis, it is certainly plausible that, given a fine-tuned notion of category of dimensioned manifolds $\textsf{DimMan}$, we can find a Hamiltonian functor relating the theories of phase spaces that directly generalizes the one presented above for the case of unit-free manifolds:
\begin{equation*}
\begin{tikzcd}[row sep=small]
 & \textsf{DimRing} \arrow[dd,"\text{Evl}"']& & & \textsf{DimPoissAlg} \arrow[dd,"\text{Evl}"] \\
 & & \textsf{DimMan} \arrow[ul, "\text{Obs}"']\arrow[dl,"\text{Dyn}"]\arrow[r,"\mathbb{H}"] & \textsf{DimPhase} \arrow[ur,"\text{Obs}"] \arrow[dr,"\text{Dyn}"'] & \\
 & \textsf{DimLie} & & & \textsf{DimLie} 
\end{tikzcd}
\end{equation*}
Which would be, indeed, the dimensioned analogue of the diagram given by the cotangent functor for conventional configuration spaces in Section \ref{CanonicalSymplectic}.

\chapter{Speculations on Dynamics, Integrability and Quantization} \label{SpeculationsDynIntQuant}

In this chapter we briefly discuss some remaining questions of a more conceptual nature regarding our discussion on the foundations of phase space and outline a few further lines of research that follow from the main mathematical results presented in this thesis.

\section{Comments on Observables, Smoothness, and Dynamics} \label{CommentsObsSmoothDyn}

The fundamental idea behind the notion of observable in a theory of phase spaces as introduced in Section \ref{FoundationsHamiltonianPhaseSpace} is that of the assignment of measurable characteristics to points in the space of states of a physical system. In all the concrete cases of theories of phase spaces presented in this thesis, starting from the conventional formulation of canonical Hamiltonian mechanics on cotangent bundles to the more general theories incorporating physical dimensions via line bundles, the spaces of observables exhibited Poisson algebra structures (or their dimensioned analogues). In fact, it was the search for such algebraic manifestation what motivated the introduction of the potentil functor for Jacobi phase spaces in the first place. This seems to point to this particular algebraic structure - Poisson algebras - as something intrinsically characteristic of classical mechanics. Indeed, this is the widely accepted viewpoint in the existing literature, as discussed in Section \ref{BriefHistoryOfPhaseSpace}. However, the general presentation of phase space theories of Section \ref{FoundationsHamiltonianPhaseSpace} doesn't seem to necessarily require the presence of commutative or Lie algebra structures to recover the overall categorical shape of classical mechanics.\newline

Hence, the present work shall be considered as an exploration of the Hamiltonian phase space theories whose observables do indeed carry Poisson-like structures. The list of principles of Section \ref{FoundationsHamiltonianPhaseSpace} provide a natural physical motivation of many of the formal aspects of conventional Hamiltonian mechanics and the several generalizations presented in subsequent chapters. Nevertheless, there are some structural components of a Hamiltonian theory of phase spaces that are left without a clear physical interpretation. Notably, the requirements that observables carry an associative multiplication and a bracket satisfying a Jacobi-like identity is hard to motivate from first principles, even when these properties are found in all the natural generalizations that seem to be compatible with a great deal of the physical interpretation of the basic motivating examples of conventional phase spaces.\newline

One way to approach this question is to consider the fact that early on in our presentation, motivated by the Principle of Refinement introduced in Section \ref{MeasurandFormalism}, any mathematical sets intended to abstractly capture physical entities were assumed to be smooth manifolds. Although this may seem a fairly arbitrary choice made to simplify, as it does, the mathematical development of the formal aspects of physical theories, from a categorical point of view, working within the class of smooth manifolds can be justified for entirely independent reasons. Recall that smooth manifolds were placed as the cornerstone of the commutative diagram of functors summarizing the theory of conventional phase spaces
\begin{equation*}
\begin{tikzcd}[row sep=small]
& \Ring \arrow[dd,"\text{Der}"] \\
\Man \arrow[ur,"\text{C}^\infty"] \arrow[dr,"\Tan"'] & \\
& \Lie_\Man 
\end{tikzcd}
\end{equation*}
where, critically, the interpretation of dynamics on manifolds as ring derivations is enabled by the isomorphism $\Sec{\Tan M}\cong\Dr{\Cin{M}}$. It should be pointed out that, far from a minor generic result, this isomorphism is a signature feature of smooth differential geometry within the larger context of algebraic geometry. Indeed, it is well known that in $\text{C}^k$-differentiable manifolds, the spaces of infinitesimal derivations at a point are infinite-dimensional and thus generically fail to coincide with the notion of directional derivative along a curve (see \cite{abraham2012manifolds}). This implies, in turn, that rings of $\text{C}^k$-differentiable functions on $\text{C}^k$-differentiable manifolds will generically fail to have non-zero global derivations, thus negating the possibility of a commutative diagram of functors similar to the one presented above for the $\text{C}^\infty$-smooth case. On the other hand, a classic result in differential topology ensures that any $\text{C}^1$-differentiable manifold contains a compatible smooth atlas (see \cite{munkres2016elementary}), which can be interpreted as the fact that smooth spaces can always approximate less regular ones with arbitrary accuracy, an thus the category of smooth manifolds is perceived to account for a large class of possible spaces, particularly when they are allowed to be possibly locally singular. On a related note, it is also an elementary fact that the derivations of rings of $\text{C}^0$-differentiable, i.e. continuous, functions of topological manifolds are all necessarily zero and that the derivations of general non-commutative algebras fail to be Lie algebras with the commutator bracket.\newline

We thus see that the choice to take smooth manifolds as phase spaces is geometrically motivated by the Principle of Refinement, involving the physical assumption that smaller scales are always measurable in principle, but it is equally motivated if we want the algebraic notion of derivation of a local ring-like structure on observables to be physically interpreted as the manifestation of time-evolution of states seen as curves on the phase space. This illustrates a general point to consider when further generalizing phase spaces: if one wishes to formulate the notion of theory of phase spaces with more exotic geometric objects, such as infinite-dimensional manifolds, singular spaces or non-commutative varieties, and ensure that the Principles of Kinematics, Dynamics and Observation can be accounted for so as to claim that some form of generalized mechanics has been defined, clear analogues of the rings of functions and the isomorphism between derivations and dynamics need to exist.\newline

Returning to the question about the algebraic structure of observables, we see that the discussion above gives a mathematical justification for the ubiquity of commutative and Lie algebra structures found in classical mechanics and the several generalizations proposed in this thesis. Turning our attention of the Jacobi identity of the Poisson-like structures on observables, in particular, we recall that The Principle of Conservation was implemented by demanding that there exist Hamiltonian maps from the ring-like structure of observables into the derivations of the ring-like commutative product that are identified with the dynamics (forcing a Leibniz-like property) that vanishes on self-images (forcing the map to be equivalently defined from an antisymmetric bracket on observables). From our discussion above, Lie-like algebraic structures seem to emerge naturally as the algebraic manifestation of dynamics, then the Jacobi identity of a the Leibniz and antisymmetric bracket induced on observables from the Principle of Conservation can be understood as the natural categorical requirement that the Hamiltonian maps are, furthermore, morphisms in the category of Lie-like algebras determined by the natural notion of dynamics in each concrete theory of phase spaces.

\section{Recovering Hamiltonian Dynamics on Unit-Free Phase Spaces} \label{RecoveringHamiltonianDyn}

Throughout the present work, the notion dynamics has been treated in a rather abstract way, focusing on how its algebraic and categorical manifestations fit in the larger picture of phase space theories rather than looking to the problem of integrating concrete differential equations on manifolds. It goes without saying that one of the core strengths of geometric mechanics is to give an intrinsic formulation to the equations of motion that have proven successful in the course of the history of physics, then, any realistic generalization of classical mechanics would need to encompass the dynamical equations of classical mechanics in a somewhat natural way.\newline

The results of sections \ref{CanonicalContact} and \ref{PotentialFunctorJacobi} led us to claim that line bundles regarded as unit-free manifolds constitute a valid generalization of the conventional notion of configuration space from a categorical perspective, however, unless we could recover the equations of motion of conventional Hamiltonian mechanics, this fact will become a mere mathematical curiosity with no immediate impact on the development of physically-realistic theories. In the general context of the theory of Jacobi phase spaces, proposition \ref{LocalUnitsJacobi} ensures that in general Jacobi manifolds we find Poisson algebras of local functions that are invariant with respect to the Hamiltonian vector field of a local unit. This means that the usual theory of Hamiltonian dynamics on Poisson manifolds can be recovered locally on Jacobi phase spaces. Furthermore, since an analogue of Darboux theorem exists for contact manifolds \cite[Proposition 4.9]{tortorella2017deformations}, we can express the dynamical equations of the Hamiltonian vector field of a choice of energy, represented as a local function $h$ dependent on the $2n+1$ Darboux coordinates $(q^i,p_i,z)$ defined on an open neighbourhood of a contact manifold $U\subset M$, as:
\begin{align*}
    \frac{dq^i}{dt} &=\frac{\partial h}{\partial p_i}\\
    \frac{dp_i}{dt} &=-\frac{\partial h}{\partial q^i}-p_i\frac{\partial h}{\partial z} \\
    \frac{dz}{dt} &=p_i\frac{\partial h}{\partial p_i}-h.
\end{align*}
Since the Hamiltonian vector field of the local unit $1\in\Cin{U}$ in the Darboux coordinates is simply $X_1=\tfrac{\partial}{\partial z}$, it is easy to see that restricting to $\tfrac{\partial}{\partial z}$-invariant functions, i.e. functions that only depend on $q^i$ and $p_i$, gives Hamilton's equations for the dynamics in those coordinates
\begin{align*}
    \frac{dq^i}{dt} &=\frac{\partial h}{\partial p_i}\\
    \frac{dp_i}{dt} &=-\frac{\partial h}{\partial q^i}.
\end{align*}
This indeed explicitly recovers the Hamiltonian dynamics found locally on Darboux charts of canonical coordinates for arbitrary symplectic phase spaces. More concretely, in the context of canonical contact phase spaces constructed as the jet bundles of unit-free configuration spaces, it is readily checked that a local trivialization on the line bundle over a configuration space turning $\lambda:L\to Q$ into $\Real_U$, with $U\subset Q$ the trivializing neighbourhood, induces a trivialization on $L_{\Jet^1 L}$ that allows to regard observables on the canonical contact phase space as real functions defined on $\Cot U \times \Real$. In this case, the Darboux coordinates $(q^i,p_i,z)$ correspond to the canonical adapted coordinates on the cotangent bundle $(q^i,p_i)$ and a global coordinate $z$ on $\Real$. Therefore, a choice of local unit in a unit-free configuration space naturally recovers the conventional notions of Hamiltonian dynamics on canonical symplectic phase spaces locally, as desired.

\section{Another Approach to Unit-Free Phase Spaces motivated by Line Bundle Dynamics} \label{AnotherApproachUnitFreeDyn}

The discussion in Section \ref{RecoveringHamiltonianDyn} above serves as a reality check that our proposal of unit-free phase spaces encompasses conventional Hamiltonian dynamics when trivializing coordinate charts are introduced on the base manifolds of the line bundles that represent the freedom of choice of unit over a set of states of a physical system. By taking a more structural point of view and setting aside, for a moment, the question of how to recover the conventional equations of motion of classical mechanics explicitly, we see that the theory of unit-free phase spaces proposed in chapter \ref{ContactPhaseSpaces} suggests a natural generalization of the fundamental notions of physical path and equations of motion.\newline

Recall that in the context of conventional phase spaces, time evolutions of a physical system were identified with smooth functions on the smooth manifold that represents the set of states of the system. Dynamics were identified with smooth vector fields $\Sec{\Tan M}$ that, in virtue of the standard results of ordinary differential equations on manifolds, can be locally integrated to give flows of smooth curves, thus connecting the algebraic notion of dynamics with the equations of motion of the physical system. If we take a similar approach in the context of unit-free phase spaces, where dynamics are identified with derivations of a line bundle $\Dr{L}$, we naturally encounter the problem of integrating derivations, regarded as infinitesimal automorphisms, to 1-parameter groups of line bundle automorphisms.\newline

Fixing an initial point $p_0\in P$ on the base manifold a unit-free phase space $\lambda: L\to P$, a solution of such integration problem gives a curve $p(t)\in P$ passing through the initial point $p(0)=p_0$, as in the case of ordinary dynamics, but also a 1-parameter family of fibre-to-fibre isomorphisms $B(t):L_{p_0}\to L_{p(t)}$ covering the curve $p(t)$. With the intuition that the line bundle over a phase space encodes the freedom of choice of unit for observables of the system of some fixed physical dimension, the 1-parameter family $B(t)$ has a clear interpretation as the information that will keep an arbitrary choice of units at a given fibre consistent throughout the time evolution of the system. In the case of a canonical contact phase space, this interpretation is apparent even when a local unit on a unit-free configuration space is introduced, since a choice of local $\tfrac{\partial}{\partial z}$-invariant energy generates non-trivial dynamics for the fibre coordinate
\begin{equation*}
    \frac{dz}{dt} =p_i\frac{\partial h}{\partial p_i}-h.
\end{equation*}

In light of the notions of Lie groupoid and Lie algebroid integrability introduced in Section \ref{LieGroupoids}, the problem of integrating a derivation of a line bundle to a 1-parameter family of line bundle automorphisms can be reformulated as an ordinary integration problem by regarding derivation $\Dr{L}=\Sec{\Der L}$ as left-invariant vector fields on the general linear groupoid $\GL{L}$ and noting that the der bundle of a line bundle is canonically isomorphic to the Lie algebroid of its general linear groupoid $A_{\GL{L}}\cong \Der L$. This calls for a broader generalization of unit-free phase spaces where the manifolds encoding the sets of states of a system and the freedom of choice of units of observables of fixed physical dimension are taken to be Lie groupoids that share basic structural features with general linear groupoids of line bundles.

\section{A Program to Investigate Line Bundle Dynamics and the possibility of Topological Effects} \label{ProgramDimDyn}

Since, following from the arguments in Section \ref{AnotherApproachUnitFreeDyn} above, line bundle dynamics can be given a broad physical interpretation as the unit-free analogues of conventional dynamics, at least locally, a natural question arises:
\begin{itemize}
    \item What is the physical significance of a line bundle representing a unit-free configuration space being generically non-trivializable?
\end{itemize}
And, directly related to this, we also pose the following:
\begin{itemize}
    \item Given the natural unit-free analogues of the dynamical theories of ordinary classical mechanics, e.g. geodesic flows, Newtonian potentials, magnetic terms, etc., what are the dynamics induced on the base manifolds? Do they reduce to the ordinary dynamics of classical mechanics on the base manifolds? Are they dependent on the fact that the line bundle may not be trivializable?
\end{itemize}

At the of Section \ref{CanonicalContact} L-metrics were introduced as the direct line bundle analogue of Riemannian metrics in an attempt to generalize Newtonian mechanics to the unit-free context. There, the notion of L-kinetic energy of a metric on a unit-free configuration space $L_Q$ gave a well-defined unit-free observable $K_\gamma\in\Sec{L_{\Jet^1 L}}$ and thus the contact Hamiltonian dynamics for such an observable give the line bundle analogue of the geodesic equation. A systematic study of such dynamical systems on line bundles appears then as a first natural approach to the question posed above.\newline

Such endeavour, being pursued by the author at the time of writing this thesis, goes beyond the scope of this work but without getting into too much detail, it is possible to illustrate that the theory of L-metrics is more general than a direct analogue of Riemannian geometry for line bundles, even in a locally trivial sense. We can see this by noting that L-metrics are, under a local trivialization giving the local description of the der bundle as $\Der L|_U\cong \Tan U \oplus \Real_U$, by non-degenerate bilinear forms
\begin{equation*}
    \gamma:(\Tan U\oplus \Real_U)\odot (\Tan U\oplus \Real_U)\to \Real_U.
\end{equation*}
Then we see that the datum of a L-metric is equivalent to a triple $(c,\alpha,g)$ with $c\in\Cin{U}$ a non-vanishing local function and $\alpha\in\Sec{\Cot U}$, $g\in \Sec{\odot^2\Cot U}$ satisfying the following non-degeneracy condition:
\begin{equation*}
    \begin{cases}
    fc+\alpha(X)=0\\
    f\alpha + g^\flat(X)=0
    \end{cases}
    \qquad \Rightarrow \qquad f=0, \quad X=0
\end{equation*}
for $f\in \Cin{U}$ and $X\in \Sec{\Tan U}$. By setting $\alpha=0$ the non-degeneracy condition holds iff $g$ is itself a non-degenerate bilinear form on $\Tan U$, which then clearly shows that L-metrics under this assumption recover the usual notion of Riemannian metric locally. However, this condition is not invariant under a change of general local trivialization and thus we see that generic L-metrics form a class of geometric structures on line bundles that encompass, in particular, Riemannian metrics on trivial line bundles. In light of the remarks made in Section \ref{AnotherApproachUnitFreeDyn} about line bundle dynamics interpreted as ordinary dynamics on Lie groupoids, we can speculate that the theory of L-metrics may correspond to the study of a particular subclass of Riemannian groupoids.\newline

Instead of considering the full generality of the problem of comparing geodesic dynamics of L-metrics to the ordinary geodesics of Riemannian metrics, it may be physically insightful to consider first a simple toy model that captures all the essential features that motivated the questions posed at the beginning of this section. On the matter of the physical significance of the trivializability of the line bundle representing the unit-free configuration space of a physical system, we can consider the simplest example of a configuration space whose topology allows for both a trivializable and non-trivializable line bundle over it. This is, of course, the circle $S^1$. On the circle we find two non-isomorphic line bundles: the cylinder $\text{Cyl}$, i.e. the trivial line bundle $S^1\times \Real$, and the M\"obius strip $\text{M\"ob}$, defined as the topological quotient of the rectangle with a twist.\newline

A first obvious difference between a mechanical theory on Cyl and on on M\"ob is that, in the cylinder, all the observables can be shifted by a constant function (that will not affect dynamics) to give non-vanishing global functions on $S^1$, whereas in the M\"obius strip all observables are global sections that vanish in at least one point. To explore the implications of this basic topological fact a consistent theory of unit-free Newtonian mechanics needs to be set up in order to endow both Cyl and M\"ob with a natural L-metric so that we can compare the trivial L-metric dynamics on Cyl, interpreted as the geodesic motion of a particle on $S^1$ equipped with the standard round metric, with the L-metric dynamics on M\"ob interpreted as the natural line bundle analogue.\newline

Even before any dynamical study is conducted, if we interpret the circle $S^1$ as the configuration space of positions of fixed physical system, and we interpret sections of the line bundles defined over it as observable physical characteristics of the system, we find a striking difference between the case of a cylinder and a M\"obius strip: since the same metric manifold $S^1$ forms the base of both Cyl and M\"ob, the smoothness (continuity) of sections implies that values of sections of Cyl repeat every integer number of loops around the circle whereas values of sections of M\"ob repeat every even number of loops around the circle. This clearly diverges from the usual physical intuitions behind configuration spaces but, provided a consistent dynamical theory is developed for line bundles, this may represent an avenue towards a classical theory of mechanics whose observables are sensitive to the topological properties of the spaces where they are defined, a feature that has been typically regarded as characteristic of quantum theories.

\section{Quantum Mechanics as a Phase Space Theory} \label{QuantumPhaseSpace}

Although we do not elaborate extensively on this topic here, we remark that the general principles of the phase space formalism outlined in Section \ref{FoundationsHamiltonianPhaseSpace}, indeed directly inspired by classical mechanics, apply to the usual formulation of quantum mechanics on Hilbert spaces in largely analogous ways except for a few notable differences. Despite the main points of contention in the discussion of classical vs. quantum being the process of measurement and the physical interpretation of measurement outcomes, trying to fit the dynamical aspects of conventional quantum physics into a theory of phase spaces in the sense of Section \ref{FoundationsHamiltonianPhaseSpace} proves illuminating in revealing the many geometric similarities between classical and quantum mechanics as well as a few of the fundamental differences that are entirely independent of the problem of measurement or interpretational issues.\newline

Phase spaces of quantum systems are commonly taken to be projective Hilbert spaces $\mathbb{P}H$, which, by a process of symplectic reduction of the imaginary part of the Hermitian form on $H$, are canonically endowed with the Fubini-Study symplectic form $\omega_H$. In this sense, canonical quantum phase spaces, i.e. projective Hilbert spaces, are subclasses of symplectic phase spaces and we can already notice a key fundamental difference with classical mechanics in the case of finite-dimensional quantum phase spaces: canonical classical phase spaces are given by cotangent bundles, always non-compact symplectic manifolds for non-trivial configuration spaces, whereas canonical quantum phase spaces are always compact, as they correspond to quotients of unit spheres on Hilbert spaces by compact groups. Observables are given by the subspace of quadratic functions on $\Cin{\mathbb{P}H}$ defined as the expectation values of Hermitian operators on $H$. These form a Lie subalgebra of the Poisson algebra induced by the Fubini-Study symplectic form on $\mathbb{P}H$ but are not closed under point-wise multiplication. It is worth noting that these functions span the cotangent spaces of projective Hilbert spaces everywhere. It is a simple computation that the Hamiltonian vector fields of these observables recover Schr\"odinger dynamics on $\mathbb{P}H$, thus we find an entirely dynamical motivation for the identification Hermitian operators of $H$ as the quantum observables. This leads to the main similarity between classical and quantum mechanics highlighted by the phase space formalism: the fundamental dynamical content of both theories corresponds to Hamiltonian dynamics in a category of canonical symplectic manifolds. Lastly, the combination product of quantum phase spaces constitutes a clear departure from the classical analogue: the correct notion of combined phase space of two quantum systems described by the Hilbert spaces $H_1$ and $H_2$ is the projective space of the tensor product $\mathbb{P}(H_1\otimes H_2)$. This construction is no longer a categorical product in a strict sense although it can be regarded as a generalization since the Segre map of projective spaces allows to give a canonical embedding $i:\mathbb{P}H_1 \times \mathbb{P}H_2 \hookrightarrow \mathbb{P}(H_1\otimes H_2)$. This is constitutes the geometric manifestation of the phenomenon of entanglement, where the points of $i(\mathbb{P}H_1 \times \mathbb{P}H_2)$ are identified as the separable states and the points of $\mathbb{P}(H_1\otimes H_2)\diagdown i(\mathbb{P}H_1 \times \mathbb{P}H_2)$ as the entangled states.

\section{Unit-Free Quantum Mechanics} \label{QuantumUnitFree}

We briefly remark that the revision of the foundations of metrology of Section \ref{MeasurandFormalism}, that later motivated the formal use of the category of lines and the potential functor to systematically account for the algebraic structure of physical dimensions and the arbitrary choices of units, is entirely independent of the concrete experimental regime where measurements will be performed. For this reason, a similar unit-free generalization of the conventional formulation of quantum mechanics can be performed by replacing the spaces of measurement outcomes (typically subsets of the real line) by subsets of measurand spaces of concrete physical dimension.\newline

A formal approach to this generalization, similar to the introduction of dimensioned manifolds as configuration spaces, that appears as a natural line of research to pursue in the future, would be to fine-tune definitions of dimensioned Hilbert spaces or dimensioned C*-algebras in such a way that the natural notion of observables defined there produces eigenvalues or evaluations with values on a dimensioned field and see what structural aspects of conventional quantum mechanics are retained or lost.

\section{Quantization as Integration of Phase Spaces by Groupoids} \label{QuantizationIntegration}

One of the great advancements in the field of Lie algebroids in recent years has been the development of the theory of integration by Lie groupoids. To the proofs of the algebroid equivalent of Lie's classic theorems introduced in Section \ref{LieGroupoids} followed many results that identified well-known geometric structures as the infinitesimal counterpart of certain multiplicative structures on an integrating Lie groupoid. In particular, and most relevant to our discussion, Poisson manifolds were shown to be the infinitesimal counterpart of symplectic groupoids (Lie groupoids with a symplectic 2-form that is compatible with the groupoid multiplication); similarly,  Dirac manifolds are shown to be the infinitesimal counterparts of presymplectic groupoids and Jacobi manifolds those of contact groupoids. From an independent stream of ideas originating in non-commutative geometry, the convolution algebras of some Lie groupoids were shown to be C*-algebras. Under the necessary technical assumptions, this presents an obvious quantization scheme in which a classical phase space is identified with a manifold with a geometric structure that is then integrated to the corresponding multiplicative structure on a Lie group and which, in turn, provides a C*-algebra via its convolution algebra. This is very much an active area of research and many fundamental problems are still open but it has been argued in the literature that this quantization scheme successfully generalizes geometric quantization and deformation quantization, at least in some well-known cases, while also giving general proofs of categorical statements such as the celebrated ``quantization commutes with reduction''.\newline

The fact there exist integration functors mapping categories of phase spaces back into categories of phase spaces, e.g. Poisson manifolds integrating to symplectic groupoids, which could, furthermore, lead to a quantization scheme, gives yet another angle to the legitimization of the proposed generalizations of phase spaces.\newline

The identification of the dimensioned Poisson algebra associated to a Jacobi manifold in theorem \ref{DimPoissonAlgebraJacobi} motivates the following integrability problem generalizing the integration of Jacobi manifolds by contact groupoids:
\begin{itemize}
    \item What are the objects integrating dimensioned Poisson algebra structures on the potential of a line bundle?
\end{itemize}
We can use the integrability result for Jacobi manifolds to safely speculate that the structure integrating the dimensioned Poisson algebra of an integrable Jacobi manifold is the multiplicative non-degenerate dimensioned Poisson algebra defined by the multiplicative contact structure on the integrating groupoid. Extending this result to non-Jacobi dimensioned Poisson algebras or, more generally, to Poisson-like structures on dimensioned manifolds is another natural line of inquiry to follow in the future.

\chapter{Conclusion} \label{Conclusion}

We began this thesis with the following question left to be made precise and resolved:
\begin{quote}
    \emph{Can a formalism that takes line bundles, the unit-free analogues of conventional configuration spaces, into contact manifolds provide a theory of mechanics that is analogous to canonical Hamiltonian mechanics and that naturally incorporates the notion of physical dimension?}
\end{quote}

In light of the results of chapter \ref{ContactPhaseSpaces} and sections \ref{UnitFreeRevisited} and \ref{JacobiManifoldsRevisited} we can affirm that, at the level of kinematics, the answer is an astounding  \emph{yes}. Firstly, our notions of theory of phase spaces and Hamiltonian functor allowed for a precise categorical characterization of the kinematic aspects of conventional canonical Hamiltonian mechanics so that any potential candidate for a generalization would be held up against that standard. The success of the introduction of our notion of unit-free manifold in Section \ref{UnitFreeManifolds} in reinterpreting Jacobi geometry as unit-free Poisson geometry, via L-vector bundle geometry and the potential functor in theorem \ref{DimPoissonAlgebraJacobi}, gave us reason to believe that unit-free manifolds will serve as configuration spaces as well. Our proofs for theorems \ref{JetLineProduct}, \ref{JetFunctor}, \ref{JetCoisotropicReduction} and \ref{JetHamiltonianReduction} confirmed that unit-free manifolds indeed fit the framework of a theory of configuration spaces with a Hamiltonian functor sending them to canonical contact manifolds.\newline

However, as we briefly discussed in Section \ref{RecoveringHamiltonianDyn} and \ref{AnotherApproachUnitFreeDyn}, there is still work do in order to prove that unit-free dynamics fit into the picture in a compatible way. A successful resolution of the problem of unit-free dynamics will likely involve some changes in the interpretation of physical paths and the question remains whether the intuitions behind the notion of configuration space can be kept unchanged after a fully dimensioned formalism of mechanics has been implemented.\newline

The breadth of results derived in sections \ref{DimRingsModules} and \ref{DimAlgebras} for dimensioned algebras, in complete analogy with the ordinary theory of rings and algebras, leads us to believe that our isolation of the axioms of dimensioned ring as the general model for the algebraic structure of physical quantities indeed constitutes a step in the right direction to incorporate physical dimension as a moving part in the mathematical formalism of scientific theories.\newline

One of the surprising results of our work was the identification of Jacobi manifolds with a special subclass of dimensioned Poisson algebras - those given by the image of the potential functor of theorem \ref{PotentialFunctorJacobi} with bracket of dimension $-1$. In a dimensioned mindset, this suggests a very natural class of geometric structures generalizing Jacobi manifolds: dimensioned Poisson algebras of arbitrary dimension on the potentil of a line bundle. It may be that some (or all) of these have already been identified in the literature, in which case it will be interesting to see whether new insights for the existing theory can be drawn from the dimensioned perspective. In a negative case the outcome might be equally interesting since the exploration of dimensioned Poisson algebras in their own right could lead to new instances of conventional Jacobi-like structures on line bundles.\newline

A first obvious context where our dimensioned perspective finds a home is, as pointed out by L. Vitagliano, in the identification of Jacobi structures with homogeneous Poisson structures, where, according to the definition found in \cite{vitagliano2019holo}, elements of the potential will be represented by polynomial functions of homogeneous degree. This offers clear avenues for future research that, as discussed above, may yield new insights for either of the two perspectives on Jacobi manifolds. A particular line of inquiry that may yield fruitful cross-pollination between standard Jacobi geometry and our dimensioned approach is the situation when we have a collection of line bundles and a dimensioned Poisson structure is present in their potential but not necessarily on each of its base line bundles. In the dimensioned picture, there is no essential difference between this kind of dimensioned Poisson algebra and the one associated to an ordinary Jacobi manifold, however, the underlying geometric data is quite different (a single line bundle versus a collection of line bundles. What would be, if any, the homogeneous Poisson equivalent of such dimensioned Poisson structures?\newline

Another natural development to follow is to adopt the formalism of positive spaces, which was introduced in \cite{janyvska2010algebraic} to give a rigorous mathematical axiomatization of conventional dimensional analysis where the notions of conversion factor and rational powers of physical units are formalized, into our dimensioned formalism. Geometrically, this will amount to the replacement of line bundles by positive bundles (fibre bundles with typical fibre $\Real^+$), and algebraically, this will amount to considering semi-vector spaces and rational tensor powers as a natural part of the mathematical theory that formalizes dimensional analysis. An interesting question would be to see whether the generalizations of dimensioned rings and algebras may require updating to incorporate dimension sets that are modelled after positive spaces. Work in this direction will also lead to a notion of positive manifold which, again, may or may not be a particular class of dimensioned manifolds as defined in Section \ref{DimManifolds}.\newline

In closing, we see that there are strong mathematical reasons to believe that our proposal for an axiomatization of the mathematical behaviour of physical quantities is on the right track for it to eventually account for all uses of conventional observables in mathematical physics. In the process of developing the formalisms of unit-free manifolds and dimensioned algebras, the central intuition has been that conventional models of physics are unit-less, in the sense that a unit that was entirely independent from the rest of the moving parts of the model was assumed for all the real variables, and we want to move to a unit-free model, in the sense that variables take real values upon freely choosing a unit. We can't help but notice a sort of irony in this transition from the unit-less to the unit-free. Indeed, we put all the mathematical technology in place to account for units to then give a description that is manifestly independent of units, very much in line with the quote proclaimed by European enlightened rulers: \emph{``Tout pour le peuple, rien par le peuple''}. Paraphrasing this motto, our working philosophy of introducing physical dimension into mechanics can be summarized as:
        \begin{center}
            Dimensioned Mechanics\\
            -\\
            All for the units, but without the units.
        \end{center}

\printbibliography

\appendix

\chapter{Supplementary Identities to be used in the proof of Theorem \ref{ProductJacobi}} \label{SymbolSquiggleProductJacobi}

Let $i=1,2$, two line bundles $\lambda_i:L_i\to M_i$ with Jacobi structures $(\Sec{L_i},\{,\}_i)$. Let $f_i,g_i\in\Cin{M}$ be smooth functions, $s_i\Sec{L_i}$ sections and $a,a'\in\Sec{L_1^\times}$, $b,b'\in\Sec{L_2^\times}$ local non-vanishing sections. Then the line product construction $L_1\utimes L_2$ allows to write the following identities
\begin{equation*}
    p_1^*f_1\tfrac{s_1}{b}=\tfrac{f_1\cdot s_1}{b} \qquad P_1^*s_1=\tfrac{s_1}{b}P^*_2b \qquad P_2^*s_2=\tfrac{s_2}{a}P^*_1a \qquad p_2^*f_2\tfrac{s_2}{a}=\tfrac{f_2\cdot s_2}{a}
\end{equation*}
The definition of the bracket is done via the following equations
\begin{equation*}
    \{P_1^*s_1,P_1^*s_1'\}_{12}:=P_1^*\{s_1,s_1'\}_1 \qquad \{P_2^*s_2,P_2^*s_2'\}_{12}:=P_2^*\{s_2,s_2'\}_2\qquad \{P_1^*s_1,P_2^*s_2\}_{12}:=0
\end{equation*}
which, in turn, imply
\begin{align*}
    X^{12}_{P_1^*s_1}[p_1^*f_1]&=p_1^*X^1_{s_1}[f] & X^{12}_{P_2^*s_2}[p_2^*f_2]&=p_2^*X^2_{s_2}[f]\\
    X^{12}_{P_1^*s_1}[p_2^*f_2]&=0 & X^{12}_{P_2^*s_2}[p_1^*f_1]&=0
\end{align*}
The definition of the symbol and squiggle is done via the following equations
\begin{align*}
    X^{12}_{P_1^*s_1}[\tfrac{a}{b}]&=\tfrac{\{s_1,a\}_1}{b}   & X^{12}_{P_2^*s_2}[\tfrac{b}{a}]&=\tfrac{\{s_2,b\}_2}{a}\\
    \Lambda^{12}(d\tfrac{a}{b}\otimes P_1^*s_1)[\tfrac{a'}{b'}]&=\tfrac{\{a,a'\}_1}{b}\tfrac{s_1}{b'} &\Lambda^{12}(d\tfrac{b}{a}\otimes P_2^*s_2)[\tfrac{b'}{a'}]&=\tfrac{\{b,b'\}_2}{a}\tfrac{s_2}{a'}
\end{align*}
which, in turn, imply
\begin{align*}
    \Lambda^{12}( dp_1^*f_1 \otimes P_1^*s_1 )[p_1^*g_1]&=p_1^*\Lambda^1(df_1\otimes s_1)[g_1] & \Lambda^{12}( dp_2^*f_2 \otimes P_2^*s_2 )[p_2^*g_2]&=p_2^*\Lambda^1(df_2\otimes s_2)[g_2]\\
    \Lambda^{12}( dp_1^*f_1 \otimes P_1^*s_1 )[p_2^*g_2]&=0 &\Lambda^{12}( dp_1^*f_1 \otimes P_1^*s_1 )[p_2^*g_2]&=0
\end{align*}
Then, using the symbol-squiggle identity inserting different combinations of spanning functions to obtain consistency relations, we obtain the following identities for the squiggle of the product Jacobi bracket
\begin{align*}
    \Lambda^{12}( dp_1^*f_1 \otimes P_1^*s_1 )[\tfrac{a}{b}] &= -p^*_1X^1_a[f_1]\tfrac{s_1}{b} \\
    \Lambda^{12}( dp_2^*f_2 \otimes P_2^*s_2 )[\tfrac{b}{a}] &= -p^*_2X^2_b[f_2]\tfrac{s_2}{a} \\
    \Lambda^{12}( dp_1^*f_1 \otimes P_2^*s_2 )[\tfrac{b}{a}] &= p_1^*X^1_a[f_1]\tfrac{s_2}{a}\tfrac{b}{a}\\
    \Lambda^{12}( dp_2^*f_2 \otimes P_1^*s_1 )[\tfrac{a}{b}] &= p_2^*X^2_b[f_2]\tfrac{s_1}{b}\tfrac{a}{b}\\
    \Lambda^{12}( dp_1^*f_1 \otimes P_2^*s_2 )[p_1^*g_1] &= -p_1^*\Lambda^1(dg_1\otimes a)[f_1]\tfrac{s_2}{a}, \quad \forall a\in\Sec{L_1^\times}\\
    \Lambda^{12}( dp_2^*f_2 \otimes P_1^*s_1 )[p_2^*g_2] &= -p_2^*\Lambda^2(dg_2\otimes b)[f_2]\tfrac{s_1}{b}, \quad \forall b\in\Sec{L_2^\times}\\
    \Lambda^{12}( dp_1^*f_1 \otimes P_2^*s_2 )[p_2^*s_2] &= 0\\
    \Lambda^{12}( dp_2^*f_2 \otimes P_1^*s_1 )[p_1^*s_1] &= 0 
\end{align*}

\chapter{Supplementary Calculations for the proof of Theorem \ref{DimPoissonAlgebraJacobi}} \label{CalculationsDimensionedPoissonBracket}

Let $\lambda: L\to M$ be a line bundle with a Jacobi structure $(\Sec{L},\{,\})$, the subindex on the dimensioned Poisson bracket $(\Sec{L}^\odot,\odot,\{,\}_{-1})$ will be omitted for simplicity. In what follows we take $a,b,c\in\Sec{L^1}=\Sec{L}$, $f,g,h\in\Sec{L^0}\cong \Cin{M}$ and $\alpha\in\Sec{L^{-1}}=\Sec{L^*}$.\newline

Consider the bracket $\{f\odot a,\alpha \odot b\}$, expanding as $\odot$-derivations we find
\begin{align*}
    \{f\odot a,\alpha \odot b\} &= a \odot \{f,\alpha\}\odot b + f\odot \{a,\alpha\}\odot b + f\odot \{a,b\}\odot \alpha + a\odot \{f,b\}\odot \alpha=\\
    &= \{f,\alpha\}(a,b) + fX_a[\alpha(b)]-f\alpha(\{a,b\})+f\alpha(\{a,b\})-X_b[f]\alpha(a).
\end{align*}
But from the definition of the $\odot$ dimensioned multiplication and the symbol-squiggle identity we have:
\begin{align*}
    \{f\odot a,\alpha \odot b\} &= \{f\cdot a,\alpha(b)\}=X_{f\cdot a}[\alpha(b)]\\
    &= fX_a[\alpha(b)]+\Lambda^\sharp(df \otimes a)[\alpha(b)],
\end{align*}
thus giving the desired bracket after simplifications:
\begin{equation*}
    \{f,\alpha\}=\Lambda^\sharp(df \otimes a)[\alpha(b)] + X_b[f]\alpha(a).
\end{equation*}

Consider the bracket $\{\alpha\odot a,\beta \odot b\}$, expanding as $\odot$-derivations we find
\begin{align*}
    \{\alpha\odot a,\beta \odot b\}(c) &= a \odot \{\alpha,\beta\}\odot b \odot c+ \alpha\odot \{a,\beta\}\odot b \odot c+ \alpha\odot \{a,b\}\odot \beta \odot c+ a\odot \{\alpha,b\}\odot \beta \odot c=\\
    &= \{\alpha,\beta\}(a,b,c) +\alpha(b)\Delta_a(\beta)(c) +\alpha(\{a,b\})\beta(c)-\beta(a)\Delta_b(\alpha)(c)=\\
    &= \{\alpha,\beta\}(a,b,c) - X_b[\alpha(a)]\beta(c)+\alpha(b)X_a[\beta(c)]-\alpha(b)\beta(\{a,c\}).
\end{align*}
But from the definition of the $\odot$ dimensioned multiplication and the definition of the bracket $\{f,g\}$ we have:
\begin{equation*}
    \{\alpha\odot a,\beta \odot b\}(c)=\{\alpha (a),\beta(b)\}(c)=\Lambda^\sharp(d\alpha(a)\otimes c)[\beta(b)],
\end{equation*}
thus giving the desired bracket after simplifications:
\begin{equation*}
     \{\alpha,\beta\}(a,b,c) = \Lambda^\sharp(d\alpha(a)\otimes c)[\beta(b)] + X_b[\alpha(a)]\beta(c)-\alpha(b)X_a[\beta(c)]+\alpha(b)\beta(\{a,c\}).
\end{equation*}

We explicitly show that the dimensioned Poisson bracket satisfies the Jacobi identity for all combinations of elements of dimension $+1$ and $0$, which, as argued in the proof of theorem \ref{DimPoissonAlgebraJacobi}, is enough to guarantee the Jacobi identity for all the other dimensions via extension as $\odot$-derivations. Firstly,
\begin{equation*}
    \{a,\{b,c\}\}=\{b,\{a,c\}\}+\{\{a,b\},c\}
\end{equation*}
follows from the Jacobi identity to the Jacobi bracket $(\Sec{L},\{,\})$ itself. Now 
\begin{align*}
    \{a,\{b,f\}\} = \{a,X_b[f]\} = X_a[X_b[f]] &= X_b[X_a[f]] + [X_a,X_b][f] =\\
    &= X_b[\{a,f\}] +X_{\{a,b\}}[f] =\\
    &=\{b,\{a,f\}\}+\{\{a,b\},f\},
\end{align*}
where we have used the fact that the symbol map $X$ is a Lie algebra morphism. Also
\begin{align*}
    \{a,\{f,g\}\}(b) = \Delta_a(\{f,g\})(b) &= X_a[\{f,g\}(b)]-\{f,g\}(\{a,b\}) =\\
    &= X_a[\Lambda^\sharp(df\otimes b)[g]]- \Lambda^\sharp(df\otimes \{a,b\})[g] =\\
    &=\Lambda^\sharp(df\otimes b)[X_a[g]] + [X_a,\Lambda^\sharp(df\otimes b)][g] - \Lambda^\sharp(df\otimes \{a,b\})[g] =\\
    &=\Lambda^\sharp(df\otimes b)[X_a[g]] + \Lambda^\sharp(dX_a[f]\otimes b)[g]+\Lambda^\sharp(df\otimes \{a,b\})[g] - \Lambda^\sharp(df\otimes \{a,b\})[g] =\\
    & =\Lambda^\sharp(df\otimes b)[\{a,g\}] + \Lambda^\sharp(d\{a,f\}\otimes b)[g]=\\
    & =\{f,\{a,g\}\}(b)+\{\{a,f\},g\}(b)
\end{align*}
where the symbol-squiggle compatibility condition 3 of proposition \ref{ExtensionBySymbol} has been used. Lastly
\begin{align*}
    \{f,\{g,h\}\}(a,b,c) &= \Lambda(df\otimes a, \Lambda(dg \otimes b , dh\otimes c)) + X_b[f]\cdot \Lambda(dg\otimes a, dh\otimes c)=\\
    &= \Lambda(dg\otimes a, \Lambda(df \otimes b , dh\otimes c)) + X_b[g]\cdot \Lambda(df\otimes a, dh\otimes c)\\
    & - \Lambda(dh\otimes a, \Lambda(df \otimes b , dg\otimes c)) - X_b[h]\cdot \Lambda(df\otimes a, dg\otimes c) =\\
    &=\Lambda(dg\otimes a, \Lambda(df \otimes b , dh\otimes c)) + X_b[g]\cdot \Lambda(df\otimes a, dh\otimes c)\\
    &+ \Lambda(\Lambda(df \otimes b , dg\otimes c), dh\otimes a) + X_b[h]\cdot \Lambda(dg\otimes a, df\otimes c) =\\
    &=\{g,\{f,h\}\}(a,b,c) + \{\{f,g\},h\}(a,b,c)
\end{align*}
where the symbol-squiggle compatibility condition 4 of proposition \ref{ExtensionBySymbol} has been used\includegraphics[width=0.069cm, height=0.0420cm]{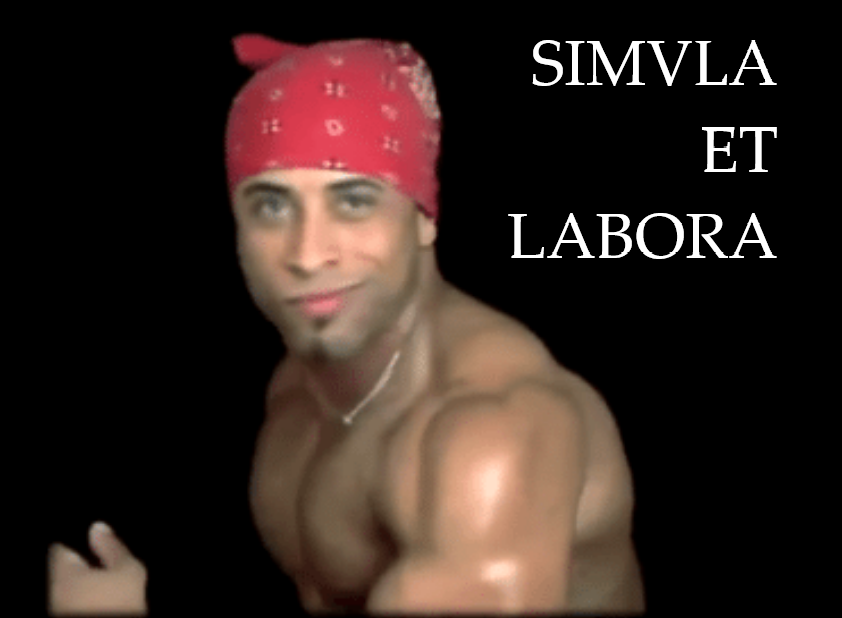}

\end{document}